%% file: 11_03_2022.tex
\documentclass[10pt, a4paper]{article}
\input{packages_notes}
\input{commandes_guillaume}
\usepackage[english]{babel}
\usepackage{braket}
\usepackage{authblk}

\usepackage{mathtools}

\begin{document}
\title{A Geometrical Analysis of Kernel Ridge Regression  and its Applications}

\author[1]{Georgios Gavrilopoulos, Guillaume Lecu{\'e}, and Zong Shang   \\ 
email: \href{mailto:georgios.gavrilopoulos@stat.math.ethz.ch}{georgios.gavrilopoulos@stat.math.ethz.ch},\\
email: \href{mailto:lecue@essec.edu}{lecue@essec.edu},\\
email: \href{mailto:zong.shang@ensae.fr}{zong.shang@ensae.fr} \\
ETH Zurich, Seminar for Statistics, Rämistrasse 101, 8092 Zurich, Switzerland.\\
ESSEC, business school, 3 avenue Bernard Hirsch, 95021 Cergy-Pontoise, France.\\ CREST, ENSAE, Institut Polytechnique de Paris, 5, avenue Henry Le Chatelier 91120 Palaiseau, France.}

\maketitle

\begin{abstract}
We obtain upper bounds for the estimation error of Kernel Ridge Regression (KRR) for all non-negative regularization parameters, offering a geometric perspective on various phenomena in KRR. As applications: 1. We address the Multiple Descents problem, unifying the proofs of \cite{liang_multiple_2020} and \cite{ghorbani_linearized_2021} for polynomial kernels in non-asymptotic regime and we establish Multiple Descents for the generalization error of KRR for polynomial kernel under sub-Gaussian design in asymptotic regimes. 2. In the non-asymptotic regime, we have established a one-sided isomorphic version of the Gaussian Equivalent Conjecture for sub-Gaussian design vectors. 3. We offer a novel perspective on the linearization of kernel matrices of non-linear kernel, extending it to the power regime for polynomial kernels. 4. Our theory is applicable to data-dependent kernels, providing a convenient and accurate tool for the feature learning regime in deep learning theory. 5. Our theory extends the results in \cite{tsigler_benign_2023} under weak moment assumption. 

Our proof is based on three mathematical tools developed in this paper that can be of independent interest: 1. Dvoretzky-Milman theorem for ellipsoids under (very) weak moment assumptions. 2. Restricted Isomorphic Property in Reproducing Kernel Hilbert Spaces with embedding index conditions. 3. A concentration inequality for finite-degree polynomial kernel functions.
\end{abstract}

\section{Introduction}
We focus on regression problems in the context of kernel learning. Let $\lambda\geq 0$ be a tuning parameter and $(\cH,\norm{\cdot}_\cH)$ be some Reproducing Kernel Hilbert space (RKHS) containing functions from the probability space $(\Omega,\mu)$ to $\bR$. Given $N$ independent design vectors $(X_i)_{i=1}^N$ distributed as the unknown probability measure $\mu$ and associated responses $(Y_i)_{i=1}^N\subset\bR$, let the Kernel Ridge Regression (KRR) estimator be
\begin{align}\label{eq:def_KRR}
    \hat f_\lambda \in \underset{f\in\cH}{\argmin}\left(\sum_{i=1}^N \left(f(X_i)-Y_i\right)^2 + \lambda\norm{f}_\cH^2\right).
\end{align}

KRR is a highly effective and adaptable method utilized in various domains, including finance, biology, natural language processing, image analysis and partial differential equations, \cite{shawe-taylor_kernel_2004,rasmussen_gaussian_2005,steinwart_support_2008,saitoh_theory_2016}. Additionally, it serves as a tool for developing mathematical foundations for deep neural networks, \cite{jacot_neural_2018,belkin_understand_2018}. This paper presents general bounds for the estimation error of KRR, aiming to investigate various phenomena in KRR and provide insights in the field of statistical deep learning.

\paragraph{From supervised learning theory to deep learning theory} 
The fundamental goal of supervised learning theory is to approximate an unknown function, denoted as $f^*$ (also referred to as the target function or signal), based on samples $(X_i,Y_i)_{i=1}^N$, where $Y_i = f^*(X_i) + \xi_i$, with $\xi_i$ representing noise terms independent of $X_i$. Our task is to select a functions class $\cF$ and construct an estimator $\hat f\in\cF$, such that, given the input of $N$ samples $(X_i, Y_i)_{i=1}^N$, $\hat f$ fits $f^*$ effectively. One metric for assessing the goodness of fit is the $L_2(\mu)$ distance $\|\hat f - f^*\|_{L_2}$; it coincides with the excess risk of $\hat f$. 

The excess risk of $\hat f$ primarily depends on two factors: approximation error and estimation error, \cite[Chapter 4]{bach_learning_2024}. In a broad sense, the approximation error describes whether the best element in $\cF$ can fit $f^*$ effectively, while the estimation error quantifies how far $\hat f$ deviates from this oracle. To reduce the approximation error, one approach is to increase the complexity of $\cF$, introducing high levels of non-linearity so that $\cF$ can describe sufficiently complex non-linear functions. Kernel methods and deep neural networks have both emerged as two approaches born from this line of thought to approach non-linear dependence between $Y$ and $X$.

\paragraph{The Core Issues of Deep Learning Theory} Introducing non-linearity increases the estimation error because as the approximation error decreases, $\cF$ becomes rich enough to contain many global minimizers of the empirical risk (even interpolant estimators, that is, functions $f$ such that $Y_i = f(X_i)$ for each $1\leq i\leq N$), leading in general to overfitting issues. This is known as the approximation-estimation trade-off \cite[section 13.3]{wainwright_high-dimensional_2019}. Classical statistical learning theory suggests that there must be a wisely-chosen regularization to learn from noisy data. As widely believed, an estimator should not interpolate the input data since it will suffer from interpolation of the noise in the data, which weakens its generalization properties. In \cite{hastie_elements_2009}, this phenomenon is called ``over-parameterization''. The motivation behind the appearance of \eqref{eq:def_KRR} is to introduce the regularization term $\|f\|_\mathcal{H}^2$ that induces some regularity in $\hat f_\lambda$ and prevents it from overfitting noise, since overfitted estimators tend to be non-smoothed.

However, as shown by the empirical results in \cite{zhang_understanding_2017}, ``over-parameterization'' may also lead to good statistical properties. In recent years, the statistical properties of interpolant estimators have attracted a lot of attention. In practical applications, deep neural networks, which possess high complexity and strong approximation capability \cite{schmidt-hieber_nonparametric_2020,suzuki_adaptivity_2018,yarotsky_error_2017}, are commonly trained as interpolant estimators (granting the global convergence of empirical risk), but with good generalization properties. Surprisingly, over-parameterized deep neural networks exhibit both a low approximation error and a low estimation error. Hence, the investigation into the phenomenon of overfitting estimators is often regarded as an essential task within the statistical theory of deep neural networks \cite{bartlett_deep_2021,belkin_fit_2021}.
A key component toward the construction of a statistical theory of deep learning is to analyze the implicit regularization exhibited by learning algorithms, specifically in terms of their selection of empirical risk minimizers. Additionally, researchers aim to investigate the statistical properties of these minimizers, and prove the the existence of benign overfitting, \cite{bartlett_benign_2020}.

\paragraph{KRR: A tool for Deep Learning Theory} Similar to deep neural networks, kernel methods exhibit a high degree of non-linearity while possessing characteristics that are relatively easier to study. This makes them a reference point for gaining insights into deep learning theory \cite{belkin_understand_2018}, as mentioned at the beginning of this paper. The connection between KRR and deep learning theory can be traced back to the concept of implicit regularization, \cite{jacot_neural_2018}.


Characterizing the implicit regularization of algorithms, such as gradient descent, stochastic gradient descent, and ADAM, \cite{kingma_adam_2017}, for the purpose of training neural networks is a complex and challenging task, \cite{jacot_neural_2018}. The parameter trajectory of even a deep \emph{linear} neural network during training exhibits a significantly non-linear behavior over time, \cite{saxe_exact_2014,chizat_infinite-width_2022}.
It is worth noting, however, that under specific initialization conditions, known as the static Neural Tangent Kernel (NTK) parameterization, the training path can be linearized. This observation is supported by the works of \cite{jacot_neural_2018,chizat_lazy_2019,yang_tensor_2020}. Under (static) NTK parameterization, it has been demonstrated (see for example \cite{boyer_living_2022}), that a deep neural network trained using the first-order gradient method, converges towards a global minimizer of the empirical risk. Furthermore, this minimizer has the lowest possible RKHS norm (thus equals to $\hat f_0$ in \eqref{eq:def_KRR}), where the RKHS is generated by the neural tangent kernel, see for example \cite{boyer_living_2022}. Through this avenue, kernel methods establish a connection with deep learning theory. However, there are numerous criticisms surrounding the concept of static NTK, for instance, \cite{woodworth_kernel_2020,donhauser_how_2021,bietti_deep_2021}.


\paragraph{Towards the feature learning regime} Many studies have already indicated that the estimation error of static NTK is significantly larger than that of deep neural networks. This is primarily due to the fact that static NTK does not exhibit alignment of its eigenfunctions with the unknown function (called lazy regime \cite{chizat_lazy_2019}). This issue is considered a significant drawback of static NTK. In fact, it is widely acknowledged that the strength of Deep Neural Networks (DNNs) lies in their powerful feature construction capability. In the context of kernel machines, this implies that the optimization algorithm for DNNs should automatically select a RKHS in such a way that the eigenfunctions of this RKHS align well with the unknown function, resulting in a smaller estimation error for the associated KRR \eqref{eq:def_KRR}. In addressing this issue, the authors in \cite{ba_high-dimensional_2022,long_properties_2021,radhakrishnan_feature_2022,damian_neural_2022,dandi_how_2023,moniri_theory_2023} examine data-dependent kernels of neural networks within the context of the \emph{feature learning} regime. 


These data-dependent kernels have eigenfunctions that rely on the sample points, which means that to some extent, they can learn from data how to align with the unknown function. This implies that such data-dependent kernels possess the characteristics of feature learning while also being firmly grounded in the theory of RKHS, making them amenable to effective theoretical research. It is for these reasons that data-dependent kernels are highly promising. However, unlike the RKHS typically studied in the past, these data-dependent kernels are neither translation-invariant nor rotationally invariant. That is why obtaining more general results (regarding both kernels and stochastic assumptions) for KRR is expected to be a promising useful tool in statistical deep learning theory.
The general foundation of KRR must be general, not dependent on specific translation invariance or rotational invariance of RKHS, and not contingent on the properties of the RKHS spectrum, such as the power decay of eigenvalues. Furthermore, this bound should be user-friendly, allowing mathematicians studying training dynamics to use it conveniently. At the same time, this bound must be accurate enough to reflect the estimation error of KRR. Finally, the bound must be general enough to be manageable for the study of DNNs.

\subsection{Our contributions}
\begin{center}
\begin{enumerate}
    \item \textbf{The primary contribution of this paper is a general, precise and user-friendly bound for the estimation error of KRR.}
\end{enumerate}
\end{center}
To illustrate the precision and convenience of our general bounds, we applied it to various problems, including Multiple Descents, Gaussian equivalent property, the linearization of kernel matrices, and the estimation error of KRR of data-dependent conjugate kernels, resulting in numerous outcomes.

\begin{enumerate}[start=2]
    \item In asymptotic scenario, we demonstrate the occurrence of the Multiple Descents phenomenon for the generalization error of KRR in RKHS with polynomial kernel, with the requirement of sub-Gaussian design, thus generalizes \cite{ghorbani_linearized_2021}. This is considered a well-known challenging problem. In the non-asymptotic case, we improve upon the conclusions of \cite{liang_multiple_2020}.
    It is worth noting that we present a \emph{unified} proof for these two fundamentally different setups, which has been considered a challenging problem, \cite{donhauser_how_2021}. 

    \item In \emph{non-asymptotic} regime, we establish the Gaussian equivalent property for the \emph{upper bound} of the estimation error of KRR in polynomial RKHS, with the stipulation of \emph{sub-Gaussian} design or spherical uniform distribution. 

    \item We obtain a geometric perspective on the linearization of kernel random matrices corresponding to polynomial kernels. Additionally, we have extended the proportional regime to the power regime (which will be explained later).

    \item In the specific case of linear ridge regression, we extend the results of \cite{tsigler_benign_2023} to a weaker set of moment assumptions.

    \item For smooth inner product kernels including NTK defined over a Euclidean sphere, we have derived upper bounds on the estimation error of KRR in the non-asymptotic regime. This upper bound is optimal in a minimax sense, which demonstrates that our general upper bounds are sharp compared to \cite{mourtada_elementary_2022}.

    \item For data-dependent conjugate kernels, we have verified that the conditions outlined in our general bounds are satisfied. This underscores the potential of our general bounds as a tool to feed a mathematical deep learning theory. We provide the excess risk for the single-neuron model when feature learning occurs.

    \item Our general bounds hold significant geometric implications, and our proof relies on three novel mathematical tools we developed:
    \begin{enumerate}
        \item A Dvoretzky-Milman theorem for ellipsoids under (very) weak moment assumptions, see Theorem~\ref{theo:DM_RKHS}.
        \item Restricted Isomorphic Property for RKHS with embedding index condition, see Proposition~\ref{prop:RIP}.
        \item A concentration inequality for finite-degree polynomial kernel function of a subgaussian vector, see Theorem~\ref{theo:polynomial_HS}.
    \end{enumerate}
We believe these new mathematical tools will be useful for addressing other statistical problems beyond KRR.
    \item We provide a general upper bound for the excess risk of KRR and sufficient conditions for benign overfitting to occur in the case of model misspecification, as discussed in Proposition~\ref{prop:dependent_noise_result} and subsequent discussions. This contributes to the examination of feature learning and generalization properties of deep neural networks. Specifically, our generalization bound takes the form of an oracle inequality, meaning that we can choose any oracle in our bound, rather than being required to select the projection of $f^*$ onto the closure on $\cH$ as the oracle (for example, \cite[Section 7.5.2]{bach_learning_2024}), thus providing an interface between the theory of approximation error (see \cite{devore_neural_2021} and references therein) and the theory of estimation error (see \cite{bartlett_deep_2021} and references therein) for deep neural networks.
\end{enumerate}

\subsection{Literature review}\label{sec:literature_review}

Thus far, motivated by the theory of DNNs, research into KRR has uncovered some remarkable and intricate phenomena. Within this context, several predominant trends can be identified:

\paragraph{Benign Overfitting} The concept of benign overfitting refers to phenomena in which the estimation error of an interpolant estimator approaches zero as the number of samples and model parameters increase infinitely when the number of parameters exceeds the number of samples. In other words, this estimator demonstrates consistency. It has been demonstrated that in the context of linear regression, many minimum norm interpolant estimators are consistent \cite{bartlett_benign_2020,tsigler_benign_2023,lecue_geometrical_2022,koehler_uniform_2021,wang_tight_2022,donhauser_fast_2022}, given certain assumptions regarding the relationship between the covariance structure of the design vector and the unknown signal.

\paragraph{Inconsistency} In contrast to the conclusions drawn about Benign Overfitting, there exist findings regarding inconsistency. These findings indicate that for certain RKHSs, including those generated by the (static) NTK defined over a sphere, the Laplace kernel and Sobolev kernel, the minimum RKHS norm interpolant estimator is inconsistent. In other words, its excess risk cannot be expected to approach zero as $N$ increases (the dimension of design vectors is fixed). Several results have reported this phenomenon, including \cite{rakhlin_consistency_2019,donhauser_how_2021,buchholz_kernel_2022,haas_mind_2023,li_asymptotic_2023}. The significance of such conclusions lies in demonstrating that for these RKHS, benign overfitting is not viable for relatively general unknown functions.

\paragraph{Multiple Descents phenomenon}  Starting from \cite{mei_generalization_2020,hastie_surprises_2022}, the double descents phenomenon of the estimation error of the minimum RKHS norm interpolant estimators $\hat f_0$ has garnered significant interest. In broad terms, the phenomena being discussed pertains to the observation that in the over-parameterized regime, the estimation error of the minimum RKHS norm interpolant estimator diminishes, despite the fact that it over-fits the training data. 
    Subsequently, \cite{ghorbani_linearized_2021,mei_generalization_2022} demonstrated the occurrence of the Multiple Descents phenomenon in asymptotic regime.  This phenomenon pertains to the existence of Multiple Descents of the estimation error of the minimum RKHS norm interpolant estimator (including static NTKs). Notably, this descent is observed when $N\sim d^\iota$ for some $\iota\in\bN_+$ and $N,d\to\infty$ (called the power regime). The significance of the Multiple Descents conclusion is to illustrate that as $N$ increases (with $d$ also increasing, but $N$ being a power of $d$), the estimation error of KRR is not monotonic. For each $\iota$, KRR learns a degree-$\iota$ polynomial approximation of the target function. Therefore, if the target function is a degree-\(\iota\) polynomial or can be well approximated by such a polynomial, over-fitting of the minimum RKHS norm interpolant estimator will be benign.

    
    \paragraph{Linearization}  In the process of proving double descents, many studies make use of a tool known as the ``linearization of kernel matrices'' developed in the studies conducted, for instance, by \cite{el_karoui_spectrum_2010,do_spectrum_2013,cheng_spectrum_2013,fan_spectral_2019}. In broad terms, this category of tools establishes the following fact: when $d$ and $N$ both tend to infinity, and their ratio remains at a constant level (referred to as the ``proportional regime''), the spectrum of the kernel matrices defined by inner product kernels and translation-invariant kernels can be approximated by the kernel matrix of a linearized kernel. Since the properties of KRR depend heavily on the kernel matrix, in a sense, in the proportional regime, KRR degenerates into a form of linear regression. This means that even with highly nonlinear kernel functions, in this scenario, KRR can only learn a linear approximation of the target function. Consequently, KRR does not offer an advantage over linear regression in this context. This observation has also been made in the context of double descents \cite{mei_generalization_2020,hastie_surprises_2022}.
    
    \paragraph{Gaussian Equivalence}
    The Gaussian Equivalent Conjecture is another novel phenomenon that has been frequently observed recently \cite{goldt_modelling_2020,aubin_generalization_2020,seddik_random_2020,dhifallah_precise_2020}. It posits that in the analysis of KRR, we can replace the feature map with a centered Gaussian random vector and replace the integral operator of the RKHS with the covariance matrix of this Gaussian random vector. The excess risk obtained from the analysis of this linear model is nearly identical to the original excess risk of KRR. From this perspective, the Gaussian Equivalent Conjecture is somewhat similar to a variant of the central limit theorem but in a non-asymptotic way. The introduction of this phenomenon is meant to facilitate the study of KRR's excess risk. This, from another angle, underscores the main theme of this paper, namely that the universality theory of KRR's excess risk is eagerly anticipated in the field of Deep Learning theory.
\paragraph{Further questions} 
\begin{enumerate}
    
    \item In \cite{ghorbani_linearized_2021,mei_generalization_2022}, Multiple Descents is proven only when both $d$ and $N$ tend to infinity. It is crucial to note that when $d$ varies, it leads to changes in the sample distribution, RKHS, and $f^*$, making it challenging to determine whether KRR in this power regime or proportional regime still approximates the initial target function we sought to approximate.
    Therefore, the question arises: \textbf{Can we prove Multiple Descents in the classical context with $d$ fixed?} Additionally, \cite{ghorbani_linearized_2021,mei_generalization_2022} only provide proofs for the case when the design vector $X$ follows a uniform distribution over the Euclidean sphere, which is highly limiting. Hence, a natural question is: \textbf{is it possible to establish Multiple descents for a fixed target function under more general stochastic assumptions?}
    This problem is widely recognized as highly challenging, \cite{haas_mind_2023,donhauser_how_2021}.

    
    \item The same question arises as well for the Gaussian Equivalent Property: \textbf{Does the Gaussian Equivalent Property hold true in  non-asympotic? for distributions other than the Gaussian or uniform over the Euclidean sphere distributions? other than the proportional regime?}
    
\end{enumerate}

\subsection{Structure of the paper}

Starting from this section, we focus on specific statistical problems with their rigorous mathematical derivations. In Section \ref{sec:notations}, we introduce various basic mathematical notations required for this paper. In Section~\ref{sec:RKHS}, we provide fundamental knowledge about RKHS. Due to the length of this paper, we provide an informal version of the main conclusions in Section~\ref{sec:RKHS}. Since the main proofs in this paper are built upon the decomposition of RKHS, we also introduce some relevant notations in this section. In Section~\ref{sec:DM_RKHS} and Section~\ref{sec:RIP}, we introduce the two primary mathematical tools used in this paper: a Dvoretzky-Milman theorem (with its proof in Section~\ref{sec:proof_DM} in supplementary material) and a Restricted Isomorphy Property (with its proof in Section~\ref{sec:proof_RIP} in supplementary material). 


In Section~\ref{sec:polynomial_kernel}, we present the upper bounds (with proofs in Sections \ref{sec:proof_main_upper_k<N} and \ref{sec:proof_main_upper_k>N} in supplementary material) for the estimation error of KRR. 
In Section \ref{sec:multiple_descent}, we apply these tools to the study of the Multiple Descents phenomenon. In Section \ref{sec:KRR_smooth}, we apply them to the study of smooth kernel functions (not necessarily finite-degree polynomial kernels), where we show that even for smooth kernels, the upper bound for KRR is nearly identical to that of linear ridge regression. This leads to the content of the next section, Section \ref{sec:gaussian_equivalence}, where we demonstrate that the upper side of the Gaussian Equivalence Property holds true in the context of finite-degree polynomial kernels. In Section~\ref{sec:conjugate_kernel}, we apply our results to the data-dependent conjugate kernel, obtaining an upper bound on the estimation error of the KRR defined by it.

In Section \ref{sec:future_directions} in supplementary material, we will list several potential research directions. For the sake of adhering to the length constraint, we place the content of non-linear random matrix linearization in Section \ref{sec:linearization} in supplementary material. Readers focused solely on statistical questions can skip this section. The proofs for other theorems, propositions, lemmas, etc., can be found in Section \ref{sec:aux_proofs} in supplementary material.

\subsection{Notations}\label{sec:notations}


We use $c,c_0,c_1,\cdots, C,C_0,C_1,\cdots$ to denote absolute constants. Usually, $C$ stands for large but finite constants and $c$ stands for small but non-zero constants. Such constants are always assumed to be positive, but may change from one instance to another. Given two quantities $A,B$, we write $A\lesssim B$ (or $A\gtrsim B$) if there exists an absolute constant $C$ such that $A\leq CB$ (or $A\geq CB$). If a constant is assumed to depend on some parameter (say $K$), we use expression $C_K$, and write $A\lesssim_K B$ (or $A\gtrsim_K B$) if there exists an absolute constant $C_K$ such that $A\leq C_KB$ (or $A\geq C_K B$). We write $A\sim B$ if $B\lesssim A\lesssim B$. We use $O_d(\cdot)$ (respectively $o_d(\cdot)$) for big-O (respectively small-O) notations, where $d$ emphasizes the asymptotic variable. Further, if $g(d) = O_d(f(d))$, we write $f = \Omega_d(g)$ and if $g(d)=o_d(f(d))$, we write $f=\omega_d(g)$. Given a sequence of random variables $(Z_d)_d$ and deterministic sequences $(s_d)_d$, we say $Z_d = o_{d,\bP}(s_d)$ if $Z_d/s_d\to 0$ in probability.


Given $r\in\bN_+$, we let $[r] := \{1,2,\cdots,r\}$. Let $(\Omega,\mu)$ be a probability space and, for $q\in\bN_+$, let $L_q(\Omega,\mu)$ be the $L_q$ space with norm $\norm{f}_{L_q}=\left(\int_{\Omega}\left|f(x)\right|^q\right)^{1/q}d\mu(x)$. When there is no ambiguity, we abbreviate $L_q(\Omega,\mu)$ as $L_q(\mu)$ or $L_q$. Given some random variables $X_1$, we write $\bE_{X_1}$ for the conditional expectation with respect to $X_1$ conditionally on all other random variables. Given probability measures $\mu_1,\mu_2$, we denote by $\mu_1\times\mu_2$ the product probability measure. We say a real-valued random variable $X$ is  sub-Gaussian, if $\norm{X}_{\psi_2}:=\inf\left(t>0: \bE\exp\left(X^2/t^2\right)\leq 2\right)<\infty$. We write $\cN(0,1)$ as the standard Gaussian random variable, and $\cN(0,I_d)$ be the standard Gaussian random vector in $\bR^d$.


Given two Hilbert spaces $\cH_1,\cH_2$, we characterize $\cH_1\otimes\cH_2$ by defining $f\otimes g\in \cH_1\otimes\cH_2$ as the mapping $(f\otimes g):h\in\cH_2\mapsto \left\langle g,h\right\rangle_{\cH_2}f\in\cH_1$. We use angle brackets $\left\langle \cdot,\cdot\right\rangle_{\cH}$ to denote an inner product in some Hilbert space $\cH$, and omit it when $\cH$ is a Euclidean space. We use $\norm{\cdot}_\cH$ to denote its Hilbert norm, and denote the Euclidean norm by $\norm{\cdot}_2$. Given a bounded linear operator $T:\cH_1\to\cH_2$, we denote $\norm{T}_{\text{op},\cH_1\to\cH_2}$ as the operator norm of $T$, that is,
$$\norm{T}_{\text{op},\cH_1\to\cH_2}=\sup\left(\norm{Tx}_{\cH_2}: \norm{x}_{\cH_1}\leq 1\right).$$
When there is no ambiguity, we abbreviate the operator norm as $\norm{T}_{\text{op}}$. For any ONB $(\varphi_j)_{j\in\bN}$ of $\cH_1$, if $\sum_{j\in\bN}\norm{T\varphi_j}_{\cH_2}^2<\infty$, we say $T$ is a Hilbert-Schmidt operator from $\cH_1$ to $\cH_2$, and denote 
$$\norm{T}_{HS,\cH_1\to\cH_2}=\sqrt{\sum_{j\in\bN}\norm{T\varphi_j}_{\cH_2}^2}.$$
One can prove that the HS norm of $T$ is independent of the choice of ONB, for instance, \cite[pp.7]{pisier_volume_1989}. When there is no ambiguity, we abbreviate the Hilbert-Schmidt norm of $T$ as $\norm{T}_{HS}$.  For any $j\in\bN$, let $\mathrm{He}_j(x)$ be the $j$-th (probabilist) Hermite polynomial, \cite[pp.16]{pisier_volume_1989}. For a Hilbert space $(\cH,\norm{\cdot}_\cH)$, we let $B_\cH=\{f\in\cH:\, \norm{f}_\cH\leq 1\}$ and $S_\cH = \{f\in\cH:\, \norm{f}_\cH=1\}$.


Let $A\in\bR^{m\times n}$ for some $m,n\in\bN_+$, let $\mathrm{spec}(A)=(\sigma_1(A),\cdots,\sigma_{m\wedge n}(A))$ and $(\sigma_i(A))_{i\leq m\wedge n}$ are singular values of $A$ and $\1_N\in\bR^N$ is the vector with all its coordinates equal to $1$. Given $N\in\bN_+$, we use $(e_i)_{i\in[N]}$ to denote an arbitrary ONB of $\ell_2^N$. Denote by $\cI$, a $d$-index, that is, given $i\in\bN$, we denote $\cI$ as a partition of $[i]=\cI_1\sqcup\cdots\sqcup\cI_d$ into $d$ (potentially empty) groups. With a little abuse of notation, we denote $\left|\cI\right|:=i$ and the cardinality of an empty set is $0$.

Notations on Reproducing Kernel Hilbert Space will be introduced in Section~\ref{sec:RKHS}.

\section{Tools}
Our analysis of KRR relies on several tools that are all exposed in this section. We start with classical tools and notation that are related to RKHS.

\subsection{Reproducing Kernel Hilbert Spaces}\label{sec:RKHS}

In this section, essential RKHS background knowledge is presented. For readers who are not familiar with RKHS, they can consider the upcoming discussion of $\phi(X)$ as if it were a Gaussian random vector (which is precisely the essence of the Gaussian Equivalent Property).

\paragraph{Structural aspect}Let $\Omega\subset\bR^d$ be a compact Hausdorff space\footnote{We emphasize that we do not actually require \(\Omega\) to be compact; it suffices for \(K\) to have a spectral decomposition.}, $\mu$ be a probability measure on $\Omega$. Let $L_2(\Omega,\mu)$ be the space of real-valued, square-integrable functions with respect to $\mu$. Suppose $K:\Omega\times\Omega\to\bR$ is a positive definite continuous function, and without loss of generality, we assume that $\norm{K}_\infty\leq 1$\footnote{We remark that we do not really need $\norm{K}_\infty\leq 1$, but need only the integral operator $\Gamma$ defined below to be a positive, compact, symmetric trace-class operator.}. In Section~\ref{sec:conjugate_kernel}, we will present a crucial example of a RKHS that does not satisfy $\norm{K}_{L_\infty(\mu\times\mu)}<\infty$. Nevertheless, our analysis remains valid.
We say a Hilbert space $\cH\subset L_2(\mu)$ of functions is an RKHS (with RKHS inner product $\left\langle \cdot,\cdot\right\rangle_\cH$ and associated RKHS norm $\norm{\cdot}_\cH$) over $\Omega$ if for every $x\in\Omega$, there exists a constant $C_x>0$ (depending on $x$), such that $\left|f(x)\right|\leq C_x\norm{f}_\cH$ for every $f\in\cH$, that is, the evaluation functional $\mathrm{ev}_x: f\mapsto f(x),\, \mathrm{ev}_x:\cH\to\bR$ is a bounded linear functional. By Riesz's representation theorem for Hilbert space, this is equivalent to saying that the inner product of $\cH$ can be characterized as follows: for every $x\in\Omega$, there exists $K(x,\cdot)\in\cH$ such that $\left\langle f,K(x,\cdot)\right\rangle_\cH = f(x) = \left\langle f,\phi(x)\right\rangle_\cH$, called the \emph{Reproducing Property}. Given a kernel $K$, the (canonical) feature map is defined as follows: $\phi:x\mapsto K(x,\cdot)\in\cH$. The RKHS can be considered as a linear model on $\cH$, where the design vector $X\in\Omega\subset\bR^d$ is embedded into $\cH$ via the feature map $\phi$; hence $\phi(X)$ plays the role of a design vector and  may therefore be called the RKHS design vector. 

By mapping $X$ from $\bR^d$ to $\cH$, we clarify the prediction function, covariance structure and the design matrix.
Given that $\phi(X)$ now serves as the design vector, its integral operator   $\Gamma = \bE\left[\phi(X)\otimes\phi(X)\right] $ (i.e. $\Gamma : f\in\cH \to \bE\left[\phi(X)\left\langle \phi(X),f\right\rangle_\cH\right]\in\cH$) will play the role of the covariance matrix. We denote the operator norm of $\Gamma$ from $\cH$ to $\cH$ as $\norm{\Gamma}_{\text{op}}$. Since $\norm{K}_\infty\leq 1$,  $\Gamma$ is compact and positive semi-definite, so it has a discrete spectrum of non-negative eigenvalues. By Mercer's theorem, see, for instance \cite[Theorem 12.20]{wainwright_high-dimensional_2019} or \cite[section 4.3]{rasmussen_gaussian_2005}, for any $\vx,\vy\in\Omega$, $K(\vx,\vy)=\left\langle \phi(\vx),\phi(\vy)\right\rangle_\cH=\sum_{j\geq1}\varphi_j(\vx)\varphi_j(\vy) = \sum_{j\geq1}\sigma_jf_j(\vx)f_j(\vy)$ where $(\varphi_j)_{j\in\bN^*}$ are eigenfunctions of $\Gamma$, and $(\sigma_j)_{j\geq1}$ are the eigenvalues of $\Gamma$ associated to $(\varphi_j)_{j\geq1}$, and where $\varphi_j = \sqrt{\sigma_j}f_j$. It is possible to check that this decomposition is unique. By \cite[section 30.5]{lax_functional_2002}, $\Gamma$ is a trace-class operator, that is, $\Tr\left(\Gamma\right) <\infty$ where $\Tr\left(\Gamma\right) = \bE\norm{\phi(X)}_\cH^2 = \bE K(X,X) = \sum_{j\in\bN^*}\sigma_j$.

It is also widely-used to embed $\cH$ into $\ell_2$ by $\phi:\vx\in\Omega\mapsto \sum_{j=1}^\infty \sqrt{\sigma_j}f_j(\vx)\ve_j$ (we use the same notation as for the feature map $\phi:x\to K(x, \cdot)$ introduced above; however, which definition of $\phi(x)$ we used is clear from the context). Therefore for every $\vv\in\ell_2$, the corresponding element of $\cH$ is $f_\vv(\vx)=\left\langle \phi(\vx),\vv\right\rangle_{\ell_2}$. Therefore, $\cH$ can be defined equivalently as the image of $\ell_2$ under the map $\vv\mapsto f_\vv$, with the inner product $\left\langle f_\vv,f_\vu\right\rangle_\cH = \left\langle \vv,\vu\right\rangle_{\ell_2}$ when $\sigma_j>0$ for all $j$ (in case $\Gamma$ is of rank $r$, $\cH$ is in a one-to-one correspondence with $\ell_2^r$).

We will frequently use a decomposition of $\cH$ so it is necessary to introduce the following notation.
Denote $\Gamma_{p:q}=\sum_{p\leq j\leq q}\sigma_j\varphi_j\otimes\varphi_j$. Given $k\in\bN\cup\{\infty\}$ and $p<q\in\bN\cup\{\infty\}$, we decompose $\Gamma = \Gamma_{1:k} + \Gamma_{k+1:\infty}$. We let $\cH_{p:q}=\Span\left(\varphi_j: p\leq j\leq q\right)$ and $P_{p:q}$ be the orthogonal projection onto $\cH_{p:q}$, and for any $f\in\cH$, denote $P_{p:q}f$ as $f_{p:q}$, for example, $\phi_{p:q}(X)=\sum_{p\leq j\leq q}\left\langle \phi(X),\varphi_j\right\rangle_\cH\varphi_j$. Consequently, we decompose $\cH = \cH_{1:k}\oplus^\perp \cH_{k+1:\infty}$. Given $\iota\in\bN_+$, we denote $f_{\leq\iota}^* := f_{1:\sum_{l\leq\iota}d^l}^*$, $f_{>\iota}^* := f^* -  f_{\leq\iota}^*$, $\Gamma_{\leq \iota}:= \Gamma_{1:\sum_{l\leq\iota}d^l}$ and $\Gamma_{> \iota}:= \Gamma - \Gamma_{\leq\iota}$, $P_{\leq l}$ the projection onto the eigen-space of $\Gamma_{\leq l}$, and we denote by $P_{>l}$ its complementary.

By reproducing property, for any $f\in\cH$,
\begin{eqnarray}\label{eq:RKHS_covariance}
    \norm{f}_{L_2(\mu)}^2 = \bE\left\langle \phi(X), f\right\rangle_\cH^2 = \norm{\Gamma^{1/2}f}_\cH^2 \leq \norm{\Gamma}_{\text{op}}\norm{f}_\cH^2.
\end{eqnarray}

Define the RKHS design matrix $\bX_{\phi}:\cH\to\bR^N$ as
\begin{eqnarray*}
    \bX_\phi = \left(\begin{matrix}
        \phi(X_1)^\top\\\vdots\\\phi(X_N)^\top
    \end{matrix}\right),\mbox{ so that }\bX_\phi f = \left(\begin{matrix}
        \left\langle \phi(X_1),f\right\rangle_\cH\\\vdots\\\left\langle \phi(X_N),f\right\rangle_\cH
    \end{matrix}\right)=\left(\begin{matrix}
        f(X_1)\\\vdots\\f(X_N)
    \end{matrix}\right),
\end{eqnarray*}
where for all $x\in\Omega$, we use $\phi^\top(x):\cH\to\bR$, that is, the operator $\phi^\top(x):f\mapsto \left<\phi(x),f\right>_\cH = f(x)$.

\paragraph{Statistical model and closed-form solution to \eqref{eq:def_KRR}}
In this paper, we always assume $f^*\in\cH$ (except  in Remark~\ref{remark:model_misspecified}). Let $\bxi = (\xi_i)_{i\in[N]}$ be the noise vector with i.i.d. zero mean and variance $\sigma_\xi^2$ coordinates that are independent with design vectors $(X_i)_{i\in[N]}$. Now the kernel ridge estimator $\hat f_\lambda$ for $(X_i,Y_i)_{i\in[N]}\subset \left(\Omega\times \bR\right)^N$ and the tuning parameter $\lambda\geq 0$ is defined as
\begin{eqnarray*}
    \hat f_\lambda \in\underset{f\in\cH}{\mathrm{argmin}}\left( \norm{\bX_{\phi}f - \vy}_2^2 +\lambda\norm{f}_\cH^2 \right),\mbox{ where }Y_i = f^*(X_i) + \xi_i,\mbox{ and } \vy = (Y_1,\cdots,Y_N)^\top,
\end{eqnarray*} which coincides with \eqref{eq:def_KRR}. By \cite[Proposition 12.33]{wainwright_high-dimensional_2019}, $\hat f_\lambda$ has an explicit solution:
\begin{align*}
    \hat f_\lambda = \bX_{\phi}^\top\left(\bX_{\phi}\bX_{\phi}^\top + \lambda I_N\right)^{-1}\vy.
\end{align*}

As $(\phi(X_i))_{i\in[N]}\subset\cH$, we write  $\bX_{\phi,p:q}=\left( (P_{p:q}\phi(X_1))|\cdots | (P_{p:q}\phi(X_N)) \right)^\top$ thus we can decompose $\bX_{\phi}$ into two parts for any $k\in\bN_+$:
\begin{align*}
    \bX_\phi = \left(\begin{matrix}
        (P_{1:k}\phi(X_1))^\top\\
        \vdots\\
        (P_{1:k}\phi(X_N))^\top
    \end{matrix} \right) + \left(\begin{matrix}
        (P_{k+1:\infty}\phi(X_1))^\top\\
        \vdots\\
        (P_{k+1:\infty}\phi(X_N))^\top
    \end{matrix} \right) =: \bX_{\phi,1:k} + \bX_{\phi,k+1:\infty}.
\end{align*}We recall that for all $f\in\cH$, $f^\top:g\in\cH\mapsto \left<f,g\right>_\cH\in\bR$, that is, $f^\top g = \left<f,g\right>_\cH$.


\paragraph{Key quantities driving the rate of convergence of KRR}
We first define some key quantities related to the convergence rate of KRR, and then provide an informal version of the main conclusions of this paper. For the formal version, please refer to Theorem~\ref{theo:upper_KRR} in Section~3 and Theorem~\ref{theo:main_upper_k>N} in the Supplementary material. Let $k\in\bN_+$, and let $\kappa_{DM}$ be some absolute constant, define
\begin{align}\label{eq:def_J1}
    J_1 := \left\{ j\in[k]: \sigma_j\geq \frac{\kappa_{DM}\left(4\lambda + \Tr\left(\Gamma_{k+1:\infty}\right)\right)}{N} \right\},\, J_2 = [k]\backslash J_1,
\end{align}
and
\begin{align*}
    \begin{aligned}
        &\tilde\Gamma_{1,\mathrm{thre}}^{-1/2} = \sum_{j=1}^k \left(\sigma_j\vee \frac{ \kappa_{DM}\left(4\lambda + \Tr\left( \Gamma_{k+1:\infty} \right)\right)}{N}\right)^{-1/2}\varphi_j\otimes\varphi_j.
    \end{aligned}
\end{align*}For the optimal choice of $k$, we will have $J_1=[k]$ and $\tilde\Gamma_{1,\mathrm{thre}}^{-1/2} = \Gamma_{1,k}^{-1/2}$.  The main conclusion of this paper is given now.
\begin{Theorembis}{theo:upper_KRR}
    Under certain conditions on the kernel, the design vector and the noise, for any \( k\leq N \) that satisfies certain conditions, the following fact holds with high probability:
    \begin{align*}
        \norm{\hat f_\lambda - f^* }_{L_2(\mu)} \lesssim r_{\lambda,k}^*,
    \end{align*}where
    \begin{align}
        \begin{aligned}\label{eq:def_rate_informal}
            r_{\lambda,k}^* &= \sigma_\xi\sqrt{\frac{\left|J_1\right|}{N}}+ \sigma_\xi\sqrt{\frac{\sum_{j\in J_2}\sigma_j}{  4\lambda + \Tr\left(\Gamma_{k+1:\infty}\right)}}+ \norm{\Gamma_{k+1:\infty}^{1/2}f_{k+1:\infty}^*}_\cH+\\
        &{\norm{\tilde\Gamma_{1,\mathrm{thre}}^{-1/2}f_{1:k}^*}_\cH \frac{2\lambda + 3\Tr\left(\Gamma_{k+1:\infty}\right)}{N} } + \sigma_\xi\frac{\sqrt{N\Tr\left(\Gamma_{k+1:\infty}^2\right)}}{\lambda + \Tr\left(\Gamma_{k+1:\infty}\right)}.
        \end{aligned}
    \end{align}
\end{Theorembis}
We provide some comments on Theorem~\ref{theo:upper_KRR}.
In general, based on the choice of \(\lambda\) and the spectrum of \(\Gamma\), KRR automatically decomposes the feature space \(\mathcal{H}\) into a direct sum of two subspaces: \(\mathcal{H} = \mathcal{H}_{1:k} \oplus^\perp \mathcal{H}_{k+1:\infty}\). In $\cH_{1:k}$, KRR estimates the target function \(f^*\), while in $\cH_{k+1:\infty}$, KRR absorbs noise. Here, \(k\) is a key quantity that determines the estimation error of KRR. For each \(k\), KRR performs an associated decomposition of \(\mathcal{H}\), utilizing different subspaces for estimation, while the remaining subspace absorbs noise, leading to varying estimation errors $r_{\lambda,k}^*$. Therefore, there exists an optimal \(k\) that minimizes the rate defined by that choice among all possible \(k\). We emphasize that \(k\) is a free parameter; that is, for KRR, statisticians do not select a priori a specific \(k\); rather, this \(k\) is needed for analysis. Thus, this \(k\) is determined by KRR itself; somehow KRR learns by itself the best features space decomposition. In the space \(\mathcal{H}_{1:k}\) (where estimation happens), \([k]\) is further divided into \(J_1\) and \(J_2\), which correspond to different parts in the final bound. In most applications, \( J_1 = [k] \) and \( J_2 = \emptyset \). In the classical case (for instance, when \(\lambda \sim N \sigma_k\), in Proposition~\ref{prop:jaouad}), the dominating term in the definition of \(\tilde{\Gamma}_{1,\mathrm{thre}}^{-1/2}\) is \(\sigma_j\), which means that KRR does not filter out the largest \( k \) eigenvalues. However, theoretically, there may exist choices of \( k \) and \(\lambda\) such that KRR filters out only certain eigenvalues, with the number of these eigenvalues being less than \( k \).
In \(\mathcal{H}_{k+1:\infty}\), \(f_{k+1:\infty}^*\) is treated as noise, which leads to the term \(\|\Gamma_{k+1:\infty}^{1/2} f_{k+1:\infty}^*\|_\mathcal{H}\) in \eqref{eq:def_rate_informal}. The space \(\mathcal{H}_{k+1:\infty}\) is responsible for noise absorption, corresponding to the term \(\sigma_\xi \frac{\sqrt{N \Tr\left( \Gamma_{k+1:\infty}^2 \right)}}{\lambda + \Tr\left( \Gamma_{k+1:\infty} \right)}\) in \eqref{eq:def_rate_informal}. Here, it is necessary for \(\Tr(\Gamma_{k+1:\infty}^2)\) to be relatively small compared to \(\lambda^2\) and \((\Tr(\Gamma_{k+1:\infty}))^2\) in order to effectively absorb noise. When \(\lambda\) dominates \(\Tr(\Gamma_{k+1:\infty}\) (which corresponds to the case of sufficiently strong regularization), this implies that \(\lambda\) needs to be large enough to absorb noise. Conversely, if \(\Tr(\Gamma_{k+1:\infty})\) dominates (which corresponds to insufficient regularization), we require the spectrum of \(\Gamma_{k+1:\infty}\) to be well spread, meaning that the ratio of its \(\ell_2\) norm to its \(\ell_1\) norm must be sufficiently small.

Next, we introduce the main tools used to prove Theorem~\ref{theo:upper_KRR}. Following that, we compare Theorem~\ref{theo:upper_KRR} with classical bounds and discuss its applications.

\subsection{A Dvoretzky-Milman theorem for RKHS}\label{sec:DM_RKHS}

Below is the probabilistic version of the Dvoretzky-Milman theorem \cite{milman_asymptotic_2023, pisier_volume_1989}.

\begin{Theorem}There are absolute constants $\kappa_{DM}\leq1$ and $c_0$ such that the following holds. Let $\vertiii{\cdot}$ be some norm on $\bR^p$ and denote by $B$ its unit ball and by $B^*$ its unit dual ball. The Dvoretsky-Milman dimension of $B$ is $d_*(B) = \left(\frac{\ell^*(B^*)}{{\rm diam}(B^*, \ell_2^p)}\right)^2$,
where $\ell^*(B^*) = \bE\vertiii{G}$, $G$ is a standard Gaussian random vector in $\bR^p$, ${\rm diam}(B^*, \ell_2^p) = \sup\left(\norm{v}_2:v\in B^*\right)$.
 Denote by $\bG:=\bG^{(N\times p)}$, the $N\times p$ standard Gaussian matrix with i.i.d. $\cN(0,1)$ Gaussian entries.  Assume that $N\leq \kappa_{DM} d_*(B)$. Then with probability at least $1-\exp(-c_0 d_*(B))$, for every $\blambda\in\bR^N$, 
\begin{equation*}
\frac{1}{\sqrt{2}}\norm{\blambda}_2 \ell^*(B^*) \leq \vertiii{\bG^\top \blambda}\leq \sqrt{\frac{3}{2}}\norm{\blambda}_2 \ell^*(B^*). 
\end{equation*}
  \end{Theorem}

The probabilistic version of the Dvoretzky-Milman theorem posits that, with high probability, the random section of a convex body (denoted as $B$), which is generated by its intersection with the image of a Gaussian random matrix, exhibits properties that are nearly Euclidean with radius of the order of the Gaussian mean width $\ell^*(B)=\bE\sup(\left<G,\vv\right>:\, \vv\in B)$ where $G$ is a standard Gaussian random vector, under the condition that the dimension of the image is smaller than the Dvoretzky-Milman dimension \cite{pisier_volume_1989}. This note focuses on the scenario where $B$ is equal to the ellipsoid $\Gamma_{k+1:\infty}^{-1/2}B_\cH$, where $B_\cH$ represents the unit ball in $\cH$. Therefore, in this paper, we are only concerned with the Dvoretzky-Milman theorem for ellipsoids and so in  that case the Dvoretzky-Milman theorem states that $\bX_{\phi,k+1:\infty}^\top$ is well-conditioned, with both smallest and largest singular values of the order of $\sqrt{\Tr(\Gamma_{k+1:\infty})}$, the Gaussian mean width of  $\Gamma_{k+1:\infty}^{-1/2}B_\cH$.

\paragraph{A warm up: Dvoretzky-Milman theorem under a strong assumption} One might question the validity of a DM theorem in a RKHS due to several challenges when compared to the Gaussian scenario. Firstly, the feature map is no longer linear in $X$, which may result in a lower concentration to counterbalance metric complexity. Secondly, the feature map is not centered. Lastly, it exhibits non-independent coordinates. The subsequent Lemma, extracted from \cite{mcrae_harmless_2022} (see also \cite{mei_generalization_2022}), demonstrates that, despite the aforementioned challenges, the Dvoretzky-Milman theorem remains applicable to $\bX_{\phi,k+1:\infty}^\top$, subject to a more stringent requirement.

\begin{Lemma}[\cite{mcrae_harmless_2022}]\label{lemma:mcrae}
Recall the notations from Section~\ref{sec:RKHS}. Then,
    \begin{align*}
        &\bE\norm{\bX_{\phi,k+1:\infty}\bX_{\phi,k+1:\infty}^\top - \Tr\left(\Gamma_{k+1:\infty}\right)I_N}_{\text{op}}^2 \\
        &\leq 2 N^2\Tr\left(\Gamma_{k+1:\infty}^2\right) + 2\sup\left( \left|K_{k+1:\infty}(\vx,\vx) - \Tr(\Gamma_{k+1:\infty}) \right|^2:\, \vx\in \Omega \right).
    \end{align*}
\end{Lemma}


It is easy to see from \cite[Lemma 4.1.5]{vershynin_high-dimensional_2018} together with Lemma~\ref{lemma:mcrae} that if $N^2\Tr\left(\Gamma_{k+1:\infty}^2\right)\lesssim \Tr^2\left(\Gamma_{k+1:\infty}\right)$, and if $K_{k+1:\infty}(\vx,\vx)$ concentrates well, the Dvoretzky-Milman theorem holds for $\Gamma_{k+1:\infty}^{-1/2}B_\cH$ with constant probability. This Lemma addresses any doubts regarding the validity of the Dvoretzky-Milman theorem for non-linear, non-centered, and dependent-coordinate features. On the other hand, we observe that $N^2\Tr\left(\Gamma_{k+1:\infty}^2\right)\lesssim \Tr^2\left(\Gamma_{k+1:\infty}\right)$ is more restrictive than the classical Dvoretzky-Milman condition $N\lesssim \Tr\left(\Gamma_{k+1:\infty}\right)/\norm{\Gamma_{k+1:\infty}}_{\text{op}}$ obtained in the Gaussian case. One may wonder whether it is possible to establish the Dvoretzky-Milman theorem under this weaker assumption.

The authors of \cite[Theorem 3]{mcrae_harmless_2022} establish the minimality of $N^2\Tr\left(\Gamma_{k+1:\infty}^2\right)\lesssim \Tr^2\left(\Gamma_{k+1:\infty}\right)$ by demonstrating that the condition number of $\bX_{\phi,k+1:\infty}\bX_{\phi,k+1:\infty}^\top$ exceeds $N^4/d^2$ with a high probability when the $\varphi_j$'s are Fourier basis functions. However, this counter-example is somewhat pessimistic; we will demonstrate that, by imposing a stronger assumption, we can broaden the applicability of the Dvoretzky-Milman theorem from $N^2\Tr\left(\Gamma_{k+1:\infty}^2\right)\lesssim \Tr^2\left(\Gamma_{k+1:\infty}\right)$ to $N\lesssim \Tr\left(\Gamma_{k+1:\infty}\right)/\norm{\Gamma_{k+1:\infty}}_{\text{op}}$, the optimal Dvoretzky-Milman dimension appearing in the Gaussian case, and we can recover the classical condition of the Dvoretzky-Milman theorem for ellipsoids in the Gaussian case. We also obtain a DM theorem that holds with large probability instead of in expectation.

\paragraph{Dvoretzky-Milman theorem for ellipsoid under an $L_{2+\eps}-L_2$ moment-equivalence assumption} For any $\lambda\geq 0$, define
    \begin{equation}\label{eq:def_DM_dimension}
        d_\lambda^*(\Gamma_{k+1:\infty}^{-1/2}B_\cH) :=
        \frac{\Tr(\Gamma_{k+1:\infty})+\lambda}{\norm{\Gamma_{k+1:\infty}}_{\text{op}}}
    \end{equation} as the modified (by $\lambda$) Dvoretzky-Milman dimension for the RKHS associated with the projected kernel $K_{k+1:\infty}$. When $\lambda=0$, $d_0^*(\Gamma_{k+1:\infty}^{-1/2}B_\cH)$ is the classical Dvoretzky-Milman dimension (up to a universal constant), \cite[pp. 42]{pisier_volume_1989}.  We grant the following assumption.
\begin{Assumption}\label{assumption:DM_L4_L2}
There are absolute constants $\nC\label{C_DM}>1$, $\oC{C_distortion_1}>1$, $0\leq\gamma<1/16$, $0\leq \delta<1/(100\sqrt{\oC{C_distortion_1}})$, $\bar\delta<\oC{C_DM}^{-1}$, $\eps>0$ and $\kappa>1$ such that
\begin{itemize}
    \item With probability at least $1-\gamma$,
    \begin{equation}\label{eq:diagonal_term_assumption}
        \max_{1\leq i\leq N}\left|\frac{\norm{\phi_{k+1:\infty}(X_i)}_\cH^2}{(\ell^*)^2} - 1 \right| \leq \delta,
    \end{equation}where we define $\ell^* = \sqrt{\bE\norm{\phi_{k+1:\infty}(X)}_\cH^2}=\sqrt{\Tr\left(\Gamma_{k+1:\infty}\right)}$.
    \item For any $f\in\cH_{k+1:\infty}$, we have
    \begin{equation}\label{eq:norm_equivalence}
    \norm{f}_{L_{2+\eps}}\leq \kappa \norm{f}_{L_2}.
    \end{equation}
    \item Depending on the choice of $\eps$, there are two cases:
    \begin{enumerate}
        \item if $\eps>2$, then no extra assumption is required.
        \item if $0<\eps\leq 2$, then
        \begin{align}\label{eq:extra_assumption}
        \kappa N^{\frac{2-\eps}{2\eps+\eps^2}}\log{(N)}\left(\frac{\sqrt{N\Tr\left(\Gamma_{k+1:\infty}^2\right)}}{\Tr\left(\Gamma_{k+1:\infty}\right)}\right)<\bar\delta.
    \end{align}
    \end{enumerate}
\end{itemize}
\end{Assumption}

The assumption described in \eqref{eq:diagonal_term_assumption} asserts that the diagonal components of the kernel matrix, denoted as
$\left(\bX_{\phi,k+1:\infty}\bX_{\phi,k+1:\infty}^\top\right)_{ii} = K_{k+1, \infty}(X_i, X_i)$,
exhibit a high level of concentration around their expected value $\ell^*$. In the context of a translation-invariant kernel, it can be observed that \eqref{eq:diagonal_term_assumption} is always satisfied for $\delta=\gamma=0$ (see Section~\ref{sec:diagonal_terms_multiple_descent} in supplementary material). Similarly, for an inner-product kernel, \eqref{eq:diagonal_term_assumption} can be validated by employing concentration inequalities of polynomials. For a data-dependent kernel, the verification of \eqref{eq:diagonal_term_assumption} relies on the specific properties of that kernel. We provide an example in Proposition~\ref{prop:conjugate_kernel} below.

When $\eps>2$, there is no need for an extra condition on $N$ and, as in the Gaussian case, everything holds under the classical DM condition $N\leq\kappa_{DM}d_\lambda^*\left(\Gamma^{-1/2}B_\cH\right)$, even though we work under a weak $L_{2+\eps}-L_2$ equivalence assumption; when $0<\eps\leq 2$, we need the extra condition \eqref{eq:extra_assumption}, which is an additional assumption on $N$ and the spectrum of $\Gamma$.

The assertion made in \eqref{eq:extra_assumption} is that the convergence rate of $\sqrt{N\Tr\left(\Gamma_{k+1:\infty}^2\right)}/\Tr\left(\Gamma_{k+1:\infty}\right)$ is greater than that of $\left(N^{\frac{2-\eps}{2\eps+\eps^2}}\log{N}\right)^{-1}$ up to a universal constant (as $\bar\delta$ and $\kappa$ will be chosen to be constant). Condition \eqref{eq:extra_assumption} can also be compared with \cite{mcrae_harmless_2022,mei_generalization_2020}, where the more stringent assumption of $\sqrt{N\Tr\left(\Gamma_{k+1:\infty}^2\right)}/\Tr\left(\Gamma_{k+1:\infty}\right)\leq 1/\sqrt{N}$ is granted. In fact, when $\eps\geq 2\sqrt{2}-2$, our result is an improvement on \cite{mcrae_harmless_2022,mei_generalization_2020}. We conjecture that \eqref{eq:extra_assumption} can be removed even when $0<\eps\leq 2$ (maybe by using a clever adaptation of the coloring technique from \cite{tikhomirov_sample_2018}).

Our main result from this section is the following, and its proof may be found in Section~\ref{sec:proof_DM}.
\begin{Theorem}\label{theo:DM_RKHS}
    Let $X$ be a random vector distributed as $\mu$ in a compact set $\Omega\subset\bR^d$, and let $X_1,\cdots,X_N$ be i.i.d. copies of $X$. Let $\phi:x\in\Omega\mapsto K(x,\cdot)\in\cH$ be the feature map of the RKHS $\cH$. Let $\oC{C_comparision_trace_lambda}$ be an absolute constant.
    
    \begin{enumerate}
        \item If $\lambda\leq \oC{C_comparision_trace_lambda}\Tr\left(\Gamma_{k+1:\infty}\right)$. Let $0<\delta,\bar\delta<1$ from Assumption~\ref{assumption:DM_L4_L2}, define
    \begin{align}\label{eq:def_tilde_delta}
    \tilde\delta = \oC{C_distortion_1}\delta^2 + \oC{C_distortion_2} \bar\delta^2 + 4\sqrt{\left(3\delta+\oC{C_distortion_3}\bar\delta\right)\left(1+\delta+\oC{C_distortion_4}\bar\delta\right)}.
    \end{align}
    Suppose that for some $\lambda\geq 0$, we have $N\leq\kappa_{DM}\bar\delta^2d_\lambda^*\left(\Gamma_{k+1:\infty}^{-1/2}B_\cH\right)$  for a sufficiently small constant $\kappa_{DM}<1$ which depends only on $\kappa$ (see \eqref{eq:bound_kappa_DM} for a precise description).
        We assume that $\phi_{k+1:\infty}$ satisfies Assumption~\ref{assumption:DM_L4_L2}. Then with probability at least 
    \begin{equation*}
        1-\gamma - \frac{1}{N^2} - \left(\frac{\kappa}{\bar\delta}\right)^{2+\eps}\left(\frac{\sqrt{N\Tr\left(\Gamma_{k+1:\infty}^2\right)}}{\Tr\left(\Gamma_{k+1:\infty}\right)}\right)^{2+\eps}\frac{\log^{2+\eps}(N)}{N^{\frac{\eps}{2}-1}} =: 1-\bar p_{DM},
    \end{equation*}
    for all $\vlambda\in \bR^N$,
    \begin{equation}\label{eq:objective_DM}
        \left(1-\tilde\delta\right)\sqrt{\Tr(\Gamma_{k+1:\infty})}\norm{\vlambda}_2 \leq \norm{\bX_{\phi,k+1:\infty}^\top\vlambda}_\cH \leq \left(1+\tilde\delta\right)\sqrt{\Tr(\Gamma_{k+1:\infty})}\norm{\vlambda}_2.
    \end{equation}
    \item If $\lambda> \oC{C_comparision_trace_lambda}\Tr\left(\Gamma_{k+1:\infty}\right)$. Suppose that $\phi_{k+1:\infty}$ satisfies the first two points of Assumption~\ref{assumption:DM_L4_L2}. Suppose that for some $\lambda\geq 0$, we have $N\leq(\kappa_{DM}/4)d_\lambda^*(\Gamma_{k+1:\infty}^{-1/2}B_\cH)$.  There then exist absolute constants $\nC\label{C_distortion_9}$ depending on $\eps,\kappa,\kappa_{DM}$, and $0<\nc\label{c_distortion_9}<1$ such that with probability at least
    \begin{align*}
        1 - \gamma - N\left(\left(\frac{\kappa_{DM}\kappa^2\log^2(N) }{ N }\right)^{1+\eps/2} N\right)^{\lceil(12+2\eps)/\eps \rceil-1} - \frac{1}{N^2}=:1-\bar p_{DM},
    \end{align*}we have
    $\norm{\bX_{\phi,k+1:\infty}\bX_{\phi,k+1:\infty}^\top + \lambda I}_{\mathrm{op}} \leq \oC{C_distortion_9}\lambda+\Tr\left(\Gamma_{k+1:\infty}\right)$ and $\sigma_N\left(\bX_{\phi,k+1:\infty}\bX_{\phi,k+1:\infty}^\top + \lambda I\right)\geq \oc{c_distortion_9}\lambda + (1-\oc{c_distortion_9})\oC{C_comparision_trace_lambda}\Tr\left(\Gamma_{k+1:\infty}\right)$.
    \end{enumerate}
\end{Theorem}

We can make some remarks comparing Theorem~\ref{theo:DM_RKHS} to Lemma~\ref{lemma:mcrae}: 
\begin{itemize}
    \item Theorem~\ref{theo:DM_RKHS} presents an ``almost-isometric'' form of the Dvoretzky-Milman theorem. In the analysis of the estimation error of KRR, it suffices to employ an isomorphic variant of the Dvoretzky-Milman theorem. Specifically, by selecting negligibly small constants $\delta$ and $\bar\delta$, we can ensure that the distortion $\tilde\delta$ is also a small constant. In fact, $\tilde\delta<1/2$ suffices for our purpose. In this case, we can speak about isomorphy instead of isometry.
    \item Assumption $N\lesssim \frac{\Tr\left(\Gamma_{k+1:\infty}\right)}{\sqrt{\Tr\left(\Gamma_{k+1:\infty}^2\right)}}$ from Lemma~\ref{lemma:mcrae} needed to obtain isomorphy is much more restrictive compared to condition $N\lesssim \frac{\Tr\left(\Gamma_{k+1:\infty}\right)}{\norm{\Gamma_{k+1:\infty}}_{\text{op}}}$ from Theorem~\ref{theo:DM_RKHS}. 
    Up to \eqref{eq:extra_assumption}, we recover the same condition as in the Gaussian case.
    \item When $\eps>2$ in Assumption~\ref{assumption:DM_L4_L2} and $N\leq\kappa_{DM}d_\lambda^*\left(\Gamma_{k+1:\infty}^{-1/2}B_\cH\right)$, \eqref{eq:extra_assumption} is automatically satisfied because $\Tr\left(\Gamma_{k+1:\infty}^2\right)\leq\norm{\Gamma_{k+1:\infty}}_{\text{op}}\Tr\left(\Gamma_{k+1:\infty}\right)$. Theorem~\ref{theo:DM_RKHS} extends the main result in \cite{guedon_interval_2017} from $L_{4+\eps}-L_2$ equivalence to $L_{2+\eps}-L_2$ equivalence.

    \item The classical version of the Dvoretzky-Milman theorem asserts that with high probability, the intersection of a convex body \(B\) with the image of a random Gaussian matrix, that is, \(B \cap \mathrm{Range}(\bG^\top)\), is almost Euclidean. In Theorem~\ref{theo:DM_RKHS}, \textit{case 1)}  extends this result to the case where \(B\) is an ellipsoid to the scenario of random matrices under weak moment assumptions. However, \textit{case 2)} of Theorem~\ref{theo:DM_RKHS} establishes that the condition number of \(\bX_{\phi,k+1:\infty}\bX_{\phi,k+1:\infty}^\top + \lambda I\) is constant, and scales with \(\lambda\) in the same order. While this does not directly imply an almost Euclidean property for \(B \cap \mathrm{Range}(\bX_{\phi,k+1:\infty}^\top)\), for convenience, we still refer to it as the Dvoretzky-Milman theorem.

    \item In \textit{case 2)} of Theorem~\ref{theo:DM_RKHS}, it is, in fact, unnecessary to compute the lower bound of the minimum (non-zero) singular values of $\bX_{\phi,k+1:\infty}$. Therefore, even the $L_4$-$L_2$ equivalence conditions of the Bai-Yin theorem is not required. Indeed, if allowing for a logarithmic factor of sub-optimality, the scenario of \textit{case 2)} of Theorem~\ref{theo:DM_RKHS} can hold for a more general kernel. This only requires replacing the upper bound in \textit{case 2)} with Proposition~\ref{prop:upper_dvoretzky_infty} and noting that the lower bound in \textit{case 2)} does not depend on the properties of $\phi(X)$, see Proposition~\ref{prop:smooth_large_regularization} in supplementary material.

\end{itemize}


\paragraph{Upper bound of the Dvoretzky-Milman theorem} In this paragraph, we prove that there exists an absolute constant $\nC\label{C_DMU}>0$ such that with high probability, for any $\vlambda\in\bR^N$, under no assumption on $N$,
\begin{equation}\label{eq:upper_dvoretzky_high_probability}
        \norm{\Gamma_{k+1:\infty}^{1/2}\bX_{\phi,k+1:\infty}^\top\vlambda}_\cH \leq\oC{C_DMU}\left(\sqrt{\Tr(\Gamma_{k+1:\infty}^2)}+\sqrt{N}\norm{\Gamma_{k+1:\infty}}_{\text{op}}\right)\norm{\vlambda}_2.
\end{equation}
This will happen to be true, under the following $L_{4+\eps}-L_2$ norm equivalent assumption and an assumption concerning diagonal terms.

\begin{Assumption}\label{assumption:upper_dvoretzky}
    There exist absolute constants $\ngamma\in \left(0,\frac{1}{16}\right)\label{gamma_DMU_L2}$, $\ndelta\geq 0\label{delta_DMU_L2}$, $\eps>0$ and $\kappa'>1$ such that
    \begin{itemize}
        \item
        \begin{align}\label{eq:diagonal_upper_dvoretzky}
            \bP\left( \max_{1\leq i\leq N} \frac{\norm{\Gamma_{k+1:\infty}^{1/2}\phi_{k+1:\infty}(X_i)}_\cH^2}{\Tr\left(\Gamma_{k+1:\infty}^2\right)}  \leq 1+ \odelta{delta_DMU_L2} \right)\geq 1-\ogamma{gamma_DMU_L2},
        \end{align}
        \item for any $f\in \cH_{k+1:\infty}$, $\norm{f}_{L_{4+\eps}}\leq\kappa'\norm{f}_{L_2}$.
    \end{itemize}
\end{Assumption}
The next result is proven in Section~\ref{sec:proof_upper_dvoretzky}. It is similar to \cite[Theorem 4]{tsigler_benign_2023} and \cite{guedon_interval_2017}.
\begin{Proposition}[\cite{tsigler_benign_2023} or \cite{guedon_interval_2017}]\label{prop:upper_dvoretzky}
Suppose Assumption~\ref{assumption:upper_dvoretzky} holds. There exists a constant $\nc\label{c_P_DMU}$ such that \eqref{eq:upper_dvoretzky_high_probability} holds with a probability of at least $1-\frac{\oc{c_P_DMU}}{N^\eps}-\ogamma{gamma_DMU_L2}$.
    The symbol $\bar p_{DMU}$ is defined as $\frac{\oc{c_P_DMU}}{N^\eps}+\ogamma{gamma_DMU_L2}$.
\end{Proposition}


The validity of Assumption~\ref{assumption:upper_dvoretzky} is not universal; therefore, we present an alternative assumption in the supplementary material, as discussed in Section~\ref{sec:omitted}. The proof of Proposition~\ref{prop:upper_dvoretzky} can be found in Section~\ref{sec:proof_upper_dvoretzky}.

\subsection{Restricted Isomorphy Property}\label{sec:RIP}


The Restricted Isomorphy Property characterizes the geometric properties of the RKHS design matrix restricted to $\cH_{1:k}$, that is, $\bX_{\phi,1:k}$. We will see later that this is the part of the space where estimation happens. It identifies the set on which, with high probability, the operator $\bX_{\phi,1:k}:\cH_{1:k}\to\ell_2^N$ forms an isomorphism when $k\lesssim N$ or restricted isomorphism when $k\gtrsim N$. This property was used for linear functionals of sub-Gaussian random vectors in \cite{lecue_geometrical_2022} in the context of benign overfitting in linear regression. In this paper, since we need to study RKHS, we must establish a corresponding version of this restricted isomorphy property.

\paragraph{When $k\lesssim N$} When $k\lesssim N$, the RKHS design operator $\bX_{\phi,1:k}$ behaves like an isomorphy over the entire space $\cH_{1:k}$ under the following assumption.

\begin{Assumption}\label{assumption:RIP}
    There exist absolute constants $\kappa''\geq 1$, $\nc\label{c_RIP}$ depending on $\kappa''$ ($\oc{c_RIP}\leq \frac{9}{1568(\kappa'')^4}$ is sufficient), $0\leq \ngamma\label{gamma_RIP}<1/16$, $\eps>0$, $\ndelta\label{delta_RIP}\geq 0$ such that
        \begin{itemize}
        \item $k\leq \oc{c_RIP}N$.
        \item with probability at least $1-\ogamma{gamma_RIP}$,
        \begin{align}\label{eq:diagonal_RIP}
            \max_{1\leq i\leq N}\norm{\Gamma_{1:k}^{-1/2}\phi_{1:k}(X_i)}_\cH^2\leq \odelta{delta_RIP} k ,
        \end{align}
        \item for any $f\in\cH_{1:k}$, $\norm{f}_{L_{4+\eps}}\leq \kappa''\norm{f}_{L_2}$;
    \end{itemize}
\end{Assumption}

The next result shows the isomorphy property of $\bX_{\phi,1:k}$ on $\cH_{1:k}$ under Assumption~\ref{assumption:RIP}.
\begin{Proposition}\label{prop:IP}
   Under Assumption~\ref{assumption:RIP}, there exist absolute constants $\nc\label{c_P_RIP}$, $\oc{c_RIP_lower}$ and $\oC{C_RIP_upper}$ such that with probability at least $1-\ogamma{gamma_RIP}-\frac{\oc{c_P_RIP} }{N^\eps}- 2\exp(-k)$,
   for all $f_{1:k}\in \cH_{1:k}$,
  \begin{align*}
      \oc{c_RIP_lower}\norm{\Gamma_{1:k}^{1/2}f_{1:k}}_\cH\leq \frac{1}{\sqrt{N}}\norm{\bX_{\phi,1:k}f_{1:k}}_2\leq \oC{C_RIP_upper}\norm{\Gamma_{1:k}^{1/2}f_{1:k}}_\cH,
  \end{align*}where $\oc{c_RIP_lower}$ can be taken as $\frac{1}{2}$ and $\oC{C_RIP_upper}$ can be taken as $\sqrt{2\oC{C_DMU}^2(1+\oc{c_RIP})}$.
  We denote $\bar p_{RIP}$ as $\ogamma{gamma_RIP}+\frac{\oc{c_P_RIP} }{N^\eps} + 2\exp(-k)$ .
\end{Proposition}

\paragraph{When $k$ is not necessarily smaller than $N$} When $k$ is not necessarily smaller than $N$, the design matrix $\bX_{\phi,1:k}$ cannot behave like an isomorphy over the entire space $\cH_{1:k}$ because it has a non-trivial kernel, but it can be an isomorphy restricted to a subset of $\cH_{1:k}$. This set can be taken to be a cone defined below in \eqref{eq:def_cone_RIP}. We refer to this property as the Restricted Isomorphy Property (RIP) as in \cite{lecue_geometrical_2022} in reminiscence to the RIP used in Compressed sensing  \cite{foucart_mathematical_2013}. Please note that we can still use the following proposition to replace Proposition~\ref{prop:IP} when $k \lesssim N$, albeit at the cost of incurring a logarithmic factor.
\begin{Proposition}\label{prop:RIP}
    For any $R>0$, let $\bar\Gamma_{1:k}^{-1/2}=\sum_{j\leq k}\min\left(\frac{1}{R},\frac{1}{\sqrt{\sigma_j}}\right)\varphi_j\otimes \varphi_j$. For some $\nc\label{c_kappa_RIP}$ sufficiently small ($\oc{c_kappa_RIP}<\frac{1}{100\oC{C_Rudelson}^2\oC{C_estimate_gamma_infty}^2}$ is sufficient), let
        \begin{align}\label{eq:def_fixed_point}
    R_N(\oc{c_kappa_RIP}) = \inf\left\{R>0:\, \norm{\max_{i\in[N]} \norm{ \bar\Gamma_{1:k}^{-1/2}\phi(X_i) }_\cH }_{L_\infty} \leq \oc{c_kappa_RIP}\frac{\sqrt{N}}{\log{N}} \right\}.
\end{align}
For any $0<\ndelta\label{delta_P_RIP}<1$, there exist absolute constants $\oc{c_RIP_lower}$, and $\oC{C_RIP_upper}$  depending on $\odelta{delta_P_RIP}$ such that when $R\geq R_N(\oc{c_kappa_RIP})$, then with probability at least $1-\odelta{delta_P_RIP}$, for all $f\in \mathrm{cone}\left(\cC(R)\right)$, where
\begin{align}\label{eq:def_cone_RIP}
    \cC(R) = R^{-1}B_{\cH_{1:k}}\cap \Gamma_{1:k}^{-1/2}S_{\cH_{1:k}} \mbox{ and } \mathrm{cone}\left(\cC(R)\right) = \left\{ f\in\cH_{1:k}:\, R\norm{f}_\cH\leq \norm{f}_{L_2} \right\},
\end{align}we have
\begin{align*}
    \oc{c_RIP_lower}\norm{\Gamma_{1:k}^{1/2}f_{1:k}}_\cH\leq \frac{1}{\sqrt{N}}\norm{\bX_{\phi,1:k}f_{1:k}}_2\leq \oC{C_RIP_upper}\norm{\Gamma_{1:k}^{1/2}f_{1:k}}_\cH.
\end{align*}
\end{Proposition}In contrast to Proposition~\ref{prop:IP}, we only have constant probability deviation in Proposition~\ref{prop:RIP}. We refer to Proposition~\ref{prop:RIP} as the ``Restricted Isomorphic Property under embedding index condition'' because the estimation of the fixed point $R_N(\oc{c_kappa_RIP})$ requires the embedding index condition. An example of the estimate of $R_N(\oc{c_kappa_RIP})$ may be found in Section~\ref{sec:estimate_RIP}. The proof of Proposition~\ref{prop:IP} and Proposition~\ref{prop:RIP} are postponed to Section~\ref{sec:proof_RIP}.

\section{Main Results}\label{sec:polynomial_kernel}

In this section, we present the upper bounds on the estimation error of KRR.

\subsection{Our results}
As remarked in several works \cite{bartlett_deep_2021, tsigler_benign_2023, lecue_geometrical_2022}, there is a fundamental parameter $k$ which is the dimension of the space endowed by the top $k$ eigenvectors of $\Gamma$ where 'estimation' happens whereas, on the orthogonal space, 'absorption of the noise (and even overfitting of the noise when $\lambda=0$)' happens. Our analysis depends on the case where this parameter is smaller or larger than the number of data $N$. 


\paragraph{When $k\leq\oc{c_RIP}N$} In this paragraph, we provide conclusions for the case of $k\lesssim N$. This scenario is precisely what linear regression and Multiple Descents problems are most concerned with. We recall that the definition of $d_\lambda^*(\Gamma_{k+1:\infty}^{-1/2}B_\cH)$ is in Equation \eqref{eq:def_DM_dimension}, the ones of  of $\bar p_{RIP}$, $\bar p_{DM}$, and $\bar p_{DMU}$ are in Proposition \ref{prop:IP}, Theorem \ref{theo:DM_RKHS}, and Proposition \ref{prop:upper_dvoretzky}, respectively. 

\begin{Theorem}\label{theo:upper_KRR}
Suppose Assumptions~\ref{assumption:DM_L4_L2},~\ref{assumption:upper_dvoretzky} and~\ref{assumption:RIP} hold. There then exist absolute constants $\nC\label{C_N_lower}$, $\oc{c_RIP}$, $\oc{c_kappa_DM}$, $\oc{c_P_bX_f_star}$, $\oC{C_noise}$ ($\oC{C_noise}$ depends on $\okappa{kappa_noise}$) and $\nC\label{C_rate_upper}$, such that the following holds. Suppose the noise $\xi$ is independent of $X$ with mean zero and variance $\sigma_\xi^2$. We assume that for some $\okappa{kappa_noise}>0$ and $r>4$,  $\norm{\xi}_{L_r}\leq \okappa{kappa_noise}\sigma_\xi$. Let $\lambda\geq 0$.
We assume that there exists $k\in\bN$ so that $\oC{C_N_lower}\leq N \leq \oc{c_kappa_DM}\kappa_{DM}d_\lambda^*\left(\Gamma_{k+1:\infty}^{-1/2}B_\cH\right)$, $\sigma_1 N > \kappa_{DM}(4\lambda + \Tr\left(\Gamma_{k+1:\infty}\right))$ and $k\leq\oc{c_RIP}N$. Let $\bar p_{\xi}$ be some probability deviation strictly less than $1$ (defined later in \eqref{eq:def_bar_p_xi}). Then with probability at least
\begin{align*}
    1 - \bar p_{RIP} - \bar p_{DM} - \bar p_{DMU} - \frac{ \oc{c_P_bX_f_star} }{N}  - \bar p_{\xi} - \left(\frac{\oC{C_noise}\Tr(\Gamma_{k+1:\infty})}{|J_1|\Tr(\Gamma_{k+1:\infty}) + N \left(\sum_{j\in J_2}\sigma_j\right)}\right)^{\frac{r}{4}},
\end{align*}
we have, for $r_{\lambda,k}^*$ is defined in \eqref{eq:def_rate_informal}, that
\begin{align*}
    \norm{\hat f_\lambda - f^*}_{L_2}\leq \oC{C_rate_upper}r_{\lambda,k}^*.
\end{align*}
\end{Theorem}

\paragraph{When $k$ is not necessarily smaller than $\oc{c_RIP}N$} In this paragraph, we provide the informal version of our conclusion for the case of $k$ is not necessarily smaller than $\oc{c_RIP}N$ (the formal version may be found in Section~\ref{sec:omitted} in supplementary material). In fact, for the context of the minimum $\norm{\cdot}_2$ norm interpolant estimator in linear regression, the authors of \cite{lecue_geometrical_2022} have already conducted research on this scenario. They have demonstrated that when the design vector is symmetric, the optimal value of $k$—the one that minimizes the estimation error among all possible $k$—falls precisely within the range of $k\lesssim N$ if one wants benign overfitting to happen. 
However, in the case of KRR, there is no lower bound indicating that the optimal $k$ must satisfy $k \lesssim N$. We therefore present the following theorem.

\begin{Theorembis}{theo:main_upper_k>N}
    Under certain conditions, for any \(\lambda \geq 0\) and \(k\) satisfying the Dvoretzky-Milman condition, with high probability, we have $\|\hat f_\lambda - f^*\|_{L_2}\lesssim r_{\lambda,k}^*$ where $r_{\lambda,k}^*$ is defined in \eqref{eq:def_rate_informal}.
\end{Theorembis}

\begin{Remark}[Misspecified model]\label{remark:model_misspecified}
    In practical applications, we often encounter cases where $f^*\notin\cH$, as exemplified by our Proposition~\ref{prop:feature_learning}. When $f^* \notin \cH$, we define $f^{**} = \arg \min (\|f^* - f\|_{L_2(\mu)}: f \in \cH)$, that is, the orthogonal projection of $f^*$ onto $\cH$ in the $L_2(\mu)$ inner product (we could also choose other oracles as we have the liberty to choose the oracle). We replace the target function in Theorem~\ref{theo:upper_KRR} and Theorem~\ref{theo:main_upper_k>N} with $f^{**}$ and replace the noise $\bxi$ with $\veps = \vr + \bxi$, where $\vr = (f^*(X_i) - f^{**}(X_i))_{i=1}^N$. This new noise $\veps$ is dependent on the kernel design matrix $\bX_\phi$. See Proposition~\ref{prop:noise_concentration_dependent}, Proposition~\ref{prop:dependent_noise_result} and the subsequent discussion of this property in supplementary material and Proposition~\ref{prop:feature_learning} for an example when $f^*\notin\cH$. This conclusion is not contradictory to the counterexamples provided in \cite{chinot_robustness_2022} and \cite{shamir_implicit_2022}, as indicated in Remark~\ref{remark:no_contradiction} in supplementary material.
    
\end{Remark}

\begin{Remark}[Uniform results in $\lambda$]
It is possible to have results equivalent to Theorem~\ref{theo:upper_KRR} and Theorem~\ref{theo:main_upper_k>N} uniform in the regularization parameter $\lambda$. This type of results is particularly useful when one wants to use a data-dependent regularization parameter in order to achieve optimal and adaptive results. It is for instance the case, when $\lambda$ is chosen according to the Lepski's method \cite{blanchard_lepskii_2019}. In that case, our results hold with the same probability and  convergence rates however, for instance in Theorem~\ref{theo:upper_KRR},  we just need to assume that  $ N \lesssim d_{\lambda_0}^*\left(\Gamma_{k+1:\infty}^{-1/2}B_\cH\right)$ for some $\lambda_0\geq 0$ and then the result of  Theorem~\ref{theo:upper_KRR} holds uniformly for all $\lambda\geq \lambda_0$. This may be particularly useful when $\lambda_0=0$.
\end{Remark}

\subsection{Comparison with previous results}

In this section, we primarily select conclusions from \cite{mcrae_harmless_2022}, and \cite{mourtada_elementary_2022}. Due to length constraints, we have placed the comparison with the results of \cite[Section 4.3.3, Theorem 4.13]{bartlett_deep_2021} in Section~\ref{sec:omitted} of the supplementary material.

\begin{Theorem}[\cite{mcrae_harmless_2022}]\label{theo:mcrae_et_al}
    Let $k\log{k}\lesssim N$. Suppose that some absolute constants $\alpha_1\geq \alpha_2>0$ exist such that $\alpha_1 I_N\succeq  \lambda I_N + \bX_{\phi,k+1:\infty}\bX_{\phi,k+1:\infty}^\top \succeq \alpha_2 I_N$ and $$\frac{\alpha_1-\alpha_2}{\alpha_1 + \alpha_2}+ \frac{2}{N}\norm{\bX_{\phi,1:k}^\top\bX_{\phi,1:k} - NI}_{\text{op},L_2\to L_2}\leq c$$ for some $c<1$ where $I$ is the identical operator on $\cH_{1:k}$. Then with constant probability
    \begin{align*}
        &\norm{\hat f_\lambda - f^*}_{L_2(\mu)}^2\lesssim \sigma_\xi^2\left(1+\frac{\alpha_1}{\alpha_2}\right)^2 \left(\frac{k}{N} + \frac{4N\Tr\left(\Gamma_{k+1:\infty}^2\right)}{\left(\alpha_1 + \alpha_2\right)^2}\right) \\
        &+ \left(\min\left\{\sqrt{\sigma_1},\frac{1}{1-c}\frac{2\alpha_1\alpha_2}{N\sqrt{\sigma_k}(\alpha_1+\alpha_2)},\frac{1}{1-c}\sqrt{\frac{2\alpha_1\alpha_2}{N(\alpha_1+\alpha_2)}}\right\}\right)^2\left( 
1+ \sqrt{\frac{N\sigma_{k+1}(\alpha_1+\alpha_2)}{2\alpha_1\alpha_2}} \right)^2\norm{f^*}_\cH^2.
    \end{align*}
\end{Theorem}

In \cite{mcrae_harmless_2022}, the authors use Lemma~\ref{lemma:mcrae} and the matrix Bernstein inequality to verify the assumptions under the conditions $k\log{k}\leq c N$, $N\sqrt{\Tr\left(\Gamma_{k+1:\infty}^2\right)}\leq c\Tr\left(\Gamma_{k+1:\infty}\right)$ for some small constant $c$ when $\cH_{1:k}$ is spanned by some bounded orthonormal system. We provide some comments on this result.

\begin{itemize}
    \item Theorem~\ref{theo:mcrae_et_al} applies only when the conditions $N\sqrt{\Tr\left(\Gamma_{k+1:\infty}^2\right)}\leq c\Tr\left(\Gamma_{k+1:\infty}\right)$ and $k \log k\lesssim N$ are valid. In contrast, thanks to Proposition~\ref{prop:IP} and Theorem~\ref{theo:DM_RKHS}, our results hold in the optimal Gaussian regime i.e. $k\lesssim N$ and $N\lesssim d_\lambda^*(\Gamma_{k+1:\infty}^{-1/2}B_\cH)$ and thus for a wider range of kernels. Please refer to the discussion after Lemma~\ref{lemma:mcrae}.
    
    \item Consequently, in order to satisfy the inequality $N\sqrt{\Tr\left(\Gamma_{k+1:\infty}^2\right)}\leq c\Tr\left(\Gamma_{k+1:\infty}\right)$, it is necessary to select a sufficiently large value for $k$. However, this choice of $k$ may (1) contradict the assumption $k\log k\lesssim N$ from Theorem~\ref{theo:mcrae_et_al} and (2) the term $k/N$ appearing in the rate in Theorem~\ref{theo:mcrae_et_al} may be sub-optimal. Due to this limitation, the findings presented by Theorem~\ref{theo:mcrae_et_al} lack the necessary precision to observe numerous phenomena that will be further examined in this study.
    
    \item Despite the utilization of Proposition~\ref{prop:IP} and Theorem~\ref{theo:DM_RKHS} in the context of Theorem~\ref{theo:mcrae_et_al}, it is evident that the other terms in the upper bound are not optimal. As demonstrated by Theorem~\ref{theo:upper_KRR} below, the term $\norm{f^*}_\cH$  can be replaced by the smaller ones $\norm{\Gamma_{k+1:\infty}^{1/2}f_{k+1:\infty}^*}_\cH$ and $\norm{\tilde\Gamma_{1,\mathrm{thre}}^{-1/2}f_{1:k}^*}_{\cH}(\lambda + \Tr\left(\Gamma_{k+1:\infty}\right))/N$. 
    \item The upper bound stated in Theorem~\ref{theo:mcrae_et_al} is dependent on the dimension $k$. The selection of the parameter $k$ is crucial as it influences the decomposition of $\cH$ into two orthogonal subspaces, namely $\cH_{1:k}$ and $\cH_{k+1:\infty}$. The former is utilized for estimating the (part of) target function $f_{1:k}^*$, while the latter is employed for accommodating the noise, see Proposition~\ref{prop:decomposition_hat_f} below for a clear statement of the 'feature space' (i.e. the RKHS here) decomposition. In a general sense, when the value of $k$ is increased, the estimator $\hat f_{\lambda}$ takes into account, a larger portion of the space as containing valuable information of $f^*$, while allocating less space for the absorption of noise. The assumption made in Theorem~\ref{theo:mcrae_et_al} is that $k\log{k}\lesssim N$, which imposes a restriction. In Section~\ref{sec:KRR_smooth}, we will explore scenarios when it is crucial for $k$ to be of at least the same order as $N$.
    \item The probability deviation appearing in Theorem~\ref{theo:mcrae_et_al} is of a constant level, whereas  Theorem~\ref{theo:upper_KRR} holds with large probability tending to $1$ when $N$ tends to infinity.
\end{itemize}

We also notice the results presented in \cite{mourtada_elementary_2022}. These results have been utilized by \cite{ba_high-dimensional_2022,bietti_learning_2022} to analyze the estimation error of KRR defined by data-dependent kernels. We present a conclusion here, demonstrating that by selecting an appropriate value of $k$, we can recover the bounds established in \cite{mourtada_elementary_2022,bach_learning_2024}. This implies that this classical result corresponds to a specific choice of $k$ in our bounds even though this choice may not be the optimal one. Indeed, many examples discussed below indicate the existence of values of $k$ for which our bounds are sharper than those associated with the classical bounds. This illustrates the benefits of our bounds.
\begin{Proposition}\label{prop:jaouad}
    Given $\lambda>0$, let $k = \min(j\in\bN:\, \lambda\gtrsim N\sigma_{j+1})$. Let $\mathrm{Proj}_{\cH}f^*$ be the orthogonal projection of $f^*$ onto the closure of $\cH$. Then
    \begin{align*}
        &(r_{\lambda,k}^*)^2 + \norm{f^* - \mathrm{Proj}_{\cH}f^* }_{L_2(\mu)}^2 \lesssim \\
        &\left( 1+ \frac{\Tr(\Gamma)}{\lambda} \right)^2\inf\left( \norm{f - f^*}_{L_2(\mu)}^2 + \frac{\lambda}{N}\norm{f}_\cH^2:\, f\in\cH\right) + \frac{\sigma_\xi^2}{N}\Tr\left( \Gamma\left(\Gamma + \frac{\lambda}{N} I\right)^{-1} \right).
    \end{align*}
\end{Proposition}The proof of Proposition~\ref{prop:jaouad} may be found in Section~\ref{sec:proof_prop_jaouad} in supplementary material. The bound on the right side of Proposition~\ref{prop:jaouad} is the upper bound for the generalization error of KRR provided by \cite{mourtada_elementary_2022,bach_learning_2024}. Here, \(\|f^* - \mathrm{Proj}_\mathcal{H} f^*\|_{L_2(\mu)}^2\) represents the approximation error. In Proposition~\ref{prop:dependent_noise_result} in supplementary material, we obtain an upper bound for the generalization error under a misspecified model, which corresponds to the left side of Proposition~\ref{prop:jaouad}. This indicates that our bound is never worse than the bound provided by \cite{mourtada_elementary_2022}. However, it's worth noting that the results from \cite{mourtada_elementary_2022} are not sharp in many cases. For instance, their variance term is only valid when $\lambda\gtrsim 1$ (hence it cannot cover interpolant estimators, where $\lambda = 0$), and their bias term is also not sharp compared to ours, as we will demonstrate in Section~\ref{sec:KRR_smooth}. Moreover, \cite{mourtada_elementary_2022} can only choose the oracle as the projection of \(f^*\) onto the closure of \(\mathcal{H}\), while our general bound allows for the selection of any oracle.

After the completion of the material preparation for this paper, we became aware of the work by \cite{barzilai_generalization_2024}. We emphasize the differences between \cite{barzilai_generalization_2024} and the current paper:
\begin{enumerate}
    \item Theorem 1 in \cite{barzilai_generalization_2024} is used to control a lower bound on the smallest eigenvalue of the kernel matrix. However, their lower bound still requires $N^2\Tr(\Gamma_{k+1:\infty}^2)<\Tr(\Gamma_{k+1:\infty})$, as seen in Equation 12 of \cite{barzilai_generalization_2024}. This poses a similar issue to Lemma~\ref{lemma:mcrae}, whereas our Theorem~\ref{theo:DM_RKHS} does not require this assumption.

    \item The lower bound in Theorem 1 of \cite{barzilai_generalization_2024} contains an additional logarithmic factor (their $\beta_k k \log(k)<n$ and $\log(k+1)$ terms), which leads to the need for an additional logarithmic factor in their statistical conclusion \cite[Theorem 2]{barzilai_generalization_2024}, similar to \cite{mcrae_harmless_2022}. In contrast, our Theorem~\ref{theo:upper_KRR} and Theorem~\ref{theo:main_upper_k>N} address both cases of $k\lesssim N$ and $k\gtrsim N$ separately, thus not imposing requirements on the choice of $k$. It's worth noting that the choice of $k$ determines the decomposition of $\hat{f}_\lambda$, hence the unnecessary logarithmic factor in the final bound due to larger $k$ in \cite[Theorem 2]{barzilai_generalization_2024}.

    \item Definition 2 in \cite{barzilai_generalization_2024} requires almost sure upper bounds on $\norm{\phi_{k+1:\infty}(X)}_\cH^2$, $\norm{\Gamma_{k+1:\infty}^{1/2}\phi_{k+1:\infty}(X)}_\cH^2$, and $\norm{\Gamma_{1:k}^{-1/2}\phi_{1:k}(X)}_\cH^2$. Our Assumption~\ref{assumption:DM_L4_L2} only requires upper or lower bounds with high probability. This renders the conclusions of \cite{barzilai_generalization_2024} inapplicable to linear regression with unbounded probability measures (for example, Gaussian/sub-Gaussian or heavy-tailed distributions). More importantly, this also renders their findings inapplicable to data-dependent conjugate kernels, as discussed in Section~\ref{sec:conjugate_kernel} and Lemma~\ref{lemma:wang_deformed}. In Section~\ref{sec:conjugate_kernel}, we will demonstrate the significance of such kernels for deep learning theory.



\end{enumerate}

\section{Applications}
In this section, we apply our main results (Theorems~\ref{theo:upper_KRR} and \ref{theo:main_upper_k>N}) to the study of four phenomenons that have been observed in empirical and theoretical studies on neural networks. Our aim is to show that sharp upper bounds on the convergence properties of KRR are relevant to the understanding of phenomenons in deep learning.

\subsection{Multiple Descents of the minimum RKHS-norm interpolant estimator}\label{sec:multiple_descent}

In this section, we apply our Theorem~\ref{theo:upper_KRR} to the inner product kernel, obtaining the phenomenon of Multiple Descents, a concept that will be explained later in this text.
Imagine that one is training an estimator, and when the number of samples increases, the test error does not decrease and, in fact, may even increase. How can we explain this phenomenon, and should we stop adding more samples? The Multiple Descents phenomenon is closely related to this scenario.

The phenomenon of Multiple Descents, as far as current knowledge indicates, was initially identified by the authors of \cite{liang_multiple_2020,ghorbani_linearized_2021}. The observation is made that the estimation error $\norm{\hat f_0 -f^*}_{L_2}$, where $\hat f_0$ is KRR for $\lambda=0$, that is, it is the minimum $\norm{\cdot}_\cH$-norm interpolant estimator, exhibits numerous decreases as the value of $N$ grows. It is important to note that in the proof of the multiple Descents phenomenon, two setups have been considered: one where $d$ is fixed and $d\ll N$ and the other where $d\to\infty$, $d\ll N$. Our results pertain to both cases where $d$ is fixed and where $d\to\infty$. In fact, \cite[Lemma 2.1]{donhauser_how_2021} demonstrated that for a fairly broad class of $\cH$ (including Laplace kernel and exponential inner product kernel), as $d\to\infty$, $\cH$ does not even contain relatively simple functions as defined in the following Lemma:
\begin{Lemma}[\cite{donhauser_how_2021}]\label{lemma:DWY21}
    Let $K:(\vx,\vy)\in\bR^d\times\bR^d\mapsto \prod_{j=1}^d q(x_j, y_j)$ where $q$ is some kernel function which may change with $d$. Assume $\int q(x,x)d\mu(x)\leq 1$. Then for any $f:\vx\in\bR^d\mapsto\prod_{j=1}^m f_j(x_j)\in\bR$ for some fixed $m\in\bN$ (in particular, $m$ can be taken to be $1$) and each $f_j$ is a function in the RKHS generated by $q$, we have $\norm{f}_\cH\to\infty$ as $d\to\infty$.
\end{Lemma} Let us once again emphasize the motivation for fixed $d$: the variation in $d$ makes it challenging to determine which $f^*$, $\hat{f}_0$ is approximating (each $d$ defines a distinct $f^*$, see Theorem~\ref{theo:MMM_21} below), and $\cH$ itself becomes very small as $d\to\infty$, as indicated by Lemma~\ref{lemma:DWY21} -- therefore restricting the possible choices of $f^*$ as $d$ increases. The specific problem of Multiple Descents pertains to the properties of the inner product kernel, that is, there exists a scalar function $h:\bR\to\bR$ such that for any $\vx,\vy\in\bR^d$, we have
\begin{align}\label{eq:polynomial_inner_product_kernel}
    K(\vx,\vy)=h\left(\frac{\left\langle \vx, \vy\right\rangle}{d}\right).
\end{align}
The phenomenon of Multiple Descents holds significance in both practical and theoretical contexts:
\begin{enumerate}
    \item In practical terms, the phenomenon of Multiple Descents informs us that there is no need to be concerned about an increase in estimation error when the number of samples is increased. This is because, once we surpass this increase, our estimator is capable of learning a more accurate approximation of the signal of higher degree. In fact, an increase in $N$ will lead to a reduction of the dimension of the space available for accommodating noise. This, in turn, results in an increase in the variance in the estimation error, which is the cause of the rise in estimation error.
    \item On the other hand, in theoretical terms, the Multiple Descents phenomenon sheds light on the ``degree of non-parametricness'' of the estimator $\hat f_0$. This refers to the effectiveness of $\hat f_0$ in learning $f^*$ when the number of parameters is infinite.
\end{enumerate}

In this section, by applying our general bounds from Theorem~\ref{theo:upper_KRR}, we extend \cite{liang_multiple_2020} (see Theorem~\ref{theo:LRZ20} in supplementary material) to the optimal interval of $d^\iota\lesssim N\lesssim d^{\iota+1}$ within the non-asymptotic regime. Moreover, we prove that Multiple Descents can still occur under the condition of sub-Gaussian random vectors. This extends the conclusions of \cite{ghorbani_linearized_2021} to the sub-Gaussian setting.

\begin{Propositionbis}{prop:multiple_descent}
    If \( h \) is a polynomial of degree \( L \), under certain conditions (particularly when \( X \) is a sub-Gaussian random vector), for any \( 0 \leq \iota \leq L-1 \), with high probability, we have $\left| \|\hat{f}_0 - f^*\|_{L_2} - \|\Gamma_{>\iota}^{1/2} f_{>\iota}^*\|_{\mathcal{H}} \right| \lesssim r_{0,k}^*$, where \( r_{0,k}^* \) is equivalent to $\max\left\{ \sigma_\xi \sqrt{\frac{d^\iota}{N}}, \frac{1}{N} \|\Gamma_{\leq \iota}^{-1/2} f_{\leq \iota}^*\|_{\mathcal{H}}, \|\Gamma_{>\iota}^{1/2} f_{>\iota}^*\|_{\mathcal{H}}, \sigma_\xi \sqrt{\frac{N}{d^{\iota+1}}} \right\}$. Moreover, as \( N, d \to \infty \) with \( \omega(d^\iota) \leq N \leq o(d^{\iota+1}) \), we have $\left| \|\hat{f}_0 - f^*\|_{L_2(\mu)} - \|\Gamma_{>\iota}^{1/2} f_{>\iota}^*\|_{\mathcal{H}} \right| = o_{d,\mathbb{P}}(1)(\sigma_\xi + \|f^*\|_{\mathcal{H}}).$
\end{Propositionbis}
The formal version of Proposition~\ref{prop:multiple_descent} can be found in Section~\ref{sec:omitted} of the supplementary material, and its proof is presented in Section~\ref{sec:proof_multiple_descent} and Section~\ref{sec:multiple_descent_asymptotic}. Actually, our bound also holds even with $o_{d,\mathbb{P}}(1)\sigma_\xi + O_{d,\mathbb{P}}(1) \|\Gamma_{>\iota}^{1/2}f^*\|_{\mathcal{H}}$ in place of $o_{d,\mathbb{P}}(1)(\sigma_\xi + \|f^*\|_{\mathcal{H}})$.

To establish Proposition~\ref{prop:multiple_descent}, we have developed a new Hanson-Wright-type concentration inequality, as presented in Theorem~\ref{theo:polynomial_HS} in supplementary material. Below is an informal version of it. This inequality may be of independent interest. For instance, Theorem~\ref{theo:polynomial_HS} verifies the assumptions in \cite[Lemma 1]{mcrae_harmless_2022} and \cite[Assumption 1(d)]{mei_generalization_2022}.
\begin{Theorembis}{theo:polynomial_HS}
    Suppose $X$ is a sub-Gaussian random vector with i.i.d. mean-zero sub-Gaussian coordinates $(x_i)_{i=1}^d$. Suppose $h$ defined in \eqref{eq:polynomial_inner_product_kernel} is a degree $L$ polynomial.
    There exists an absolute constant $\oC{C_non_linear_HS}$ depending on $L$ and $\norm{x_1}_{\psi_2}$, such that for any $t>0$, we have: with high probability (which is exponential on $-t^{1/L}$), we have $\left| K(X,X)- \bE K(X,X) \right|< t$.
\end{Theorembis}To our knowledge, Theorem~\ref{theo:polynomial_HS} is the first result proving the concentration of kernel functions for sub-Gaussian random vectors with inner product kernels in the non-asymptotic regime. In fact, for uniform distributions on spheres, this is trivial (see Section~\ref{sec:diagonal_terms_multiple_descent} of the supplementary material), as the diagonal concentration in this case is equivalent to the diagonal concentration of translation-invariant kernels. This poses a challenge that restricts the applicability of a series of works by \cite{ghorbani_linearized_2021,mei_generalization_2022} to general probability measures. In the asymptotic case, such as in \cite{el_karoui_spectrum_2010}, verifying diagonal concentration is also trivial (see Section~\ref{sec:linearization} of the supplementary material) and only requires the use of Lagrange's theorem. However, in the non-asymptotic case, the situation is more complex, which is why we present Theorem~\ref{theo:polynomial_HS}.

Several remarks on Proposition~\ref{prop:multiple_descent} are presented subsequently. Due to length constraints, we have included the results of \cite{liang_multiple_2020} and \cite{ghorbani_linearized_2021} in Section~\ref{sec:omitted} of the supplementary material, see Theorem~\ref{theo:LRZ20} and Theorem~\ref{theo:MMM_21} there.
\begin{enumerate}
\item In the high-dimensional asymptotic limit,  Proposition~\ref{prop:multiple_descent} extends the results of \cite{ghorbani_linearized_2021} to sub-Gaussian designs. To our knowledge, this is the first such extension in the literature, addressing a well-known challenging problem. This establishes the universality of Multiple Descents under sub-Gaussian designs and in the high-dimensional asymptotic limit. Moreover, the interval in which Multiple Descents occurs is characterized by \( \omega_d(d^\iota) \leq N \leq o_d(d^{\iota+1}) \), which improves upon the results of \cite{liang_multiple_2020} and \cite{ghorbani_linearized_2021}. Furthermore, we do not require the assumptions \( \mathbb{E} h_d(X) = \mathbb{E} f_d^*(X) = 0 \) as in \cite{ghorbani_linearized_2021}.
\item When $L=2$, Proposition~\ref{prop:multiple_descent} provides evidence that the upper bound of the absolute value of $\|\hat f_0 - f^*\|_{L_2} - \|\Gamma_{>1}^{1/2}f_{>1}^*\|_{\cH}$ exhibits a descent when $d \lesssim N\lesssim d^2 $. When $d,N\to\infty$ with $\omega(d)\leq N\leq o(d^2)$, Proposition~\ref{prop:multiple_descent} provides sharp asymptotics, not just an upper bound, and this descent being referred to as the second descent within the context of the ``double descents'' phenomenon, as demonstrated in \cite{mei_generalization_2020}. Nevertheless, it is important to highlight the superiority of Proposition~\ref{prop:multiple_descent} in comparison to the findings presented in \cite{mei_generalization_2020}:
\begin{enumerate}
    \item  The underlying assumption in our analysis is that of a sub-Gaussian design. However, it is worth noting that in the work by \cite{mei_generalization_2020}, they need the random variable $X$ to be uniformly distributed over $\sqrt{d}S_2^{d-1}$, which is more restrictive and unrealistic.
    \item 
    In \cite{mei_generalization_2020}, it is demonstrated that KRR defined by the ``random feature kernel'' exhibits the  double descents phenomenon, whereas what we establish here is the double descents phenomenon for KRR defined by the polynomial kernel. We emphasize that our approach is also applicable to the random feature kernel (even data-dependent kernels), as seen in Proposition~\ref{prop:conjugate_kernel}, but we leave this for future work.
\end{enumerate}
    \item \cite{ghorbani_linearized_2021,mei_generalization_2022} proved that when $X$ is uniformly distributed over $\sqrt{d}S_2^{d-1}$ and under some assumptions on $h$, $\|\hat f_0 - f_d^*\|_{L_2} - \norm{P_{>\iota}f_d^*}_{L_2}$ is infinitesimal compared to $\norm{f_d^*}_{L_2}+\sigma_\xi$ when $\omega_d\left(d^\iota\log{d}\right)\leq N\leq O_d\left(d^{\iota+1-\delta_0}\right)$ when $N,d\to\infty$. This result is recently improved in the asymptotic regime. In fact, \cite{misiakiewicz_spectrum_2022} proved that in asymptotic regime, that is, when both $N,d\to\infty$ and when $X$ is uniformly distributed over $\sqrt{d}S_2^{d-1}$ and under some assumptions on $h$ and when $f_d^*$ satisfies a kind of randomness property, the Multiple Descents happens when $d^\iota\lesssim N\lesssim d^{\iota+1}$ for each $0\leq\iota\leq L$. We place significant emphasis on the superiority of our non-asymptotic results, particularly in the context of the ``fixed $d$'' situation, as it allows us to avoid the ambiguity that arises from altering $f^*$ and $\cH$ as $d$ grows, as indicated by Lemma~\ref{lemma:DWY21} and the follow-up discussion.

    \item We unify the two theories developed by \cite{ghorbani_linearized_2021,mei_generalization_2020,mei_generalization_2022,misiakiewicz_spectrum_2022} and \cite{liang_multiple_2020}, whose proof methods are very different. Proposition~\ref{prop:multiple_descent} indicates that they are in fact of the same nature. This is believed to be a challenging task \cite{donhauser_how_2021}. This proves that the geometrical viewpoint underlying our approach captures the phenomenon of Multiple Descents appearing in different scenarii that needed before specific tools.  

    \item We have only provided results for \( L < \infty \), while in the case where \( h \) is a smooth function, it is necessary to investigate how to approximate it using finite-degree polynomials. Specifically, if we denote \(\bK^{[\leq L]}_N = (K^{[\leq L]}(X_i,X_j))_{1 \leq i,j \leq N}\) as the kernel matrix defined by the kernel function \( K^{[\leq L]}(\vx, \vy) = h^{[\leq L]}(\langle \vx, \vy \rangle / d) \), where \( h^{[\leq L]} \) is a polynomial of degree \( L \), then to prove Multiple Descents, we need to examine whether \(\|\bK_N - \bK_N^{[\leq L]}\|_{\text{op}} \to 0\) as $L,d,N\to\infty$, where \(\bK_N\) is the kernel matrix defined by \( K \). Since this paper primarily focuses on non-asymptotic theory, this issue is beyond the scope of our current work.
\end{enumerate}

\subsection{Estimation error of KRR with smooth kernel in the non-asymptotic regime}\label{sec:KRR_smooth}


In Section~\ref{sec:multiple_descent}, we considered the polynomial kernels with finite degree, that is, when there exists $h:\bR\to\bR$ being a finite-degree polynomial function such that for any $\vx,\vy\in\bR^d$, $K(\vx,\vy)=h\left(\left\langle \vx,\vy\right\rangle/d\right)$. One may question whether our general bounds are limited to these specific kernels. In this section, we will demonstrate that this is not the case.
In this section, we study the upper bound for the estimation error of KRR under the assumption that the kernel function $K$ is smooth in the non-asymptotic regime. In \cite{bietti_inductive_2019,haas_mind_2023}, it is proven that when $X$ is uniform over the Euclidean sphere, the eigenvalues of NTK have a power decay. Motivated by this observation, in this section, we grant the following assumption:

\begin{Assumption}\label{assumption:smooth}
    \begin{itemize}
        \item $X\in\bR^d$ is distributed uniformly over $\sqrt{d}S_2^{d-1}$, the noise $\xi$ is distributed according to a centered Gaussian variable with variance $\sigma_\xi^2$ and is independent of $X$'s.
        \item $h\in C^\infty(\bR)$, $K(\vx,\vy)=h(\left\langle \vx,\vy\right\rangle/d)$ for any $\vx,\vy\in\bR^d$ and $\norm{K}_{L_\infty(\mu)}<\infty$. Moreover, $h(t)=\sum_{i=0}^\infty \alpha_i t^i$ with $\alpha_i>0$ for any $i$.
        \item For any $j\in\bN_+$, $\sigma_j\sim j^{-\alpha}$ for some $\alpha>1$\footnote{We allow \(\alpha\) to depend on \(d\), for example, \(\alpha = \frac{2r}{d}\) when \(\mathcal{H}\) is an \(r\)-th order Sobolev space (where \(r > d/2\)), as discussed in \cite[Chapter 7]{bach_learning_2024}.}.
    \end{itemize}
\end{Assumption}

We further grant the source condition and the embedding index condition. The source condition is often assumed in the context of inverse problems and in the study of RKHS from the perspective of classical integral operators, and is often verified in application, for example, \cite[Section 4.2]{pillaud-vivien_kernelized_2023}. It is employed to provide insights into the ``smoothness'' of the unknown function $f^*$, as seen in, for instance, \cite[Equation 11]{bauer_regularization_2007},  \cite[section 2]{blanchard_kernel_2016} and \cite{zhang_optimality_2023}. We assume the following assumption holds: 
\begin{Assumption}\label{assumption:source_condition_2}
    There exist $s\geq 1$ and a sequence $(a_j)_{j\in\bN}\in\ell_2$ such that $f^*=\sum_{j=1}^\infty a_j\sigma_j^{1/2}f_j$ and $\sum_{j=1}^\infty \sigma_j^{1-s}a_j^2<\infty$.
\end{Assumption}
The relationship between Assumption~\ref{assumption:smooth} and the source condition in the language of interpolation spaces can be found in Section~\ref{sec:omitted}.

The embedding index condition in the study of RKHS and spectral algorithms, to the best of our knowledge, originated from \cite{mendelson_regularization_2010} and was later refined by \cite{steinwart_optimal_2009,fischer_sobolev_2020} and others. Specifically, we make the assumption that there exist absolute constants $\oC{C_embedding}$ and $0<\otheta{theta_embedding}<1$ such that for any $f\in\cH$, we have
\begin{align}\label{eq:def_embedding_index}
    \norm{f}_{L_\infty}\leq\oC{C_embedding}\norm{f}_\cH^{\otheta{theta_embedding}}\norm{f}_{L_2}^{1-\otheta{theta_embedding}}.
\end{align}
Under Assumption~\ref{assumption:smooth}, it is proven that $\otheta{theta_embedding}$ can be taken to be $1/\alpha$, \cite{mendelson_regularization_2010}.

\begin{Proposition}\label{prop:KRR_smooth_applied}
    Suppose that Assumption~\ref{assumption:smooth}, Assumption~\ref{assumption:source_condition_2} and \eqref{eq:def_embedding_index} hold. 
If $(s\wedge 2)\alpha>1$, then for $k = N^{\frac{1}{1+(s\wedge 2)\alpha}}$ and $\lambda\sim N\sigma_k$, with constant probability,
\begin{align}\label{eq:result_power_decay}
\norm{\hat f_\lambda  - f^*}_{L_2}\lesssim_{\sigma_\xi} \log^3(N) N^{-\frac{\alpha (s\wedge 2)}{2(1+\alpha (s\wedge 2))}}.    
\end{align}
\end{Proposition}
Proposition~\ref{prop:KRR_smooth_applied} is an unpublished argument by Zhifan Li \cite{li_personal_nodate}. The authors would like to express their gratitude to Zhifan Li for allowing them to cite his conclusion. The proof of Proposition~\ref{prop:KRR_smooth_applied} can be found in Section~\ref{sec:proof_KRR_smooth_applied}.

Some comments are in order.

\begin{itemize}
    \item Up to logarithmic factors, Proposition~\ref{prop:KRR_smooth_applied} achieves the minimax optimal rate established in \cite{caponnetto_optimal_2007} under the source condition and embedding index conditions. \cite{mourtada_elementary_2022} also provides a general upper bound on the estimation error of KRR. However, their results fall short of achieving the minimax optimal rate, as discussed in \cite[Section 3]{mourtada_elementary_2022}. Therefore, our general bounds are more accurate compared to that of \cite{mourtada_elementary_2022}.
    \item As $\alpha$ increases, $\cH$ contains increasingly smoother functions. Similarly, as $s$ increases, the smoothness of the optimal function $f^*$ also increases, resulting in a faster convergence rate. When $s \geq 2$, KRR is unable to further capture the smoothness of $f^*$, a phenomenon referred to as the saturation effect, \cite{bauer_regularization_2007}. Our geometric perspective further indicates that the saturation effect occurs within the space $\cH_{1:k}$.
    \item Additionally, we highlight that, due to the geometric nature of our general bounds, we can clearly demonstrate that for any sufficiently large $N$, $\hat f_\lambda$ captures a degree $\lfloor \frac{1}{1+(s\wedge 2)\alpha} \frac{\log{N}}{\log{d}}\rfloor$-polynomial approximation of $f^*$ while treating the higher-order terms as noise. \cite{donhauser_how_2021} proved that for a more general rotational invariant kernel, $f^*$ can be approximated by a polynomial of degree $\lfloor 2\frac{\log{N}}{\log{d}}\rfloor$, while \cite{ghorbani_linearized_2021} provided a result for degree $\lfloor \frac{\log{N}}{\log{d}}\rfloor$. Their results are asymptotic. In contrast, our conclusion is not only non-asymptotic, but also requires lower-degree polynomials to approximate $f^*$. Our results are not contradictory to \cite{ghorbani_linearized_2021} because we assume eigenvalue decay and the source condition of $f^*$. This implies that $f^*$ is smoother under these conditions, hence it can be approximated by lower-degree polynomials. Moreover, to our knowledge, we are the first to demonstrate how the eigenvalue decay and the source condition of $f^*$ affect the degree of polynomial approximation.
\end{itemize}

\paragraph{The polynomial barrier} In a rough sense, the polynomial barrier phenomenon asserts that in the power regime where $N\sim d^\iota$ for some $\iota\in\bN_+$ and as both $d$ and $N$ approach infinity, regardless of the choice of $\lambda$, $\hat f_\lambda$ can only approximate a degree-$\iota$ polynomial approximation of $f^*$, \cite{donhauser_how_2021}.
We note that, as discussed in Section~\ref{sec:proof_multiple_descent} in supplementary material, $\cH_{1:k}\subset\oplus_{l\in [k]}^\perp V_{d,l}$ for some spaces $V_{d,l}$ spanned by spherical harmonics of dimension $d$ and whose degree is less than $l$ with $\dim(V_{d,l})\sim d^l$. Since \(\mathcal{H}_{1:k}\) is \(k\)-dimensional, it is contained in the space spanned by \(\log_d(k)\) order spherical harmonics. Hence, our choice of $k=N^{\frac{1}{1+(s\wedge 2)\alpha}}$ indicates that $\hat f_\lambda$ with optimal choice of $\lambda$ approximates $f_{1:N^{\frac{1}{1+(s\wedge 2)\alpha}}}^*$ and treats $f_{N^{\frac{1}{1+(s\wedge 2)\alpha}}+1:\infty}^*$ as part of the noise, which corresponds to approximating $f^*$ as a degree-$\lfloor\frac{1}{1+(s\wedge 2)\alpha}\frac{\log{N}}{\log{d}}\rfloor$ polynomial, with the remaining higher-order terms considered as noise. 


Breaking the polynomial barrier is a good motivation for studying non-rotationally invariant kernels. For instance, one such kernel function is the conjugate kernel employed in the context of feature learning, as discussed in \cite{ba_high-dimensional_2022,moniri_theory_2023}, see also Section~\ref{sec:conjugate_kernel}. We leave this question opens even though we have no doubt that Theorem~\ref{theo:upper_KRR} and \ref{theo:main_upper_k>N} could be useful tools to answer this question.

\subsection{The Gaussian Equivalence Property}\label{sec:gaussian_equivalence}
Up to this point, we apply our general bounds to some specific choices of kernels such as polynomial or smooth ones. However, as mentioned in Section~\ref{sec:RKHS}, ``readers unfamiliar with RKHS theory can think of the feature map $\phi(X)$ as a Gaussian random vector, even though their nature is fundamentally different.'' One might have noticed a phenomenon: when treating $\phi(X)$ as a centered Gaussian random vector with covariance operator $\Gamma$ and calculating the estimation error by substituting it into Theorem~\ref{theo:upper_KRR} and Theorem~\ref{theo:main_upper_k>N}, the results seem to be almost identical to what we obtained through our ``non-Gaussian'' theory. This observation is known as the Gaussian equivalence property. 

In order to prove such a property, we first need a sharp, i.e. matching upper and lower bounds for KRR in the Gaussian setup. In order to do so, we apply Theorem~\ref{theo:upper_KRR} to $\phi:x\mapsto x$ and $X$ is a centered Gaussian random vector with covariance operator $\Gamma$. This provides an upper bound that is given in the result below. For the matching lower bound (which is also presented in this result), see Section~\ref{sec:omitted} in supplementary material for its formal version; its proof can be found in Section~\ref{sec:proof_linear_case} in supplementary material.

\begin{Propositionbis}{prop:linear}
    For the Gaussian linear regression problem, if we choose \( k \) as \( k_{b,\lambda}^* = \min\left(k\in\bN:\, \frac{\lambda + \Tr\left(\Gamma_{k+1:\infty}\right)}{\norm{\Gamma_{k+1:\infty}}_{\text{op}}}\geq bN\right) \), then, when $r_{\lambda,k_{b,\lambda}^*}^*\lesssim N$,  we have \(\mathbb{E}\|\Gamma^{1/2}(\hat{f}_\lambda - f^*)\|_2 \gtrsim r_{\lambda,k_{b,\lambda}^*}^*\) and with high probability, we have $\|\Gamma^{1/2}(\hat{f}_\lambda - f^*)\|_2 \lesssim r_{\lambda,k_{b,\lambda}^*}^*$.
\end{Propositionbis}

From Proposition~\ref{prop:linear} we can see that the estimation error of KRR in the linear Gaussian case is (up to a multiplicative constant) equivalent to that of Theorem~\ref{theo:upper_KRR}, i.e., of KRR for a large class of kernels and design vectors. It therefore shows the Gaussian equivalence property: we can indeed 'replace' $\phi(X)$ by a centered Gaussian vector with covariance operator identical to the integral operator of $\phi(X)$, i.e., $\Gamma$. 

This is a counter-intuitive and practical property that should hold in the non-asymptotic regime, called the (one-sided) \emph{Gaussian Equivalent Property} that we can formalize in the following definition.
\begin{Definition}\label{def:GEP}
Let $K$ be a kernel on $\bR^d\times \bR^d$ (with associated RKHS $\cH$ and feature map $\phi$), $X$ be a random variable with values in $\bR^d$, $\xi$ be a real valued random variable independent of $X$ and $Y=\inr{f^*, \phi(X)}_\cH + \xi$ for some $f^*\in\cH$. We say that the triple $(K, X, \xi)$ satisfies the (one-sided) \textbf{Gaussian equivalence property (GEP)} when the convergence rate of the KRR  constructed from $N$ i.i.d. copies of $(X,Y)$ and with regularization parameter $\lambda$ is smaller than $r_{\lambda,k_{b,\lambda}^*}^*$, up to some absolute multiplicative constant, for $b\sim \max\{\oC{C_lower_kappa_DM}\kappa_{DM}^{-1}, 2\oc{c_lower_7}^{-1}\}$ where $k_{b,\lambda}^*$ is defined in  \eqref{eq:def_k_b_star} and $r^*_{\lambda, k}$ is defined in \eqref{eq:def_rate}.
\end{Definition}

\begin{Remark}
    \begin{enumerate}
        \item Our Definition~\ref{def:GEP} only provides that if we replace the feature map with a Gaussian random vector, the resulting excess risk can serve as an upper bound for the excess risk of the original problem. In \cite{montanari_universality_2022,hu_universality_2022}, the GEP is typically regarded as a two-sided bound. However, in this paper, we focus solely on the upper bound.
        \item We notice that in some papers, the requirements for the Gaussian Equivalent Conjecture are stricter; they necessitate the explicit construction of an operator $\Gamma_0$ to approximate $\Gamma$ (in the sense of operator norm), as seen in \cite{moniri_theory_2023,ba_high-dimensional_2022,hu_universality_2022}. This approach often involves a specific analysis of the kernel function, for example, the random feature kernel. In the next section, we aim to provide a more universally applicable method, so we only provide the following insight: for the conjugate kernel defined later, potentially data-dependent, \cite[Lemma 5.2]{wang_deformed_2023} offers a way to construct such an approximating operator $\Gamma_0$.
    \end{enumerate}
\end{Remark}

This property was rigorously proven by \cite{montanari_universality_2022,hu_universality_2022,ba_high-dimensional_2022} but in the limited range of the proportional asymptotic regime and a Gaussian design assumption (i.e. $X$ is assumed to be Gaussian). Given the inherent appeal and practical value of this property, natural questions arise: \textbf{does the Gaussian Equivalence Property hold in non-asymptotic scenarios? Or for sub-Gaussian designs?} We emphasize that the significance of non-asymptotic theory and power regimes in the context of deep learning theory has been thoroughly expounded upon in Section~\ref{sec:multiple_descent}. Thanks to Proposition~\ref{prop:gaussian_equivalent_conjecture} and Theorem~\ref{theo:upper_KRR}, we are in a position to show the (one-sided) Gaussian equivalence property as defined above for a large class of kernels, design vectors $X$ and noises $\xi$.

\begin{Proposition}\label{prop:gaussian_equivalent_conjecture}
If Assumption~\ref{assumption:DM_L4_L2}, Assumption~\ref{assumption:upper_dvoretzky} and Assumption~\ref{assumption:RIP} are valid, $(K, X, \xi)$ satisfies the (one-sided) Gaussian Equivalence Property.

\end{Proposition}

This finding is applicable in a  broad range of kernels, extending to non-asymptotic scenarios, and allowing for general design vectors and noises. This represents a significant improvement compared to previous results \cite{montanari_universality_2022,hu_universality_2022,ba_high-dimensional_2022}. In particular, the non-asymptotic nature of Proposition~\ref{prop:gaussian_equivalent_conjecture} is important to us because the GEP may be seen as a 'non-asymptotic CLT' for KRR which is one aspect of the 'surprising' nature of this property.


\subsection{Data-dependent conjugate kernel}\label{sec:conjugate_kernel}
In recent years, research in Deep Learning Theory has increasingly focused on data-dependent kernels due to the fact that rotationally invariant kernels exhibit a polynomial barrier, \cite{ghorbani_linearized_2021,donhauser_how_2021}. One approach involves training a neural network on a subset of samples using gradient descent and using the obtained ``representation feature'' (which will be defined later) as a feature map to define a (data-dependent) kernel. The advantage of the associated RKHS lies in its data-dependent feature map, enabling the alignment of its eigenfunctions with $f^*$ for a better approximation of $f^*$. This alignment can significantly reduce the estimation error and is sometimes call feature learning. 

For instance, \cite[section 4.3]{ba_high-dimensional_2022} established that under certain conditions, performing a single, sufficiently large step in gradient descent followed by data-dependent kernel KRR, results in an estimation error smaller than that of any inner product kernel at any tuning parameter. Moreover, this alignment can also learn hidden subspaces, thereby mitigating the ``curse of dimensionality'', as discussed in \cite{dandi_how_2023,damian_neural_2022,bach_breaking_2017}. This field has seen rapid development in the recent years, as evidenced by papers such as \cite{ba_high-dimensional_2022,dandi_how_2023,damian_neural_2022,moniri_theory_2023,ben_arous_online_2021,bietti_learning_2022,bietti_learning_2023}. However, due to the complexity of these data-dependent kernels, especially their lack of rotational invariance and the absence of specific spectral decay, researchers in Deep Learning Theory often find it challenging to study the estimation error of KRR defined over these RKHS, resorting to rough bounds (\cite[Theorem 11]{ba_high-dimensional_2022}, \cite{damian_neural_2022,mousavi-hosseini_neural_2023,bietti_learning_2022}) that does not lead to benign overfitting, or limited investigation in asymptotic regimes \cite{ba_high-dimensional_2022}, or they may have to assume the validity of a certain Gaussian Equivalent Conjecture, \cite{moniri_theory_2023,dandi_how_2023}. This poses difficulties for the development of DNN approximation and feature learning via data-dependent kernels. In this section, we apply Theorem~\ref{theo:upper_KRR} to such data-dependent kernels to obtain upper bounds on the estimation error of the associated KRR estimator $\hat f_\lambda$.

In this section, we are concerned with a fully-connected shallow neural network with $m$ neurons (for deep neural networks, we can similarly define its conjugate kernel),
\begin{align}\label{eq:def_shallow_NN}
    f:\vx\in\bR^d\mapsto \frac{1}{\sqrt{m}}\left\langle \vw, \sigma\left(W \vx\right) \right\rangle_{\ell_2^m}\in\bR,
\end{align}where $\vw\in\bR^m$, $W\in\bR^{m\times d}$, and $\sigma:\bR\to\bR$ is the non-linearity, applied coordinate-wisely, that is, $\sigma(W\cdot)=(\sigma(\left\langle W_j, \cdot\right\rangle_{\ell_2^d}))_{j\in[m]}$, where $(W_i)$'s are row vectors of $W$. The RKHS $\cH$ generated by the feature map $\phi:\vx\mapsto \frac{1}{\sqrt{m}}\sigma\left(W\vx\right)$ is called the conjugate kernel of shallow neural network with non-linearity $\sigma$ and first layer $W$.
We identify $\cH$ with $\bR^m$ and $\phi$ with $m^{-1/2}\sum_{j=1}^m \sigma(\left<W_j,\vx\right>)\ve_j$ where $(\ve_j)_{j=1}^m$ is the canonical basis of $\bR^m$.
In other words, $\cH = \{f_\vw(\cdot)=\left<\vw,\phi(\cdot)\right>: \vw\in\bR^m\}$ is composed of all shallow neural networks with $\sigma$ as the activation function and $W$ as the weight matrix of the first layer and each $f_\vw$ is identified by $\vw\in\bR^m$.  At this point, \(\Gamma = \bE[\phi(X)\otimes\phi(X)]\) is identified as a linear operator on \(\mathbb{R}^m\), and its eigenfunctions are identified as the vectors in \(\mathbb{R}^m\) given by \(\varphi_j\).
Using the spectral decomposition of $\Gamma$, there exists a unitary operator $U = [\varphi_1|\varphi_2|\cdots]$ such that for any $j\in[m]$, $\varphi_j = U\ve_j$.

As the analysis of KRR depends on the decomposition of $\cH$ into two orthogonal subspaces $\cH_{1:k}\oplus^\perp \cH_{k+1:m}$ for some $k\in[m]$, we need to introduce the decomposition of $\phi$. Recall that $\cH_{1:k}=\Span(\varphi_j:1\leq j\leq k)$, we rewrite $\phi(\vx)$ as $\phi(\vx)=\sum_{j=1}^m \frac{1}{\sqrt{m}}\sigma\left(\left\langle W_j,\vx\right\rangle\right)\ve_j = \sum_{j\in J} \frac{1}{\sqrt{m}}\sigma\left(\left\langle W_j,\vx\right\rangle\right)\ve_j + \sum_{j\in J^c} \frac{1}{\sqrt{m}}\sigma\left(\left\langle W_j,\vx\right\rangle\right)\ve_j=: \phi_J(\vx)+\phi_{J^c}(\vx)$ for $J=[k]$ and $\phi_J=P_{\cH_{1:k}}\phi$ and $\phi_{J^c}=P_{\cH_{k+1:m}}\phi$. Then $\cH_{1:k}=U\Span(\ve_j:j\in J)$ and $\cH_{k+1:m}=U\Span(\ve_j:j\in J^c)$. One can check that they are orthogonal spaces.
We write the restriction of $W$ onto $J$ as $W_J$, that is, the $j$-th row of $W_J$, denoted as $(W_J)_j$, equals $W_j$, if $j\in J$, and $(W_J)_j=0$, if $j\in J^c$. Intuitively, this can be understood as using the neurons associated with the subset $J$ to estimate $f^*$, while the neurons corresponding to $J^c$ are utilized to absorb noise.

Given the samples $(X_i, Y_i)_{i\in[N_1+N_2]}$ for some $N_1,N_2\in\bN_+$, we use $(X_i, Y_i)_{i\in[N_1]}$ to perform gradient descent with respect to the squared loss and with step size $\eta\sqrt{m}$ for training the weight matrix $W$. 
After the completion of this gradient descent, we employ $(X_i, Y_i)_{N_1+1\leq i\leq N_1+N_2}$ to perform KRR on the obtained conjugate kernel with non-linearity $\sigma$ and first layer $W$ obtained by gradient descent in the previous step -- so that $W$ depends on $(X_i, Y_i)_{i=1}^{N_1}$. We call this RKHS \emph{data-dependent conjugate kernel}. Due to the independence of the samples, in the second stage, we can treat $W$ as a deterministic matrix (though $W$ has undergone alignment with respect to $f^*$ through gradient descent). 
In the remaining part of this section, we always take condition on random vectors $(X_i, Y_i)_{i\in[N_1]}$ (that is, trained $W$) because we want to study the features learning property of gradient descent. Hence, we assume that the gradient descent step is finished and we then want to study the features that have been learned during that step. Before introducing our results, we emphasize that through similar analysis, our results can be extended to neural networks with bias, or even deep neural networks. Here, we only present the simplest case.
In this section, we acknowledge the following assumption.

\begin{Assumption}\label{assumption:conjugate_kernel}
There exist $B,D>0$ depending on the gradient descent procedure, such that the following holds.
    \begin{enumerate}
    \item For each $i\in[N_1+N_2]$, $X_i$ are i.i.d. $\cN(0,I_d)$, the noise $\xi$ is independent of $X$ with mean zero and variance $\sigma_\xi^2$. We assume that there exist absolute constants $\okappa{kappa_noise}>0$ and $r>4$ such that  $\norm{\xi}_{L_r}\leq \okappa{kappa_noise}\sigma_\xi$.
        \item $\sigma$ is Lipschitz on $\bR$ with Lipschitz constant $\norm{\sigma}_{Lip}$. Moreover, $\sup_{j\in[m]}\norm{\sigma(\left\langle W_j,\cdot\right\rangle)}_{L_2}<D$.
        \item The first layer $W$ after gradient descent satisfies the following equations: for any $i,j\in J^c$, $\left<W_i,W_j\right>=\delta_{ij}$, $\norm{W_J}_{\text{op}}\leq B$ and $\sum_{j\in J}(\norm{W_j}_2^2-1)^2\leq B^2$.
        \item The choice of $k$ for $k\lesssim N_2$ and tuning parameter $\lambda\geq 0$ satisfies the following equations (we recall that $\Gamma = \bE \phi(X)\phi(X)^\top\in\bR^{m\times m}$):
        \begin{enumerate}
            \item $m\Tr\left(\Gamma_{k+1:m}^2\right)\gtrsim \lambda_\sigma^2 \norm{\Gamma_{k+1:m}}_{\text{op}}$,
            \item $(\Tr\left(\Gamma_{k+1:m}\right)+\lambda)\norm{\sigma}_{L_2(\gamma_d)}^2\gtrsim \lambda_\sigma^2N_2\|\Gamma_{k+1:m}\|_{\text{op}}$,
            \item $m\sigma_k\sqrt{N_2}\geq 2\lambda_\sigma^2 B^2 $.
        \end{enumerate}
    \end{enumerate}
\end{Assumption}

Points \textit{1.} and \textit{2.} in the context of feature learning and conjugate kernel research are the most common settings \cite{wang_deformed_2023,fan_spectra_2020,mei_generalization_2022}. When $X$ has covariance matrix $\Sigma$, the hidden layer matrix $W$ can be replaced by $W\Sigma^{1/2}$. For a general probability distribution, such as the sub-Gaussian distribution, refer to Remark~\ref{remark:general_distribution} in the supplementary material. Points \textit{3.} and \textit{4.}, on the other hand, have been relatively rare in prior studies, primarily because we characterize which subspace of the RKHS absorbs noise when using $\hat{f}_\lambda$, a perspective that previous work did not investigate (except \cite[Appendix D]{ba_high-dimensional_2022}). In \textit{3.}, there is a partitioning of $[m]$. Neurons belonging to $J$ are used for prediction, while neurons belonging to $J^c$ are used to absorb noise. In many cases, $\|W_j\|_2=1$ holds for all $j\in[m]$, for example, when using Stiefel gradient flow to learn the hidden layer (in this case, the hidden layer matrix $W$ consists of orthonormal row vectors of $\bR^d$), as seen in \cite{bietti_learning_2023}. This model yields a typical weight matrix produced by gradient descent reported in many works, such as \cite{boursier_early_2024} and references therein.
However, upper bounds on $\|W_J\|_{\text{op}}$ and $\sum_{j\in J}(\|W_j\|_2^2-1)^2$ as well as verifying condition \textit{4.} require a more detailed examination of the optimization algorithm, particularly with regard to the feature learning it induces. We provide some examples when the spectrum and eigenvectors of $\Gamma$ may be obtained in Section~\ref{sec:example_Gamma} in supplementary material. We assume $k\lesssim N_2$ for simplicity. 
We provide an example satisfying Assumption~\ref{assumption:conjugate_kernel} in Proposition~\ref{prop:feature_learning} after the statement on the performance of KRR for this data-dependent kernel learned by gradient descent that follows from our general bounds.

\begin{Propositionbis}{prop:conjugate_kernel}
    Grant Assumption~\ref{assumption:conjugate_kernel}. If $k\lesssim N$, for any $\lambda\geq 0$, the same conclusions as in Theorem~\ref{theo:upper_KRR} hold.
\end{Propositionbis}

The formal version of Proposition~\ref{prop:conjugate_kernel} may be found in Section~\ref{sec:formal_data_dependent} in supplementary material, and the proof of Proposition~\ref{prop:conjugate_kernel} can be found in Section~\ref{sec:proof_conjugate_kernel} in the supplementary material. As far as we know, Proposition~\ref{prop:conjugate_kernel} is the first result on the estimation error of KRR on data-dependent conjugate kernels.

Below, we provide a toy example to validate Assumption~\ref{assumption:conjugate_kernel} and utilize this example to illustrate how gradient flow performs feature learning.

\begin{Propositionbis}{prop:feature_learning}
    Under certain conditions, there exists a learning algorithm such that after training the first layer with \((X_i, Y_i)_{i=1}^{N_1}\), minimum data-dependent conjugate RKHS norm interpolant estimator applied to the samples \((X_i, Y_i)_{i=N_1+1}^{N_1+N_2}\) yields \(\hat{f}_0\) that satisfies, with probability at least $1/2$,
    \begin{align}\label{eq:result_feature_learning}
    \begin{aligned}
        \norm{\hat f_0 - f^*}_{L_2}&\lesssim \left|a^*\right|\norm{\sigma}_{Lip}\tilde O\left(\max\left\{\sqrt{\frac{d+m}{N_1}},\frac{d^2}{N_1}\right\}\right)+ \frac{\sigma_\xi}{\sqrt{N_2}} + \sigma_\xi\sqrt{\frac{N_2}{m}} + \frac{\left|a^*\right|}{\|\sigma\|_{L_2(\gamma_d)}}\frac{1}{N_2}
    \end{aligned}
\end{align}in the single neuron model, i.e. when $f^*:x\to a^*\sigma(\inr{w^*, x})$ for $a^*\in\bR$ and $w^*\in S_2^{d-1}$.
\end{Propositionbis}

The formal version of Proposition~\ref{prop:feature_learning} may be found in Section~\ref{sec:omitted} in supplementary material and the proof of it may be found in Section~\ref{sec:proof_feature_learning}. Some comments are in order.

We utilize the conclusions of feature learning provided in \cite{bietti_learning_2022}. The probability deviation $1/2$ in Proposition~\ref{prop:feature_learning} arises from the special ``random feature approximation'' technique utilized in \cite{bietti_learning_2022}, as discussed in \cite[Appendix D.4]{bietti_learning_2022}. Since feature learning is beyond the scope of this paper, we directly rely on the conclusions from \cite{bietti_learning_2022} instead of attempting to improve upon them.
In Proposition~\ref{prop:feature_learning}, the excess risk is the approximation error (the first term in the right-hand-side of Equation~\eqref{eq:result_feature_learning}) plus estimation error (other terms). The latter term becomes very small when the layer width $m$ is sufficiently large relative to the number of samples $N_2$ used to train the second layer. In \cite[Section 6]{bietti_learning_2022}, the authors conjecture that $\sqrt{d+m}$ can be improved to $\sqrt{d}$. If this conjecture holds true, then in our Proposition~\ref{prop:feature_learning}, after learning the hidden layer through gradient descent, the excess risk is decreasing with the layer width $m$. This implies that even when the parameters of the shallow neural network are numerous (indicating a large layer width $m$), resulting in small approximation error for the shallow neural network, the estimation error obtained through stochastic gradient descent (where feature learning happens) and the minimum data-dependent conjugate RKHS norm interpolant estimator remains small. This contradicts the intuition of the classical approximation-estimation trade-off. However, we emphasize that a large value of $m$ will lead to computational difficulties.

\section{Supplementary material}

The supplementary material is organized as follows.
\begin{itemize}
 \item Section~\ref{sec:linearization} contains our result of the linearization of kernel matrices.
 \item Section~\ref{sec:future_directions} contains some future directions.
 \item Section~\ref{sec:omitted} contains the parts that were omitted from the main text due to length constraints. This includes the following content:
 \begin{itemize}
     \item Section~\ref{sec:formal_k_large}: formal version of Theorem~\ref{theo:main_upper_k>N} on the case when $k$ is not necessarily smaller than $N$.
     \item Section~\ref{sec:alternative_assumption}: alternative assumptions of Assumption~\ref{assumption:upper_dvoretzky} and Assumption~\ref{assumption:DMU_used_for_RIP}.
     \item Section~\ref{sec:comparision_BMR}: comparison between our Theorem~\ref{theo:upper_KRR} and \cite[section 4.3.3, Theorem 4.13]{bartlett_deep_2021}.
     \item Section~\ref{sec:previous_results_multiple_descent}: previous results on multiple descent.
     \item Section~\ref{sec:formal_multiple_descent}: formal version of Proposition~\ref{prop:multiple_descent} on multiple descent.
     \item Section~\ref{sec:formal_data_dependent}: formal version of Proposition~\ref{prop:conjugate_kernel} on data-dependent conjugate kernel.
     \item Section~\ref{sec:example_feature_learning}: formal version of Proposition~\ref{prop:feature_learning} on feature learning.
     \item Section~\ref{sec:matcing_lower_bound}: formal version of Proposition~\ref{prop:linear} on Gaussian linear regression.
     \item Section~\ref{sec:smooth_omitted}: estimation error of KRR with smooth kernel in non-asymptotic regime.
     \item Section~\ref{sec:source_condition_omitted}: source condition in the language of interpolation space.
 \end{itemize}
 \item Section~\ref{sec:proof_main_upper_k<N} contains the proof of Theorem~\ref{theo:upper_KRR}.
 \item Section~\ref{sec:proof_main_upper_k>N} contains the proof of Theorem~\ref{theo:main_upper_k>N}.
 \item Section~\ref{sec:proof_DM} contains the proof of Theorem~\ref{theo:DM_RKHS}.
\item Section~\ref{sec:proof_smooth_kernel} contains the proof of Theorem~\ref{prop:KRR_smooth_upper}.
\item Section~\ref{sec:aux_proofs} contains the proof of all the other theorems and lemmas.
\end{itemize}

\subsection{Linearization of non-linear kernel matrix.} \label{sec:linearization}
Despite Theorem~\ref{theo:DM_RKHS} being of a non-asymptotic nature, one can question the insights that Theorem~\ref{theo:DM_RKHS} can provide in the asymptotic regime. 
A series of papers \cite{el_karoui_spectrum_2010,do_spectrum_2013,cheng_spectrum_2013,fan_spectral_2019} investigated the spectral properties of $(K_{k+1:\infty}(X_i,X_j))_{1\leq i,j\leq N}$ as the dimension $d$ (recall that $K_{k+1:\infty}$ is a kernel defined on $\bR^d\times\bR^d$) and number of samples $N$ tend to infinity, while maintaining the ratio $d/N$ to be a constant level. As $(K_{k+1:\infty}(X_i,X_j))_{1\leq i,j\leq N} = \bX_{\phi,k+1:\infty}\bX_{\phi,k+1:\infty}^\top$, our Theorem~\ref{theo:DM_RKHS} also produces results on the spectrum of $(K_{k+1:\infty}(X_i,X_j))_{1\leq i,j\leq N}$. The reader may find it of interest to compare Theorem~\ref{theo:DM_RKHS} with the theorems presented in the works of \cite{el_karoui_spectrum_2010,do_spectrum_2013,cheng_spectrum_2013,fan_spectral_2019}. In the following, we consider a general $(K(X_i,X_j))_{1\leq i,j\leq N}$ instead of $(K_{k+1:\infty}(X_i,X_j))_{1\leq i,j\leq N}$. The following theorem is the main result of \cite{el_karoui_spectrum_2010} applied to isotropic design vector.
\begin{Theorem}[\cite{el_karoui_spectrum_2010}]\label{theo:El_Karoui}
    Assume that there exist some absolute constants $c,C$ such that $c\leq d/N\leq C$ and $(X_i)_{i\in\bN}\subset\bR^d$ are i.i.d. isotropic random vectors. Assume the coordinates of $X_i$ are $(x_{ij})_{j\leq d}$, and assume that $x_{ij}$'s have mean $0$, variance $1$ and $(4+\varepsilon)$ moment for some $\varepsilon>0$. Let $K$ be an inner product kernel as in \eqref{eq:polynomial_inner_product_kernel} with $h:\bR\to\bR$ is a $C^1$ function in a neighborhood of $1$ and a $C^3$ function in a neighborhood of $0$. There exists a matrix $M$ defined as
    \begin{align*}
        M = \left(h(0)+h''(0)\frac{1}{2d}\right) \1_N\1_N^\top + \frac{h'(0)}{d}\bX\bX^\top + \left(h(1)-h(0)-h'(0)\right)I_N,
    \end{align*}where $\bX$ is the design matrix $[X_1|X_2|\cdots|X_N]^\top$, such that as $N,d\to\infty$ while $d/N$ remain bounded, we have $\norm{\left(K(X_i,X_j)\right)_{1\leq i,j\leq N}-M}_{\text{op}}\to 0$ in probability.
\end{Theorem} 

Theorem~\ref{theo:El_Karoui} suggests that when both $N$ and $d$ tend to infinity, while $d/N$ remains at a constant level, the spectrum of the kernel matrix $(K(X_i,X_j))_{1\leq i,j\leq N}$ can be approximated by that of a linearized kernel, referred to as the linearization of the kernel matrix with non-linear kernel. The regularity of $h$ in Theorem~\ref{theo:El_Karoui} is relaxed by \cite{cheng_spectrum_2013} to $C^1$ near $0$ and by \cite{do_spectrum_2013} to be differentiable at $0$ and continuous at $1$.
The results from \cite{cheng_spectrum_2013,fan_spectral_2019} demand $X$ to be a Gaussian random vector. Moreover, \cite[Theorem 2.1]{el_karoui_spectrum_2010} establishes convergence in probability, while \cite{cheng_spectrum_2013,do_spectrum_2013,fan_spectral_2019} only achieve convergence in distribution using the Stieltjes transform or Free Probability techniques.

One might ask the question: what does the linearization of non-linear kernel random matrices entail? From a statistical perspective, linearization offers a means to analyze KRR through a kernel matrix of a linear kernel. By examining the spectrum of the limit matrix $M$ defined in Theorem~\ref{theo:El_Karoui}, and applying Weyl's inequality, we can gain insights into the spectrum of the kernel matrix. This, in turn, allows us to determine the estimation error of KRR. In fact, this approach is commonly adopted in many works within the field of Deep Learning theory, \cite{liang_just_2020,hastie_surprises_2022}. 

However, as we have observed in Section~\ref{sec:multiple_descent}, Theorem~\ref{theo:El_Karoui} only holds in the proportional regime, specifically when $d\sim N$. This limitation implies that we cannot observe 1) the multiple descent phenomenon through Theorem~\ref{theo:El_Karoui} and 2) the higher degree approximation of $f^*$ that KRR performs -- an important property of deep neural networks . Natural questions therefore arise: can we encounter the linearization of  kernel matrices with non-linear kernel beyond the proportional regime? can we get quantitative (i.e. non-asymptotic) results concerning this approximation?

Our subsequent conclusions reveals that although we have not been able to establish the approximation of a kernel matrix $(K(X_i, X_j))_{1\leq i,j\leq N}$ by the kernel matrix of a linearized kernel in the topology generated by the operator norm as in Theorem~\ref{theo:El_Karoui}, we have demonstrated that the spectrum of the kernel matrix converges in probability to a considerably smaller region (this convergence pertains to the topology on $\bR^N$). 
However, given the interest in the theoretical foundations of deep learning community in results such as Theorem~\ref{theo:El_Karoui} (despite its non-quantitative and proportional regimes range of applications), we expect that these findings will prove to be useful in future research within the domain of statistical deep learning Theory, for instance, \cite{donhauser_how_2021}.

The following is the assumption that will be granted in this section.


\begin{Assumption}\label{assumption:linearization}
        \begin{enumerate}
        \item $N,d\to\infty$, while $d^\iota=o(N)$, $N=o(d^{\iota+1})$ for some $\iota\in\bN_+$.
        \item $X\in\bR^d$ is a mean-zero random vector with i.i.d. sub-Gaussian coordinates. Moreover, for any finite set $S$, the first coordinate (denoted as $x_1$) of $X$ satisfies $\bP(x_1\in S)<1$ and $\bP(x_1 = 0)>0$. We denote $\mu_1$ as the distribution of $x_1$, thus $\mu=\otimes_{j=1}^d \mu_1$.
        \item For some $\iota,L\in\bN_+$, there exists a function $h:t\in\bR\mapsto \sum_{i=\iota}^L \alpha_i t^i\in\bR$ for some $\alpha\in\bR_+$, $\alpha_i\sim 1$ for all $i\in[L]$, such that \eqref{eq:polynomial_inner_product_kernel} holds.
        
    \end{enumerate}
\end{Assumption}

The following proposition is derived from a combination of Theorem~\ref{theo:DM_RKHS} and Lemma~\ref{lemma:liang_et_al}.
\begin{Proposition}\label{prop:linearization}
    Suppose Assumption~\ref{assumption:linearization} holds. There exist absolute constant $\nc\label{c_distortion_asymptotic}$ such that the following holds. Then
    \begin{align}\label{eq:objective_linearization}
        \begin{aligned}
            \bP\bigg\{&(1-\oc{c_distortion_asymptotic})h(1)\leq \sigma_{N}\left({(K(X_i,X_j))_{1\leq i,j\leq N}}\right)\leq\\
        &\sigma_{1}\left({(K(X_i,X_j))_{1\leq i,j\leq N}}\right)\leq (1+\oc{c_distortion_asymptotic})h(1)\bigg\} \to 1
        \end{aligned}
    \end{align}as $N,d\to\infty$, where we recall that $\sigma_N(A)$ and $\sigma_1(A)$ represents the smallest and the largest singular value of some matrix $A$.
    
\end{Proposition}

The proof of Proposition~\ref{prop:linearization} can be found in Section~\ref{sec:proof_linearization}. We now provide some comments on Proposition~\ref{prop:linearization}.

\begin{enumerate}
    \item Theorem~\ref{theo:El_Karoui} and Proposition~\ref{prop:linearization} both characterize the properties of the spectrum of the kernel matrix as $d$ and $N$ approach infinity. However, these two conclusions exhibit several differences. Firstly, Theorem~\ref{theo:El_Karoui} states that the kernel matrix can be approximated by a matrix in the operator norm on $(\bR^N, \ell_2) \to (\bR^N, \ell_2)$, but it does not provide information about whether the condition number of the kernel matrix is bounded. In fact, if one uses the conclusion of Theorem~\ref{theo:El_Karoui}, further analysis of the spectrum of the matrix $M$ is necessary. On the other hand, Proposition~\ref{prop:linearization} directly indicates that the condition number of the kernel matrix is nearly 1, and all eigenvalues are concentrated around $h(1)$. This type of results are useful for the analysis of KRR via the geometric method.

    \item Furthermore, Proposition~\ref{prop:linearization} holds beyond the proportional regime $N\sim d$ in the power regime $N\sim d^\iota$, whereas Theorem~\ref{theo:El_Karoui} only applies to the proportional regime. Even if we set $\iota=1$, the range of validity of Proposition~\ref{prop:linearization} is approximately $d\lesssim N\lesssim d^2$, which is notably broader than the range of validity of Theorem~\ref{theo:El_Karoui} (which is $N\sim d$).
    
    \item However, it must be emphasized that the aforementioned advantages of Proposition~\ref{prop:linearization} come at a cost. It's worth noting that Proposition~\ref{prop:linearization} is applicable only to finite-degree polynomial kernels, and the Taylor coefficients of the function $h$ must be positive. In contrast, Theorem~\ref{theo:El_Karoui} has a much broader range of applicability regarding the choice of kernels. One may wish to extend Proposition~\ref{prop:linearization} using the Stone-Weierstrass theorem; however, this is not feasible. This is because in Proposition~\ref{prop:linearization}, the constants associated with the condition $d^\iota\lesssim N\lesssim d^{\iota+1}$ conceal information about the Taylor coefficients of $h$. Readers can verify for themselves that when $\iota=1$, we require $N\leq 4\frac{\sum_{i\leq L}\alpha_i}{a_1}\kappa_{DM}d$. We conjecture that methods similar to those in \cite{lu_equivalence_2023,dubova_universality_2023} may help extend Proposition~\ref{prop:linearization} to smooth kernel functions.

    \item One may observe that Proposition~\ref{prop:linearization} aligns with the ``diagonal'' component as discussed in \cite[section 2.2 (B)]{el_karoui_spectrum_2010}. By ``diagonal components,'' we are referring to the elements along the diagonal of the kernel matrix $(K(X_i,X_j))_{1\leq i,j\leq N}$. One may notice that this precisely corresponds to \eqref{eq:diagonal_term_assumption} in Assumption~\ref{assumption:DM_L4_L2}.
    This phenomenon occurs in the Dvoretzky-Milman regime, specifically when the value of $N$ is less than $\kappa_{DM}d_0^*$. In this regime, the off-diagonal or cross terms are overshadowed by the diagonal part. This observation is also evident in the proof of Theorem~\ref{theo:DM_RKHS}. The aforementioned phenomena suggests that the Dvoretzky-Milman theorem serves as a characterization of the ``diagonal dominance'' in the context of asymptotic non-linear random matrix theory. This observation has not been   discerned using the Stieltjes transform or Free Probability approaches.
\end{enumerate}


\subsection{Further work}\label{sec:future_directions}

\paragraph{Beyond the rotational invariant kernel} In Section~\ref{sec:conjugate_kernel}, we provide an upper bound on the estimation error of KRR on data-dependent conjugate kernels and a toy example. However, for more general cases, such as the multi-index model, computing the eigenvalues and eigenvectors of $\Gamma_{k+1:\infty}$ remains a challenge.
An intriguing avenue of research involves investigating the approximation and estimation error of KRR with the kernel discussed above, specifically when the weight matrix $W$ is obtained from one \cite{ba_high-dimensional_2022,damian_neural_2022,moniri_theory_2023,ba_learning_2023} or many gradient steps \cite{bietti_learning_2023,mousavi-hosseini_neural_2023}. The objective is to understand how the data-dependent kernel aligns with the target function and how its spectrum and eigenvectors change, ultimately leading to a reduction in the estimation error even when overfitting happens ($\lambda=0$) or near overfitting happens ($\lambda$ is close to $0$), especially when the width of the neural network exceeds the sample size. In a word, we would like to understand the feature learning for this kernel after some gradient descent steps.

\paragraph{The Dvoretzky-Milman theorem generated by an RKHS feature map} In order to satisfy the classical Dvoretzky-Milman criterion, Theorem~\ref{theo:DM_RKHS} requires $\eps$ to be greater than 2 if one wants to avoid the extra assumption \eqref{eq:extra_assumption}. A question that poses a significant math challenge is the appropriate application of the coloring technique proposed in the work of \cite{tikhomirov_sample_2018} in order to enhance the value of $\eps$ from being greater than 2 to $\eps>0$.

\paragraph{Point-wise lower bound} The lower bound proven in Proposition~\ref{prop:linear} for the linear Gaussian case holds for every target function, hence it is referred to as a ``point-wise'' lower bound, and sometimes also known as an optimistic lower bound (in contrast to the minimax lower bound), \cite{bellec_optimistic_2017,kur_minimal_2021}.

A challenging problem is to obtain a point-wise lower bound for the estimation error of a general KRR. In this context, we aim to derive results that do not rely on the spectrum of the RKHS (such as power decay) or specific properties of the kernel function (e.g., inner-product or translation-invariant kernels).


\subsection{The parts omitted in the main text}\label{sec:omitted}

\subsubsection{Upper bound of the estimation error of KRR when $k$ is not necessarily smaller than $N$}\label{sec:formal_k_large}

In this section, we introduce the formal version of Theorem~\ref{theo:main_upper_k>N}, which was omitted in the main text. First, we need to introduce another assumption.

\begin{Assumption}\label{assumption:DMU_used_for_RIP}
    There exists absolute constants $0\leq\ngamma\label{gamma_RIP_k>N}<1/16$, $\ndelta\label{delta_RIP_k>N}\geq 0$, $\eps>0$, $\nkappa\label{kappa_DMU_for_RIP}\geq 1$ such that
    \begin{itemize}
    \item with probability at least $1-\ogamma{gamma_RIP_k>N}$, $\max_{i\in[N]}\norm{\tilde\Gamma_{1:k}^{-1/2}\phi_{1:k}(X_i)}_\cH^2\leq (1+\odelta{delta_RIP_k>N})\Tr\left(\tilde\Gamma_{1:k}\Gamma_{1:k}\right)$,
    \item for any $f\in\cH_{1:k}$, $\norm{f}_{L_{4+\eps}}\leq\okappa{kappa_DMU_for_RIP}\norm{f}_{L_2}$.
    \end{itemize}
\end{Assumption}The validity of Assumption~\ref{assumption:DMU_used_for_RIP} is not always true. When it is not valid, we provide below in the next section (Section~\ref{sec:alternative_assumption}) Assumption~\ref{assumption:upper_dvoretzky_infty} as an alternative assumption.

\begin{Theorem}\label{theo:main_upper_k>N}
Suppose Assumptions~\ref{assumption:DM_L4_L2},~\ref{assumption:upper_dvoretzky} and~\ref{assumption:DMU_used_for_RIP} hold. Suppose the noise $\xi$ is independent of $X$ with mean zero and variance $\sigma_\xi^2$. We assume that for some $\okappa{kappa_noise}>0$ and $r>4$,  $\norm{\xi}_{L_r}\leq \okappa{kappa_noise}\sigma_\xi$.
There then exist absolute constants $\oc{c_kappa_DM}$, $\oc{c_P_bX_f_star}$, $\oC{C_noise}$ ($\oC{C_noise}$ depends on $\okappa{kappa_noise}$) and $\nC\label{C_rate_2}$ such that the following holds. Let $\lambda\geq0$. Assume that there exists $k\in\bN\cup\{\infty\}$ such that $N\leq\oc{c_kappa_DM}\kappa_{DM}d_\lambda^*\left(\Gamma_{k+1:\infty}^{-1/2}B_\cH\right)$, and such that the the following equation holds:
    \begin{align}\label{eq:extra_condition_k>N}
    \begin{cases}
        \sum_{j\in J_2}\sigma_j \leq \kappa_{DM}(4\lambda + \Tr\left(\Gamma_{k+1:\infty}\right))\left(1-\frac{\left|J_1\right|}{N}\right), & \mbox{ when } N\sigma_1 \geq \kappa_{DM}(4\lambda + \Tr\left(\Gamma_{k+1:\infty}\right))\\
        \Tr\left(\Gamma_{1:k}\right)\leq N\sigma_1 , &\mbox{otherwise}.
    \end{cases}
\end{align}
    Let $\bar p_{\xi}$ be some probability deviation strictly less than $1$ (defined in \eqref{eq:def_bar_p_xi} later). Suppose that
    \begin{align}\label{eq:RIP_condition_k>N}
        \kappa_{DM}\left(4\lambda + \Tr\left(\Gamma_{k+1:\infty}\right)\right)\geq N\left(R_N^*(\oc{c_kappa_RIP})\right)^2.
    \end{align}
    Then, for all such $k$'s, with probability at least
    \begin{align*}
    1 - \bar p_{RIP} - \bar p_{DM} - 2\bar p_{DMU} - \frac{ \oc{c_P_bX_f_star} }{N}  - \bar p_{\xi}-\left(\frac{\oC{C_noise}\Tr(\Gamma_{k+1:\infty})}{|J_1|\Tr(\Gamma_{k+1:\infty}) + N \left(\sum_{j\in J_2}\sigma_j\right)}\right)^{\frac{r}{4}},
\end{align*}we have
\begin{align*}
    \norm{\hat f_\lambda - f^*}_{L_2}\leq \oC{C_rate_2} r_{\lambda,k}^*.
\end{align*}

\end{Theorem}

\subsubsection{Alternative assumption of Assumption~\ref{assumption:upper_dvoretzky} and Assumption~\ref{assumption:DMU_used_for_RIP}}\label{sec:alternative_assumption}

We provide an alternative setup for cases where Assumption~\ref{assumption:upper_dvoretzky} and Assumption~\ref{assumption:DMU_used_for_RIP} are not valid. This setup is given by the next assumption where there is no norm-equivalence assumption on $\cH$.

\begin{Assumption}\label{assumption:upper_dvoretzky_infty}
    There exist absolute constants $\ngamma\label{gamma_DMU_infty_1}$, $\ngamma\label{gamma_DMU_infty_2}$ and $\ngamma\label{gamma_DMU_infty_3}$ in $(0,1/16)$, and $\ndelta\label{delta_DMU_infty_1}$, $\ndelta\label{delta_DMU_infty_2}$ and $\ndelta\label{delta_DMU_infty_3}$ that are at least $1$ such that
    \begin{equation*}
        \bP\left(\max_{i\in[N]}\frac{\norm{\phi_{k+1:\infty}(X)}_\cH}{\sqrt{\Tr\left(\Gamma_{k+1:\infty}\right)}}\leq \odelta{delta_DMU_infty_1}\right)\geq 1-\ogamma{gamma_DMU_infty_1},\quad \bP\left(\max_{i\in[N]}\frac{\norm{\Gamma_{k+1:\infty}^{1/2}\phi_{k+1:\infty}(X)}_\cH}{\sqrt{\Tr\left(\Gamma_{k+1:\infty}^2\right)}}\leq \odelta{delta_DMU_infty_2}\right)\geq 1-\ogamma{gamma_DMU_infty_2}
        \end{equation*}and 
        \begin{equation*}
            \bP\left(\max_{i\in[N]}\frac{\norm{\tilde\Gamma_{1:k}^{-1/2}\phi_{1:k}(X)}_\cH}{\sqrt{\Tr\left(\tilde\Gamma_{1:k}^{-1}\Gamma_{1:k}\right)}}\leq \odelta{delta_DMU_infty_3}\right)\geq 1-\ogamma{gamma_DMU_infty_3}.
        \end{equation*}
\end{Assumption}

\begin{Proposition}\label{prop:upper_dvoretzky_infty}
    Suppose Assumption~\ref{assumption:upper_dvoretzky_infty} holds. There exists an absolute constant $\nC\label{C_DMU_pre}$, such that for any $0<\ndelta\label{delta_P_upper_dvoretzky}<1$, $0<\ndelta\label{delta_P_upper_dvoretzky_2}<1$ and  $0<\ndelta\label{delta_P_upper_dvoretzky_3}<1$, with probability at least $1-\odelta{delta_P_upper_dvoretzky} - \odelta{delta_P_upper_dvoretzky_2} - \odelta{delta_P_upper_dvoretzky_3} - \ogamma{gamma_DMU_infty_1} - \ogamma{gamma_DMU_infty_2} - \ogamma{gamma_DMU_infty_3}$, for any $\vlambda\in\bR^N$,
    \begin{align*}
        \norm{\bX_{\phi,k+1:\infty}^\top\vlambda}_\cH&\leq \odelta{delta_P_upper_dvoretzky}^{-1/2}\oC{C_DMU_pre}\log{N}\sqrt{\odelta{delta_DMU_infty_1}\Tr\left(\Gamma_{k+1:\infty}\right)}\norm{\vlambda}_2,\\
        \norm{\Gamma_{k+1:\infty}^{1/2}\bX_{\phi,k+1:\infty}^\top\vlambda}_{\cH}&\leq \odelta{delta_P_upper_dvoretzky_2}^{-1/2}\oC{C_DMU_pre}\left(\sqrt{N}\norm{\Gamma_{k+1:\infty}}_{\text{op}} + \log{N}\sqrt{\odelta{gamma_DMU_infty_2}\Tr\left(\Gamma_{k+1:\infty}^2\right)}\right)\norm{\vlambda}_2,\\
        \norm{\tilde\Gamma_{1:k}^{-1/2} \bX_{\phi,1:k}^\top \vlambda }_\cH &\leq \odelta{delta_P_upper_dvoretzky_3}^{-1/2}\oC{C_DMU_pre}\left(\log(N)\sqrt{\odelta{gamma_DMU_infty_3}}\sqrt{ \square^2 \left|J_1\right| + \triangle^2\sum_{j\in J_2}\sigma_j } + \sqrt{N}\sigma(\square,\triangle)\right)\norm{\vlambda}_2
    \end{align*}In this case, the symbol $\bar p_{DMU}$ is defined as $\odelta{delta_P_upper_dvoretzky_2} + \ogamma{gamma_DMU_infty_2}$.
\end{Proposition}
When the norm-equivalence condition is not satisfied, an additional logarithmic factor needs to be taken into account, as observed by comparing the upper side of \eqref{eq:objective_DM}, \eqref{eq:upper_dvoretzky_high_probability} and the corresponding expressions in Proposition~\ref{prop:upper_dvoretzky_infty} and the probability deviation is of constant level. The proof of Proposition~\ref{prop:upper_dvoretzky_infty} can be found in Section~\ref{sec:proof_upper_dvoretzky}.

\subsubsection{Comparison with \cite{bartlett_deep_2021}}\label{sec:comparision_BMR}

The following Theorem is taken from \cite[section 4.3.3, Theorem 4.13]{bartlett_deep_2021}. Theorem~\ref{theo:HMMT} is non-asymptotic, because they use the non-asymptotic version of \cite{el_karoui_spectrum_2010} developed in \cite{liang_just_2020}, see \cite[Lemma A.6]{bartlett_deep_2021}.

\begin{Theorem}[\cite{bartlett_deep_2021}]\label{theo:HMMT}
    Suppose that $d\sim N$. Let $\Sigma$ be a semi-positive definite matrix and let $X=\Sigma^{1/2}Z$ where $Z$ has i.i.d., variance $1$, centered sub-Gaussian coordinates. Denote by $M$ an upper bound on the sub-Gaussian norm of the coordinates of $Z$, on $\norm{\Sigma}_{\text{op}}$, and on $(1/d)\sum_{i=1}^d \sigma_i^{-1}(\Sigma)$. Assume that $k$ is a continuous function over $\bR$ and is smooth in a neighborhood of $0$ with $k(0),k'(0)>0$,  and let $K(\vx,\vy)=k\left(\left\langle \vx,\vy\right\rangle/d\right)$. Let $Y = f^*(X) + \xi$, where $\norm{f^*}_{L_{4+\eta}}\leq M$ for some $\eta>0$ and $\xi\sim\cN(0,\sigma_\xi^2)$. Denote $\bbeta_0 = \Sigma^{-1}\bE\left[Xf^*(X)\right]$ and $\lambda_*>0$ be the unique positive solution of
    \begin{align*}
        N\left(1-\frac{\gamma}{\lambda_*}\right) = \Tr\left(\Sigma\left(\Sigma+\lambda_* I \right)^{-1}\right),\mbox{ where }\gamma = k\left(\frac{\Tr\left(\Sigma\right)}{d}\right)-k(0)-k'(0)\left( \frac{\Tr\left(\Sigma\right)}{d} \right).
    \end{align*}Then, there exist absolute constants $c,C$ such that with probability at least $1-CN^{-1/4}$, the bias term $\mathrm{bias} = \big(\bE_X[(f^*(X)-\bE_Y \hat f_0(X))^2\big|\bX]\big)^{1/2}$ and the variance term $\mathrm{var} = \bE_{X,\xi_1,\cdots,\xi_N}[(\hat f_0(X)-\bE_Y \hat f_0(X))^2\big| \bX]$ satisfies
    \begin{align*}
        \left|\mathrm{bias}^2 - \cB(\Sigma,\bbeta_0) - \norm{P_{>1}f^*}_{L_2}^2(1+\cV(\Sigma)) \right|\lesssim N^{-c},\quad \left|\mathrm{var} - \sigma_\xi^2\cV(\Sigma)\right|\lesssim N^{-c},
    \end{align*}where
    \begin{align*}
        \cV(\Sigma) = \frac{\Tr\left( \Sigma^2\left( \Sigma + \lambda_* I \right)^{-2}\right) }{ N - \Tr\left( \Sigma^2\left( \Sigma + \lambda_* I \right)^{-2}\right)},\quad \cB(\Sigma,\bbeta_0) = \frac{\lambda_*^2\left\langle \bbeta_0, \Sigma^2\left( \Sigma + \lambda_* I \right)^{-2}\Sigma\bbeta_0\right\rangle}{ 1-N^{-1}\Tr\left( \Sigma^2\left( \Sigma + \lambda_* I \right)^{-2}\right) }.
    \end{align*}
\end{Theorem}

\begin{itemize}
    \item In the scenario where $k(t)=t$ (in which KRR reduces to a linear ridge regression), see \cite[Corollary 4.14]{bartlett_deep_2021}, Theorem~\ref{theo:HMMT} provides a unified bound of the estimation error in the context of linear regression \cite[section 4.3.3]{bartlett_deep_2021}, as discussed in the works of \cite{tsigler_benign_2023,lecue_geometrical_2022}. Nevertheless, the estimation error is exacerbated due to the inclusion of an additional term $N^{-c}$. Moreover, there is a deterioration of the probability deviation. In contrast, our Theorem~\ref{theo:upper_KRR} demonstrates that our estimation error is equivalent to that of the Gaussian case. In other words, we prove that the kernel ridge regression has the (one-sided) Gaussian Equivalence Property, see Section~\ref{sec:gaussian_equivalence} for more details. Moreover, if we substitute Proposition~\ref{prop:IP} and Theorem~\ref{theo:DM_RKHS} with the Gaussian counterparts proposed in \cite{lecue_geometrical_2022}, we achieve the optimal probability deviation in the linear case under a Gaussian assumption concerning the design vector and noise.
    \item Theorem~\ref{theo:HMMT} is limited to the case where $d\sim N$, and hence it is unable to demonstrate the occurrence of the multiple descent of the upper bound for the estimation error of KRR. This phenomenon will be established as a corollary of our Theorem~\ref{theo:upper_KRR} in subsequent sections of this work. It is crucial to observe that Theorem~\ref{theo:HMMT} highlights the learning of solely the linear approximation of $f^*$, whereas in the context of multiple descent, $\hat f_0$ encompasses the higher-degree approximation of $f^*$, hence resulting in a significant decrease in the estimation error.
    \item The fixed point $\lambda_*$ mentioned in Theorem~\ref{theo:HMMT} presents challenges in terms of computation. The analysis of the two terms $\cV(\Sigma)$ and $\cB(\Sigma,\bbeta_0)$ is also challenging for this reason. In contrast, Theorem~\ref{theo:upper_KRR} only relies on the spectrum of $\Gamma$ as well as the decomposition of $f^*$ in its eigenbasis, thereby rendering it more feasible in practical applications. The reader is directed to Section~\ref{sec:multiple_descent} for specific examples.
\end{itemize}

\subsubsection{Previous results on multiple descent}\label{sec:previous_results_multiple_descent}

\begin{Theorem}[\cite{liang_multiple_2020}]\label{theo:LRZ20}
    Let $\iota\in\bN_+$ and consider $d = N^\alpha$ where $(\iota+1)^{-1}\leq \alpha < \iota^{-1}$. Consider a general function $h(t)=\sum_{i=0}^\infty\alpha_i t^i$ with corresponding Taylor coefficients $\{\alpha_i\}_{i\in\bN}$ and a kernel $K$ such that \eqref{eq:polynomial_inner_product_kernel} holds. Suppose the first bullet of Assumption~\ref{assumption:polynomial_feature_sub_gaussian} below holds. Suppose $Y = f^*(X) + \xi$, where $f^*\in\cH$, $\xi$ is independent of $X$, and the variance of \(\xi\) is \(\sigma_{\xi}^2\). Suppose there exists an absolute constant $C>0$ such that $\norm{f^*}_{L_4(\mu)}\leq C$. Then
    \begin{enumerate}
        \item if $\alpha_1,\cdots,\alpha_\iota>0$ and there exists $\iota'\geq 2\iota+3$ such that $\alpha_{\iota'}>0$. Suppose $d^\iota\log{d}\lesssim N\lesssim d^{\iota+1}$, then with constant probability ,
        \begin{align*}
            \norm{\hat f_0 - f^*}_{L_2}^2 \lesssim \frac{d^\iota}{N} + \frac{N}{d^{\iota+1}},
        \end{align*}
        \item if for some $\iota>0$, $\alpha_1,\cdots,\alpha_\iota>0$ and for all $\iota'>\iota$, $\alpha_{\iota'}=0$. Suppose $d^\iota\log{d}\lesssim N$. Then with constant probability,
        \begin{align*}
            \norm{\hat f_0 - f^*}_{L_2}^2 \lesssim \frac{d^\iota}{N} + \frac{1}{N}.
        \end{align*}
    \end{enumerate}
\end{Theorem}

Let us provide some comments on this theorem.
\begin{itemize}
    \item The prediction for the interval of $N$, that is, $d^\iota\log{d}\lesssim N\lesssim d^{\iota+1}$, where multiple descent happens provided by Theorem~\ref{theo:LRZ20} is not sharp, as proved by our Proposition~\ref{prop:multiple_descent} which provides a sharp interval of $N$, that is, $d^\iota\lesssim N\lesssim d^{\iota+1}$.
    \item The dependence of the estimation error on $\sigma_\xi$ and $f^*$ is not explicitly stated in Theorem~\ref{theo:LRZ20}. On the contrary, our Proposition~\ref{prop:multiple_descent} stated the explicit dependence on $\sigma_\xi$ and $f^*$.
    \item The probability deviation in Theorem~\ref{theo:LRZ20} is a constant probability, which is because Proposition 4 in \cite{liang_multiple_2020} only provides a constant-level probability deviation (that is, this deviation does not converge to $0$ as $N$ increases). Our Proposition~\ref{prop:multiple_descent} provides a better probability deviation, which converges to $0$ as $N$ increases.
\end{itemize}

Another result on the multiple descent phenomenon comes from \cite{ghorbani_linearized_2021}. 
We emphasize that the following result, as well as other results such as in \cite{mei_generalization_2022,misiakiewicz_spectrum_2022,xiao_precise_2022}, holds as $d\to\infty$.

\begin{Theorem}[Theorem 4 of \cite{ghorbani_linearized_2021}]\label{theo:MMM_21}
    Let $\{f_d^*\in L_2(\Omega_d,\mu_d)\}_{d\in\bN_+}$ be a sequence of target functions. Let $(X_i)_{i\in[N]}\sim \mu_d$ be a sequence of i.i.d.\ random vectors, which are uniformly distributed over $\sqrt{d}S_2^{d-1}$. Let $Y_i = f_d^*(X_i)+\xi_i$ where $\xi_i\sim \cN(0,\sigma_\xi^2)$ for some $\sigma_\xi>0$. For some $\iota\in\bN_+$ and $\delta_0>0$, assume that $\omega_d\left(d^\iota\log{d}\right)\leq N\leq O_d\left(d^{\iota+1-\delta_0}\right)$
    .
    Let $\{h_d\}_{d\in\bN_+}$, $h_d:\bR\to\bR$ be a sequence of real-valued functions such that
    \begin{enumerate}
        \item $h_d(\cdot/\sqrt{d})\in L_2([-\sqrt{d},\sqrt{d}],\mu_{d-1}^1)$ where $\mu_{d-1}^1$ is the distribution of $\left\langle X_1,\ve_1\right\rangle$.
        \item Let $u(\iota)=\sum_{l\leq\iota}d^l$. There exists a constant $c_\iota>0$ such that $\lambda_*(d,\iota) := d^\iota \min_{k\leq \iota}\sigma_{u(\iota)}(\Gamma)\geq c_\iota \Tr\left(\Gamma_{u(\iota)+1:\infty}\right)$.
    \end{enumerate}Suppose $\bE h_d(X) = \bE f_d^*(X) = 0$. Denote the KRR with tuning parameter $\lambda$ by $\hat f_{d,\lambda}$. Then for any $0<\lambda<\lambda_*(d,\iota)$,
    \begin{align*}
    \left| \norm{ f_d^* - \hat f_{d,\lambda} }_{L_2(\mu_d)}^2 - \norm{P_{>\iota}f_d^*}_{L_2(\mu_d)}^2 \right|=o_{d,\bP}(1)\left(\norm{f_d^*}_{L_2(\mu_d)}^2 + \sigma_\xi^2\right).
    \end{align*}
\end{Theorem}

Theorem~\ref{theo:MMM_21} requires that the design vector $X$ is uniformly distributed over a Euclidean sphere $\sqrt{d}S_2^{d-1}$, which is very restrictive. On the contrary, our Proposition~\ref{prop:multiple_descent} holds for sub-Gaussian design. Moreover,  the validity of Theorem~\ref{theo:MMM_21} hinges on the assumption that $d$ tends towards infinity. As we saw at the beginning of Section~\ref{sec:multiple_descent}, this assumption is not convenient (it requires to consider a sequence of target functions $f_d^*$). Our result holds for fixed $d$.
The multiple descent appears when $\omega_d\left(d^\iota\log{d}\right)\leq N\leq O_d\left(d^{\iota+1-\delta_0}\right)$. On the contrary, we only need $Cd^\iota\leq N \leq cd^{\iota+1}$ in Proposition~\ref{prop:multiple_descent}.

Recently, \cite{misiakiewicz_spectrum_2022} proves multiple descent when the target function $f_d^*$ is a certain type of random function (see Assumption 2 in \cite{misiakiewicz_spectrum_2022}), when $X$ follows a uniform distribution over $\sqrt{d}S_2^{d-1}$, and when both the dimension $d$ and the sample size $N$ go to infinity. In a more precise manner, it can be shown that varying values of $\lfloor\log{N}/\log{d}\rfloor$ result in distinct descents. Within each descent, the function $\hat f_0$ acquires a polynomial approximation of degree $\lfloor\log{N}/\log{d}\rfloor$ of the true function $f^*$.

\subsubsection{Formal version of Proposition~\ref{prop:multiple_descent} (multiple descent)}\label{sec:formal_multiple_descent}

To align with the setups of  \cite{liang_multiple_2020} and \cite{ghorbani_linearized_2021}, we introduce two assumptions:

\begin{Assumption}\label{assumption:polynomial_feature_sub_gaussian}
\begin{enumerate}
    \item $X\in\bR^d$ is a random vector with i.i.d. zero-mean sub-Gaussian coordinates. Moreover, for any finite set $S$, the first coordinate (denoted as $x_1$) of $X$ satisfies $\bP(x_1\in S)<1$ and $\bP(x_1 = 0)>0$. We denote $\mu_1$ as the distribution of $x_1$, thus $\mu=\otimes_{j=1}^d \mu_1$.
    \item Suppose the noise $\xi$ is independent of $X$ with mean zero and variance $\sigma_\xi^2$. We assume that there exist some $\okappa{kappa_noise}>0$ and $r>4$ such that for all $i$'s,  $\norm{\xi}_{L_r}\leq \okappa{kappa_noise}\sigma_\xi$.
    
    \item There exists a polynomial function $h:t\in\bR\mapsto \sum_{i=0}^L \alpha_i t^i$ for some $\alpha_i\in\bR_+$, $\alpha_i\sim 1$ for all $i\in [L]\cup\{0\}$ and $1\leq L\log(2L)< \log{d}$, such that \eqref{eq:polynomial_inner_product_kernel} holds.

    \item $d\gtrsim_{\mu_1} (\log{N})^{2L}L^{6L}$, where $N$ is the number of samples.
\end{enumerate}
\end{Assumption}

\begin{Assumption}\label{assumption:polynomial_feature_gaussian}
\begin{enumerate}
    \item $X$ is distributed uniformly over $\sqrt{d}S_2^{d-1}$, whose distribution is denoted as $\mu$. $(X_i)_{i\in[N]}$ are i.i.d. copies of $X$.
    \item Suppose the noise $\xi$ is independent of $X$ with mean zero and variance $\sigma_\xi^2$. We assume that there exist some absolute constants $\okappa{kappa_noise}>0$ and $r>4$, such that for all $i$'s,  $\norm{\xi}_{L_r}\leq \okappa{kappa_noise}\sigma_\xi$.
    \item There exists a function $h:t\in\bR\mapsto \sum_{i=0}^L \alpha_i t^i$ for some $\alpha_i\in\bR_+$, $\alpha_i\sim 1$ for all $i\in [L]\cup\{0\}$, $L!=O(d)$ such that \eqref{eq:polynomial_inner_product_kernel} holds. Moreover, $\bE h(G)\mathrm{He}_\iota(G)\neq 0$ for all $\iota\in[L]$, where $G$ is a standard Gaussian random variable.
\end{enumerate}
\end{Assumption}

Applying Theorem~\ref{theo:upper_KRR}, we obtain the following proposition.
\begin{Proposition}
    \label{prop:multiple_descent}
    There exist absolute constants $c$, $\oC{C_N_lower},\oC{C_multiple_lower},\oC{C_noise},\oc{c_multiple_upper}$, $\oc{c_P_DMU}$, $\oc{c_P_RIP}$ and $\oc{c_P_bX_f_star}$ such that the following holds. Suppose $N\geq \oC{C_N_lower}$.
    \begin{itemize} 
        \item Grant Assumption~\ref{assumption:polynomial_feature_sub_gaussian}. For any $0\leq \iota\leq L-1$, $r_{0,k}^*$ for some well-chosen $k$ (see Section~\ref{sec:proof_multiple_descent}), is equivalent to
        \begin{align*}
            \max\left\{ \sigma_\xi\sqrt{\frac{d^\iota}{N}}, \frac{1}{N}\norm{\Gamma_{\leq \iota}^{-1/2}f_{\leq \iota}^*}_{\cH}, \norm{\Gamma_{>\iota}^{1/2}f_{>\iota}^*}_{\cH}, \sigma_\xi\sqrt{\frac{N}{d^{\iota+1}}} \right\},
        \end{align*}when $\oC{C_multiple_lower} d^\iota\leq N \leq \oc{c_multiple_upper}d^{\iota+1}$. Moreover, with probability at least 
        \begin{align}\label{eq:prob_devi_multiple_descent_sub_Gaussian}
            1-\gamma - \ogamma{gamma_RIP} - \ogamma{gamma_DMU_L2} - 2\exp\left(-cd^\iota\right) - \frac{1}{N} - c\left(\frac{\oC{C_noise}}{N}\right)^{r/4},
        \end{align}we have $\norm{\hat f_0 - f^*}_{L_2}\lesssim r_{0,k}^*$ and $\left|\|\hat f_0 - f^*\|_{L_2} - \norm{\Gamma_{>\iota}^{1/2}f_{>\iota}^*}_{\cH}\right|\lesssim r_{0,k}^*$, where the value of $\gamma$, $\ogamma{gamma_RIP}$ and $\ogamma{gamma_DMU_L2}$ for each case can be found in Section~\ref{sec:proof_diagonal_concentration}.

        \item Under Assumption~\ref{assumption:polynomial_feature_gaussian}, for every $0\leq\iota\leq L-1$, $r_{0,k}^*$ for some well-chosen $k$ (see Section~\ref{sec:proof_multiple_descent}) is equivalent to
            \begin{align}\label{eq:def_r_star_multiple_descent}
                \max\left\{ \sigma_\xi\sqrt{\frac{d^\iota}{N}}, \frac{1}{N}\norm{\Gamma_{\leq\iota}^{-1/2}f_{\leq\iota}^*}_{\cH}, \norm{\Gamma_{>\iota}^{1/2}f_{>\iota}^*}_{\cH}, \sigma_\xi\sqrt{\frac{N}{d^{\iota+1}}} \right\},
            \end{align}
            when $\oC{C_multiple_lower} d^\iota \leq N\leq \oc{c_multiple_upper} d^{\iota+1}$. Moreover, with probability at least 
            \begin{align*}
                1- 2\exp\left(-cd^\iota\right) - \frac{1}{N}- c\left(\frac{\oC{C_noise}}{N}\right)^{r/4},
            \end{align*}we have $\norm{\hat f_0 - f^*}_{L_2}\lesssim r_{0,k}^*$  and $\left|\|\hat f_0 - f^*\|_{L_2} - \norm{\Gamma_{>\iota}^{1/2}f_{>\iota}^*}_{\cH}\right|\lesssim r_{0,k}^*$.
        \item Under the first 2 items of Assumption~\ref{assumption:polynomial_feature_sub_gaussian}, if we further have that $\{f_d^{*}\in L_4(\bR^d,\mu)\}_{d\in\bN_+}$ is a sequence of target functions, such that $Y_i = f_d^{*}(X_i)+\xi_i$.
        Then when $N,d\to\infty$ with $d^\iota=o(N)$ and $N=o(d^{\iota+1})$, we have $\left|\|\hat f_0 - f_d^*\|_{L_2(\mu)}- \|\Gamma_{>\iota}^{1/2}(f_d)_{>\iota}^*\|_\cH\right| = o_{d,\bP}(1)(\sigma_\xi+\|f_d^*\|_\cH)$.
    \end{itemize}

Moreover, the above results are still valid when $\hat f_0$ is replaced by $\hat f_\lambda$ where $\lambda\lesssim 1$.
\end{Proposition}

\begin{Remark}
    Under Assumption~\ref{assumption:polynomial_feature_gaussian} and Assumption~\ref{assumption:polynomial_feature_sub_gaussian}, we assume $\alpha_i > 0$. If there exists $\gamma < L-1$ such that $\alpha_{\gamma},\alpha_{\gamma+2} \neq 0$ but $\alpha_{\gamma+1} = 0$, then by recalculating \eqref{eq:quantities_sub_gaussian} and \eqref{eq:quantities_uniform} below, Proposition~\ref{prop:multiple_descent} still holds. However, the range for $N$ will change from $\oC{C_multiple_lower} d^\gamma \leq N \leq \oc{c_multiple_upper}d^{\gamma+1}$ to $\oC{C_multiple_lower} d^\gamma \leq N \leq \oc{c_multiple_upper}d^{\gamma+2}$. Moreover, in $r_{0,k}^*$, $d^{\gamma+1}$ will be replaced by $d^{\gamma+2}$, and $\norm{\Gamma_{>\gamma}^{1/2}f_{>\gamma}^*}_\cH$ will be replaced by $\norm{\Gamma_{>\gamma+1}^{1/2}f_{>\gamma+1}^*}_\cH$. Similarly, we can handle cases where several $\alpha_i$ are equal to 0.
\end{Remark}

\subsubsection{Conclusions for data-dependent conjugate kernel}\label{sec:formal_data_dependent}

\begin{Proposition}\label{prop:conjugate_kernel}
    Suppose Assumption~\ref{assumption:conjugate_kernel} holds. For some $\bar p_{DM}$, $\bar p_{DMU}$, $\bar p_{RIP}$ that are defined in Section~\ref{sec:proof_conjugate_kernel} later and $\bar p_{\xi}$ defined in \eqref{eq:def_bar_p_xi} later. If $k\leq\oc{c_kappa_RIP}N_2$, then, for any $\lambda\geq 0$, with probability at least
    \begin{align*}
        &1-\bar p_{DM}-\bar p_{DMU}-\bar p_{RIP}-\exp\left(-\left|J_1\right|-\frac{N_2\sum_{j\in J_2}\sigma_j}{2\kappa_{DM}\left(4\lambda + \Tr\left(\Gamma_{k+1:m}\right)\right)}\right)\\
        &-\exp\left(-\oc{c_P_Bernstein}N_2\min\left\{\left(\frac{m\norm{f_{k+1:m}^*}_{L_2}^2}{\lambda_\sigma^2 B_c^2 \norm{f_{k+1:m}^*}_\cH^2}\right)^2, \frac{m\norm{f_{k+1:m}^*}_{L_2}^2}{\lambda_\sigma^2 B_c^2 \norm{f_{k+1:m}^*}_\cH^2} \right\}\right)-\bar p_{\xi},
    \end{align*}the same conclusions as in Theorem~\ref{theo:upper_KRR} hold.
\end{Proposition}
The proof of Proposition~\ref{prop:conjugate_kernel} may be found in Section~\ref{sec:proof_conjugate_kernel}.

\subsubsection{An example of feature learning}\label{sec:example_feature_learning}

The following assumption is from \cite{bietti_learning_2022}. It is easy to verify that the activation functions $\sigma: t\in\bR \mapsto \tanh(t)$ and $\sigma: t\in\bR \mapsto 1/(1+\exp(-t)) - 1/2$ satisfy \cite[Assumption 5.2]{bietti_learning_2022}, with $s=1$ in \cite[Assumption 5.2]{bietti_learning_2022}.

\begin{Assumption}\label{assumption:Bietti}
    There exists an activation function $\sigma$ and a vector $\vw^*\in S_2^{d-1}$ such that $y = f^*(\vx)+\xi$, where $f^*(\vx) = a^*\sigma(\left<\vw^*,\vx\right>)$,  $a^*\in\bR\backslash\{0\}$, $\xi\sim\cN(0,\sigma_\xi^2)$ for some $\sigma_\xi>0$, and is independent with $\vx$, and $\vx\sim\cN(\vzero,I_d)$. Suppose $\sigma:t\in\bR\mapsto \tanh(t)$ or $\sigma:t\in\bR\mapsto 1/(1+\exp(-t))-1/2$.
    
Let $X_1,\cdots,X_{N_1+N_2}$ be independent copies of $\vx$, and $Y_1,\cdots,Y_{N_1+N_2}$ be independent copies of $y$. Let $\varepsilon_j$'s are Rademacher random variables independent with $X_i$ and $Y_i$'s, $b_j\sim\cN(0,\tau^2)$ for some $\tau>0$ independent with other random variables. Define
\begin{align*}
    L_N(\va,\vw)=\frac{1}{N_1}\sum_{i=1}^{N_1}\left(\left<\va, \Phi(\left<\vw,X_i\right>)\right>-Y_i\right)^2+\eta\norm{\va}_2^2,\mbox{ where }\Phi(u)=\frac{1}{\sqrt{m}}(\sigma(\varepsilon_j u - b_j))_{j\in[m]}.
\end{align*}
Suppose $N_1\gtrsim m+d$. Let $\eta=O(1)$ and $\eta = \Omega(\sqrt{\Delta_{crit}})$ where 
\begin{align*}
    \Delta_{crit}:=\max\bigg\{\sqrt{(d+m)/N_1},\left(\frac{d^2}{N_1}\right)^2\bigg\}.
\end{align*}
Let $N_1=\tilde\Omega(\max\{ \frac{(d+m)}{\eta^4}, \frac{d^2}{\eta^2} \})$. Let $m=\Omega(\frac{1}{\eta}\log(\frac{1}{\eta\delta}))$ and $m=\tilde O(\eta\Delta_{crit}^{-1})$. Let $m_0 = \Theta(\log(\frac{1}{\delta}))$ and $\rho = \Theta(\sqrt{m}m_0^{-\frac{3}{2}}( \tau^2 + \frac{\eta m}{m_0} )^{-1})$. Let $T_0 = \tilde\Theta(1)$, and $T_1 = \tilde\Theta(\frac{\eta^4 N}{d+m})$.
\end{Assumption}
In Assumption~\ref{assumption:Bietti}, we assume that $f^*$ is a single-neuron function, a class of target functions that has been extensively studied in the field of deep learning theory, as seen in works such as \cite{mei_landscape_2018,oymak_overparameterized_2019,vardi_implicit_2021,tan_online_2023} and references therein.

Given $W(T)$ be some vector in $\bR^d$, we construct a weight matrix $W$ by performing Gram-Schmidt process under the inner product $\ell_2^d$ starting from $W(T)$. That is, we obtain $\vw_1,\cdots,\vw_m$ such that for any $i,j\in[m]$, $\left<\vw_i,\vw_j\right>=\delta_{ij}$, and $\vw_1=W(T)$. Take $\ndelta\label{delta_feature_learning}<1$ such that $\|\sigma(\odelta{delta_feature_learning}\cdot)\|_{L_2(\gamma_d)}^2\sim \frac{1}{m}$. Such $\odelta{delta_feature_learning}$ always exists since $\delta\in\{0\leq\delta<1\}\mapsto \|\sigma(\delta\cdot)\|_{L_2(\gamma)}$ is a continuously increasing function, and $\|\sigma(0\cdot)\|_{L_2(\gamma)}=0$. We let $W = [\vw_1|\odelta{delta_feature_learning}\vw_2|\odelta{delta_feature_learning}\vw_3|\cdots|\odelta{delta_feature_learning}\vw_m]^\top=:[W_1|W_2|\cdots|W_m]^\top\in\bR^{m\times d}$.

\begin{Proposition}\label{prop:feature_learning}
Grant Assumption~\ref{assumption:Bietti}. Initialize $\vw(0)$ to be uniformly distributed over $S_2^{d-1}$, and $\va(0)$ uniformly distributed over $\{\va\in\rho S_2^{m-1}: \norm{\va}_0 = m_0\}$, where $\norm{\va}_0 = \left|\{j\in [m]: a_j\neq 0\}\right|$.
Run the following gradient flow algorithm up to time $T=T_0+T_1$:
\begin{align*}
    \dot{\va}(t) = -\1_{\{t>T_0\}}\nabla_{\va}L_N(\va,\vw),\mbox{ and } \dot{\vw}(t) = -\Pi_{\vw^\perp} \nabla_{\vw} L_N(\va,\vw),
\end{align*}where $\Pi_{\vw^\perp}:\vv\in\bR^d\mapsto \vv-\left<\vv,\vw\right>\vw\in\bR^d$.

Let $W(T) = \mathrm{sgn}(\left<\vw(T),\vw^*\right>)\vw(T)$ where where $\mathrm{sgn}(t)=1 $ if $t\geq 0$ and $-1$ otherwise. Let $\phi:\vx\in\bR^d\mapsto \frac{1}{\sqrt{m}}\sigma(W\vx)$. Denote the data-dependent conjugate RKHS by $(\cH_T,\norm{\cdot}_{\cH_T})$. For any $\delta>0$, with probability at least $\frac{1}{2}-\delta$, the minimum $\norm{\cdot}_{\cH_T}$ norm interpolant estimator (interpolating $(X_i,Y_i)_{i=N_1+1}^{N_1+N_2}$) $\hat f_0$ constructed on sample $(X_i,Y_i)_{i=N_1+1}^{N_1+N_2}$ satisfies
\begin{align*}
    \norm{\hat f_0 - f^*}_{L_2}&\lesssim \left|a^*\right|\norm{\sigma}_{Lip}\tilde O\left(\eta^{-2}\max\left\{\sqrt{\frac{d+m}{N_1}},\frac{d^2}{N_1}\right\}\right)\\
    &+ \frac{\sigma_\xi}{\sqrt{N}} + \sigma_\xi\sqrt{\frac{N_2}{m}} + \frac{\left|a^*\right|}{\|\sigma\|_{L_2(\gamma_d)}}\frac{1}{N_2}.
\end{align*}
\end{Proposition}

The proof of Proposition~\ref{prop:feature_learning} may be found in Section~\ref{sec:proof_feature_learning}.

\subsubsection{A matching lower bound in Gaussian linear regression}\label{sec:matcing_lower_bound}

\begin{Proposition}\label{prop:linear}
There are absolute constants $\nc\label{c_linear_1}$, $\nc\label{c_linear_2}$, and $\nC\label{C_linear_3}$ such that the following holds. Let $G$ be a mean $0$ Gaussian random vector in $\ell_2$ with covariance matrix $\Gamma$ and let $\xi\sim\cN(0, \sigma_\xi^2)$ be independent of $G$. Let $Y = f^*(G)+\xi$ and $(X_i, Y_i)_{i=1}^N$ be $N$ i.i.d. copies of $(X,Y)$. We write $f^*(G):=\left<\bbeta^*,G\right>$ for some unknown $\bbeta^*\in\ell_2$. We identify $\cH$ with $\ell_2$ and so $f^*$ is identified with $\bbeta^*$. For any $\lambda\geq 0$, let $\hat\bbeta_\lambda\in\argmin\left(\sum_{i=1}^N\left(Y_i - \left\langle \bbeta,G_i\right\rangle\right)^2 + \lambda\norm{\bbeta}_2^2\right)$.
We assume that
\begin{enumerate}
    \item either there exists $k\leq \oc{c_RIP} N$,
    \item or $k$ is not necessarily smaller than $N$ but \eqref{eq:extra_condition_k>N} and \eqref{eq:RIP_condition_k>N} hold,
\end{enumerate}
such that $N\leq \kappa_{DM}d_\lambda^*(\Gamma_{k+1:\infty}^{-1/2}B_2^\infty)$ where $B_2^\infty$ is the unit ball of $\ell_2$. The following then holds for all such $k$'s and $\lambda$'s.
With probability at least
\begin{align*}
    1-\oc{c_linear_1}\exp\left(-\oc{c_linear_2}\left(|J_1| + N \left(\sum_{j\in J_2}\sigma_j\right)/\left(\Tr(\Gamma_{k+1:\infty})\right)\right)\right),
\end{align*}we have
\begin{equation*}
\norm{\Gamma^{1/2}(\hat\bbeta_\lambda-\bbeta^*)}_{2}\leq \oC{C_linear_3} r_{\lambda,k}^*,
\end{equation*}where $J_1, J_2$ have been defined in Section~\ref{sec:notations}.
Moreover, given $b>\max\{\oC{C_lower_kappa_DM}\kappa_{DM}^{-1}, 2\oc{c_lower_7}^{-1}\}$ for some absolute constants $\nC\label{C_lower_kappa_DM}$ and $\oc{c_lower_7}$, we define
\begin{align}\label{eq:def_k_b_star}
    k_{b,\lambda}^* = \min\left(k\in\bN:\, \frac{\lambda + \Tr\left(\Gamma_{k+1:\infty}\right)}{\norm{\Gamma_{k+1:\infty}}_{\text{op}}}\geq bN\right),
\end{align}where the infimum of empty set is defined as $\infty$.
If $N$ is large enough such that $k_{b,\lambda}^*\leq N$, then 
\begin{align*}
    \bE \norm{\Gamma^{1/2}(\hat\bbeta_\lambda-\bbeta^*)}_2\gtrsim r_{\lambda,k_{b,\lambda}^*}^*.
\end{align*}
\end{Proposition}

\subsubsection{Estimation error of KRR with smooth kernel in non-asymptotic regime}\label{sec:smooth_omitted}

In this paragraph, we apply our general bounds from Theorem~\ref{theo:upper_KRR} and Theorem~\ref{theo:main_upper_k>N}. Recall the definition on $J_1$ and $J_2$ in \eqref{eq:def_J1}.

\begin{Proposition}\label{prop:KRR_smooth_upper}
Suppose that $k\in\bN$ and $\lambda\geq 0$ satisfy $N\leq\oc{c_kappa_DM}\kappa_{DM}d_\lambda^*(\Gamma_{k+1:\infty}^{-1/2}B_\cH)$, $\sum_{j\in J_2}\sigma_j \leq \kappa_{DM}(4\lambda + \Tr\left(\Gamma_{k+1:\infty}\right))\left(1-\frac{\left|J_1\right|}{N}\right)$, and $\kappa_{DM}\left(4\lambda + \Tr\left(\Gamma_{k+1:\infty}\right)\right)\geq N\left(R_N^*(\oc{c_kappa_RIP})\right)^2$.
Under Assumption~\ref{assumption:smooth}, if $(2k)^{\frac{1}{\log{d}}}\lesssim N/\log^8(N)$. Then, with constant probability, $\norm{ \hat f_\lambda - f^*}_{L_2}$ is smaller than (up to an absolute multiplicative constant) 
    \begin{align}\label{eq:result_KRR_smooth}
         \begin{aligned}
             \log^3(N)&\max\bigg\{ \sigma_\xi\sqrt{\frac{\left|J_1\right|}{N}}, \sigma_\xi\sqrt{\frac{\sum_{j\in J_2}\sigma_j}{4\lambda + k^{1-\alpha}}}, \norm{ \Gamma_{k+1:\infty}^{1/2}f_{k+1:\infty}^* }_\cH,\\
             &{\norm{\tilde\Gamma_{1,\mathrm{thre}}^{-1/2}f_{1:k}^*}_\cH\frac{\lambda + k^{1-\alpha}}{N}}, \sigma_\xi\frac{\sqrt{Nk^{1-2\alpha}}}{\lambda + k^{1-\alpha}} \bigg\}.
         \end{aligned}
    \end{align}
\end{Proposition}

The proof of Proposition~\ref{prop:KRR_smooth_upper} can be found in Section~\ref{sec:proof_smooth_kernel}. The constant probability in Proposition~\ref{prop:KRR_smooth_upper} arises because we aim to allow for the case where $k\gtrsim N$. If we assume $k\lesssim N$ (Proposition~\ref{prop:KRR_smooth_applied} below precisely satisfies this assumption), then we can use the conclusions from \cite[Lemma 3]{mcrae_harmless_2022} or \cite[Lemma 2]{barzilai_generalization_2024} instead of Proposition~\ref{prop:RIP}, thereby obtaining the result with high probability.

\subsubsection{Source condition in the language of interpolation space}\label{sec:source_condition_omitted}

Below, we present the description of this condition in the language of interpolation spaces, which is taken from \cite[section 2]{zhang_optimality_2023}. For any $s\geq 0$, let $L^s:f\in L_2(\mu)\mapsto \sum_{j=1}^\infty \sigma_j^s\left<f,f_j\right>_{L_2}f_j$, where we recall that $(f_j)_{j=1}^\infty\in L_2(\mu)$ form an ONB of $L_2(\mu)$. Let $[\cH]^s$ be a Hilbert space defined as
\begin{align*}
    [\cH]^s := \left\{\sum_{j=1}^\infty a_j \sigma_j^{s/2}f_j:\, (a_j)_{j=1}^\infty\in \ell_2 \right\},\quad \left<f,g\right>_{[\cH]^s}:=\left<L^{-s/2}f, L^{-s/2}g\right>_{L^2(\mu)}.
\end{align*}
We have:
\begin{align}\label{eq:equivalent_source_condition}
    \norm{f^*}_{[\cH]^s}^2 & = \norm{\sum_{j=1}^\infty a_j \sigma_j^{1/2}\sigma_j^{-s/2}f_j}_{L^2(\mu)}^2 = \sum_{j=1}^\infty a_j^2 \sigma_j^{1-s}\leq 1.
\end{align}

\section{Proofs}

We first provide some definitions. Define
\begin{align}\label{eq:def_square_1}
    \begin{aligned}
        \square &= 
        \max\bigg\{ \sigma_\xi\sqrt{\frac{\Tr\left(\Gamma_{1:k}\right)}{4\lambda + \Tr\left(\Gamma_{k+1:\infty}\right)}}, \sqrt{\frac{\sigma_1 N}{4\lambda + \Tr\left(\Gamma_{k+1:\infty}\right)}}\norm{\Gamma_{k+1:\infty}^{1/2}f_{k+1:\infty}^*}_\cH,\\
        &\norm{f_{1:k}^*}_\cH\sqrt{\frac{4\lambda + \Tr\left(\Gamma_{k+1:\infty}\right)}{N}} \bigg\},
    \end{aligned}
\end{align}if $\sigma_1 N \leq \kappa_{DM}(4\lambda + \Tr\left(\Gamma_{k+1:\infty}\right))$; and
\begin{align}\label{eq:def_square_2}
    \begin{aligned}
        \square &= \max\bigg\{ \sigma_\xi\sqrt{\frac{\left|J_1\right|}{N}}, \sigma_\xi\sqrt{\frac{\sum_{j\in J_2}\sigma_j}{  4\lambda + \Tr\left(\Gamma_{k+1:\infty}\right)}}, \norm{\Gamma_{k+1:\infty}^{1/2}f_{k+1:\infty}^*}_\cH,\\
        &{\norm{\tilde\Gamma_{1,\mathrm{thre}}^{-1/2}f_{1:k}^*}_\cH \frac{2\lambda + 3\Tr\left(\Gamma_{k+1:\infty}\right)}{N} } \bigg\},
    \end{aligned}
\end{align}if $\sigma_1 N > \kappa_{DM}(4\lambda + \Tr\left(\Gamma_{k+1:\infty}\right))$. Let $\triangle = \square\sqrt{N}/\sqrt{\kappa_{DM}\left(4\lambda + \Tr\left(\Gamma_{k+1:\infty}\right)\right)}$.

For any $\lambda\geq 0$ and $k\in\bN$, define
\begin{align}\label{eq:def_rate}
    r_{\lambda,k}^* := \max\left\{\square, \sigma_\xi\frac{\sqrt{N\Tr\left(\Gamma_{k+1:\infty}^2\right)}}{\lambda + \Tr\left(\Gamma_{k+1:\infty}\right)}\right\}.
\end{align}

Define
\begin{align}\label{eq:def_tilde_Gamma}
    \tilde\Gamma_{1:k}^{1/2} = \sum_{j=1}^k \max\left(\frac{\sqrt{\sigma_j}}{\square}, \frac{1}{\triangle}\right)\varphi_j\otimes\varphi_j.
\end{align}

For the sake of simplicity, we denote by $A:\bR^N\to\cH$, a random matrix with i.i.d. column vectors denoted by $\phi(X_1),\cdots,\phi(X_N)$: $A=\left[\phi(X_1)|\cdots|\phi(X_N)\right]$ such that for any $\vlambda\in\bR^N$, $A\vlambda = \sum_{i=1}^N\lambda_i\phi(X_i)$. In other words, $A$ is the adjoint of $\bX_\phi$, i.e., $A = \bX_\phi^\top$. We denote $\ell^* = \sqrt{\bE\norm{\phi(X)}_\cH^2} = \sqrt{\bE K(X,X)} = \sqrt{\Tr\left(\Gamma\right)}$.

\subsection{Proof of Theorem~\ref{theo:upper_KRR} (the $k\lesssim N$ case)}\label{sec:proof_main_upper_k<N}

In this section, we establish the proof of Theorem~\ref{theo:upper_KRR}. The proof is generally divided into two main parts: the Stochastic Argument and the Deterministic Argument. In the following subsection, we commence with the Stochastic Argument.

\subsubsection{Stochastic event behind Theorem~\ref{theo:upper_KRR}.} Let $\oC{C_DMU},\nC\label{C_bX_f_star},\nC\label{C_sum_Gamma_phi},\nc\label{c_RIP_lower}$, and $\nC\label{C_RIP_upper}$ be absolute constants. We denote by $\Omega_0$ the event which we have:
\begin{itemize}
  \item for all $\blambda\in\bR^N$,
  \begin{align}\label{eq:DM_applied}
      \begin{aligned}
          \left(\frac{1}{2}\Tr(\Gamma_{k+1:\infty})+\lambda\right)\norm{\blambda}_2&\leq \norm{\left(\bX_{\phi,k+1:\infty}\bX_{\phi,k+1:\infty}^\top+\lambda I\right)\blambda}_2\\
          &\leq \left(\frac{3}{2}\Tr(\Gamma_{k+1:\infty})+\lambda\right)\norm{\blambda}_2
      \end{aligned}
  \end{align}
  \item for all $f_{1:k}\in \cH_{1:k}$,
  \begin{align}\label{eq:RIP_applied}
      \oc{c_RIP_lower}\norm{\Gamma_{1:k}^{1/2}f_{1:k}}_\cH\leq (1/\sqrt{N})\norm{\bX_{\phi,1:k}f_{1:k}}_2\leq \oC{C_RIP_upper}\norm{\Gamma_{1:k}^{1/2}f_{1:k}}_\cH
  \end{align}
  \item for all $\blambda\in\bR^N$,
  \begin{align}\label{eq:DM_upper_applied}
      \norm{\Gamma_{k+1:\infty}^{1/2}\bX_{\phi,k+1:\infty}^\top\blambda}_\cH\leq \oC{C_DMU}\left(\sqrt{\Tr\left(\Gamma_{k+1:\infty}^2\right)} + \sqrt{N}\norm{\Gamma_{k+1:\infty}}_{\text{op}}\right)\norm{\blambda}_2
  \end{align}
  \item \begin{align}\label{eq:bX_f_star_applied}
      \norm{\bX_{\phi,k+1:\infty}f_{k+1:\infty}^*}_2\leq \oC{C_bX_f_star}\kappa\sqrt{N}\norm{\Gamma_{k+1:\infty}^{1/2}f_{k+1:\infty}^*}_\cH
  \end{align}
  \item \begin{align}\label{eq:Gamma_phi_applied}
      \sum_{i=1}^N \norm{\left(\Gamma_{k+1:\infty}^{1/2}\phi_{k+1:\infty}\right)(X_i)}_\cH^2\leq \oC{C_sum_Gamma_phi} N\Tr\left(\Gamma_{k+1:\infty}^2\right).
  \end{align}
\end{itemize}

By the definition of $\tilde\delta$, see \eqref{eq:def_tilde_delta}, there exists an absolute constant $\nc\label{c_kappa_DM}$ such that if $\delta^2,\bar\delta^2<\oc{c_kappa_DM}$, we have $\tilde\delta<1/2$. When $\lambda>\oC{C_comparision_trace_lambda}\Tr\left(\Gamma_{k+1:\infty}\right)$, one may replace the absolute constants $\frac{1}{2}$ and $\frac{3}{2}$ in \eqref{eq:DM_applied} by those in Theorem~\ref{theo:DM_RKHS}.  It follows from Theorem~\ref{theo:DM_RKHS}, Proposition~\ref{prop:upper_dvoretzky} and Proposition~\ref{prop:IP} that if $N\leq \oc{c_kappa_DM}\kappa_{DM}d_\lambda^*(\Gamma_{k+1:\infty}^{-1/2}B_\cH)$, then with probability larger than $1-\bar p_{RIP} - \bar p_{DM} - \bar p_{DMU},$
\eqref{eq:DM_applied}, \eqref{eq:RIP_applied} and \eqref{eq:DM_upper_applied} hold.

For \eqref{eq:bX_f_star_applied}, we use \cite[Lemma 3.2]{mendelson_upper_2016} on the $L_r$-norm of a sum of i.i.d. random variables to deduce the following result.

\begin{Lemma}\label{lem:moment-sum-variid-positive}
    There exist absolute constants $\nc\label{c_Men_1},\nc\label{c_Men_2},\nc\label{c_Men_3}$ such that the following holds. Let $1\leq r<q$, set $Z\in L_q$ and put $Z_1,\cdots,Z_N$ to be independent copies of $Z$. Fix $1\leq p\leq N$, let $j_0 = \lceil (\oc{c_Men_1} p)/\left( ((q/r)-1)\log{(4+eN/p)} \right) \rceil$ and $t>2$.
    If $j_0 = 1$ and $0<\beta<(q/r)-1$ then with probability at least $1-\oc{c_Men_3}t^{-q}N^{-\beta}$,
    \begin{align*}
        \left(\sum_{j=1}^N \left|Z_i\right|^r\right)^{1/r} \leq \oc{c_Men_2}\left(\frac{q}{q-(\beta+1)r}\right)^{1/r}t\norm{Z}_{L_q}N^{1/r}.
    \end{align*}
\end{Lemma}
Without loss of generality, we take $\oc{c_Men_2}>1$.
Let $Z_i = f_{k+1:\infty}^*(X_i)$, $r=2$, $p=1$, $q=4+\eps$ (where $\eps$ is from Assumption~\ref{assumption:DM_L4_L2}, Assumption~\ref{assumption:upper_dvoretzky} and Assumption~\ref{assumption:RIP}) and $\beta=1$ in Lemma~\ref{lem:moment-sum-variid-positive}, and by the fact that $\norm{f_{k+1:\infty}^*}_{L_{4+\eps}}\leq\kappa\norm{f_{k+1:\infty}^*}_{L_2}$, Lemma~\ref{lem:moment-sum-variid-positive} indicates that if $N\geq (e^{\oc{c_Men_1}}-4)/e\vee\oc{c_Men_3}^2$, then there exists absolute constant $\nc\label{c_P_bX_f_star}$ with probability at least $1-\oc{c_P_bX_f_star}/N$,
\begin{align*}
    \norm{\bX_{\phi,k+1:\infty}f_{k+1:\infty}^*}_2 &= \left(\sum_{i=1}^N \left(f_{k+1:\infty}^*(X_i)\right)^2 \right)^{1/2} \leq \oC{C_bX_f_star}\sqrt{N}\norm{f_{k+1:\infty}^*}_{L_{4+\eps}}\\
    &\leq \oC{C_bX_f_star}\kappa\sqrt{N}\norm{\Gamma_{k+1:\infty}^{1/2}f_{k+1:\infty}^*}_{\cH},
\end{align*}where $\oc{c_P_bX_f_star} = \oc{c_Men_3}$, $\oC{C_bX_f_star} = \oc{c_Men_2}\sqrt{(4+\eps)/\eps}.$

We are left with checking \eqref{eq:Gamma_phi_applied}. However, this is a simple consequence of Assumption~\ref{assumption:upper_dvoretzky}. By Assumption~\ref{assumption:upper_dvoretzky}, \eqref{eq:Gamma_phi_applied} holds with probability at least $1-\ogamma{gamma_DMU_L2}$ with constant $\oC{C_sum_Gamma_phi} = 1+ \odelta{delta_DMU_L2}.$


Combining the above probabilistic estimates, we have the following Proposition
\begin{Proposition}\label{prop:stochastic_argument_upper_bound}
    Suppose Assumption~\ref{assumption:DM_L4_L2}, Assumption~\ref{assumption:upper_dvoretzky} and Assumption~\ref{assumption:RIP} hold. There exist absolute constants $\oc{c_Men_1}$, $\oc{c_Men_3}$, $\oc{c_RIP}$, $\oc{c_kappa_DM}$ and $\oc{c_P_bX_f_star}$, such that if $N\geq (e^{\oc{c_Men_1}}-4)/e\vee\oc{c_Men_3}^2$, and if there exists $k\leq \oc{c_RIP} N,$ such that $N\leq   \oc{c_kappa_DM}\kappa_{DM}d_\lambda^*\left(\Gamma_{k+1:\infty}^{-1/2}B_\cH\right)$, then
\begin{align*}
    \bP(\Omega_0)&\geq 1 - \bar p_{RIP} - \bar p_{DM} - \bar p_{DMU} - \frac{ \oc{c_P_bX_f_star} }{N} - \ogamma{gamma_DMU_L2}.
\end{align*}
\end{Proposition}



We now place ourselves on the event $\Omega_0$ up to the end of the proof of Theorem~\ref{theo:upper_KRR}. All the remaining material does not rely on any stochastic arguments since they all have been collectd in $\Omega_0$.

\subsubsection{Decomposition of $\hat f_\lambda$}

As in the linear situation, the KRR estimator $\hat f_\lambda$ is decomposed into two components: $\hat f_{1:k}\in \cH_{1:k}=\Span\left(\varphi_j:1\leq j\leq k\right)$ and $\hat f_{k+1:\infty}\in \cH_{k+1:\infty}=\Span\left(\varphi_j:\, j>k\right)$. The two components have their own role in estimating $f^*$: $\hat f_{1:k}$ is used as a ridge estimator of $P_{1:k}f^*$ whereas $\hat f_{k+1:\infty}$ is used to absorb noise, thus is not expected to be a good estimator of $P_{k+1:\infty}f^*$.
\begin{Proposition}\label{prop:decomposition_hat_f}
    For any $k\in\bN_+$, the KRR defined by \eqref{eq:def_KRR} can be written as $\hat f_\lambda = \hat f_{1:k}+\hat f_{k+1:\infty}$, where
    \begin{align}
        \hat f_{1:k} \in \underset{f_{1:k}\in \cH_{1:k}}{\argmin}\left(\norm{Q\left(\vy - \bX_{\phi,1:k}f_{1:k}\right)}_\cH^2 + \norm{f_{1:k}}_\cH^2\right)\label{eq:objective_f_1_decompose},
    \end{align}and
    \begin{align}
    \hat f_{k+1:\infty} = \bX_{\phi,k+1:\infty}^\top \left(\bX_{\phi,k+1:\infty}\bX_{\phi,k+1:\infty}^\top  + \lambda I_N\right)^{-1} \left(\vy - \bX_{\phi,1:k}\hat f_{1:k}\right),\label{eq:objective_f_2_decompose}
    \end{align}
    where $Q:\bR^N\to \cH_{k+1:\infty}$ is a bounded linear operator such that $Q^\top Q = \left(\bX_{\phi,k+1:\infty}\bX_{\phi,k+1:\infty}^\top  + \lambda I_N\right)^{-1}.$ Such an operator exists because $Q^\top Q$ is semi positive-definite.
\end{Proposition}
\beginproof The empirical regularized loss functional is defined as $L:f\in\cH\mapsto\norm{\vy - \bX_{\phi}f}_2^2 + \lambda\norm{f}_\cH^2$. As $\hat f_\lambda$ is a minimizer of $L(f)$, we decompose $\hat f_\lambda=\hat f_{1:k}+\hat f_{k+1:\infty}$ and take the derivative of $L$ with respect to the canonical inner product on $\cH$ at $\hat f_{1:k}$ and $\hat f_{k+1:\infty}$ and set them to $0$, as a result, we obtain
\begin{align}
    \left(\bX_{\phi,1:k}^\top\bX_{\phi,1:k}+\lambda I\right)\hat f_{1:k} + \bX_{\phi,1:k}^\top\bX_{\phi,k+1:\infty}\hat f_{k+1:\infty}&= \bX_{\phi,1:k}^\top \vy\label{eq:derivative_1}\\
    \bX_{\phi,k+1:\infty}^\top \bX_{\phi,1:k}\hat f_{1:k}+ \left(\bX_{\phi,k+1:\infty}^\top\bX_{\phi,k+1:\infty}+\lambda I\right)\hat f_{k+1:\infty} &= \bX_{\phi,k+1:\infty}^\top\vy,\label{eq:derivative_2}
\end{align}where $I:\cH\to\cH$ is identity operator. Solving \eqref{eq:derivative_2} gives
\begin{align*}
    \hat f_{k+1:\infty} = \left(\bX_{\phi,k+1:\infty}^\top\bX_{\phi,k+1:\infty}+\lambda I\right)^{-1}\bX_{\phi,k+1:\infty}^\top\left(\vy-\bX_{\phi,1:k}\hat f_{1:k}\right),
\end{align*}which coincides with \eqref{eq:objective_f_2_decompose} because of the Woodbury formula
\[\bX_{\phi,k+1:\infty}^\top\left(\bX_{\phi,k+1:\infty}\bX_{\phi,k+1:\infty}^\top + \lambda I_N\right)^{-1} = \left(\bX_{\phi,k+1:\infty}^\top\bX_{\phi,k+1:\infty}+\lambda I\right)^{-1}\bX_{\phi,k+1:\infty}^\top.\]
For \eqref{eq:objective_f_1_decompose}, we plug \eqref{eq:objective_f_2_decompose} into \eqref{eq:derivative_1} to obtain
\begin{align*}
    &\left( \bX_{\phi,1:k}^\top\bX_{\phi,1:k} + \lambda I - \bX_{\phi,1:k}^\top\bX_{\phi,k+1:\infty}\left( \bX_{\phi,k+1:\infty}^\top \bX_{\phi,k+1:\infty} + \lambda I \right)^{-1} \bX_{\phi,k+1:\infty}^\top\bX_{\phi,1:k} \right)\hat f_{1:k}\\
    &= \bX_{\phi,1:k}^\top\left(I - \bX_{\phi,k+1:\infty}\left( \bX_{\phi,k+1:\infty}^\top \bX_{\phi,k+1:\infty} + \lambda I \right)^{-1}\bX_{\phi,k+1:\infty}^\top\right)\vy.
\end{align*}Let $F = I - \bX_{\phi,k+1:\infty}\left( \bX_{\phi,k+1:\infty}^\top \bX_{\phi,k+1:\infty} + \lambda I \right)^{-1}\bX_{\phi,k+1:\infty}^\top$. The above equation is then equivalent to
\begin{align*}
    \left(\bX_{\phi,1:k}^\top F \bX_{\phi,1:k} + \lambda I \right)\hat f_{1:k} = \bX_{\phi,1:k}^\top F \vy.
\end{align*}Applying the Woodbury formula to $F$ gives $F = \lambda\left(\lambda I + \bX_{\phi,k+1:\infty}\bX_{\phi,k+1:\infty}^\top\right)^{-1}$. As $\bX_{\phi,1:k}^\top F \bX_{\phi,1:k} + \lambda I$ is invertible (because $F\succeq 0$), we have
\begin{align}\label{eq:solution_hat_f_1}
    \begin{aligned}
        \hat f_{1:k}&=\left( \bX_{\phi,1:k}^\top \left(\lambda I + \bX_{\phi,k+1:\infty}\bX_{\phi,k+1:\infty}^\top\right)^{-1}\bX_{\phi,1:k} + I \right)^{-1} \bX_{\phi,1:k}^\top\left(\lambda I + \bX_{\phi,k+1:\infty}\bX_{\phi,k+1:\infty}^\top\right)^{-1}\vy.
    \end{aligned}
\end{align}To check that \eqref{eq:solution_hat_f_1} is equivalent to \eqref{eq:objective_f_1_decompose}, we take the gradient of the convex objective function $f_{1:k}\mapsto\norm{Q\left(\vy - \bX_{\phi,1:k}f_{1:k}\right)}_\cH^2 + \norm{f_{1:k}}_\cH^2$ from \eqref{eq:objective_f_1_decompose} and set it to $0$. This gives $\hat f_{1:k} =   \bX_{\phi,1:k}^\top Q^\top Q (\vy - \bX_{\phi,1:k}\hat f_{1:k})$. Recall that $Q^\top Q=\left(\bX_{\phi,k+1:\infty}\bX_{\phi,k+1:\infty}^\top  + \lambda I_N\right)^{-1}$, and $\bX_{\phi,1:k}^\top Q^\top Q \bX_{\phi,1:k}+I$ is invertible. Hence $\hat f_{1:k}$ from \eqref{eq:solution_hat_f_1} is the unique solution to the optimization problem from \eqref{eq:objective_f_1_decompose} and so \eqref{eq:objective_f_1_decompose} holds.
\endproof

\subsubsection{Estimation properties of the ``ridge estimator'' $\hat f_{1:k}$}

For any $f_{1:k}\in \cH_{1:k}$, we define its (empirical) excess risk as follows (note that KRR's empirical excess risk includes a regularization term):
\begin{eqnarray}
    \cL_{f_{1:k}} &=&\norm{Q\left(\vy -\bX_{\phi,1:k}f_{1:k}\right)}_\cH^2 + \norm{f_{1:k}}_\cH^2 - \left( \norm{Q\left(\vy - \bX_{\phi,1:k}f_{1:k}^*\right)}_\cH^2 + \norm{f_{1:k}^*}_\cH^2 \right)\notag\\
    &=& \norm{(\bX_{\phi,k+1:\infty}\bX_{\phi,k+1:\infty}^\top+ \lambda I_N)^{-1/2}\bX_{\phi,1:k}( f_{1:k} - f_{1:k}^*)}_2^2\notag\\
    &+&2 \inr{\bX_{\phi,1:k}^\top (\bX_{\phi,k+1:\infty}\bX_{\phi,k+1:\infty}^\top+ \lambda I_N)^{-1}(\bX_{\phi,k+1:\infty} f_{k+1:\infty}^*+\bxi) -  f_{1:k}^*, f_{1:k} - f_{1:k}^*}_\cH \notag\\
    &+& \norm{ f_{1:k} - f_{1:k}^*}_\cH^2,\label{eq:new_Q+M+R}
\end{eqnarray}where we have used the fact that $Q^\top Q = \left(\bX_{\phi,k+1:\infty}\bX_{\phi,k+1:\infty}^\top + \lambda I_N\right)^{-1}$, $\norm{Q\vlambda}_\cH = \norm{\left(\bX_{\phi,k+1:\infty}\bX_{\phi,k+1:\infty}^\top + \lambda I_N\right)^{-1/2}\vlambda}_2$ from Proposition~\ref{prop:decomposition_hat_f} and $\norm{f_{1:k}}_\cH^2 - \norm{f_{1:k}^*}_\cH^2 = \norm{f_{1:k} - f_{1:k}^*}_\cH^2 - 2\left\langle f_{1:k}^*,\, f_{1:k}^* - f_{1:k}\right\rangle_\cH$.

We denote the three terms of the decomposition \eqref{eq:new_Q+M+R} by $\cQ_{f_{1:k}}$, $\cM_{f_{1:k}}$ and $\cR_{f_{1:k}}$ respectively:
\begin{align}
    \cQ_{f_{1:k}} &= \norm{(\bX_{\phi,k+1:\infty}\bX_{\phi,k+1:\infty}^\top+ \lambda I_N)^{-1/2}\bX_{\phi,1:k}( f_{1:k} - f_{1:k}^*)}_2^2,\label{eq:def_quadratic_term}\\
    \cM_{f_{1:k}} &= 2 \inr{\bX_{\phi,1:k}^\top (\bX_{\phi,k+1:\infty}\bX_{\phi,k+1:\infty}^\top+ \lambda I_N)^{-1}(\bX_{\phi,k+1:\infty} f_{k+1:\infty}^*+\bxi) -  f_{1:k}^*, f_{1:k} - f_{1:k}^*}_\cH, \label{eq:def_multipler_term}\\
    \cR_{f_{1:k}} &= \norm{ f_{1:k} - f_{1:k}^*}_\cH^2 \label{eq:def_regularization_term}.
\end{align}
We notice that of these three terms, only the multiplier term $\cM_{f_{1:k}}$ can take negative values, whereas the quadratic term $\cQ_{f_{1:k}}$ and the regularization term $\cR_{f_{1:k}}$ are always positive.

We will show that with high probability, $\norm{\Gamma_{1:k}^{1/2}(\hat f_{1:k} - f_{1:k}^*)}_\cH\leq\square$ and $\norm{\hat f_{1:k} - f_{1:k}^*}_\cH\leq\triangle$, where $\square, \triangle>0$ will be defined later. In other words, we want to show that $\hat f_{1:k} \in f_{1:k}^* + B$ where $B$ is the unit ball of the norm $\vertiii{\cdot}$ defined as
\begin{align*}
    \vertiii{f} :=\max\left\{\frac{\norm{\Gamma_{1:k}^{1/2}f}_\cH}{\square}, \frac{\norm{f}_\cH}{\triangle} \right\}.
\end{align*}

From the definition of $\hat f_{1:k}$ in \eqref{eq:objective_f_1_decompose}, we know that $\cL_{\hat f_{1:k}}\leq 0$ so it suffices to show that for all $f_{1:k}\notin f_{1:k}^* + B$ we have $\cL_{f_{1:k}}>0$. We denote the border of $B$ in $V_{1:k}$ by $\partial  B$. Let  $ f_{1:k}  \in V_{1:k}$ be such that $ f_{1:k} \notin  f_{1:k} ^*+ B$.  There exists $f_0\in \partial  B$ and $\theta>1$ such that $ f_{1:k}  -  f_{1:k} ^* = \theta(f_0 -  f_{1:k} ^*)$. Using \eqref{eq:new_Q+M+R}, it follows from the convexity that $\cL_{ f_{1:k} }\geq \theta \cL_{f_0}$. As a consequence, if we prove that $\cL_{ f_{1:k} }>0$ for all $ f_{1:k} \in f_{1:k} ^*+\partial B$, this will imply that  $\cL_{ f_{1:k} }>0$ for all $ f_{1:k} \notin f_{1:k} ^*+ B$. Hence, we only need to show the positivity of the excess regularized risk $\cL_{ f_{1:k} }$ on the border $ f_{1:k} ^*+\partial B$.

For $f_{1:k}\in\partial B$, there are two cases:
\begin{enumerate}
    \item $\norm{\Gamma_{1:k}^{1/2}( f_{1:k} - f_{1:k}^*)}_\cH=\square$ and $\norm{ f_{1:k} - f_{1:k}^*}_\cH\leq \triangle$, or
    \item $\norm{\Gamma_{1:k}^{1/2}( f_{1:k} - f_{1:k}^*)}_\cH\leq \square$ and $\norm{ f_{1:k} - f_{1:k}^*}_\cH= \triangle$.
\end{enumerate}

We will prove that either we have $\cQ_{f_{1:k}}>\cM_{f_{1:k}}$ (this occurs in case [1]), or $\cR_{f_{1:k}}>\cM_{f_{1:k}}$ (this occurs in case [2]). Combined with \eqref{eq:new_Q+M+R}, this will show that $\cL_{f_{1:k}}>0$.
To achieve this goal, we need to obtain a lower bound for $\cQ_{f_{1:k}}$ in case [1] and an upper bound for $\cM_{f_{1:k}}$ in case [1] and [2]. The lower bound for $\cR_{f_{1:k}}$ is straightforward because this term is not random and is positive.

\paragraph{Bound of the multiplier term} We show an upper bound on $\cM_{f_{1:k}}$ when $f_{1:k}\in f_{1:k}^* + \partial B$.

Observe that
\begin{align*}
    \left|\cM_{f_{1:k}}\right|\leq 2\underset{f\in B}{\sup}\left|\left\langle \bX_{\phi,1:k}^\top (\bX_{\phi,k+1:\infty}\bX_{\phi,k+1:\infty}^\top + \lambda I_N)^{-1}(\bX_{\phi,k+1:\infty} f_{k+1:\infty}^*+\bxi) -  f_{1:k}^*,f \right\rangle_\cH\right|.
\end{align*}Also, for $f\in \cH_{1:k}$, we have $\vertiii{f}\leq\norm{\tilde\Gamma_{1:k}^{1/2}f}_\cH\leq \sqrt{2}\vertiii{f}$ where $\tilde\Gamma_{1:k}$ is defined in \eqref{eq:def_tilde_Gamma}. Therefore, $\vertiii{\cdot}$'s dual norm $\vertiii{\cdot}_*$ is also equivalent to $\norm{\tilde{\Gamma}_{1:k}^{1/2}\cdot}_\cH$'s dual norm  which is given by $\norm{\tilde{\Gamma}_{1:k}^{-1/2}\cdot}_\cH$: for all $f\in \cH_{1:k}$, $(1/\sqrt{2})\vertiii{f}_*\leq \norm{\tilde{\Gamma}_{1:k}^{-1/2}f}_\cH\leq \vertiii{f}_*$. Hence, for all $f_{1:k}\in f_{1:k}^* + \partial B$, we have
 \begin{align}\label{eq:multiplier_decomp}
     \left|\cM_{f_{1:k}}\right|&\leq 2 \sqrt{2}\norm{\tilde \Gamma_{1:k}^{-1/2}\left(\bX_{\phi,1:k}^\top (\bX_{\phi,k+1:\infty}\bX_{\phi,k+1:\infty}^\top + \lambda I_N)^{-1}(\bX_{\phi,k+1:\infty} f_{k+1:\infty}^*+\bxi) -  f_{1:k}^*\right)}_\cH\notag\\
 &\leq 2 \sqrt{2}\bigg(\norm{\tilde \Gamma_{1:k}^{-1/2}\bX_{\phi,1:k}^\top (\bX_{\phi,k+1:\infty}\bX_{\phi,k+1:\infty}^\top+ \lambda I_N)^{-1}\bX_{\phi,k+1:\infty} f_{k+1:\infty}^*}_\cH  \\
 &\quad+ \norm{\tilde \Gamma_{1:k}^{-1/2}\bX_{\phi,1:k}^\top (\bX_{\phi,k+1:\infty}\bX_{\phi,k+1:\infty}^\top+ \lambda I_N)^{-1}\bxi}_\cH + \norm{\tilde \Gamma_{1:k}^{-1/2} f_{1:k}^*}_\cH\bigg)\notag
 \end{align}

We next handle the first two terms in \eqref{eq:multiplier_decomp} in the next two lemmas.

\begin{Lemma}\label{lem:first_beta_2_star_term}
Under the assumptions of Proposition~\ref{prop:stochastic_argument_upper_bound},
\begin{align}\label{eq:multiplier_bX_f_star_k<N}
\begin{aligned}
    &\norm{ \tilde\Gamma_{1:k}^{-1/2} \bX_{\phi,1:k}^\top\left( \bX_{\phi,k+1:\infty}\bX_{\phi,k+1:\infty}^\top + \lambda I_N \right)^{-1}\bX_{\phi,k+1:\infty}f_{k+1:\infty}^* }_\cH\\
    &\leq \frac{4\oC{C_bX_f_star}\oC{C_RIP_upper}\kappa N\norm{\Gamma_{k+1:\infty}^{1/2}f_{k+1:\infty}^*}_\cH}{4\lambda+ \Tr\left(\Gamma_{k+1:\infty}\right)}\sigma\left(\square,\triangle\right),
\end{aligned}
\end{align}where
\begin{align}\label{eq:def_sigma_square_triangle}
    \sigma(\square,\triangle) := \begin{cases}
        \square, & \mbox{if } \triangle\sqrt{\sigma_1}\geq \square\\
        \triangle\sqrt{\sigma_1}, & \mbox{otherwise}.
    \end{cases}
\end{align}
\end{Lemma}

\beginproof 
On the event $\Omega_0$ we have
\begin{align*}
    \norm{\tilde\Gamma_{1:k}^{-1/2}\bX_{\phi,1:k}^\top }_{\text{op}} = \norm{\bX_{\phi,1:k}\tilde\Gamma_{1:k}^{-1/2} }_{\text{op}}\leq \oC{C_RIP_upper}\sqrt{N}\norm{\Gamma_{1:k}^{1/2}\tilde\Gamma_{1:k}^{-1/2}}_{\text{op}},
\end{align*}because of the isomorphic property of $\bX_{\phi,1:k}$ and
\begin{align*}
    \norm{\left(\bX_{\phi,k+1:\infty}\bX_{\phi,k+1:\infty}^\top + \lambda I_N\right)^{-1}}_{\text{op}} \leq \frac{4}{4\lambda+\Tr\left(\Gamma_{k+1:\infty}\right)}.
\end{align*}Hence, on $\Omega_0$,
\begin{align*}
    &\norm{ \tilde\Gamma_{1:k}^{-1/2} \bX_{\phi,1:k}^\top\left( \bX_{\phi,k+1:\infty}\bX_{\phi,k+1:\infty}^\top + \lambda I_N \right)^{-1}\bX_{\phi,k+1:\infty}f_{k+1:\infty}^* }_\cH \\
    &\leq \norm{\tilde\Gamma_{1:k}^{-1/2}\bX_{\phi,1:k}^\top}_{\text{op}} \norm{\left( \bX_{\phi,k+1:\infty}\bX_{\phi,k+1:\infty}^\top + \lambda I_N \right)^{-1}}_{\text{op}}\norm{\bX_{\phi,k+1:\infty}f_{k+1:\infty}^*}_2 \\
    &\leq \frac{4\oC{C_bX_f_star}\oC{C_RIP_upper}\kappa N}{4\lambda+\Tr\left(\Gamma_{k+1:\infty}\right)} \norm{\Gamma_{1:k}^{1/2}\tilde\Gamma_{1:k}^{-1/2}}_{\text{op}}\norm{\Gamma_{k+1:\infty}^{1/2}f_{k+1:\infty}^*}_\cH\leq \frac{4\oC{C_bX_f_star}\oC{C_RIP_upper} \kappa N\norm{\Gamma_{k+1:\infty}^{1/2}f_{k+1:\infty}^*}_\cH}{4\lambda+ \Tr\left(\Gamma_{k+1:\infty}\right)}\sigma\left(\square,\triangle\right),
\end{align*}where the last inequality follows from the definition of $\Gamma_{1:k}$ and $\tilde\Gamma_{1:k}$.
\endproof

We define the sets
\begin{align*}
    J_1 := \left\{ j\in[k]:\, \sigma_j\geq \left(\frac{\square}{\triangle}\right)^2 \right\} ,\quad J_2:=[k]\backslash J_1.
\end{align*}
They will match the definition of $J_1$ and $J_2$ in \eqref{eq:def_J1} once $\triangle$ and $\square$ have been chosen.

We prove the following lemma:
\begin{Lemma}\label{lem:noise_term}
    Under the assumptions of Proposition~\ref{prop:stochastic_argument_upper_bound}. We define
    \begin{align}\label{eq:def_t_square_triangle}
        t(\square,\triangle) := \frac{1}{\sigma^2\left(\square,\triangle\right)}\left( \left|J_1\right|\square^2 + \triangle^2\sum_{j\in J_2}\sigma_j \right).
    \end{align} Recall $r$ and $\okappa{kappa_noise}$ from Theorem~\ref{theo:upper_KRR}.
    There then exists an absolute constant $\oC{C_noise}$ depending only on $\okappa{kappa_noise}$ and there exists an absolute constant $\nC\label{C_noise_term}$ such that with probability at least $1- (\oC{C_noise}/\lfloor t(\square,\triangle)\rfloor)^{r/4} - \bP(\Omega_0^c)$,
    \begin{align}\label{eq:multiplier_bxi_k<N}
        \norm{\tilde\Gamma_{1:k}^{-1/2}\bX_{\phi,1:k}^\top\left( \bX_{\phi,k+1:\infty}\bX_{\phi,k+1:\infty}^\top + \lambda I_N \right)^{-1}\bxi}_\cH \leq \frac{8\oC{C_noise_term}\sigma_\xi\sqrt{N}}{4\lambda + \Tr\left(\Gamma_{k+1:\infty}\right)}\sqrt{\left|J_1\right|\square^2 + \triangle^2\sum_{j\in J_2}\sigma_j}.
    \end{align}
\end{Lemma}
\beginproof Let $D = \tilde\Gamma_{1:k}^{-1/2}\bX_{\phi,1:k}^\top\left( \bX_{\phi,k+1:\infty}\bX_{\phi,k+1:\infty}^\top + \lambda I_N \right)^{-1}$. We calculate separately the upper bounds for $\sqrt{\Tr\left(DD^\top\right)}$ and $\norm{D}_{\text{op}}$. For the basis $(\varphi_j)_{j\in\bN}$ of eigenfunctions $\Gamma$, and since $\Tr$ is independent with the choice of the basis, on $\Omega_0$ we have
\begin{align*}
    \Tr\left(DD^\top\right) &= \sum_{j\in\bN}\left\langle DD^\top \varphi_j, \varphi_j\right\rangle_\cH = \sum_{j\in\bN}\norm{D^\top\varphi_j}_2^2 = \sum_{j\in[k]}\norm{D^\top \varphi_j}_2^2 \\
    &= \sum_{j=1}^k \left(\frac{\sqrt{\sigma_j}}{\square}\vee \frac{1}{\triangle}\right)^{-2} \norm{\left(\bX_{\phi,k+1:\infty}\bX_{\phi,k+1:\infty}^\top + \lambda I_N \right)^{-1}\bX_{\phi,1:k}\varphi_j}_2^2\\
    &\leq \sum_{j=1}^k \left(\frac{\sqrt{\sigma_j}}{\square}\vee \frac{1}{\triangle}\right)^{-2} \norm{\left(\bX_{\phi,k+1:\infty}\bX_{\phi,k+1:\infty}^\top + \lambda I_N \right)^{-1}}_{\text{op}}^2\norm{\bX_{\phi,1:k}\varphi_j}_2^2\\
    &\leq \sum_{j=1}^k \left(\frac{\sqrt{\sigma_j}}{\square}\vee \frac{1}{\triangle}\right)^{-2} \left(\lambda+\frac{\Tr\left(\Gamma_{k+1:\infty}\right)}{4}\right)^{-2}\oC{C_RIP_upper}^2N\sigma_j,
\end{align*}where the last inequality follows from \eqref{eq:DM_applied}. Hence
\begin{align}\label{eq:trace_DD_top_1}
    \sqrt{\Tr\left(DD^\top\right)}\leq \frac{4\oC{C_RIP_upper}\sqrt{N}}{4\lambda + \Tr\left(\Gamma_{k+1:\infty}\right)}\sqrt{\left|J_1\right|\square^2 + \triangle^2\sum_{j\in J_2}\sigma_j}.
\end{align}This implies that $D$ is a Hilbert-Schmidt operator.
Using the inequality $\norm{\Gamma_{1:k}^{1/2}\tilde\Gamma_{1:k}^{-1/2}}_{\text{op}}\leq\sigma(\square,\triangle)$ we obtain
\begin{align}\label{eq:op_D_1}
\begin{aligned}
    \norm{D}_{\text{op}} &= \norm{D^\top}_{\text{op}} \leq \norm{\left(\bX_{\phi,k+1:\infty}\bX_{\phi,k+1:\infty}^\top + \lambda I_N \right)^{-1}}_{\text{op}} \norm{\bX_{\phi,1:k}\tilde\Gamma_{1:k}^{-1/2}}_{\text{op}}\\
    &\leq\oC{C_RIP_upper}\sqrt{N}\norm{\Gamma_{1:k}^{1/2}\tilde\Gamma_{1:k}^{-1/2}}_{\text{op}}\cdot\frac{4}{4\lambda+\Tr\left(\Gamma_{k+1:\infty}\right)}\\
    &\leq \frac{4\oC{C_RIP_upper}\sqrt{N}\sigma\left(\square,\triangle\right)}{4\lambda + \Tr\left(\Gamma_{k+1:\infty}\right)}. 
\end{aligned}
\end{align}
We finish the proof by Proposition~\ref{prop:noise_concentration} with the $k$ from Proposition~\ref{prop:noise_concentration} set as
\begin{align*}
    k = \biggl\lfloor\frac{\left|J_1\right|\square^2 + \triangle^2\sum_{j\in J_2}\sigma_j}{\sigma^2\left(\square,\triangle\right)} \biggr\rfloor
\end{align*}and with $\oC{C_noise_term}=\frac{3}{2}\oC{C_RIP_upper}$.
\endproof

\paragraph{Bound of the quadratic term and choice of $\square$ and $\triangle$}

In the previous section, we obtained an upper bound on $\cM_{f_{1:k}}$. Our main approach, as outlined in the previous section, is to separately prove $\cQ_{f_{1:k}} > \cM_{f_{1:k}}$ in case [1], and $\cR_{f_{1:k}} > \cM_{f_{1:k}}$ in case [2]. Now that we have the upper bound for $\cM_{f_{1:k}}$, it only remains to bound $\cQ_{f_{1:k}}$ in case [1].

Before we begin, we need to make another classification. This time, the classification is based on the values of $\sigma(\square, \triangle)$. In the upcoming proof, we will \emph{firstly} start by classifying based on $\sigma(\square, \triangle)$, and \emph{then} proceed to prove the desired propositions separately in cases [1] and [2]. This parameter is crucial in the analysis as it determines whether the regularization is too strong, potentially completely submerging the signal. One can revisit the classification discussion regarding $\sigma_1 N$ and $4\lambda + \Tr\left(\Gamma_{k+1:\infty}\right)$ in Theorem~\ref{theo:upper_KRR}. Doing so will reveal that this corresponds to the classification based on the values of $\sigma(\square, \triangle)$. When $\sigma_1 N$ is too small, it signifies excessive regularization that drowns out the signal.


\paragraph{If $\sigma\left(\square,\triangle\right)=\square$} 
Let us first study case [1]. Consider $f_{1:k}\in \cH_{1:k}$ such that $\norm{\Gamma_{1:k}^{1/2}( f_{1:k} - f_{1:k}^*)}_\cH = \square$ and $\norm{ f_{1:k} - f_{1:k}^*}_\cH\leq\triangle$. In this case, we show that $\cQ_{f_{1:k}}>\cM_{f_{1:k}}$. Notice that on $\Omega_0$ we have
\begin{align}\label{eq:quadratic_k<N}
\begin{aligned}
    \cQ_{f_{1:k}} &= \norm{\left(\bX_{\phi,k+1:\infty}\bX_{\phi,k+1:\infty}^\top + \lambda I_N\right)^{-1/2}\bX_{\phi,1:k}(f_{1:k}-f_{1:k}^*)}_2^2\\
    &\geq \left(\lambda + \frac{3\Tr\left(\Gamma_{k+1:\infty}\right)}{2}\right)^{-1}\oc{c_RIP_lower}^2N\norm{\Gamma_{1:k}^{1/2}(f_{1:k}-f_{1:k}^*)}_\cH^2\\
    &=\frac{2\oc{c_RIP_lower}^2N\square^2}{2\lambda + 3\Tr\left(\Gamma_{k+1:\infty}\right)}.
\end{aligned}
\end{align}To prove that $\cQ_{f_{1:k}}>\cM_{f_{1:k}}$, it suffices to show that
\begin{align*}
    \frac{\oc{c_RIP_lower}^2N\square^2}{\sqrt{2}\left(2\lambda + 3\Tr\left(\Gamma_{k+1:\infty}\right)\right)} &> \norm{\tilde \Gamma_{1:k}^{-1/2}\bX_{\phi,1:k}^\top (\bX_{\phi,k+1:\infty}\bX_{\phi,k+1:\infty}^\top+ \lambda I_N)^{-1}\bX_{\phi,k+1:\infty} f_{k+1:\infty}^*}_\cH \\
    &+ \norm{\tilde \Gamma_{1:k}^{-1/2}\bX_{\phi,1:k}^\top (\bX_{\phi,k+1:\infty}\bX_{\phi,k+1:\infty}^\top+ \lambda I_N)^{-1}\bxi}_\cH + \norm{\tilde \Gamma_{1:k}^{-1/2} f_{1:k}^*}_\cH
\end{align*}
Lemma~\ref{lem:first_beta_2_star_term} and Lemma~\ref{lem:noise_term} then make clear that is suffices to prove that the following conditions hold for well-chosen $\square$ and $\triangle$.
\begin{itemize}
    \item $$\frac{\oc{c_RIP_lower}^2N\square^2}{\sqrt{2}\left(2\lambda + 3\Tr\left(\Gamma_{k+1:\infty}\right)\right)}> \frac{4\oC{C_bX_f_star}\oC{C_RIP_upper}\kappa N\norm{\Gamma_{k+1:\infty}^{1/2}f_{k+1:\infty}^*}_\cH}{4\lambda+ \Tr\left(\Gamma_{k+1:\infty}\right)}\sigma\left(\square,\triangle\right).$$ This is equivalent to $\square > \frac{12\sqrt{2}\oC{C_bX_f_star}\oC{C_RIP_upper}}{\oc{c_RIP_lower}^2}\kappa\norm{\Gamma_{k+1:\infty}^{1/2}f_{k+1:\infty}^*}_\cH.$
    \item $$\frac{\oc{c_RIP_lower}^2N\square^2}{\sqrt{2}\left(2\lambda + 3\Tr\left(\Gamma_{k+1:\infty}\right)\right)}> \frac{8\oC{C_noise_term}\sqrt{N}\sigma_\xi}{4\lambda+\Tr\left(\Gamma_{k+1:\infty}\right)}\sqrt{\left|J_1\right|\square^2},$$ which holds if $\square > \frac{24\sqrt{2}\oC{C_noise_term}}{\oc{c_RIP_lower}^2}\sigma_\xi\sqrt{\frac{\left|J_1\right|}{N}}.$
    \item $$\frac{\oc{c_RIP_lower}^2N\square^2}{\sqrt{2}\left(2\lambda + 3\Tr\left(\Gamma_{k+1:\infty}\right)\right)}>\frac{8\oC{C_noise_term}\sqrt{N}\sigma_\xi}{4\lambda+\Tr\left(\Gamma_{k+1:\infty}\right)}\sqrt{\triangle^2\sum_{j\in J_2}\sigma_j},$$ which holds if $\square > \left( \frac{24\oC{C_noise_term}}{\oc{c_RIP_lower}^2}\triangle \sigma_\xi \sqrt{\frac{2}{N}\sum_{j\in J_2}\sigma_j} \right)^{1/2}.$
    \item $$\frac{\oc{c_RIP_lower}^2N\square^2}{\sqrt{2}\left(2\lambda + 3\Tr\left(\Gamma_{k+1:\infty}\right)\right)}> \norm{\tilde\Gamma_{1:k}^{-1/2}f_{1:k}^*}_\cH,$$ which is equivalent to $\square > \sqrt{\norm{\tilde\Gamma_{1:k}^{-1/2}f_{1:k}^*}_\cH \frac{\sqrt{2}\left(2\lambda + 3\Tr\left(\Gamma_{k+1:\infty}\right)\right)}{\oc{c_RIP_lower}^2 N} }.$
\end{itemize}
In conclusion, there exists an absolute constant $\nC\label{C_square_1}$ (for example, $\oC{C_square_1} = \frac{24\sqrt{2}\oC{C_bX_f_star}\oC{C_noise_term}}{\oc{c_RIP_lower}^2}$) so that if we have
    \begin{align}\label{eq:square}
    \begin{aligned}
        \square &> \oC{C_square_1}\kappa\max\bigg\{ \sigma_\xi\sqrt{\frac{\left|J_1\right|}{N}}, \left(\triangle\sigma_\xi\sqrt{\frac{1}{N}\sum_{j\in J_2}\sigma_j}\right)^{1/2}, \norm{\Gamma_{k+1:\infty}^{1/2}f_{k+1:\infty}^*}_\cH,\\
        &\sqrt{\norm{\tilde\Gamma_{1:k}^{-1/2}f_{1:k}^*}_\cH \frac{2\lambda + 3\Tr\left(\Gamma_{k+1:\infty}\right)}{N} } \bigg\},
    \end{aligned}
    \end{align}then $\cQ_{f_{1:k}}>\cM_{f_{1:k}}$.

In case [2]. We consider a function $f_{1:k}\in \cH_{1:k}$ such that $\norm{\Gamma_{1:k}^{1/2}( f_{1:k} - f_{1:k}^*)}_\cH\leq \square$ and $\norm{ f_{1:k} - f_{1:k}^*}_\cH=\triangle$. In this case, we show that $\cR_{f_{1:k}}>\cM_{f_{1:k}}$. Since $\cR_{f_{1:k}}=\triangle^2$, this amounts to showing that
\begin{align*}
    \triangle^2 &> 2\sqrt{2}\norm{\tilde \Gamma_{1:k}^{-1/2}\bX_{\phi,1:k}^\top (\bX_{\phi,k+1:\infty}\bX_{\phi,k+1:\infty}^\top+ \lambda I_N)^{-1}\bX_{\phi,k+1:\infty} f_{k+1:\infty}^*}_\cH \\
    &+ 2\sqrt{2}\left(\norm{\tilde \Gamma_{1:k}^{-1/2}\bX_{\phi,1:k}^\top (\bX_{\phi,k+1:\infty}\bX_{\phi,k+1:\infty}^\top+ \lambda I_N)^{-1}\bxi}_\cH + \norm{\tilde \Gamma_{1:k}^{-1/2} f_{1:k}^*}_\cH\right).
\end{align*}By Lemma~\ref{lem:first_beta_2_star_term} and Lemma~\ref{lem:noise_term}, $\triangle$ must satisfy the following conditions:
\begin{itemize}
    \item $$\triangle^2 > \frac{4\oC{C_bX_f_star}\oC{C_RIP_upper}\kappa N\norm{\Gamma_{k+1:\infty}^{1/2}f_{k+1:\infty}^*}_\cH}{4\lambda+ \Tr\left(\Gamma_{k+1:\infty}\right)}\sigma\left(\square,\triangle\right).$$
    \item $$\triangle^2 > \frac{8\oC{C_noise_term}\sqrt{N}\sigma_\xi}{4\lambda + \Tr\left(\Gamma_{k+1:\infty}\right)}\sqrt{\left|J_1\right|\square^2}.$$
    \item $$\triangle^2> \frac{8\oC{C_noise_term}\sqrt{N}\sigma_\xi}{4\lambda + \Tr\left(\Gamma_{k+1:\infty}\right)}\sqrt{\triangle^2\sum_{j\in J_2}\sigma_j},$$ which is equivalent to $\triangle^2 > \frac{64\oC{C_noise_term}^2N\sigma_\xi^2}{\left(4\lambda + \Tr\left(\Gamma_{k+1:\infty}\right)\right)^2} \left(\sum_{j\in J_2}\sigma_j\right).$
    \item $$\triangle^2 > \norm{\tilde\Gamma_{1:k}^{-1/2}f_{1:k}^*}_\cH.$$
\end{itemize}Hence, there exists an absolute constant $\nC\label{C_triangle_1}$ (for example, $\oC{C_triangle_1} = 64\oC{C_bX_f_star}\oC{C_noise_term}^2$). We need to choose
\begin{align}\label{eq:triangle}
    \triangle^2  > \oC{C_triangle_1}\kappa\max\left\{ \frac{\sigma_\xi\square\sqrt{\left|J_1\right|N}}{4\lambda + \Tr\left(\Gamma_{k+1:\infty}\right)}, \frac{\sigma_\xi^2 N \sum_{j\in J_2}\sigma_j}{\left(4\lambda + \Tr\left(\Gamma_{k+1:\infty}\right)\right)^2}, \frac{N\square\norm{\Gamma_{k+1:\infty}^{1/2}f_{k+1:\infty}^*}_\cH}{4\lambda + \Tr\left(\Gamma_{k+1:\infty}\right)}, \norm{\tilde\Gamma_{1:k}^{-1/2}f_{1:k}^*}_\cH \right\}.
\end{align}Then $\cR_{f_{1:k}}>\cM_{f_{1:k}}$.

While we have shown that $\cQ_{f_{1:k}} > \cM_{f_{1:k}}$ in case [1] and $\cR_{f_{1:k}} > \cM_{f_{1:k}}$ in case [2] if \eqref{eq:square} and \eqref{eq:triangle} hold, our task is not yet complete because \eqref{eq:triangle} and \eqref{eq:square} do not explicitly define for $\triangle$ and $\square$. We next derive explicit definition for $\square$ and $\triangle$ through these two equations.

We fix $\triangle$ so that
\begin{align}\label{eq:quotien_square_triangle}
    \frac{\square}{\triangle} = \sqrt{\frac{\kappa_{DM} \left(4\lambda + \Tr\left(\Gamma_{k+1:\infty}\right)\right) }{N}}
\end{align}and take $\nC\label{C_square_2} = \oC{C_square_1}^2\kappa^2\kappa_{DM}^{-1/2}\vee 2\oC{C_triangle_1}\kappa$ and
\begin{align}\label{eq:def_square_3}
    \begin{aligned}
        \square &> \oC{C_square_2}\max\bigg\{ \sigma_\xi\sqrt{\frac{\left|J_1\right|}{N}}, \sigma_\xi\left(\frac{\sum_{j\in J_2}\sigma_j}{  4\lambda + \Tr\left(\Gamma_{k+1:\infty}\right)}\right)^{1/2}, \norm{\Gamma_{k+1:\infty}^{1/2}f_{k+1:\infty}^*}_\cH,\\
        &\sqrt{\norm{\tilde\Gamma_{1:k}^{-1/2}f_{1:k}^*}_\cH \frac{2\lambda + 3\Tr\left(\Gamma_{k+1:\infty}\right)}{N} } \bigg\}.
    \end{aligned}
\end{align}One may observe that the second term inside the max is different from the corresponding term in \eqref{eq:square}. However, if $\square > \oC{C_square_2}\sigma_\xi\left(\frac{\sum_{j\in J_2}\sigma_j}{  4\lambda + \Tr\left(\Gamma_{k+1:\infty}\right)}\right)^{1/2}$, and $\square,\triangle$ satisfy \eqref{eq:quotien_square_triangle}, it follows that
\begin{align*}
    \frac{\square}{\sqrt{\triangle \sigma_\xi \sqrt{\frac{1}{N} \sum_{j\in J_2} \sigma_j } }} &= \sqrt{\frac{\square}{\triangle} \frac{\square}{\sigma_\xi \sqrt{\frac{1}{N} \sum_{j\in J_2}\sigma_j } } }\\
    &> \left( \sqrt{\frac{\kappa_{DM} \left(4\lambda + \Tr\left(\Gamma_{k+1:\infty}\right)\right) }{N}} \cdot \frac{ \oC{C_square_2}\sigma_\xi\left(\frac{\sum_{j\in J_2}\sigma_j}{  4\lambda + \Tr\left(\Gamma_{k+1:\infty}\right)}\right)^{1/2} }{\sigma_\xi \sqrt{\frac{1}{N} \sum_{j\in J_2}\sigma_j } } \right)^{1/2}\\
    &= \kappa_{DM}^{1/4}\sqrt{\oC{C_square_2}} \geq \oC{C_square_1}\kappa.
\end{align*}Hence the new choice of $\square$ in \eqref{eq:def_square_3} satisfies \eqref{eq:square}.

We now need to check that for this choice of $\square$, \eqref{eq:triangle} is also satisfied.
\begin{itemize}
    \item 
    \begin{align*}
        \frac{\triangle^2}{\frac{\sigma_\xi\square\sqrt{\left|J_1\right|N}}{ 4\lambda + \Tr\left(\Gamma_{k+1:\infty}\right) }} &= \square^2 \frac{N}{\kappa_{DM} \left(4\lambda + \Tr\left(\Gamma_{k+1:\infty}\right)\right)} \cdot\frac{ 4\lambda + \Tr\left(\Gamma_{k+1:\infty}\right) }{ \sigma_\xi \square \sqrt{\left|J_1\right|N} } > \frac{\oC{C_square_2}}{\kappa_{DM}}>\oC{C_triangle_1}\kappa.
    \end{align*}
    \item 
    \begin{align*}
        \frac{\triangle^2}{ \frac{\sigma_\xi^2N\sum_{j\in J_2}\sigma_j}{\left( 4\lambda + \Tr\left(\Gamma_{k+1:\infty}\right)\right)^2 } } &= \square^2 \frac{N}{ \kappa_{DM} \left(4\lambda + \Tr\left(\Gamma_{k+1:\infty}\right)\right) }\cdot\frac{\left( 4\lambda + \Tr\left(\Gamma_{k+1:\infty}\right)\right)^2}{\sigma_\xi^2N\sum_{j\in J_2}\sigma_j}\\
        &>\oC{C_square_2}^2  \frac{\sigma_\xi^2\sum_{j\in J_2}\sigma_j}{4\lambda + \Tr\left(\Gamma_{k+1:\infty}\right)} \frac{N}{ \kappa_{DM} \left(4\lambda + \Tr\left(\Gamma_{k+1:\infty}\right)\right) }\cdot\frac{\left( 4\lambda + \Tr\left(\Gamma_{k+1:\infty}\right)\right)^2}{\sigma_\xi^2N\sum_{j\in J_2}\sigma_j}>\oC{C_triangle_1}\kappa.
    \end{align*}
    \item 
    \begin{align*}
        \frac{\triangle^2}{ \frac{N\square \norm{\Gamma_{k+1:\infty}^{1/2}f_{k+1:\infty}^*}_\cH}{ 4\lambda + \Tr\left(\Gamma_{k+1:\infty}\right) } } &= \square^2 \frac{N}{ \kappa_{DM} \left(4\lambda + \Tr\left(\Gamma_{k+1:\infty}\right)\right) }\cdot \frac{ 4\lambda + \Tr\left(\Gamma_{k+1:\infty}\right) }{ N\square\norm{\Gamma_{k+1:\infty}^{1/2}f_{k+1:\infty}^*}_\cH }>\frac{\oC{C_square_2}}{\kappa_{DM}}>\oC{C_triangle_1}\kappa.
    \end{align*}
    \item 
    \begin{align*}
        \frac{\triangle^2}{ \norm{\tilde\Gamma_{1:k}^{-1/2}f_{1:k}^*}_\cH } &= \square^2 \frac{N}{ \kappa_{DM} \left(4\lambda + \Tr\left(\Gamma_{k+1:\infty}\right)\right) }\cdot \frac{1}{\norm{\tilde\Gamma_{1:k}^{-1/2}f_{1:k}^*}_\cH}>\frac{\oC{C_square_2}}{2\kappa_{DM}}>\oC{C_triangle_1}\kappa.
    \end{align*}
\end{itemize}We deduce that with the right choice of the absolute constants, such a choice of $\square,\triangle$ satisfies \eqref{eq:square} and \eqref{eq:triangle}.

We have established that by selecting appropriate values for $\square$ and $\triangle$, we can conclude the following: if $f_{1:k} - f_{1:k}^* \in \partial B$, then we necessarily have $\cL_{f_{1:k}} > 0$, and thanks to a homogeneity argument, it follows that for all $f_{1:k}\notin f_{1:k}^*+B$, we have $\cL_{f_{1:k}}\leq 0$ hence $\hat f_{1:k}\in f_{1:k}^*+B$.

In the beginning of the analysis, we assumed that $\sigma(\square,\triangle)=\square$, which is true if and only if $\sigma_1\geq \kappa_{DM}\frac{4\lambda+\Tr\left(\Gamma_{k+1:\infty}\right)}{N}$. Hence if this inequality is satisfied, $\square$ is an upper bound on the estimation error $\norm{\Gamma_{1:k}^{1/2}(\hat f_{1:k} - f_{1:k}^*)}_\cH$ and $\triangle$ is an upper bound on $\norm{\hat f_{1:k}-f_{1:k}^*}_\cH$. Notice also that $\tilde\Gamma_{1:k}^{-1/2} = U\tilde D_{1}^{-1/2}U^\top$ where $\tilde D_1^{-1/2}=:\square D_{1,\mathrm{thre}}^{-1/2}$. Hence we can express $\square$ as in Equation~\eqref{eq:def_square_2}.
\begin{align*}
    \square &= \oC{C_square_2}\max\bigg\{ \sigma_\xi\sqrt{\frac{\left|J_1\right|}{N}}, \sigma_\xi\sqrt{\frac{\sum_{j\in J_2}\sigma_j}{  4\lambda + \Tr\left(\Gamma_{k+1:\infty}\right)}}, \norm{\Gamma_{k+1:\infty}^{1/2}f_{k+1:\infty}^*}_\cH,\\
    &{\norm{\tilde\Gamma_{1,\mathrm{thre}}^{-1/2}f_{1:k}^*}_\cH \frac{2\lambda + 3\Tr\left(\Gamma_{k+1:\infty}\right)}{N} } \bigg\}.
\end{align*}

\paragraph{If $\sigma(\square,\triangle)=\triangle\sqrt{\sigma_1}$} In this case, it follows by definition that $J_1=\emptyset$, $J_2=[k]$,
\begin{align*}
    t\left(\square,\triangle\right) = \frac{\Tr\left(\Gamma_{1:k}\right)}{\sigma_1},\mbox{ and }\tilde D_1^{1/2}=\frac{1}{\triangle}\diag\left(1,\cdots,1,0,\cdots\right),
\end{align*}where there are $k$ ones in the definition of $\tilde D_1^{1/2}$. 
Since we have completed a similar proof in the previous paragraph, we will expedite the presentation in this paragraph.

Suppose that $\norm{\Gamma_{1:k}^{1/2}( f_{1:k} - f_{1:k}^*)}_\cH=\square$ and $\norm{ f_{1:k} - f_{1:k}^*}_{\cH}\leq\triangle$. As we discussed in the previous subsections, on $\Omega_0$, $\cQ_{f_{1:k}} \geq \frac{N\square^2}{4\lambda + 6\Tr\left(\Gamma_{k+1:\infty}\right)}$. To show that $\cQ_{f_{1:k}}>\cM_{f_{1:k}}$, it suffices to show that
\begin{align*}
    \frac{\oc{c_RIP_lower}^2N\square^2}{\sqrt{2}\left(2\lambda + 3\Tr\left(\Gamma_{k+1:\infty}\right)\right)} &> \norm{\tilde \Gamma_{1:k}^{-1/2}\bX_{\phi,1:k}^\top (\bX_{\phi,k+1:\infty}\bX_{\phi,k+1:\infty}^\top+ \lambda I_N)^{-1}\bX_{\phi,k+1:\infty} f_{k+1:\infty}^*}_\cH \\
    &+ \norm{\tilde \Gamma_{1:k}^{-1/2}\bX_{\phi,1:k}^\top (\bX_{\phi,k+1:\infty}\bX_{\phi,k+1:\infty}^\top+ \lambda I_N)^{-1}\bxi}_\cH + \norm{\tilde \Gamma_{1:k}^{-1/2} f_{1:k}^*}_\cH.
\end{align*}Recall that Lemma~\ref{lem:first_beta_2_star_term} and Lemma~\ref{lem:noise_term} hold true for all possible values of $\square,\triangle$, so we can still use them in the current setting. Hence:
\begin{itemize}
    \item
    \begin{align*}
        \frac{\oc{c_RIP_lower}^2N\square^2}{\sqrt{2}\left(2\lambda + 3\Tr\left(\Gamma_{k+1:\infty}\right)\right)} > \frac{4\oC{C_bX_f_star}\oC{C_RIP_upper}\kappa N\norm{\Gamma_{k+1:\infty}^{1/2}f_{k+1:\infty}^*}_\cH}{4\lambda+ \Tr\left(\Gamma_{k+1:\infty}\right)}\sigma\left(\square,\triangle\right),
    \end{align*}which is equivalent to $\square^2 > \frac{12\sqrt{2}\oC{C_bX_f_star}\oC{C_RIP_upper}}{\oc{c_RIP_lower}^2}\kappa\triangle\sqrt{\sigma_1}\norm{\Gamma_{k+1:\infty}^{1/2}f_{k+1:\infty}^*}_\cH.$
    \item $$\frac{\oc{c_RIP_lower}^2N\square^2}{\sqrt{2}\left(2\lambda + 3\Tr\left(\Gamma_{k+1:\infty}\right)\right)}>\frac{8\oC{C_noise_term}\sqrt{N}\sigma_\xi}{4\lambda+\Tr\left(\Gamma_{k+1:\infty}\right)}\sqrt{\left|J_1\right|\square^2},$$ which is true since $\left|J_1\right|=0$.
    \item $$\frac{\oc{c_RIP_lower}^2N\square^2}{\sqrt{2}\left(2\lambda + 3\Tr\left(\Gamma_{k+1:\infty}\right)\right)}>\frac{8\oC{C_noise_term}\sqrt{N}\sigma_\xi}{4\lambda+\Tr\left(\Gamma_{k+1:\infty}\right)}\sqrt{\triangle^2\sum_{j\in J_2}\sigma_j},$$ which is true if $\square^2 > \frac{24\sqrt{2}\oC{C_noise_term}}{\oc{c_RIP_lower}^2}\sigma_\xi\triangle\sqrt{\frac{\Tr\left(\Gamma_{1:k}\right)}{N}}.$
    \item $$\frac{\oc{c_RIP_lower}^2N\square^2}{\sqrt{2}\left(2\lambda + 3\Tr\left(\Gamma_{k+1:\infty}\right)\right)}>\norm{\tilde\Gamma_{1:k}^{-1/2}f_{1:k}^*}_\cH,$$ which is true if $\square^2 > \frac{\sqrt{2}}{\oc{c_RIP_lower}^2}\norm{\tilde\Gamma_{1:k}^{-1/2}f_{1:k}^*}_\cH\frac{2\lambda + \Tr\left(\Gamma_{k+1:\infty}\right)}{N}$.
\end{itemize}
We conclude that there exists an absolute constant $\nC\label{C_square_3} = \frac{24\sqrt{2}\oC{C_bX_f_star}\oC{C_noise_term}}{\oc{c_RIP_lower}^2}\kappa$, such that we can take
\begin{align*}
    \square^2 > \oC{C_square_3}\max\left\{ \triangle\sqrt{\sigma_1}\norm{\Gamma_{k+1:\infty}^{1/2}f_{k+1:\infty}^*}_\cH, \sigma_\xi\triangle\sqrt{\frac{\Tr\left(\Gamma_{1:k}\right)}{N}}, \frac{2\lambda + 3\Tr\left(\Gamma_{k+1:\infty}\right)}{N}\norm{\tilde\Gamma_{1:k}^{-1/2}f_{1:k}^*}_\cH \right\}.
\end{align*}Since $\norm{\tilde\Gamma_{1:k}^{-1/2}}_{\text{op}}=\triangle$, it suffices to choose
\begin{align*}
    \square^2 > \oC{C_square_3}\max\left\{ \triangle\sqrt{\sigma_1}\norm{\Gamma_{k+1:\infty}^{1/2}f_{k+1:\infty}^*}_\cH, \sigma_\xi\triangle\sqrt{\frac{\Tr\left(\Gamma_{1:k}\right)}{N}}, \frac{2\lambda + 3\Tr\left(\Gamma_{k+1:\infty}\right)}{N}\triangle\norm{f_{1:k}^*}_\cH \right\}.
\end{align*}
In the case where $\norm{\Gamma_{1:k}^{1/2}(\hat f_{1:k} - f_{1:k}^*)}_\cH\leq\square$ and $\norm{\hat f_{1:k} - f_{1:k}^*}_\cH=\triangle$, a similar analysis gives us that there exists an absolute constant $\nC\label{C_triangle_3}$ depending on $\kappa$ such that we can take
\begin{align}\label{eq:triangle_large_regularization}
    \triangle^2 > \oC{C_triangle_3}\max\left\{ \norm{f_{1:k}^*}_\cH^2, \frac{N\sigma_\xi^2\Tr\left(\Gamma_{1:k}\right)}{\left(4\lambda + \Tr\left(\Gamma_{k+1:\infty}\right)\right)^2}, \frac{\sigma_1N^2\norm{\Gamma_{k+1:\infty}^{1/2}f_{k+1:\infty}^*}_\cH^2}{\left(4\lambda + \Tr\left(\Gamma_{k+1:\infty}\right)\right)^2}\right\}.
\end{align}
Again, we choose that $\triangle = \square\sqrt{N/\left(\kappa_{DM}(4\lambda+\Tr\left(\Gamma_{k+1:\infty}\right))\right)}$. There exists an absolute constant $\nC\label{C_square_4}$ such that we can express $\square$ as in \eqref{eq:def_square_1}.
\begin{align*}
    \square &= \oC{C_square_4}\max\bigg\{ \sigma_\xi\sqrt{\frac{\Tr\left(\Gamma_{1:k}\right)}{4\lambda + \Tr\left(\Gamma_{k+1:\infty}\right)}}, \sqrt{\frac{\sigma_1 N}{4\lambda + \Tr\left(\Gamma_{k+1:\infty}\right)}}\norm{\Gamma_{k+1:\infty}^{1/2}f_{k+1:\infty}^*}_\cH,\\
    &\norm{f_{1:k}^*}_\cH\sqrt{\frac{4\lambda + \Tr\left(\Gamma_{k+1:\infty}\right)}{N}} \bigg\}.
\end{align*}
In particular, we check that for this choice, $\triangle$ satisfies \eqref{eq:triangle_large_regularization}.

\begin{Proposition}\label{prop:KRR_estimate}
    Under the assumption of Theorem~\ref{theo:upper_KRR}, there exist absolute constants $\oC{C_square_2}$, $\oC{C_noise}$ and $\oC{C_square_4}$ such that the following holds for all such $k$'s and all $\lambda\geq 0$. Recall the definition of $t(\square,\triangle)$ from \eqref{eq:def_t_square_triangle}, with probability at least $1-(\oC{C_noise}/\lfloor t(\square,\triangle)\rfloor)^{r/4} - \bP\left(\Omega_0^c\right)$ we have
    \begin{align*}
        \norm{\Gamma_{1:k}^{1/2}\left(\hat f_{1:k}-f_{1:k}^*\right)}_\cH\leq \square,\quad \norm{\hat f_{1:k}-f_{1:k}^*}_\cH\leq \square\sqrt{\frac{N}{\kappa_{DM}(4\lambda+\Tr\left(\Gamma_{k+1:\infty}\right))}},
    \end{align*}where
    \begin{enumerate}
        \item If $\sigma_1 N \leq \kappa_{DM}(4\lambda + \Tr\left(\Gamma_{k+1:\infty}\right))$,
        \begin{align*}
            \square &= \oC{C_square_2}\max\bigg\{ \sigma_\xi\sqrt{\frac{\Tr\left(\Gamma_{1:k}\right)}{4\lambda + \Tr\left(\Gamma_{k+1:\infty}\right)}}, \sqrt{\frac{\sigma_1 N}{4\lambda + \Tr\left(\Gamma_{k+1:\infty}\right)}}\norm{\Gamma_{k+1:\infty}^{1/2}f_{k+1:\infty}^*}_\cH,\\
            &\norm{f_{1:k}^*}_\cH\sqrt{\frac{4\lambda + \Tr\left(\Gamma_{k+1:\infty}\right)}{N}} \bigg\}.
        \end{align*}
        \item If $\sigma_1 N > \kappa_{DM}(4\lambda + \Tr\left(\Gamma_{k+1:\infty}\right))$,
        \begin{align*}
            \square &= \oC{C_square_4}\max\bigg\{ \sigma_\xi\sqrt{\frac{\left|J_1\right|}{N}}, \sigma_\xi\sqrt{\frac{\sum_{j\in J_2}\sigma_j}{  4\lambda + \Tr\left(\Gamma_{k+1:\infty}\right)}}, \norm{\Gamma_{k+1:\infty}^{1/2}f_{k+1:\infty}^*}_\cH,\\
            &{\norm{\tilde\Gamma_{1,\mathrm{thre}}^{-1/2}f_{1:k}^*}_\cH \frac{2\lambda + 3\Tr\left(\Gamma_{k+1:\infty}\right)}{N} } \bigg\}.
        \end{align*}
    \end{enumerate}
\end{Proposition}

\begin{Remark}
    Since we are using a concentration inequality rather than a deviation inequality, we cannot determine absolute constants, such as \(\oC{C_square_2}\) (for a discussion on concentration vs. deviation, see, for instance, \cite{klochkov_uniform_2020}). However, if we were to replace the concentration inequality with a deviation inequality, our proof could also establish absolute constants. For example, in our proof of \eqref{eq:bX_f_star_applied}, we utilized the concentration provided by Lemma~\ref{lem:moment-sum-variid-positive}. If we replace Lemma~\ref{lem:moment-sum-variid-positive} with Markov's inequality, then \eqref{eq:bX_f_star_applied} holds with a probability of at least \(1 - t^{-1}\), and the coefficient \(\oC{C_bX_f_star}\kappa\) on the right side of \eqref{eq:bX_f_star_applied} will be replaced by some \(t > 0\). Consequently, this additional deviation parameter \(t\) can be used to absorb the absolute constants accumulated during the proof, transferring them into the deviation probability and ensuring that the coefficients in the rate are explicit.
\end{Remark}

\subsubsection{Upper bound on $\norm{\Gamma_{k+1:\infty}^{1/2}(\hat f_{k+1:\infty}-f_{k+1:\infty}^*)}_\cH$}

We do not expect $\hat f _{k+1:\infty}$ to be a good estimator of $ f ^*_{k+1:\infty}$ because the minimum $\norm{\cdot}_\cH$-norm estimator $\hat f $ is using the 'remaining part' of $\cH$ endowed by the eigenfunctions $(\varphi_j)_{j\geq k+1}$ of $\Gamma$ (we denoted this space by $\cH_{k+1:\infty}$) to absorb the influence of noise introduced by $\bxi$ and not to estimate $ f_{k+1:\infty}^*$ which is why we call the error term $\norm{\Gamma_{k+1:\infty}^{1/2}(\hat f _{k+1:\infty} -  f_{k+1:\infty}^*)}_\cH$ a price for noise absorption instead of an estimation error. A consequence is that we can only upper bound this term by
\begin{equation*}
\norm{\Gamma_{k+1:\infty}^{1/2}(\hat f _{k+1:\infty} -  f_{k+1:\infty}^*)}_\cH\leq \norm{\Gamma_{k+1:\infty}^{1/2}\hat f _{k+1:\infty}}_\cH + \norm{\Gamma_{k+1:\infty}^{1/2} f_{k+1:\infty}^*}_\cH.
\end{equation*}We then just need to find a high probability upper bound on $\norm{\Gamma_{k+1:\infty}^{1/2}\hat f _{k+1:\infty}}_\cH$.

We have for $A:=\bX_{\phi,k+1:\infty}^\top (\bX_{\phi,k+1:\infty} \bX_{\phi,k+1:\infty}^\top + \lambda I_N)^{-1}$,
\begin{align}\label{eq:decomp_second_term}
\notag &\norm{\Gamma_{k+1:\infty}^{1/2}\hat f _{k+1:\infty}}_\cH= \norm{\Gamma_{k+1:\infty}^{1/2}A(\vy-\bX_{\phi,1:k}\hat f _{1:k})}_\cH\\ 
&\leq \norm{\Gamma_{k+1:\infty}^{1/2}A\bX_{\phi,1:k}( f_{1:k}^*-\hat f _{1:k})}_\cH + \norm{\Gamma_{k+1:\infty}^{1/2}A\bX_{\phi,k+1:\infty} f_{k+1:\infty}^*}_\cH + \norm{\Gamma_{k+1:\infty}^{1/2}A\bxi}_\cH
\end{align}and now we obtain high probability upper bounds on the three terms in \eqref{eq:decomp_second_term}.

On $\Omega_0$, for all $\blambda\in\bR^N$, 
\begin{equation}\label{eq:upper_dvoretzky}
\norm{\Gamma_{k+1:\infty}^{1/2} \bX_{\phi,k+1:\infty}^\top \blambda}_\cH\leq  \oC{C_DMU}\left(\sqrt{\Tr(\Gamma_{k+1:\infty}^2)}+\sqrt{N}\norm{\Gamma_{k+1:\infty}}_{\text{op}}\right)\norm{\blambda}_2.
\end{equation}Notice that this result holds without any extra assumption on $N$. We have
\begin{align}\label{eq:price_noise_inter_X1}
\begin{aligned}
     &\norm{\Gamma_{k+1:\infty}^{1/2}A\bX_{\phi,1:k}( f_{1:k}^*-\hat f _{1:k})}_\cH\\
     &= \norm{\Gamma_{k+1:\infty}^{1/2} \bX_{\phi,k+1:\infty}^\top (\bX_{\phi,k+1:\infty} \bX_{\phi,k+1:\infty}^\top + \lambda I_N)^{-1} \bX_{\phi,1:k}( f_{1:k}^*-\hat f _{1:k})}_\cH\\
&\leq \norm{\Gamma_{k+1:\infty}^{1/2} \bX_{\phi,k+1:\infty}^\top}_{\text{op}} \norm{(\bX_{\phi,k+1:\infty} \bX_{\phi,k+1:\infty}^\top + \lambda I_N)^{-1}}_{\text{op}} \norm{\bX_{\phi,1:k}( f_{1:k}^*-\hat f _{1:k})}_2\\
&\leq 4\oC{C_DMU}\oC{C_RIP_upper} \frac{\left(\sqrt{N\Tr(\Gamma_{k+1:\infty}^2)}+N \norm{\Gamma_{k+1:\infty}}_{\text{op}}\right)}{4\lambda + \Tr(\Gamma_{k+1:\infty})}\norm{\Gamma_{1:k}^{1/2}( f_{1:k}^*-\hat f _{1:k})}_\cH.
\end{aligned}
\end{align}

On $\Omega_0$,
\begin{align}\label{eq:upper_dvoretzky_signal}
\begin{aligned}
    &\norm{\Gamma_{k+1:\infty}^{1/2}A\bX_{\phi,k+1:\infty} f_{k+1:\infty}^*}_\cH\\
    &\leq \norm{\Gamma_{k+1:\infty}^{1/2} \bX_{\phi,k+1:\infty}^\top}_{\text{op}} \norm{(\bX_{\phi,k+1:\infty} \bX_{\phi,k+1:\infty}^\top+ \lambda I_N)^{-1}}_{\text{op}} \norm{\bX_{\phi,k+1:\infty} f_{k+1:\infty}^*}_2\\
&\leq 4\oC{C_DMU}\oC{C_bX_f_star} \kappa \frac{\sqrt{N\Tr(\Gamma_{k+1:\infty}^2)} + N \norm{\Gamma_{k+1:\infty}}_{\text{op}}}{4\lambda + \Tr(\Gamma_{k+1:\infty})}\norm{\Gamma_{k+1:\infty}^{1/2} f_{k+1:\infty}^*}_\cH.
\end{aligned}
\end{align}

Finally, let $D = \Gamma_{k+1:\infty}^{1/2}A$. As $\bX_{\phi,k+1:\infty}\Gamma_{k+1:\infty}\bX_{\phi,k+1:\infty}^\top:\bR^N\to\bR^N$, 
\[\Tr\left(\bX_{\phi,k+1:\infty}\Gamma_{k+1:\infty}\bX_{\phi,k+1:\infty}^\top\right) = \norm{\bX_{\phi,k+1:\infty}\Gamma_{k+1:\infty}^{1/2}}_{HS}^2 = \sum_{i=1}^N \norm{\left(\Gamma_{k+1:\infty}^{1/2}\phi_{k+1:\infty}\right)(X_i)}_\cH^2\]
is the sum of $N$ i.i.d. random variables appearing in \eqref{eq:Gamma_phi_applied}. On $\Omega_0$, $\sum_{i=1}^N \norm{\left(\Gamma_{k+1:\infty}^{1/2}\phi_{k+1:\infty}\right)(X_i)}_\cH^2\leq \oC{C_sum_Gamma_phi} N\Tr\left(\Gamma_{k+1:\infty}^2\right)$.
So
\begin{align}
&\Tr(DD^\top) = \Tr\left(D^\top D\right) \leq \frac{\Tr\left(\bX_{\phi,k+1:\infty}\Gamma_{k+1:\infty}\bX_{\phi,k+1:\infty}^\top\right)}{\norm{\bX_{\phi,k+1:\infty}\bX_{\phi,k+1:\infty}^\top+ \lambda I_N}_{\text{op}}^2} \leq \frac{16\oC{C_sum_Gamma_phi}N \Tr(\Gamma_{k+1:\infty}^2)}{\left(4\lambda + \Tr(\Gamma_{k+1:\infty})\right)^2}\label{eq:trace_DD_top_2}
\end{align}and
\begin{align}
\norm{D}_{\text{op}}
&= \norm{\Gamma_{k+1:\infty}^{1/2}\bX_{\phi,k+1:\infty}^\top (\bX_{\phi,k+1:\infty} \bX_{\phi,k+1:\infty}^\top+ \lambda I_N)^{-1}}_{\text{op}}\notag\\
&\leq \norm{\Gamma_{k+1:\infty}^{1/2}\bX_{\phi,k+1:\infty}^\top }_{\text{op}} \norm{(\bX_{\phi,k+1:\infty} \bX_{\phi,k+1:\infty}^\top + \lambda I_N)^{-1}}_{\text{op}}\notag\\
&\leq \frac{4\oC{C_DMU}}{4\lambda + \Tr(\Gamma_{k+1:\infty})}\left(\sqrt{\Tr(\Gamma_{k+1:\infty}^2)}+\sqrt{N}\norm{\Gamma_{k+1:\infty}}_{\text{op}}\right). \label{eq:op_D_2}
\end{align}
Set $k$ from Proposition~\ref{prop:noise_concentration} as
\begin{align}\label{eq:def_bar_p_xi}
    k = \left\lfloor\frac{\sqrt{\oC{C_sum_Gamma_phi}}}{\oC{C_DMU}}\frac{\sqrt{N \Tr(\Gamma_{k+1:\infty}^2)}}{\sqrt{\Tr(\Gamma_{k+1:\infty}^2)}+\sqrt{N}\norm{\Gamma_{k+1:\infty}}_{\text{op}}} \right\rfloor^2=: \oC{C_noise}(\bar p_{\xi})^{-4/r}.
\end{align}Then by Proposition~\ref{prop:noise_concentration}, with probability at least $1-\bar p_{\xi}-\bP(\Omega_0^c)$, we have
\begin{equation}\label{eq:upper_A_xi}
\begin{aligned}
\norm{\Gamma_{k+1:\infty}^{1/2}A\bxi}_\cH\leq \frac{3}{2}\sigma_\xi\frac{\sqrt{16\oC{C_sum_Gamma_phi}N \Tr(\Gamma_{k+1:\infty}^2)}}{4\lambda + \Tr(\Gamma_{k+1:\infty})}.
\end{aligned}
\end{equation}

Let us summarize the above discussion into the following Proposition:
\begin{Proposition}\label{prop:KRR_noise_absorption}
    Under the assumption of Theorem~\ref{theo:upper_KRR}. The following then holds for all such $k$'s and all $\lambda\geq 0$. With probability at least $1-\bar p_{\xi}-\bP(\Omega_0^c)$ we have
    \begin{align}
        &\norm{\hat f_{k+1:\infty}-f_{k+1:\infty}^*}_{L_2}\leq \oC{C_DMU}\frac{\left(\sqrt{N\Tr(\Gamma_{k+1:\infty}^2)}+N \norm{\Gamma_{k+1:\infty}}_{\text{op}}\right)}{4\lambda + \Tr(\Gamma_{k+1:\infty})}\norm{\Gamma_{1:k}^{1/2}( f_{1:k}^*-\hat f _{1:k})}_\cH \notag\\
        &+ 4\oC{C_DMU}\oC{C_sum_Gamma_phi}\kappa \frac{\sqrt{N\Tr(\Gamma_{k+1:\infty}^2)} + N \norm{\Gamma_{k+1:\infty}}_{\text{op}}}{4\lambda + \Tr(\Gamma_{k+1:\infty})}\norm{\Gamma_{k+1:\infty}^{1/2} f_{k+1:\infty}^*}_\cH + \norm{\Gamma_{k+1:\infty}^{1/2}f_{k+1:\infty}^*}_{\cH}\notag \\
        &+ \frac{3}{2}\sigma_\xi\frac{\sqrt{16\oC{C_sum_Gamma_phi}N \Tr(\Gamma_{k+1:\infty}^2)}}{4\lambda + \Tr(\Gamma_{k+1:\infty})}.\label{eq:KRR_noise_absorption_sigma}
    \end{align}
\end{Proposition}

\subsection{Proof of Theorem~\ref{theo:main_upper_k>N} (the case $k$ is not necessarily smaller than $N$)}\label{sec:proof_main_upper_k>N}



We divide the proof into two cases: either $\sigma(\square, \triangle) = \square$ or $\sigma(\square, \triangle) = \triangle\sqrt{\sigma_1}$. Within each case, we further examine two specific scenarios, as discussed in Section~\ref{sec:proof_main_upper_k<N}, labeled as case [1] and case [2].

In each of these two cases, we will prove separately that $\cQ_{f_{1:k}} > \cM_{f_{1:k}}$ and $\cR_{f_{1:k}} > \cM_{f_{1:k}}$, thereby showing that $\cL_{f_{1:k}} > 0$ when $f_{1:k}\notin f_{1:k}^*+B$. This implies that $\hat f_{1:k}$ cannot be outside of $f_{1:k}^*+B$. See Section~\ref{sec:proof_main_upper_k<N} for a more detailed explanation. The main difference with the proof of Theorem~\ref{theo:upper_KRR} is that $\bX_{\phi,1:k}$ no longer acts as an isomorphism on $\cH_{1:k}$ but it may act as an RIP on a cone defined in \eqref{eq:def_cone_RIP}.

\subsubsection{Stochastic Argument}

The stochastic argument of the proof of Theorem~\ref{theo:main_upper_k>N} is almost the same as that of Theorem~\ref{theo:upper_KRR}, except that we make use of an \emph{restricted} isomorphy property on the $\mathrm{cone}\left(\cC(R_N(\oc{c_kappa_RIP}))\right)$ (see \eqref{eq:def_cone_RIP}) instead of the isomorphy property on $\cH_{1:k}$. As a consequence, we shall prove that a vector (to be defined later) belongs to this cone. This requires some extra work in addition to the proof of Theorem~\ref{theo:upper_KRR}.

In this subsection, we prove that there exist absolute constants $\oc{c_RIP_lower}$, $\oC{C_RIP_upper}$, $\oc{c_kappa_RIP}$, $\nC\label{C_DMU_k>M}$, $\oC{C_bX_f_star}$ and $\nC\label{C_sum_Gamma_phi_infty}$, such that the following random event $\Omega'$ happens with high probability:
\begin{itemize}
    \item for any $\vlambda\in\bR^N$,
    \begin{align}\label{eq:DM_k>M_applied}
        \left(\frac{1}{2}\Tr(\Gamma_{k+1:\infty})+\lambda\right)\norm{\blambda}_2\leq \norm{\left(\bX_{\phi,k+1:\infty}\bX_{\phi,k+1:\infty}^\top+\lambda I\right)\blambda}_2\leq \left(\frac{3}{2}\Tr(\Gamma_{k+1:\infty})+\lambda\right)\norm{\blambda}_2,
    \end{align}
    \item for any $f_{1:k}\in\mathrm{cone}\left(\cC(R_N(\oc{c_kappa_RIP}))\right)$,
    \begin{align}\label{eq:RIP_k>M_applied}
        \oc{c_RIP_lower}\norm{\Gamma_{1:k}^{1/2}f}_\cH \leq \frac{1}{\sqrt{N}}\norm{\bX_{\phi,1:k}f}_2 \leq \oC{C_RIP_upper}\norm{\Gamma_{1:k}^{1/2}f}_\cH,
    \end{align}
    \item for any $\vlambda\in\bR^N$,
    \begin{align}\label{eq:Dvreotzky_upper_k>M_applied}
        \norm{\tilde\Gamma_{1:k}^{-1/2} \bX_{\phi,1:k}^\top \vlambda }_\cH \leq \oC{C_DMU_k>M}\left(\sqrt{ \square^2 \left|J_1\right| + \triangle^2\sum_{j\in J_2}\sigma_j } + \sqrt{N}\sigma(\square,\triangle)\right)\norm{\vlambda}_2,
    \end{align}where $\triangle = \square\sqrt{\frac{N}{\kappa_{DM}(4\lambda+\Tr\left(\Gamma_{k+1:\infty}\right))}}$,
      \item for all $\blambda\in\bR^N$,
  \begin{align}\label{eq:DM_upper_k>N_applied}
      \norm{\Gamma_{k+1:\infty}^{1/2}\bX_{\phi,k+1:\infty}^\top\blambda}_\cH\leq \oC{C_DMU}\left(\sqrt{\Tr\left(\Gamma_{k+1:\infty}^2\right)} + \sqrt{N}\norm{\Gamma_{k+1:\infty}}_{\text{op}}\right)\norm{\blambda}_2,
  \end{align}
    \item \begin{align}\label{eq:bX_f_star_applied_k>N}
        \norm{\bX_{\phi,k+1:\infty}f_{k+1:\infty}^*}_2\leq \oC{C_bX_f_star}\kappa\sqrt{N}\norm{\Gamma_{k+1:\infty}^{1/2}f_{k+1:\infty}^*}_\cH,
    \end{align}
    \item \begin{align}\label{eq:Gamma_phi_applied_k>N}
        \sum_{i=1}^N \norm{\left(\Gamma_{k+1:\infty}^{1/2}\phi_{k+1:\infty}\right)(X_i)}_\cH^2\leq \oC{C_sum_Gamma_phi_infty}N\Tr\left(\Gamma_{k+1:\infty}^2\right).
    \end{align}
\end{itemize}

\begin{Proposition}\label{prop:stochastic_argument_k>N}
    Under the assumption of Theorem~\ref{theo:main_upper_k>N}, there exist absolute constants $\oc{c_RIP_lower}$, $\oC{C_RIP_upper}$, $\oc{c_kappa_RIP}$, $\oC{C_DMU_k>M}$, $\oC{C_bX_f_star}$ and $\oC{C_sum_Gamma_phi_infty}$ such that we have $\bP(\Omega')\geq 1-\bar p_{DM}-2\bar p_{DMU}-\bar p_{RIP}-\ogamma{gamma_DMU_L2} - \frac{ \oc{c_P_bX_f_star} }{N}$.
\end{Proposition}
\beginproof
\begin{itemize}
    \item Notice that Assumption~\ref{assumption:DM_L4_L2} is granted in Theorem~\ref{theo:main_upper_k>N}. As in the proof of \eqref{eq:DM_applied}, there exists an absolute constant $\oc{c_kappa_DM}$ such that when $N\leq\oc{c_kappa_DM}\kappa_{DM}d_\lambda^*\left(\Gamma_{k+1:\infty}^{-1/2}B_\cH\right)$ (which is exactly one of the assumptions of Theorem~\ref{theo:main_upper_k>N}), then $\tilde\delta$ defined in \eqref{eq:def_tilde_delta} is strictly smaller than $1/2$. \eqref{eq:DM_k>M_applied} thus follows by applying Theorem~\ref{theo:DM_RKHS}.

    \item \eqref{eq:RIP_k>M_applied} follows by Proposition~\ref{prop:RIP}. In fact, we can take $\odelta{delta_P_RIP}$ from Proposition~\ref{prop:RIP} to be $1/100$ (thus $\bar p_{RIP}=1/100$) and $\oc{c_kappa_RIP}<(10^6\oC{C_Rudelson}^2\oC{C_estimate_gamma_infty}^2)^{-1}$, then \eqref{eq:RIP_k>M_applied} follows with $\oc{c_RIP_lower}=1-\frac{\sqrt{2}}{10}$, $\oC{C_RIP_upper}=1+\frac{\sqrt{2}}{10}$ with probability at least $1-\bar p_{RIP}$.
    
    \item Proving \eqref{eq:Dvreotzky_upper_k>M_applied} makes use of Assumption~\ref{assumption:DMU_used_for_RIP}.    Replacing $\Gamma_{k+1:\infty}$ by $\tilde\Gamma_{1:k}$ and $\bX_{\phi,k+1:\infty}$ by $\bX_{\phi,1:k}$ in \eqref{eq:proof_upper_DM} from Section~\ref{sec:proof_upper_dvoretzky} leads to the fact that \eqref{eq:Dvreotzky_upper_k>M_applied} with probability at least $1-\bar p_{DMU}$ with $\bar p_{DMU}=\frac{\oc{c_P_DMU}}{N} + \ogamma{gamma_RIP_k>N}$, and $\oC{C_DMU_k>M} = \oC{C_DMU}$. Moreover, recalling the definition of $\tilde\Gamma_{1:k}^{1/2}$ (see \eqref{eq:def_tilde_Gamma}) the spectrum of $\tilde\Gamma_{1:k}^{-1/2}\Gamma_{1:k}^{1/2}$ is
\begin{align*}
    \frac{\square\sqrt{\sigma_1}}{\max\{\sqrt{\sigma_1}, \square/\triangle\}},\cdots,\frac{\square\sqrt{\sigma_k}}{\max\{\sqrt{\sigma_k}, \square/\triangle\}},0,\cdots.
\end{align*}Squaring everything and summing them together gives us the trace of $\tilde\Gamma_{1:k}^{-1}\Gamma_{1:k}$. Therefore,
\begin{align}\label{eq:trace_tilde_Gamma_Gamma}
    \Tr\left(\tilde\Gamma_{1:k}^{-1}\Gamma_{1:k}\right) = \sum_{j\in J_1}\square^2 + \triangle^2 \sum_{j\in J_2}\sigma_j = \square^2\left|J_1\right| + \triangle^2\sum_{j\in J_2}\sigma_j.
\end{align}

    \item \eqref{eq:DM_upper_k>N_applied} and \eqref{eq:bX_f_star_applied_k>N} follow from the same idea as in the proof of Proposition~\ref{prop:stochastic_argument_upper_bound} as \eqref{eq:DM_k>M_applied}, \eqref{eq:DM_upper_k>N_applied} and \eqref{eq:bX_f_star_applied_k>N} do not depend on $\cH_{1:k}$. As a result, \eqref{eq:DM_upper_k>N_applied} and \eqref{eq:bX_f_star_applied_k>N} hold with probability at least $1 -\bar p_{DMU}- \frac{ \oc{c_P_bX_f_star} }{N}$.
\item By Assumption~\ref{assumption:upper_dvoretzky}, \eqref{eq:Gamma_phi_applied_k>N} holds with probability at least $1-\ogamma{gamma_DMU_L2}$ with constant $\oC{C_sum_Gamma_phi_infty} = 1+ \odelta{delta_DMU_L2}.$
\end{itemize}

\endproof

Up to the end of the proof, we put ourselves on the event $\Omega'$.

\subsubsection{Deterministic Argument}

As in the proof of Theorem~\ref{theo:upper_KRR}, the overall idea is to decompose $\hat f_\lambda$ into two components: $\hat f_{1:k}+\hat f_{k+1:\infty}$, used respectively to estimate $f_{1:k}^*$ and to absorb noise. Let us remind the ideas briefly.
\begin{itemize}
    \item The estimation property of $\hat f_{1:k}$: We employ the idea from Section~\ref{sec:proof_main_upper_k<N}, that is, to decompose $\partial B$ into two categories: $\norm{\Gamma_{1:k}^{1/2}(f_{1:k}-f_{1:k}^*)}_\cH=\square, \norm{f_{1:k}-f_{1:k}^*}_\cH\leq\triangle$ and $\norm{\Gamma_{1:k}^{1/2}(f_{1:k}-f_{1:k}^*)}_\cH\leq\square, \norm{f_{1:k}-f_{1:k}^*}_\cH=\triangle$. In each case, we prove that either we have $\cQ_{f_{1:k}}>\cM_{f_{1:k}}$ or we have $\cR_{f_{1:k}}>\cM_{f_{1:k}}$ as long as $\square,\triangle$ are larger than its correct level. These indicate that $\cL_{f_{1:k}}>0$ (refer to \eqref{eq:new_Q+M+R}). By the definition of $\hat f_{1:k}$, this means that $\hat f_{1:k}\notin f_{1:k}^* +B$. Necessarily, we must have $\norm{\Gamma_{1:k}^{1/2}(f_{1:k}-f_{1:k}^*)}_\cH\leq\square, $ and $ \norm{f_{1:k}-f_{1:k}^*}_\cH\leq\triangle$.
    \item The price for noise absorption by $\hat f_{k+1:\infty}$ can be quantified. We leave this part for the corresponding paragraph.
\end{itemize}

\paragraph{Estimation property of $\hat f_{1:k}$} 

For proving $\cQ_{f_{1:k}}>\cM_{f_{1:k}}$ or $\cR_{f_{1:k}}>\cM_{f_{1:k}}$, we must bound $\cM_{f_{1:k}}$ from above and $\cQ_{f_{1:k}}$ from below. As discussed in Section~\ref{sec:proof_main_upper_k<N} (refer to \eqref{eq:multiplier_decomp}), we need to establish the analog of Lemma~\ref{lem:first_beta_2_star_term} and Lemma~\ref{lem:noise_term} in the case where $k\gtrsim N$, that is, to find high-probability upper bounds for the following two quantities:
\begin{align*}
    &\norm{\tilde\Gamma_{1:k}^{-1/2}\bX_{\phi,1:k}^\top\left( \bX_{\phi,k+1:\infty}\bX_{\phi,k+1:\infty}^\top + \lambda I_N \right)^{-1}\bxi}_\cH,\mbox{ and }\\
&\norm{ \tilde\Gamma_{1:k}^{-1/2} \bX_{\phi,1:k}^\top\left( \bX_{\phi,k+1:\infty}\bX_{\phi,k+1:\infty}^\top + \lambda I_N \right)^{-1}\bX_{\phi,k+1:\infty}f_{k+1:\infty}^* }_\cH.
\end{align*}
The difference between the case where $k\lesssim N$ is that we no longer have the isomorphy property (\eqref{eq:RIP_applied}), thus we have to deal with $\norm{\tilde\Gamma_{1:k}^{-1/2}\bX_{\phi,1:k}^\top}_{\text{op}}$ by using another method. Fortunately, \eqref{eq:Dvreotzky_upper_k>M_applied} exactly serves as the tool that we want. One issue remains: there are two terms in the upper bound of \eqref{eq:Dvreotzky_upper_k>M_applied}, namely, $\sqrt{ \square^2 \left|J_1\right| + \triangle^2\sum_{j\in J_2}\sigma_j }$ and $\sqrt{N}\sigma(\square,\triangle)$. In fact, we have already discussed in \cite[Remark 2]{lecue_geometrical_2022} that the only interesting case is when the latter term dominates, that is,
\begin{align}\label{eq:another_assumption_k>N}
    \sqrt{ \square^2 \left|J_1\right| + \triangle^2\sum_{j\in J_2}\sigma_j }\leq \sqrt{N}\sigma(\square,\triangle).
\end{align}Further, recall the definition of $\square/\triangle$ from Section~\ref{sec:proof_main_upper_k<N}, that is,
\begin{align*}
    \frac{\square}{\triangle} = \sqrt{\frac{ \kappa_{DM}\left(4\lambda + \Tr\left(\Gamma_{k+1:\infty}\right)\right) }{N}}.
\end{align*} \eqref{eq:another_assumption_k>N} is equivalent to \eqref{eq:extra_condition_k>N} in the assumption of Theorem~\ref{theo:main_upper_k>N}. 

\begin{itemize}
    \item Lower bound of $\cQ_{f_{1:k}}$ (Analog of \eqref{eq:quadratic_k<N}.)

    As explained in Section~\ref{sec:RIP}, the kernel design matrix $\bX_{\phi,1:k}$ does not behave like an isomorphy on the entire RKHS $\cH_{1:k}$ where $k\gtrsim N$, but only on a restricted cone $\cC$ (for the sake of simplicity, we denote $\cC$ for $\mathrm{cone}(\cC(R_N(\oc{c_kappa_RIP}))$) as defined in \eqref{eq:def_cone_RIP}. We therefore have no choice but to assume that $f_{1:k}-f_{1:k}^*\in\cC$. To ensure $f_{1:k}-f_{1:k}^*\in\cC$, a sufficient condition is: $\square/\triangle\geq R_N^*(\oc{c_kappa_RIP})$. This is because under this assumption, and when $\norm{\Gamma_{1:k}^{1/2}( f_{1:k}-f_{1:k}^*)}_\cH=\square$, $\norm{ f_{1:k}-f_{1:k}^*}_\cH\leq\triangle$, we have
\begin{align*}
    R_N^*(\oc{c_kappa_RIP})\norm{f_{1:k}-f_{1:k}^*}_\cH\leq R_N^*(\oc{c_kappa_RIP})\triangle = R_N^*(\oc{c_kappa_RIP})\frac{\triangle}{\square}\norm{\Gamma_{1:k}^{1/2}(f_{1:k} - f_{1:k}^*)}_\cH \leq \norm{\Gamma_{1:k}^{1/2}(f_{1:k} - f_{1:k}^*)}_\cH.
\end{align*}This indicates that $f_{1:k}-f_{1:k}^*\in\cC$. Due to this fact, on $\Omega'$ we have
    \begin{align}\label{eq:quadratic_k>N}
    \begin{aligned}
        \cQ_{f_{1:k}} &= \norm{ \left(\bX_{\phi,k+1:\infty}\bX_{\phi,k+1:\infty}^\top + \lambda I_N\right)^{-1/2}\bX_{\phi,1:k}(f_{1:k}-f_{1:k}^*) }_2^2 \\
        &\geq \frac{4}{4\lambda + \Tr\left(\Gamma_{k+1:\infty}\right)}\norm{\bX_{\phi,1:k}(f_{1:k}-f_{1:k}^*)}_2^2\notag\\
    &\geq \frac{4\oc{c_RIP_lower}^2N\square^2}{4\lambda + \Tr\left(\Gamma_{k+1:\infty}\right)},
    \end{aligned}
\end{align}which is the same as \eqref{eq:quadratic_k<N} in Section~\ref{sec:proof_main_upper_k<N}.
\item Analog of \eqref{eq:multiplier_bX_f_star_k<N}.
On $\Omega'$ together with \eqref{eq:another_assumption_k>N}, we have
\begin{align*}
    &\norm{\tilde\Gamma_{1:k}^{-1/2} \bX_{\phi,1:k}^\top\left( \bX_{\phi,k+1:\infty}\bX_{\phi,k+1:\infty}^\top + \lambda I_N\right)^{-1}\bX_{\phi,k+1:\infty}f_{k+1:\infty}^* }_\cH \\
    &\leq \norm{\tilde\Gamma_{1:k}^{-1/2} \bX_{\phi,1:k}^\top}_{\text{op}} \norm{\left( \bX_{\phi,k+1:\infty}\bX_{\phi,k+1:\infty}^\top + \lambda I_N\right)^{-1}}_{\text{op}}\norm{\bX_{\phi,k+1:\infty}f_{k+1:\infty}^*}_\cH \\
    &\leq 2\oC{C_DMU_k>M}\sqrt{N}\sigma(\square,\triangle)\frac{4}{4\lambda + \Tr\left(\Gamma_{k+1:\infty}\right)}\oC{C_bX_f_star}\norm{\Gamma_{k+1:\infty}^{1/2}f_{k+1:\infty}^*}_\cH,
\end{align*} which is also equivalent to Lemma~\ref{lem:first_beta_2_star_term} up to universal constants.
\item Analog of \eqref{eq:multiplier_bxi_k<N}. Let $$D = \tilde\Gamma_{1:k}^{-1/2}\bX_{\phi,1:k}^\top\left( \bX_{\phi,k+1:\infty}\bX_{\phi,k+1:\infty}^\top + \lambda I_N \right)^{-1},$$ thus $DD^\top:\bR^N\to\bR^N$. Let us recall the following assumption in Theorem~\ref{theo:main_upper_k>N}: with probability at least $1-\ogamma{gamma_RIP_k>N}$, $\max_{i\in[N]}\norm{\tilde\Gamma_{1:k}^{-1/2}\phi_{1:k}(X_i)}_\cH^2\leq (1+\odelta{delta_RIP_k>N})\Tr\left(\tilde\Gamma_{1:k}\Gamma_{1:k}\right)$. This indicates that with probability at least $1-\bP((\Omega'')^c)-\ogamma{gamma_RIP_k>N}$, $\sum_{i=1}^N\norm{\tilde\Gamma_{1:k}^{-1/2}\phi_{1:k}(X_i)}_\cH^2\leq (1+\odelta{delta_RIP_k>N})N\Tr\left(\tilde\Gamma_{1:k}\Gamma_{1:k}\right)$. Combining this observation with \eqref{eq:trace_tilde_Gamma_Gamma}, we have:
\begin{align*}
    \sqrt{\Tr\left(DD^\top\right)}&= \sqrt{\Tr\left(D^\top D\right)}\leq  \norm{\left(\bX_{\phi,k+1:\infty}\bX_{\phi,k+1:\infty}^\top + \lambda I_N \right)^{-1}}_{\text{op}}\sqrt{\sum_{i=1}^N\norm{\tilde\Gamma_{1:k}^{-1/2}\phi_{1:k}(X_i)}_\cH^2} \\
    &\leq\frac{4\sqrt{1+\odelta{delta_RIP_k>N}}\sqrt{N}}{4\lambda + \Tr\left(\Gamma_{k+1:\infty}\right)}\sqrt{\left|J_1\right|\square^2 + \triangle^2\sum_{j\in J_2}\sigma_j}.
\end{align*}Moreover,
\begin{align*}
    \norm{D}_{\text{op}} &= \norm{D^\top}_{\text{op}} \leq \norm{\left(\bX_{\phi,k+1:\infty}\bX_{\phi,k+1:\infty}^\top + \lambda I_N \right)^{-1}}_{\text{op}} \norm{\bX_{\phi,1:k}\tilde\Gamma_{1:k}^{-1/2}}_{\text{op}}\\
    &\leq 2\oC{C_DMU_k>M}\sqrt{N}\sigma(\square,\triangle)\cdot\frac{4}{4\lambda+\Tr\left(\Gamma_{k+1:\infty}\right)}\\
    &\leq \frac{8\oC{C_DMU_k>M}\sqrt{N}\sigma\left(\square,\triangle\right)}{4\lambda + \Tr\left(\Gamma_{k+1:\infty}\right)}.
\end{align*}The upper bounds on $\sqrt{\Tr\left(DD^\top\right)}$ and $\norm{D}_{\text{op}}$ are equivalent to those in Lemma~\ref{lem:noise_term} up to universal constants. Therefore,  by Proposition~\ref{prop:noise_concentration} with the $k$ from Proposition~\ref{prop:noise_concentration} set as
\begin{align*}
    k = \biggl\lfloor\frac{\left(1+\odelta{delta_RIP_k>N}\right)^2\left(\left|J_1\right| \square^2 + \triangle^2 \sum_{j\in J_2}\sigma_j\right)}{\oC{C_DMU_k>M}^2\sigma^2(\square, \triangle)}\biggr\rfloor,
\end{align*}
with probability at least $1-\bP\left((\Omega')^c\right)-\left(\oC{C_noise}/\lfloor t(\square,\triangle)\rfloor\right)^{r/4}$,
\begin{align*}
    \norm{\tilde\Gamma_{1:k}^{-1/2} \bX_{\phi,1:k}^\top\left( \bX_{\phi,k+1:\infty}\bX_{\phi,k+1:\infty}^\top + \lambda I_N\right)^{-1} \bxi }_\cH \leq \frac{3}{2}\frac{12\sqrt{1+\odelta{delta_RIP_k>N}}\sqrt{N}\sigma_\xi}{4\lambda + \Tr\left(\Gamma_{k+1:\infty}\right)}\sqrt{\left|J_1\right|\square^2 + \triangle^2\sum_{j\in J_2}\sigma_j}.
\end{align*}
\end{itemize}

The above analysis indicates that the analysis of the estimation property of $\hat f_{1:k}$ in the case $k\leq\oc{c_RIP} N$ also holds in the case where $k\gtrsim N$. We therefore have the following proposition:
\begin{Proposition}\label{prop:estimation_k>N}
    Suppose the assumptions of Theorem~\ref{theo:main_upper_k>N} hold. There exist absolute constants $\nC\label{C_square_k>N}$, $\nC\label{C_square_k>N_2}$ such that the following holds. Assume that there exists $k\in\bN\cup\{\infty\}$ such that $N\leq\oc{c_kappa_DM}\kappa_{DM}d_\lambda^*\left(\Gamma_{k+1:\infty}^{-1/2}B_\cH\right)$, \eqref{eq:extra_condition_k>N} holds, and $\kappa_{DM}\left(4\lambda + \Tr\left(\Gamma_{k+1:\infty}\right)\right)\geq N\left(R_N^*(\oc{c_kappa_RIP})\right)^2$, The following then holds for all such $k$'s. With probability at least $1-\bP\left((\Omega')^c\right)-\left(\oC{C_noise}/\lfloor t(\square,\triangle)\rfloor\right)^{r/4}$,
    \begin{align*}
        \norm{\Gamma_{1:k}^{1/2}(\hat f_{1:k} - f_{1:k}^*)}_\cH \leq \square \mbox{ and } \norm{\hat f_{1:k} - f_{1:k}^*}_\cH \leq \square\sqrt{\frac{\kappa_{DM}\left(4\lambda + \Tr\left(\Gamma_{k+1:\infty}\right)\right)}{N}},
    \end{align*}where
    \begin{itemize}
        \item if $\sigma_1 N \leq \kappa_{DM}\left(4\lambda+\Tr\left(\Gamma_{k+1:\infty}\right)\right)$,
        \begin{align*}
            \square &= \oC{C_square_k>N}\max\bigg\{ \sigma_\xi\sqrt{\frac{\Tr\left(\Gamma_{1:k}\right)}{4\lambda + \Tr\left(\Gamma_{k+1:\infty}\right)}}, \sqrt{\frac{\sigma_1 N}{4\lambda + \Tr\left(\Gamma_{k+1:\infty}\right)}}\norm{\Gamma_{k+1:\infty}^{1/2}f_{k+1:\infty}^*}_\cH,\\
            &\norm{f_{1:k}^*}_\cH\sqrt{\frac{4\lambda + \Tr\left(\Gamma_{k+1:\infty}\right)}{N}} \bigg\}
        \end{align*}
        \item if $\sigma_1 N < \kappa_{DM}\left(4\lambda+\Tr\left(\Gamma_{k+1:\infty}\right)\right)$,
        \begin{align*}
            \square &= \oC{C_square_k>N_2} \max\bigg\{ \sigma_\xi\sqrt{\frac{\left|J_1\right|}{N}}, \sigma_\xi\sqrt{\frac{\sum_{j\in J_2}\sigma_j}{  4\lambda + \Tr\left(\Gamma_{k+1:\infty}\right)}}, \norm{\Gamma_{k+1:\infty}^{1/2}f_{k+1:\infty}^*}_\cH,\\
            &{\norm{\tilde\Gamma_{1,\mathrm{thre}}^{-1/2}f_{1:k}^*}_\cH \frac{2\lambda + 3\Tr\left(\Gamma_{k+1:\infty}\right)}{N} } \bigg\}.
        \end{align*}
    \end{itemize}

\end{Proposition}

\paragraph{Price for noise absorption} As in Section~\ref{sec:proof_main_upper_k<N}, we do not expect $\hat f_{k+1:\infty}$ to be a good estimator of $f_{k+1:\infty}^*$. We therefore use the simple decomposition:
\begin{equation*}
\norm{\Gamma_{k+1:\infty}^{1/2}(\hat f _{k+1:\infty} -  f_{k+1:\infty}^*)}_\cH\leq \norm{\Gamma_{k+1:\infty}^{1/2}\hat f _{k+1:\infty}}_\cH + \norm{\Gamma_{k+1:\infty}^{1/2} f_{k+1:\infty}^*}_\cH.
\end{equation*}We then just need to find a high probability upper bound on $\norm{\Gamma_{k+1:\infty}^{1/2}\hat f _{k+1:\infty}}_\cH$.

We have for $A:=\bX_{\phi,k+1:\infty}^\top (\bX_{\phi,k+1:\infty} \bX_{\phi,k+1:\infty}^\top + \lambda I_N)^{-1}$,
\begin{align*}
&\norm{\Gamma_{k+1:\infty}^{1/2}\hat f _{k+1:\infty}}_\cH= \norm{\Gamma_{k+1:\infty}^{1/2}A(\vy-\bX_{\phi,1:k}\hat f _{1:k})}_\cH\\ 
&\leq \norm{\Gamma_{k+1:\infty}^{1/2}A\bX_{\phi,1:k}( f_{1:k}^*-\hat f _{1:k})}_\cH + \norm{\Gamma_{k+1:\infty}^{1/2}A\bX_{\phi,k+1:\infty} f_{k+1:\infty}^*}_\cH + \norm{\Gamma_{k+1:\infty}^{1/2}A\bxi}_\cH
\end{align*}and now we obtain high probability upper bounds on the three terms. The second and the third term are dealt with in exactly the same way as before, as they do not depend on $\bX_{\phi,1:k}$. Our approach with the first term has to be different, however. From the decomposition of the excess risk and from the fact that $\cL_{\hat f_{1:k}}\leq 0$, it follows that $\cQ_{\hat f_{1:k}} + \cR_{\hat f_{1:k}}\leq\left|\cM_{\hat f_{1:k}}\right|$. We have also seen that
\begin{align}\label{eq:lower_bound_of_quadratic_k>N}
    \cQ_{\hat f_{1:k}} \geq \frac{4}{4\lambda + \Tr\left(\Gamma_{k+1:\infty}\right)}\norm{\bX_{\phi,1:k}(\hat f_{1:k}- f_{1:k}^*)}_2^2.
\end{align}We would now like to obtain a bound for $\left|\cM_{\hat f_{1:k}}\right|$, so for all $f\in \cH_{1:k}$, we define
\begin{align}\label{eq:def_vertiiii}
    \vertiiii{f} := \max\left\{ \norm{\Gamma_{1:k}^{1/2}f}_\cH,\, \sqrt{\frac{\kappa_{DM}\left(4\lambda + \Tr\left(\Gamma_{k+1:\infty}\right)\right)}{N}}\norm{f}_\cH \right\}.
\end{align}By convention, we set $\vertiiii{f}=0$ for $f\in \cH_{k+1:\infty}$. By Proposition~\ref{prop:estimation_k>N}, it follows that on $\Omega'$, we have $\vertiiii{\hat f_{1:k}-f_{1:k}^*}\leq\square$. Therefore if we define
\begin{align}\label{eq:def_theta_k>N}
    \theta:= \sup\left( \left|\cM_{f_{1:k}}\right|:\, \vertiiii{f_{1:k} - f_{1:k}^*}\leq 1 \right),
\end{align}then it follows from the inequality $\cQ_{\hat f_{1:k}}\leq \cQ_{\hat f_{1:k}} + \cR_{\hat f_{1:k}}\leq\left|\cM_{\hat f_{1:k}}\right|$ that
\begin{align*}
    &\frac{4}{4\lambda + \Tr\left(\Gamma_{k+1:\infty}\right)}\norm{\bX_{\phi,1:k}\left(\hat f_{1:k} - f_{1:k}^*\right)}_2^2\\
    &\leq \norm{\hat f_{1:k}-f_{1:k}^*}_\cH^2 + \frac{4}{4\lambda + \Tr\left(\Gamma_{k+1:\infty}\right)}\norm{\bX_{\phi,1:k}\left(\hat f_{1:k} - f_{1:k}^*\right)}_2^2\\
    &\leq \theta\square.
\end{align*}Therefore, $\norm{\bX_{\phi,1:k}(\hat f_{1:k} - f_{1:k}^*)}_2\leq \frac{1}{2}\sqrt{(4\lambda +\Tr\left(\Gamma_{k+1:\infty}\right))\theta\square}$.
Since
\begin{align*}
    \frac{\square}{\triangle} = \sqrt{\frac{ \kappa_{DM}\left(4\lambda + \Tr\left(\Gamma_{k+1:\infty}\right)\right) }{N}},
\end{align*}it follows that $\square\vertiii{f} = \vertiiii{f}$. Hence if $\vertiiii{\cdot}^*$ is the dual norm of $\vertiiii{\cdot}$ and $\vertiii{\cdot}^*$ the dual norm of $\vertiii{\cdot}$, then $\vertiii{\cdot}^*=\square\vertiiii{\cdot}^*$. This implies that $\vertiiii{\cdot}^*$ is equivalent to $\square^{-1}\norm{\tilde\Gamma_{1:k}^{-1/2}\cdot}_\cH$. From the definition of $\theta$ it follows that
\begin{align*}
    &\frac{\theta}{2} = \sup\bigg( \left|\left\langle \bX_{\phi,1:k}\left(\bX_{\phi,k+1:\infty}\bX_{\phi,k+1:\infty}^\top + \lambda I_N\right)^{-1}\left(\bX_{\phi,k+1:\infty}f_{k+1:\infty}^*+\bxi\right)-f_{1:k}^*, f_{1:k}-f_{1:k}^*\right\rangle_\cH\right|\\
    &:\, \vertiiii{f_{1:k}-f_{1:k}^*}\leq 1 \bigg)\\
    &=\vertiiii{ \bX_{\phi,1:k}\left(\bX_{\phi,k+1:\infty}\bX_{\phi,k+1:\infty}^\top + \lambda I_N\right)^{-1}\left(\bX_{\phi,k+1:\infty}f_{k+1:\infty}^*+\bxi\right)-f_{1:k}^* }^*\\
    &\leq \vertiiii{f_{1:k}^*}^* + \vertiiii{ \bX_{\phi,1:k}\left(\bX_{\phi,k+1:\infty}\bX_{\phi,k+1:\infty}^\top + \lambda I_N\right)^{-1} \bX_{\phi,k+1:\infty}f_{k+1:\infty}^*}^*\\
    &+ \vertiiii{\bX_{\phi,1:k}\left(\bX_{\phi,k+1:\infty}\bX_{\phi,k+1:\infty}^\top + \lambda I_N\right)^{-1}\bxi }^*\\
    &\lesssim \square^{-1}\norm{\tilde\Gamma_{1:k}^{-1/2}f_{1:k}^*}_\cH + \square^{-1}\norm{\tilde\Gamma_{1:k}^{-1/2}\bX_{\phi,1:k}\left(\bX_{\phi,k+1:\infty}\bX_{\phi,k+1:\infty}^\top + \lambda I_N\right)^{-1} \bX_{\phi,k+1:\infty}f_{k+1:\infty}^* }_\cH \\
    &+ \square^{-1}\norm{\tilde\Gamma_{1:k}^{-1/2} \bX_{\phi,1:k}\left(\bX_{\phi,k+1:\infty}\bX_{\phi,k+1:\infty}^\top + \lambda I_N\right)^{-1}\bxi }_\cH,
\end{align*} which are exactly the quantities studied in the last paragraph, and we have shown that they are bounded from above by the right-hand-side of \eqref{eq:lower_bound_of_quadratic_k>N}, and $\norm{\bX_{\phi,1:k}(\hat f_{1:k} - f_{1:k}^*)}_2^2\leq \oC{C_RIP_upper}N\square^2$ due to the fact that $\hat f_{1:k}-f_{1:k}^*\in\cC$. As a result, with probability at least $1-\bP\left((\Omega')^c\right)-\exp\left(-t(\square,\triangle)/2\right)$,  $\norm{\bX_{\phi,1:k}(\hat f_{1:k}-f_{1:k}^*)}_2^2\lesssim \sqrt{N}\square.$ The above analysis leads to the following proposition:
\begin{Proposition}
    Suppose that Assumption~\ref{assumption:DM_L4_L2},  Assumption~\ref{assumption:upper_dvoretzky} and Assumption~\ref{assumption:DMU_used_for_RIP} hold. There exist absolute constants $\nC\label{C_price_overfitting_1}$, $\nC\label{C_price_overfitting_2}$ and $\nC\label{C_price_overfitting_3}$ such that the following holds. Assume that there exists $k\in\bN\cup\{\infty\}$ such that $N\leq\oc{c_kappa_DM}\kappa_{DM}d_\lambda^*\left(\Gamma_{k+1:\infty}^{-1/2}B_\cH\right)$, \eqref{eq:extra_condition_k>N} holds, and $\kappa_{DM}\left(4\lambda + \Tr\left(\Gamma_{k+1:\infty}\right)\right)\geq N\left(R_N^*(\oc{c_kappa_RIP})\right)^2$, then the following holds for all such $k$'s and for all $t>0$. With probability at least $1-\bP\left((\Omega')^c\right)-\exp\left(-t(\square,\triangle)/2\right) - \bar p_\xi$ (see \eqref{eq:def_bar_p_xi}),
    \begin{align*}
        &\norm{\Gamma_{k+1:\infty}^{1/2}\left(\hat f_{k+1:\infty}  -f_{k+1:\infty}^*\right)}_\cH\\
        &\leq \oC{C_price_overfitting_1}\frac{\sqrt{N\Tr\left(\Gamma_{k+1:\infty}^2\right)} + N\norm{\Gamma_{k+1:\infty}}_{\text{op}}}{4\lambda + \Tr\left(\Gamma_{k+1:\infty}\right)}\left(\norm{\Gamma_{k+1:\infty}^{1/2}f_{k+1:\infty}^*}_\cH + \square\right)\\
        &+ \norm{\Gamma_{k+1:\infty}^{1/2}f_{k+1:\infty}^*}_\cH+\frac{3}{2}\sigma_\xi \frac{\sqrt{N\Tr\left(\Gamma_{k+1:\infty}^2\right)} + N\norm{\Gamma_{k+1:\infty}}_{\text{op}}}{4\lambda + \Tr\left(\Gamma_{k+1:\infty}\right)} .
    \end{align*}
\end{Proposition}

\subsection{Proof of Theorem~\ref{theo:DM_RKHS}}\label{sec:proof_DM}

\beginproof We prove Theorem~\ref{theo:DM_RKHS}, and for the sake of generality, we replace the covariance matrix by $\Gamma$, that is, $\Gamma = \bE\left[\phi(X)\otimes\phi(X)\right]$. For some $2<p\leq 2+\eps$, let
\begin{equation}\label{eq:def_B}
    B=\sup\left(\bE\left|\left\langle \phi(X),f\right\rangle_\cH\right|^p:\, \norm{f}_\cH=1\right).
\end{equation}
Since we have $L_{2+\eps}-L_2$ norm equivalence of marginals (recall \eqref{eq:norm_equivalence}), $B\leq \kappa^p\norm{\Gamma}_{\text{op}}^{p/2}$.
To obtain a high-probability upper bound for Equation \eqref{eq:objective_DM}, we first use a one-scale net argument. Set $V_{\eps_0} = \left\{\pi \lambda:\, \lambda\in S_2^{N-1} \right\}$ to be an $\epsilon_0$-net of $S_2^{N-1}$ so that for all $\lambda\in S_2^{N-1}$, $\norm{\pi \lambda-\lambda}_2\leq \epsilon_0$ and $\left|V_{\epsilon_0}\right|\leq(5/\epsilon_0)^N$. Unlike the situation in \cite{lecue_geometrical_2022}, the choice of $\eps_0$ here does not affect the probability deviation because we will use a different discretization procedure to control $\bE_\eta V_{I_\eta,\vu}$(which will be defined later) uniformly on $S_2^{N-1}$, which does not depend on $\eps_0$. Therefore, one may choose an $\eps_0$ as small as possible(say, $1/10$) to compensate for absolute constants at the end of the proof.

It follows from the triangle inequality that for every $\vlambda\in  S_2^{N-1}$,
\begin{eqnarray*}
    \left|\frac{\norm{A \vlambda}_\cH^2}{\left(\ell^*\right)^2} - 1 \right| &=& \left| \frac{\norm{A\left(\vlambda-\pi\vlambda\right)}_\cH^2}{\left(\ell^*\right)^2} + \frac{\norm{A(\pi\vlambda)}_\cH^2}{\left(\ell^*\right)^2} + \frac{2\left\langle A \left(\vlambda-\pi\vlambda \right), A (\pi \vlambda)\right\rangle_\cH}{\left(\ell^*\right)^2}- 1 \right|\\
    &\leq& \frac{\norm{A\left(\vlambda-\pi\vlambda\right)}_\cH^2}{\left(\ell^*\right)^2} + \left|\frac{\norm{A(\pi\vlambda)}_\cH^2}{\left(\ell^*\right)^2} - 1\right| + \left|\frac{2\left\langle A \left(\vlambda-\pi \vlambda \right), A (\pi  \vlambda)\right\rangle_\cH}{\left(\ell^*\right)^2}\right|
\end{eqnarray*}and so, from the Cauchy-Schwarz inequality,
\begin{equation*}
    \underset{\vlambda\in  S_2^{N-1}}{\mathrm{sup}}\left|\frac{\norm{A \vlambda}_\cH^2}{\left(\ell^*\right)^2} - 1 \right| \leq \Phi^2 + \Psi^2 + 2\Phi\sqrt{\Psi^2 + 1}
\end{equation*}where
\begin{equation*}
    \Phi^2 := \underset{\vlambda\in S_2^{N-1}}{\mathrm{sup}} \frac{\norm{A\left(\vlambda-\pi\vlambda\right)}_\cH^2}{\left(\ell^*\right)^2} \mbox{ and } \Psi^2 := \sup_{\vlambda\in S_2^{N-1}}\left|\frac{\norm{A(\pi\vlambda)}_\cH^2}{\left(\ell^*\right)^2} - 1\right|.
\end{equation*}

Thus, we only need to bound $\Phi$ and $\Psi$ from above. To that end, we start with a decoupling argument to deal with the cross terms. 

\paragraph{A decoupling argument}
Let $(\eta_i)_{1\leq i\leq N}$ be $N$ i.i.d. selectors (i.e. $\bP(\eta_i=0) = \bP(\eta_i = 1) = 1/2$) and define $ I_\eta := \left\{1\leq i\leq N:\, \eta_i = 1 \right\}$. For all $\vu\in\bR^N$ and $I\subset[N]$, we also define
\begin{equation*}
    V_{I,\vu} := \frac{1}{\left(\ell^*\right)^2}\left\langle  \left(\sum_{i\in I} u_i \phi(X_i) \right), \, \left(\sum_{j\in I^c} u_j \phi(X_j) \right)\right\rangle_\cH.
\end{equation*}
For every $\vu\in B_2^{N}$, a decoupling technique (see for instance Chapter 6 of \cite{vershynin_high-dimensional_2018}) leads to the following result
\begin{eqnarray}\label{eq:decoupling}
    \nonumber\frac{\norm{A \vu}_\cH^2}{\left(\ell^*\right)^2} &=& \sum_{i=1}^N\frac{\norm{\phi(X_i)}_\cH^2u_i^2}{\left(\ell^*\right)^2}  + \sum_{i\neq j} u_i u_j \frac{\left\langle \phi(X_i), \phi(X_j)\right\rangle_\cH}{\left(\ell^*\right)^2}\\
    \nonumber&=&\sum_{i=1}^N\frac{\norm{\phi(X_i)}_\cH^2u_i^2}{\left(\ell^*\right)^2} + \sum_{i,j=1}^N \bE_\eta 4\eta_i {(1-\eta_j)} u_i u_j \frac{\left\langle \phi(X_i), \phi(X_j)\right\rangle_\cH}{\left(\ell^*\right)^2}\\
    \nonumber&=&\sum_{i=1}^N\frac{\norm{\phi(X_i)}_\cH^2u_i^2}{\left(\ell^*\right)^2} + \frac{4}{\left(\ell^*\right)^2}\bE_{\eta}{\left\langle \left(\sum_{i\in I_\eta} u_i \phi(X_i)\right),\, \left(\sum_{j\in I_\eta^c} u_j \phi(X_j) \right) \right\rangle_\cH}\\
    &=&\sum_{i=1}^N\frac{\norm{\phi(X_i)}_\cH^2u_i^2}{\left(\ell^*\right)^2} +4\bE_{\eta}V_{I_\eta,\vu},
\end{eqnarray}
where $\bE_\eta$ is the expectation with respect to $(\eta_i)_{i\in[N]}$ conditionally on all other random variables. With the decoupling technique, we are going to estimate $\Phi^2$ and $\Psi^2$ from above. We start with $\Psi^2$ and get 
\begin{align}\label{eq:PSI_1}
     \nonumber \Psi^2 =&  \underset{\vlambda\in S_2^{N-1}}{\mathrm{sup}}\left|\frac{\norm{A(\pi\vlambda)}_\cH^2}{\left(\ell^*\right)^2} - 1\right| \leq \underset{\lambda\in S_2^{N-1}}{\mathrm{sup}}\left|\sum_{i=1}^N\frac{\norm{\phi(X_i)}_\cH^2}{\left(\ell^*\right)^2} \left(\pi   \lambda\right)_i^2 -1 \right|+ 4\underset{\vlambda\in S_2^{N-1}}{\mathrm{sup}}\left|\bE_\eta V_{I_\eta,\pi   \vlambda}\right|\\
     &\leq \max_{1\leq i\leq N}\left|\frac{\norm{\phi(X_i)}_\cH^2}{\left(\ell^*\right)^2} -1 \right|+ 4\underset{\vlambda\in S_2^{N-1}}{\mathrm{sup}}\left|\bE_\eta V_{I_\eta,\pi   \vlambda}\right|.
\end{align}

We next have  $\Phi^2=\sup\left(\norm{A(\vlambda-\pi\vlambda)}_\cH^2/\left(\ell^*\right)^2:\, \vlambda\in S_2^{N-1}\right)\leq (\epsilon_0/\ell^*)^2\norm{A}_{\text{op}}^2$, where $\norm{A}_{\text{op}}$ is the operator norm of $A:(\bR^N,\ell_2)\to (\cH,\norm{\cdot}_\cH)$, thus it remains to prove a high probability upper bound on $\norm{A}_{\text{op}}$. From \eqref{eq:decoupling} we know that
\begin{align}\label{eq:upper_bound_operator_norm_A}
    \frac{\norm{A}_{\text{op}}^2}{\left(\ell^*\right)^2} \leq \underset{i\in[N]}{\max}\frac{\norm{\phi(X_i)}_\cH^2}{\left(\ell^*\right)^2} + 4\underset{\norm{\vmu}_2=1}{\sup}\bE_\eta V_{I_\eta, \vmu}.
\end{align}As a result,
\begin{equation}\label{eq:PHI_1}
    \Phi^2 \leq \frac{4\epsilon_0^2}{\left(\ell^*\right)^2} \underset{\pi\vmu\in V_{1/2}}{\mathrm{sup}}\norm{A(\pi\vmu)}_\cH^2 
    \leq 4\epsilon_0^2\left(\underset{1\leq i\leq N}{\mathrm{max}}\frac{\norm{\phi(X_i)}_\cH^2}{\left(\ell^*\right)^2} + 4\underset{\pi\vmu\in V_{1/2}}{\mathrm{sup}}\bE_\eta V_{I_\eta, \pi\vmu}\right).
\end{equation}
Therefore, we only need to find a high probability upper bound on $\left|\bE_\eta V_{I_\eta,\vu}\right|$ uniformly for all $\vu$ in $V_{1/2}$ and $V_{\eps_0}$. In contrast to the method employed in \cite{lecue_geometrical_2022}, we immediately derive this upper bound uniformly over $B_2^N$. This will result in $\epsilon_0$ being a free parameter, as we will show at the end of the proof (\eqref{eq:upper_bound_V} holds on $S_2^{N-1}$, instead of $V_{\eps_0}$).

To achieve this, we adapt \cite{tikhomirov_sample_2018}'s argument. We recall that Theorem~\ref{theo:DM_RKHS} contains two aspects:
\begin{enumerate}
    \item when $\Tr(\Gamma)$ is the dominating term in $\lambda+\Tr(\Gamma)$,
    \item when $\lambda$ is dominating.
\end{enumerate}
In case [1], we need not only the upper bound of $\sup\left(\norm{\bX_\phi^\top\vlambda}_\cH:\, \vlambda\in S_2^{N-1}\right)$, but also the lower bound for $\inf\left(\norm{\bX_\phi^\top\vlambda}_\cH:\, \vlambda\in S_2^{N-1}\right)$. However, in case [2], we only need the upper bound for $\sup\left(\norm{\bX_\phi^\top\vlambda}_\cH:\, \vlambda\in S_2^{N-1}\right)$. This will be clear in Section~\ref{sec:proof_DM_lambda_dominating}.

We first introduce some notation.

For all $I,J\subset[N]$ and $\ell\in[N]$, let $S_2^J$ be the unit sphere in Euclidean norm of $\bR^J$(see as a subspace of $\bR^N$ endowed by the canonical vectors indexed by $J$), and let
\begin{eqnarray*}
    S_I^J := \left\{\lambda\in S_2^J:\, \lambda_i=0,\forall i\in I^c \right\} \mbox{ and }
    S_{I,\ell}^J := \left\{\lambda\in S_I^J:\, \left|\left\{ i\in I:\, \lambda_i\neq 0  \right\}\right|\leq \ell \right\}.
\end{eqnarray*}For the sake of simplicity, we let $S_I^N:=S_I^{[N]}$ and $S_{I,\ell}^N := S_{I,\ell}^{[N]}$.
For each $k\leq N$ and $I,\cC\subset[N]$, we denote
\begin{align}\label{eq:definition_g}
\begin{aligned}
    g(k,\cC,I) &:= \sup\left( \left|\left\langle \sum_{i\in I\cap\cC}v_i\phi(X_i),\, \sum_{j\in I^c\cap\cC}u_j\phi(X_j) \right\rangle_\cH\right| ,\, \vu,\vv\in S_{[N],k}^N\right), \mbox{ and }\\
    g(N,I)&:=g(k,[N],I).
\end{aligned}
\end{align}
Further, for any vector $\vu\in\bR^N$, and any $i\leq N$, we set
\begin{equation}\label{eq:definition_W}
    W_{\vu,i} := \left\langle \phi(X_i),\, \sum_{j=1}^N u_j\phi(X_j)\right\rangle_\cH,
\end{equation} and for some $J\subset[N]$, we denote $(\left|W_{\vu,J}\right|)_j^*$ as the $j$-th largest absolute value of the coordinates of $(W_{\vu,j})_{j\in J}$. Given $1\leq k\leq N$ and $I\subset[N]$, we define
\begin{equation*}
    \cM_{I, k} = \max\left( \norm{\sum_{i\in I}\lambda_i\phi(X_i)}_\cH:\, \vlambda\in S_{I,k}^N \right),\quad \mbox{ and } \cM_{k} = \cM_{[N],k}.
\end{equation*}
In particular, $\norm{A}_{\text{op}}=\cM_N$, and by the decoupling argument again, see \eqref{eq:upper_bound_operator_norm_A}
\begin{equation}\label{eq:upper_bound_M}
    \cM_N ^2 \leq \max\left(\norm{\phi(X_i)}_\cH^2:\, 1\leq i\leq N\right) + 4(\ell^*)^2\sup\left(\bE_\eta V_{I_\eta,\vu}:\, \vu\in S_2^{N-1} \right).
\end{equation}

Once \eqref{eq:upper_bound_M} has been obtained, the next step is to derive an upper bound for the supremum of $\left|\bE_\eta V_{I_\eta,\vu}\right|$ over all $\vu\in S_2^{N-1}$ in relation to $\cM_N$. The purpose of the remaining part of the proof will be to establish this point. In order to accomplish this, we begin with
\begin{align}\label{eq:E_V_and_g_N_I}
    \begin{aligned}
        &4\left(\ell^*\right)^2\sup\left(\left|\bE_\eta V_{I_\eta,\vu}\right|:\, \vu\in S_2^{N-1}\right)\\
    &=4\sup\left(\bE_{\eta}{\left|\left\langle \left(\sum_{i\in I_\eta} u_i \phi(X_i)\right),\, \left(\sum_{j\in I_\eta^c} u_j \phi(X_j) \right) \right\rangle_\cH\right|}:\, \vu\in S_2^{N-1}\right)\\
    &= \frac{4}{2^N}\sup\left(\sum_{I\subset[N]}\left|\left\langle  \sum_{i\in I}u_i\phi(X_i),\, \sum_{j\in I^c}u_j\phi(X_j) \right\rangle_\cH\right|:\, \vu\in S_2^{N-1}\right)\\
    &\leq \frac{4}{2^N}\sum_{I\subset[N]}g\left(N,I\right), 
    \end{aligned}
\end{align}where the last step is via Jensen's inequality and by the fact that optimizing over $(\vv,\vu)\in S_{[N],k}^N\times S_{[N],k}^N$ leads to a larger supreme.

\paragraph{A sparsifying argument}We begin by the sparsifying lemma for $g\left(N,\cC,I\right)$.

\begin{Lemma}\label{lemma:sparsifying_tikhomirov}
    Let $0<\varepsilon\leq 1$, $k\geq 12/\varepsilon^2$, $4\leq m \leq k$, and let $I,\cC\subset[N]$. There then exists an absolute constant $\nC\label{C_DM_1}>0$ such that
    \begin{align}\label{eq:reduction_last_step}
    \begin{aligned}
        g(k,\cC,I)&\leq g(m,\cC,I) + 2 \oC{C_DM_1}\varepsilon^{-2}\max\left(\left|\left\langle \phi(X_i),\, \phi(X_j)\right\rangle_\cH\right|:\, i\neq j\in\cC\right)\\
        &+ \oC{C_DM_1}\sqrt{k}\varepsilon^{-2}\bigg(\sup\left( \left|W_{\vy,I^c\cap\cC}\right|_{\lfloor m/4 \rfloor }^* :\, \vy\in S_{I\cap\cC,\varepsilon k}^N\right)\\
        &+ \sup\left( \left|W_{\vz,I\cap\cC}\right|_{\lfloor m/4 \rfloor }^* :\, \vz\in S_{I^c\cap\cC,\varepsilon k}^N\right) \bigg).
    \end{aligned}
    \end{align}
\end{Lemma}
As mentioned in \cite{tikhomirov_sample_2018}, when compared to the findings of \cite{bartl_random_2022}, this method of sparsifying yields a lower cardinality when combined with the subsequent discretization argument. This gives rise to a broad condition about the tail of kernel features (from $L_{\log{N}}-L_2$ to $L_{2+\eps}-L_2$).

The next lemma is a combination of \cite[Lemma 13]{tikhomirov_sample_2018} and \cite[Lemma 4]{tikhomirov_sample_2018}:
\begin{Lemma}\label{lemma:discretization}
    Let $0<\rho\leq 1$, $r,h,p,q\in\bN$ and $r\geq 2$, let $V_\rho$ be a support-preserving Euclidean $\rho$-net of $S_{[q],h}^q$. Further, let $T$ be a $p\times q$ matrix. Then
    \begin{equation*}
        \sup\left( \left|T\vu\right|_r^* :\, \vu\in S_{[q],h}^q \right) \leq 2\sup\left( \left| T\vv \right|_{\lfloor r/2\rfloor}^* :\, \vv\in V_\rho\right) + \frac{4\rho}{\sqrt{r}}\sup\left( \sqrt{\sum_{i=1}^r \left(\left|T\vu\right|_{i}^*\right)^2 } :\, \vu\in S_{[q],h}^q\right).
    \end{equation*}
    Here, the support-preserving Euclidean $\rho$-net means that $V_\rho$ is a $\rho$-net, such that for all $\vx\in S_{[q],k}^q$, there exists $\vy\in V_\rho$ with $\mathrm{supp}(\vy)\subset\mathrm{supp}(\vx)$ and $\norm{\vx-\vy}_2\leq \rho$. By \cite[Lemma 4]{tikhomirov_sample_2018}, there exists $\nC\label{C_DM_2}>0$ and $V_\rho\subset S_{[q],k}^q$ such that $\left|V_\rho\right|\leq (\oC{C_DM_2} q / \rho h)^h$, and $V_\rho$ is a support-preserving Euclidean $\rho$-Net of $S_{[q],h}^q$.
\end{Lemma}

\paragraph{A dimension reduction argument}
Combining Lemma~\ref{lemma:sparsifying_tikhomirov} and Lemma~\ref{lemma:discretization} with an induction argument, we obtain the following proposition, as a deterministic argument. Compared to \cite[Proposition 14]{tikhomirov_sample_2018}, we remove the condition $N\geq 128 \oC{C_DM_1} \varepsilon^{-2}k$. This assumption $N\geq 128 \oC{C_DM_1} \varepsilon^{-2}k$ prevents us from being in the Dvoretzky-Milman regime, because when $N=k$, we need $\varepsilon$ sufficiently large, however, Proposition~\ref{prop:g_stochastic} needs $\varepsilon$ to be sufficiently small. Removing this assumption is possible because we are going to choose $k=N$ in Proposition~\ref{prop:g_stochastic}, which keeps $128\oC{C_DM_1}k/(\varepsilon^2N)$ as a constant in the proof.

\begin{Proposition}\label{prop:induction}
Let $I,\cC\subset[N]$, and let $0<\varepsilon<1/3$ and $k\geq 24/\varepsilon^2\vee N$.
Denote $t:=\lfloor \log_2\left(\varepsilon^2 k/24\right)\rfloor$ and define $k_j:=\lfloor k/2^j\rfloor$, $0\leq j\leq t$.
There are then subsets $V_j\subset S_{I,\varepsilon k_j}^I$ and $V_j'\subset S_{I^c,\varepsilon k_j}^{I^c}$ for all $0\leq j\leq t-1$ such that for the absolute constants from Lemma~\ref{lemma:sparsifying_tikhomirov} and Lemma~\ref{lemma:discretization},
\begin{itemize}
    \item $\left|V_j\right|,\left|V_j'\right|\leq \left(\oC{C_DM_2} N/\varepsilon k_j\right)^{2\varepsilon k_j}$ for all $0\leq j\leq t-1$.
    \item We have
    \begin{eqnarray*}
        g(k,I) &\leq& \exp\left(\frac{128\oC{C_DM_1}k}{\varepsilon^2 N}\right)\left(48+2\oC{C_DM_1}\log_2(k)\right)\varepsilon^{-2}\max\left(\left|\left\langle \phi(X_i),\,\phi(X_j)\right\rangle_\cH\right|:\, i\neq j\in\cC\right)\\
        &+& 2\oc{C_DM_1}\exp\left(\frac{128C_{1}k}{\varepsilon^2 N}\right)\varepsilon^{-2}\sum_{j=0}^{t-1}\sqrt{k_j}\bigg(\sup\left( \left| W_{\vu, I^c} \right|_{\lfloor k_{j+1}/16\rfloor}^* :\,\vu\in V_j\right)\\
        &+& \sup\left( \left| W_{\vv, I} \right|_{\lfloor k_{j+1}/16\rfloor}^* :\,\vv\in V_j'\right)\bigg).
    \end{eqnarray*}
\end{itemize}
\end{Proposition}

\beginproof First, we let $j<t$ and consider the quantity $g(k_j,\cC,I)$. 
By Lemma~\ref{lemma:sparsifying_tikhomirov}, for $k=k_j$ and $m=k_{j+1}$ (where we have $k_j\geq 12/\varepsilon^2$ and $k_{j+1}\geq 4$),
\begin{eqnarray*}
    g(k_j,\cC,I) &\leq& g(k_{j+1},\cC,I) + 2\oC{C_DM_1}\varepsilon^{-2} \max\left(\left|\left\langle \phi(X_i),\phi(X_j)\right\rangle_\cH\right|:\, i\neq j\in\cC\right) \\
    &+&\oC{C_DM_1}\sqrt{k_j}\varepsilon^{-2}\bigg(\sup\left( \left|W_{\vy,I^c\cap\cC}\right|_{\lfloor k_{j+1}/4\rfloor}^* :\, \vy\in S_{I\cap\cC,\varepsilon k_j}^N\right)\\
    &+& \sup\left( \left|W_{\vz,I\cap\cC}\right|_{\lfloor k_{j+1}/4\rfloor}^* :\, \vz\in S_{I^c\cap\cC,\varepsilon k_j}^N\right) \bigg).
\end{eqnarray*}We first find upper bounds for
\begin{align*}
    &\sup\left( \left|W_{\vy,I^c\cap\cC}\right|_{\lfloor k_{j+1}/4\rfloor}^* :\, \vy\in S_{I\cap\cC,\varepsilon k_j}^N\right),\mbox{ and }\\
    &\sup\left( \left|W_{\vz,I\cap\cC}\right|_{\lfloor k_{j+1}/4\rfloor}^* :\, \vy\in S_{I^c\cap\cC,\varepsilon k_j}^N\right).
\end{align*}
We only obtain the former, as the latter follows from the same ideas.

From the definition of $\vy\in S_{I\cap\cC,\varepsilon k_j}^N$ we know that $\vy$ is supported on $I\cap\cC$. By the definition of $W$ (see \eqref{eq:definition_W}), $W_{\vy, I^c\cap\cC} = \left(\left\langle \phi(X_i),\, \sum_{j\in I\cap\cC}y_j\phi(X_j)\right\rangle_\cH\right)_{i\in I^c\cap\cC}$. Therefore, taking supreme over $S_{I\cap\cC,\varepsilon k_j}^N$ is equivalent to taking it over $\vy\in S_{I\cap\cC,\varepsilon k_j}^{I\cap\cC}$.

Apply Lemma~\ref{lemma:discretization} with $r=\lfloor k_{j+1}/4\rfloor$, $\rho = k_j/N$, $h=\lfloor \varepsilon k_j\rfloor$, and matrix $T = \left(\left\langle \phi(X_i),\, \phi(X_j)\right\rangle_\cH\right)_{i,j}$ for $(i,j)\in\left((I^c\cap\cC)\times (I\cap\cC)\right)$, then $T\vy = \left(\left\langle \phi(X_i),\, \sum_{j\in I\cap\cC}y_j\phi(X_j)\right\rangle_\cH\right)_{i\in I^c\cap\cC} = W_{\vy, I^c\cap\cC}$ and thus there exists $V_j\subset S_{I\cap\cC,\varepsilon k_j}^{I\cap\cC}\subset S_{I,\varepsilon k_j}^N$, which is a support-preserving $k_j/N$-net of cardinality at most
\begin{equation*}
    \left(\frac{\oC{C_DM_2}\left|I\cap\cC\right|}{((k_j/N)\varepsilon k_j)}\right)^{\varepsilon k_j}\leq \left(\frac{\oC{C_DM_2}}{\varepsilon}\right)^{\varepsilon k_j}\left(\frac{N}{k_j}\right)^{2\varepsilon k_j}\leq \left(\frac{\oC{C_DM_2}N}{\varepsilon k_j}\right)^{2\varepsilon k_j}
\end{equation*}since we can choose $\oC{C_DM_2}>1\vee \varepsilon$, such that
\begin{eqnarray*}
    &&\sup\left( \left|W_{\vy,I^c\cap\cC}\right|_{\lfloor k_{j+1}/4\rfloor}^* :\, \vy\in S_{I\cap\cC,\varepsilon k_j}^N\right)\leq 2\sup\left( \left|W_{\vu, I^c\cap\cC}\right|_{\lfloor \lfloor k_{j+1}/4\rfloor /2 \rfloor}^* :\, \vu\in V_j\right) \\
    &+& \frac{4k_j/N}{\sqrt{\lfloor k_{j+1}/4\rfloor}}\sup\left( \sqrt{\sum_{i=1}^{\lfloor k_{j+1}/4\rfloor} \left(\left| W_{\vy,I^c\cap\cC} \right|_{i}^*\right)^2 } :\, \vy\in S_{I\cap\cC,\varepsilon k_j}^N\right)\\
    &\leq& 2\sup\left( \left|W_{\vu, I^c}\right|_{\lfloor \lfloor k_{j+1}/4\rfloor /2 \rfloor}^* :\, \vu\in V_j\right) + \frac{4k_j/N}{\sqrt{\lfloor k_{j+1}/4\rfloor}} g(k_{j+1},\cC,I),
\end{eqnarray*}where we have used $k_{j+1}/4\leq k_{j+1}$ and the fact that for any $\vu\in\bR^N$ and $\ell\in[N]$,
\begin{equation*}
    \sup\left(\sum_{i\in I}u_i v_i:\, \vv\in S_{I,\ell}^I\right) = \sqrt{\sum_{i=1}^\ell \left((\vu_I)_i^*\right)^2}.
\end{equation*}
The analysis for $\sup\left( \left|W_{\vz,I\cap\cC}\right|_{\lfloor m/4 \rfloor }^* :\, \vz\in S_{I^c\cap\cC,\varepsilon k}^N\right)$ is similar. We obtain
\begin{eqnarray*}
    g(k_j,\cC,I) &\leq& \left(1+\frac{32\oC{C_DM_1}k_j}{\varepsilon^2 N}\right)g(k_{j+1},\cC,I) + 2\oC{C_DM_1}\varepsilon^{-2}\max\left(\left|\left\langle \phi(X_i),\phi(X_j)\right\rangle_\cH\right|:\, i\neq j\in\cC\right)\\
    &+& 2\oC{C_DM_1}\sqrt{k_j}\varepsilon^{-2}\left(\sup\left( \left| W_{\vu, I^c} \right|_{\lfloor k_{j+1}/16\rfloor}^* :\,\vu\in V_j\right) + \sup\left( \left| W_{\vv, I} \right|_{\lfloor k_{j+1}/16\rfloor}^* :\,\vv\in V_j'\right)\right).
\end{eqnarray*}
Notice that for any $t'<t$,
\begin{equation*}
    \prod_{j=0}^{t'} \left(1+\frac{32\oC{C_DM_1}k_j}{\varepsilon^2 N}\right) \leq \prod_{j=0}^t \left(1+\frac{32\oC{C_DM_1}k_j}{\varepsilon^2 N}\right) \leq \exp\left(\sum_{j=0}^t \frac{32\oC{C_DM_1}k_j}{\varepsilon^2 N}\right) \leq \exp\left(\frac{128\oC{C_DM_1}k}{\varepsilon^2 N}\right).
\end{equation*}By induction over $0\leq j<t$, we obtain
\begin{eqnarray*}
    &&\exp\left(-\frac{128\oC{C_DM_1}k}{\varepsilon^2 N}\right)g(k,\cC,I) \leq g(k_t,\cC,I) + 2\oC{C_DM_1}\varepsilon^{-2}\log_2(k)\max\left(\left|\left\langle \phi(X_i),\phi(X_j)\right\rangle_\cH\right|:\, i\neq j\in\cC\right) \\
    &+& 2\oC{C_DM_1}\varepsilon^{-2}\sum_{j=0}^{t-1}\sqrt{k_j}\left(\sup\left( \left| W_{\vu, I^c} \right|_{\lfloor k_{j+1}/16\rfloor}^* :\,\vu\in V_j\right) + \sup\left( \left| W_{\vv, I} \right|_{\lfloor k_{j+1}/16\rfloor}^* :\,\vv\in V_j'\right)\right).
\end{eqnarray*}
Finally, we bound $g(k_t,\cC,I)$ from above:
\begin{eqnarray*}
    g(k_t,\cC,I) &\leq& \sup\left(\sum_{i,j=1}^N \left|y_i z_j\left\langle \phi(X_i),\phi(X_j)\right\rangle_\cH\right| :\, \vy\in S_{I\cap\cC,k_t}^N ,\vz\in S_{I^c\cap\cC,k_t}^N \right)\\
    &\leq& \max\left(\left|\left\langle \phi(X_i),\,\phi(X_j)\right\rangle_\cH\right|:\, i\neq j\in\cC\right)\sup\left(\sum_{i,j=1}^N \left|y_i z_j\right| :\, \vy\in S_{I\cap\cC,k_t}^N ,\vz\in S_{I^c\cap\cC,k_t}^N \right)\\
    &=&\max\left(\left|\left\langle \phi(X_i),\,\phi(X_j)\right\rangle_\cH\right|:\, i\neq j\in\cC\right)\cdot k_t \\
    &\leq& 48\varepsilon^{-2}\max\left(\left|\left\langle \phi(X_i),\,\phi(X_j)\right\rangle_\cH\right|:\, i\neq j\in\cC\right).
\end{eqnarray*}
\endproof

\paragraph{Stochastic arguments}We next introduce the randomness of $\phi(X)$. We will see in the next proposition that choosing $k=N$ does not destroy the proof, but leads to a larger constant compared with \cite[Proposition 15]{tikhomirov_sample_2018}.

\begin{Proposition}\label{prop:g_stochastic}
    For the absolute constant $\oC{C_DM_1}$ from Lemma~\ref{lemma:sparsifying_tikhomirov} and $\oC{C_DM_2}$ from Lemma~\ref{lemma:discretization}, there are sufficiently large universal constants $\nC\label{C_DM_3} = (12\oC{C_DM_2}/\varepsilon^3)^{36/\varepsilon}$ for $\varepsilon<1/256$ , $\nC\label{C_DM_4}$ and $\nC\label{C_DM_5}$ depending on $p$ such that: Suppose $N\geq (12\oC{C_DM_2}/\varepsilon^3)^6 \oC{C_DM_3}^{-\varepsilon/6}\vee \exp\left(48/\oC{C_DM_1}\right)$, let $I\subset[N]$. With probability at least $1-\oC{C_DM_4}/N^3$, for all $\cC\subset[N]$,
    \begin{equation*}
        g(N,\cC,I) \leq \oC{C_DM_5} \log_2(N)\max\left(\left|\left\langle \phi(X_i),\,\phi(X_j)\right\rangle_\cH\right|:\, i\neq j\in\cC\right) + \oC{C_DM_5} B^{1/p} \sqrt{N}\cM_{N}.
    \end{equation*}
\end{Proposition}
\beginproof Let $t=\lfloor \log_2\left(\varepsilon^2 N/24\right)\rfloor$ and $k_j=\lfloor N/2^j\rfloor$ for each $j\in\{0,1,\cdots,t\}$. Fix $j\in\{0,1,\cdots,t-1\}$, we study the term $\sup\left( \left| W_{\vu, I^c} \right|_{\lfloor k_{j+1}/16\rfloor}^* :\,\vu\in V_j\right)$ appearing by Proposition~\ref{prop:induction} for $k=N$. In the whole proof below $k=N$, in particular, $k_j = \lfloor N/2^j\rfloor$, $j=0,\cdots,t$ and $t = \lfloor \log_2{(\varepsilon^2 N/24)}\rfloor$. At this time, $48<\oC{C_DM_1}\log_2(N)$. For each $\vu\in V_j$, condition on $(X_i)_{i\in I}$, $\left(\left\langle \phi(X_j),\sum_{i\in I} u_i\phi(X_i)\right\rangle_\cH\right)_{j\in I^c} = \left(W_{\vu,j}\right)_{j\in I^c}$ are i.i.d. random variables. Recalling the definition of $B$, see \eqref{eq:def_B}, the conditional expectation of $\left|W_{\vu,j}\right|^p$ given $(X_i)_{i\in I}$ satisfies
\begin{equation*}
    \bE\left[ \left|W_{\vu,j}\right|^p\big| (X_i)_{i\in I} \right] \leq B \norm{\sum_{i\in I} u_i\phi(X_i)}_{\cH}^p \leq B \cM_{N}^p.
\end{equation*}Let $\vu\in V_j$(in particular, $\mathrm{supp}(\vu)\subset I$). Condition on $(X_i)_{i\in I}$, $\bE\left[ \left|W_{\vu,j} \right|^p/\norm{ \sum_{i\in I} u_i\phi(X_i) }_\cH \big| (X_i)_{i\in I}\right]\leq B$. For $\tau_j^p = 32eB\left(N/k_{j+1}\right)^{1+256\varepsilon}$,
\begin{align*}
    &\bP\left( \frac{\left| W_{\vu, I^c} \right|_{\lfloor k_{j+1}/16\rfloor}^*}{\norm{ \sum_{i\in I} u_i\phi(X_i) }_\cH}\geq \tau_j\bigg| \left(X_i\right)_{i\in I} \right)\\
    &= \bP\left(\exists J\subset I^c,\, \left|J\right|\geq \lfloor k_{j+1}/16\rfloor,\, \mbox{such that } \forall j\in J,\, \frac{\left| W_{\vu,j} \right|}{\norm{ \sum_{i\in I} u_i\phi(X_i) }_\cH} \geq \tau_j\bigg| \left(X_i\right)_{i\in I}\right)\\
    &\leq \binom{\left|I^c\right|}{\lfloor k_{j+1}/16\rfloor}\left(\bP\left(\frac{\left| W_{\vu,j} \right|}{\norm{ \sum_{i\in I} u_i\phi(X_i) }_\cH} \geq \tau_j\bigg| \left(X_i\right)_{i\in I}\right)\right)^{\lfloor k_{j+1}/16\rfloor}\\
    &\leq \left(\frac{e\left|I^c\right|}{\lfloor k_{j+1}/16\rfloor}\right)^{\lfloor k_{j+1}/16\rfloor}\left(\frac{B}{\tau_j^p}\right)^{\lfloor k_{j+1}/16\rfloor}\\
    &\leq \left(\frac{\left|I^c\right|}{N}\left(\frac{k_{j+1}}{N}\right)^{256\varepsilon}\right)^{\lfloor k_{j+1}/16\rfloor}\leq \left(\frac{k_{j+1}}{N}\right)^{4\varepsilon k_j}.
\end{align*}Hence, conditionally on $(X_i)_{i\in I}$, with probability at least $1-(k_{j+1}/N)^{4\varepsilon k_j}$,
\begin{equation*}
    \left| W_{\vu, I^c} \right|_{\lfloor k_{j+1}/16\rfloor}^*\leq \tau_j\norm{ \sum_{i\in I} u_i\phi(X_i) }_\cH\leq \tau_j\cM_N.
\end{equation*}Therefore, by Fubini's theorem, with probability at least $1-(k_{j+1}/N)^{4\varepsilon k_j}$, $\left| W_{\vu, I^c} \right|_{\lfloor k_{j+1}/16\rfloor}^*\leq \tau_j\cM_N$.
Taking the union bound over all $\vu\in V_j$(note that the cardinality of $V_j$ is given in Proposition~\ref{prop:induction}) and using $k_{j+1}\leq k_j$,
\begin{equation*}
    \bP\left(\sup\left(\left|W_{\vu,I^c}\right|_{\lfloor k_{j+1}/16\rfloor}^*:\, \vu\in V_j \right) \geq \tau_j \right) \leq \left(\frac{k_{j+1}^2}{N^2}\cdot\frac{\oC{C_DM_2} N}{\varepsilon k_j}\right)^{2\varepsilon k_j} \leq \left(\frac{k_{j+1}}{N}\cdot\frac{\oC{C_DM_2}}{\varepsilon}\right)^{2\varepsilon k_j}.
\end{equation*} Since $t\mapsto t\log{(eN/t)}$ is increasing on $\{t:\, 0<t\leq N\}$ and $k_t\geq 12/\varepsilon^2$, we have for all $j=0,\cdots,t-1$,
\begin{eqnarray*}
    2\varepsilon k_j\log\left(\frac{\varepsilon N}{\oC{C_DM_2} k_{j+1}}\right) &\geq& 2\varepsilon k_j\log\left(\frac{\varepsilon N}{\oC{C_DM_2} k_j}\right) \geq 2\varepsilon k_t\log\left(\frac{\varepsilon N}{\oC{C_DM_2} k_t}\right) \geq 2\varepsilon \frac{12}{\varepsilon^2}\log\left(\frac{\varepsilon N}{\oC{C_DM_2}\frac{12}{\varepsilon^2}}\right),\\
    &=&\frac{24}{\varepsilon}\log\left(\frac{\varepsilon^3 N}{12\oC{C_DM_2}}\right) \geq \log\left(\frac{N^4}{\oC{C_DM_3}}\right),
\end{eqnarray*}when $\left(\varepsilon^3 N/12 \oC{C_DM_2}\right)^{24/\varepsilon}\geq N^4/\oC{C_DM_3}$, i.e., $N\geq (12\oC{C_DM_2}/\varepsilon^3)^6\cdot \oC{C_DM_3}^{-\varepsilon/6}$.
The polynomial rate can therefore balance the union bound over $j\in\{0,1,\cdots,t-1\}$, which is a logarithmic rate with respect to $N$. With probability at least $1-\oC{C_DM_4}/N^3$, for all $j\in\{0,1,\cdots,t-1\}$,
\begin{equation*}
    \sup\left(\left|W_{\vu,I^c}\right|_{\lfloor k_{j+1}/16\rfloor}^*:\, \vu\in V_j \right) < (32eB)^{1/p}\cM_{N}\left(\frac{N}{k_{j+1}}\right)^{p^{-1}(1+256\varepsilon)}.
\end{equation*}
Finally, notice that
\begin{equation*}
    \sum_{j=0}^{t-1} \sqrt{k_j}\sup\left(\left|W_{\vu,I^c}\right|_{\lfloor k_{j+1}/16\rfloor}^*:\, \vu\in V_j \right) \leq (32eB)^{1/p}\cM_{N}\sqrt{N} \sum_{j=0}^{t-1}2^{j/p+256\varepsilon j/p-j/2}.
\end{equation*}
There exists an absolute constant $\oC{C_DM_5}$ such that $\sum_{j=0}^{t-1}2^{\left(\frac{1+256\varepsilon}{p}-\frac{1}{2}\right)j}\leq \oC{C_DM_5}$ for $1/2>(1+256\varepsilon)/p$. We finish the proof of Proposition~\ref{prop:g_stochastic} by applying Proposition~\ref{prop:induction}.
\endproof

\subsubsection{Case [1]: when $\Tr(\Gamma)$ is dominating.}\label{sec:proof_DM_trace_dominate}

In this case, we apply all the aforementioned results (Proposition~\ref{prop:g_stochastic}, Proposition~\ref{prop:induction} and Lemma~\ref{lemma:sparsifying_tikhomirov}) to $\cC=[N]$.

Let us now apply Proposition~\ref{prop:g_stochastic},
\begin{align*}
    &\bE\left|\left\{I\subset[N]:\, g(N,I) > \oC{C_DM_5} \log_2(N)\max\left(\left|\left\langle \phi(X_i),\,\phi(X_j)\right\rangle_\cH\right|:\, i\neq j\in[N]\right) + \oC{C_DM_5} B^{1/p} \sqrt{N}\cM_{N} \right\}\right|\\
    &\leq \frac{2^N}{N^3},
\end{align*}
so with probability at least $1-1/N^2$, there are at most $2^N/N$ subsets $I\subset[N]$, such that
\begin{equation}\label{eq:spiky_subsets}
    g(N,I) > \oC{C_DM_5} \log_2(N)\max\left(\left|\left\langle \phi(X_i),\,\phi(X_j)\right\rangle_\cH\right|:\, i\neq j\in[N]\right) + \oC{C_DM_5} B^{1/p} \sqrt{N}\cM_{N}.
\end{equation}
For these $2^N/N$ ``spiky'' subsets, we simply use a deterministic argument: for all $I\subset[N]$, we have
\begin{align*}
    g(N,I)&\leq N\max\left(\left|\left\langle \phi(X_i),\, \phi(X_j)\right\rangle_\cH\right|:\, i\neq j\in[N]\right).
\end{align*}
As a consequence, if we denote by $\cI$ the set of all subsets $I\subset[N]$ satisfying \eqref{eq:spiky_subsets}, with probability at least $1-1/N^2$,
\begin{align}
    \sum_{I\subset[N]} g(N,I) &= \sum_{I\in\cI} g(N,I) + \sum_{I\notin\cI} g(N,I) \notag\\
    &\leq 2\oC{C_DM_5} 2^N\log_2(N)\max\left(\left|\left\langle \phi(X_i),\, \phi(X_j)\right\rangle_\cH\right|:\, i\neq j\in[N]\right) + 2^N \oC{C_DM_5} B^{1/p} \sqrt{N}\cM_{N}, \label{eq:uniform_over_I}
\end{align}

We are left with an upper bound that has a high probability of $\max\left(\left|\left\langle \phi(X_i),\,\phi(X_j)\right\rangle_\cH\right|:\, i\neq j\in[N]\right)$. We emphasize again that we do not use the sample coloring technique developed by \cite{tikhomirov_sample_2018}, but instead, we use the strong concentration of $\left|\left\langle \phi(X_i),\,\phi(X_j)\right\rangle_\cH\right|$ to absorb this $\log{N}$ factor, because of \eqref{eq:extra_assumption}.

\paragraph{Upper bound for $\max\left(\left|\left\langle \phi(X_i),\,\phi(X_j)\right\rangle_\cH\right|:\, i\neq j\in[N]\right)$}  Let $p=2+\eps$.


For any $i\neq j\in[N]$,
\begin{align*}
    &\bE\left|\left\langle \phi(X_i), \phi(X_j)\right\rangle_\cH\right|^{2+\eps} = \bE\left[\bE\left[ \left|\left\langle \phi(X_i), \phi(X_j)\right\rangle_\cH\right|^p \big| X_i\right]\right] \leq \kappa^p \left( \bE\left|\left\langle \phi(X_i), \phi(X_j)\right\rangle_\cH\right|^2 \right)^{p/2}\\
    &=\kappa^p \left(\bE K_{k+1:\infty}(X_i, X_j)\right)^{p/2} = \kappa^p \left(\Tr\left(\Gamma^2\right)\right)^{p/2},
\end{align*}where we used \eqref{eq:norm_equivalence} to obtain the inequality.
By union bound, for any $\tau>0$,
\begin{align*}
    \bP\left(\max\left(\left|\left\langle \phi(X_i),\phi(X_j)\right\rangle\right|_\cH:\, i\neq j\in[N]\right)>\tau\right)&\leq N^2\frac{\bE\left|\left\langle \phi(X_i), \phi(X_j)\right\rangle_\cH\right|^p}{\tau^p}\\
    &\leq N^2\frac{\kappa^p \left(\Tr\left(\Gamma^2\right)\right)^{p/2}}{\tau^p}.
\end{align*}Let $\tau = \bar\delta\Tr\left(\Gamma\right)/\log{N}$, with probability at least
\begin{align*}
    1- N^2\frac{\kappa^p \left(\Tr\left(\Gamma^2\right)\right)^{p/2}}{\left(\frac{\bar\delta\Tr\left(\Gamma\right)}{\log{N}}\right)^p} = 1 - \left( \frac{\kappa}{\bar\delta} \right)^p (\log{N})^p N^2\left(\frac{\sqrt{\Tr\left(\Gamma^2\right)}}{\Tr\left(\Gamma\right)}\right)^p =: 1-\bar p,
\end{align*}we have
\begin{align*}
    \max\left(\left|\left\langle \phi(X_i),\phi(X_j)\right\rangle\right|_\cH:\, i\neq j\in[N]\right)\leq \frac{\bar\delta \Tr\left(\Gamma\right)}{\log{N}},
\end{align*}where
\begin{align*}
    \bar p = \left(\frac{\kappa}{\bar\delta}\right)^{2+\eps}\left(\log{N}\right)^{2+\eps}N^{1-\frac{\eps}{2}}\left(\frac{\sqrt{N\Tr\left(\Gamma^2\right)}}{\Tr\left(\Gamma\right)}\right)^{2+\eps}.
\end{align*}
Together with \eqref{eq:uniform_over_I}, with probability at least $1 - \bar p - N^{-2}$,
\begin{align*}
    \sum_{I\subset[N]} g\left(N,I\right) \leq 2\oC{C_DM_5}2^N \bar\delta\Tr\left(\Gamma\right) + 2^N \oC{C_DM_5} B^{1/(2+\eps)}\sqrt{N}\cM_N.
\end{align*}By \eqref{eq:E_V_and_g_N_I}, we obtain that
\begin{align}\label{eq:upper_bound_V}
    4(\ell^*)^2\sup\left(\left|\bE_{\eta}V_{I_\eta,\vu}\right| :\, \vu\in S_2^{N-1}\right) \lesssim \bar\delta \Tr\left(\Gamma\right) + \sqrt{N}B^{1/(2+\eps)}\cM_N.
\end{align}

Combining Eq.\eqref{eq:upper_bound_M}, Eq.\eqref{eq:upper_bound_V} and $B\leq\kappa^{2+\eps}\norm{\Gamma}_{\text{op}}^{\frac{2+\eps}{2}}$, there exists an absolute constant $\nC\label{C_DM_6}>1$ such that with probability at least $1-\gamma-\bar p - N^{-2}$,
\begin{align}
    &\sup\left(\left|\bE_{\eta}V_{I_\eta,\vu}\right| :\, \vu\in S_{[N],N}^N\right) \leq \oC{C_DM_6}\bar\delta \vee \oC{C_DM_6} \kappa^2\frac{N\norm{\Gamma}_{\text{op}}}{\Tr\left(\Gamma\right)}\vee \oC{C_DM_6}\kappa\sqrt{\frac{N\norm{\Gamma}_{\text{op}}}{\Tr\left(\Gamma\right)}}\cdot\frac{\max\left(\norm{\phi(X_i)}_\cH:\, i\in[N]\right)}{\sqrt{\Tr\left(\Gamma\right)}}\notag\\
    &\vee \oC{C_DM_6}\kappa\sqrt{\bar\delta}\sqrt{\frac{N\norm{\Gamma}_{\text{op}}}{\Tr\left(\Gamma\right)}}\leq \oC{C_DM_6}\bar\delta\left(1 + \kappa^2\kappa_{DM}\frac{\Tr\left(\Gamma\right)+\lambda}{\Tr\left(\Gamma\right)}+\kappa\sqrt{\kappa_{DM}\frac{\Tr\left(\Gamma\right)+\lambda}{\Tr\left(\Gamma\right)}}\left(1+\delta+\sqrt{\bar\delta}\right)\right),\label{eq:final_upper_bound_V}
\end{align}where we used that $N\leq\kappa_{DM}\bar\delta^2d_\lambda^*\left(\Gamma_{k+1:\infty}^{-1/2}B_\cH\right)$, $\bar\delta<1$ and \eqref{eq:diagonal_term_assumption} from Assumption~\ref{assumption:DM_L4_L2}.

Note that $S_2^{N-1} = S_{[N],N}^N$, and we plug \eqref{eq:final_upper_bound_V} into \eqref{eq:PSI_1} and \eqref{eq:PHI_1} and take $\eps_0 = \sqrt{\delta}$ in \eqref{eq:PHI_1},
\begin{align*}
    &\Psi^2\leq \delta + 4\oC{C_DM_6}\bar\delta\left(1 + \kappa^2\kappa_{DM}\frac{\Tr\left(\Gamma\right)+\lambda}{\Tr\left(\Gamma\right)}+\kappa\sqrt{\kappa_{DM}\frac{\Tr\left(\Gamma\right)+\lambda}{\Tr\left(\Gamma\right)}}\left(1+\delta+\sqrt{\bar\delta}\right)\right),\\
    &\Phi^2 \leq \delta\left(1+\delta+4\oC{C_DM_6}\bar\delta\left(1 + \kappa^2\kappa_{DM}\frac{\Tr\left(\Gamma\right)+\lambda}{\Tr\left(\Gamma\right)}+\kappa\sqrt{\kappa_{DM}\frac{\Tr\left(\Gamma\right)+\lambda}{\Tr\left(\Gamma\right)}}\left(1+\delta+\sqrt{\bar\delta}\right)\right)\right)
\end{align*}
Since we have the right to choose sufficiently small $\bar\delta$ and $\kappa_{DM}$ as long as \eqref{eq:extra_assumption} holds, we can set
\begin{align}\label{eq:bound_kappa_DM}
    \kappa_{DM}\leq \left(\frac{1}{12\oC{C_DM_6}\kappa}\right)^{2}<\frac{1}{4\oC{C_DM_6}\kappa^2}.
\end{align}Because $\delta,\bar\delta<1$,
\begin{align*}
    &\Psi^2< \delta + 4\oC{C_DM_6}\bar\delta + 2\bar\delta\frac{\Tr\left(\Gamma\right)+\lambda}{\Tr\left(\Gamma\right)}\leq \delta + (4\oC{C_DM_6}+2)\bar\delta\frac{\Tr\left(\Gamma\right)+\lambda}{\Tr\left(\Gamma\right)},\\
    &\Phi^2< 4\left(\delta + \delta + \delta\left(4\oC{C_DM_6} + 2\right)\bar\delta \frac{\Tr\left(\Gamma\right)+\lambda}{\Tr\left(\Gamma\right)}\right).
    \end{align*}
Recall that in this subsection, we assume that there exists an absolute constant $\nC\label{C_comparision_trace_lambda}$ such that $\lambda\leq \oC{C_comparision_trace_lambda} \Tr\left(\Gamma\right)$. In this case, there exist absolute constants $\nC\label{C_distortion_1}$, $\nC\label{C_distortion_2}$, $\nC\label{C_distortion_3}$, $\nC\label{C_distortion_4}$ and $\nC\label{C_distortion_5}$ such that \begin{align*}
    \Phi^2+\Psi^2 + 2\Phi\sqrt{\Psi^2+1} < \oC{C_distortion_1}\delta^2 + \oC{C_distortion_2} \bar\delta^2 + 4\sqrt{\left(3\delta+\oC{C_distortion_3}\bar\delta\right)\left(1+\delta+\oC{C_distortion_4}\bar\delta\right)}=:\tilde\delta<1,
\end{align*}provided that $\delta<1/(100\sqrt{\oC{C_distortion_1}})$ and $\bar\delta<1/\oC{C_distortion_5}$ (thus we can take $\oC{C_DM}$ in Assumption~\ref{assumption:DM_L4_L2} as $\oC{C_distortion_5}$), and where
\begin{align*}
    \tilde\delta = \oC{C_distortion_1}\delta^2 + \oC{C_distortion_2} \bar\delta^2 + 4\sqrt{\left(3\delta+\oC{C_distortion_3}\bar\delta\right)\left(1+\delta+\oC{C_distortion_4}\bar\delta\right)}.
\end{align*}
This proves that with probability at least $1-\gamma - \bar p - N^{-2}$, for all $\vlambda\in\bR^N$,
\begin{align*}
    \left(1- \tilde\delta\right)\ell^*\norm{\vlambda}_2\leq \norm{\bX_{\phi,k+1:\infty}^\top\vlambda}_\cH\leq \left(1+  \tilde\delta\right)\ell^*\norm{\vlambda}_2,
\end{align*}provided that $N\leq\kappa_{DM}\bar\delta^2d_\lambda^*\left(\Gamma_{k+1:\infty}^{-1/2}B_\cH\right)$.

\subsubsection{Case [2]: when $\lambda$ is dominating.}\label{sec:proof_DM_lambda_dominating}

When $\lambda>\oC{C_comparision_trace_lambda}\Tr\left(\Gamma\right)$. In this case, we only make use of the fact that $\bX_{\phi}\bX_{\phi}^\top$ is of rank-$N$ and is positive semi-definite, hence $\sigma_N\left(\bX_{\phi}\bX_{\phi}^\top + \lambda I_N\right)\geq \lambda + \sigma_N\left(\bX_{\phi}\bX_{\phi}^\top\right)\geq \lambda\geq \oc{c_distortion_9}\lambda + (1-\oc{c_distortion_9})\oC{C_comparision_trace_lambda}\Tr\left(\Gamma\right)$, where we recall that $0<\oc{c_distortion_9}<1$ is some absolute constant. Hence our objective is to prove that there exists an absolute  constant $\nC\label{C_DM_lambda_dominate_1}$ such that with high probability, we have $\norm{\bX_{\phi}\bX_{\phi}^\top + \lambda I_N}_{\mathrm{op}}\leq \oC{C_DM_lambda_dominate_1}\lambda$. We prove this by proving that there exists an absolute constant $\nC\label{C_DM_lambda_dominate_2}$ such that $\oC{C_DM_lambda_dominate_2}^2\leq \oC{C_DM_lambda_dominate_1}-1$, and with high probability we have $\norm{\bX_{\phi}^\top}_{\mathrm{op}}\leq\oC{C_DM_lambda_dominate_2}\sqrt{\lambda}$.

Let $\{\cC_m\}_{m\leq\chi}$ for some $\chi\in\bN_+$ be a partition of $[N]$, by Jensen's inequality, for any $\vlambda\in S_2^{N-1}$,
\begin{align*}
    \norm{\bX_\phi^\top\vlambda}_\cH^2\leq \chi\sum_{m=1}^\chi \norm{\sum_{i\in \cC_m}\lambda_i \phi(X_i)}_\cH^2.
\end{align*}
Applying \eqref{eq:upper_bound_operator_norm_A} and \eqref{eq:E_V_and_g_N_I} but with $A$ replaced by its restriction onto $\cC_m$ for each $m\leq \chi$, we obtain that for any $\vlambda\in S_2^{N-1}$,
\begin{align}\label{eq:DM_lambda_dominate_bridge}
    \frac{\norm{\bX_\phi^\top\vlambda}_\cH^2}{(\ell^*)^2}\leq \chi^2 \underset{i\in[N]}{\max}\frac{\norm{\phi(X_i)}_\cH^2}{(\ell^*)^2} +\frac{\chi}{(\ell^*)^2}\sum_{m=1}^\chi \frac{4}{2^N}\sum_{I\subset[N]}g(N,\cC_m,I).
\end{align}

At this time, we can make use of the sample coloring technique in \cite{tikhomirov_sample_2018}. It is a technique used to truncate the inner products $\left(\left<\phi(X_i),\phi(X_j)\right>_\cH\right)_{i\neq j\in\cC}$. Given i.i.d. random vectors $\left(\phi(X_i)\right)_{i\leq N}$ and $H>0$, there exists an undirected graph $\cG_H$ whose vertex set is $[N]$, and its edge set is:
\begin{equation*}
    \left\{ (i,j):\, 1\leq i<j\leq N,\, \left|\left<\phi(X_i),\phi(X_j)\right>_\cH\right|> H\max\left(\norm{\phi(X_h)}_\cH:\, h\leq N\right) \right\}.
\end{equation*} The coloring of $\cG_H$ is an assignment of ``colors'' to all vertices such that no adjacent vertices share the same color. The smallest possible number of colors sufficient to assign such a coloring is called the chromatic number of $\cG_H$, denoted as $\chi(\cG_H)$, and the collection $\{\cC_m^H\}_{m\leq \chi(\cG_H)}$ is the associated partition by colors of $[N]$. That is to say, for any $m\leq \chi(\cG_H)$, and $i\neq j\in \cC_m^H$, the vertices $i,j$ are not adjacent, so $\left|\left<\phi(X_i),\phi(X_j)\right>_\cH\right|\leq H \max\left(\norm{\phi(X_h)}_\cH:\, h\leq N\right)$. Since $(\phi(X_i))_{i\leq N}$ are random, $\cG_H$ is a random graph, and the following lemma is a high probability estimate of $\chi(\cG_H)$. The following lemma is a weaker version of \cite[Proposition 10]{tikhomirov_sample_2018}, which is sufficient for our purpose.

\begin{Lemma}\label{lemma:sample_coloring}
    Assume that for some $p>2$ we have $\sup\left(\bE\left|\left<\phi(X),f\right>_\cH\right|^p:\, \norm{f}_\cH=1\right)= B$. Then for any $H>0$ and any integer $m>1$, the chromatic number of $\cG_H$ satisfies $\chi(\cG_H)\leq m$ with probability at least $1-(BNH^{-p})^{m-1}N$.
\end{Lemma}
\beginproof Let us introduce an auxiliary random process $(Y_i)_{i\in[N]}$ with values in $\bN$, where $Y_1:=1$ as a constant, and for all $i=2,\cdots,N$,
\begin{equation*}
    Y_i := \min\left( r\in\bN_+:\, \forall j<i,\, j\in\bN_+,\, \mbox{ with } Y_j=r, \mbox{ we have }\left|\left<\phi(X_i),\phi(X_j)\right>_\cH\right|\leq H\norm{\phi(X_j)}_\cH \right).
\end{equation*}The process $(Y_i)_{i\in[N]}$ is ``classifying'' each $\phi(X_i)$ is if $\left|\left<\phi(X_1),\phi(X_2)\right>_\cH\right|>H\norm{\phi(X_1)}_\cH$, $\left|\left<\phi(X_1),\phi(X_3)\right>_\cH\right|\leq H\norm{\phi(X_1)}_\cH$, $\left|\left<\phi(X_2),\phi(X_3)\right>_\cH\right|\leq H\norm{\phi(X_2)}_\cH$, then $Y_2 = 2$(because $(1,2)$ is adjacent in $\cG_H$), $Y_3 = 1$, because either $(1,3)$ or $(2,3)$ is not adjacent in $\cG_H$. Such a $\cG_H$ has chromatic number $2$.

By the definition of $Y_i$, we have that any two numbers $i\neq j\in[N]$ such that $Y_i=Y_j$ are not adjacent in $\cG_H$, and $Y_i = Y_j$ is a sufficient but not necessary condition for adjacency of $(i,j)$. In particular, $\chi(\cG_H)\leq\max\left(Y_i:\, i\in[N]\right)$. Next for each $i>1$ and $m\geq 1$, we have
\begin{eqnarray*}
    \bP\left(Y_i = m+1\right) &\leq& \bP\left( \exists j\leq i-1 \mbox{ s.t. } \left|\left<\phi(X_i),\phi(X_j)\right>_\cH\right|>H\norm{\phi(X_j)}_\cH,\mbox{ and }Y_{j}=m \right)\\
    &\leq& \sum_{j=1}^{i-1} \bP\left(\left|\left<\phi(X_i),\phi(X_{j})\right>_\cH\right|>H\norm{\phi(X_j)}_\cH,\mbox{ and }Y_{j}=m \right).
\end{eqnarray*}For all $j=0,\cdots,i-1$, since $Y_j$ is $\sigma\left(X_1,\cdots,X_j\right)$-measurable, it is independent of $X_i$, hence
\begin{align*}
    &\bP\left( Y_j = m,\, \left|\left< \phi(X_i),\,\phi(X_j)\right>_\cH\right|>H\norm{\phi(X_j)}_\cH \big| \left(X_\ell\right)_{\ell=1}^{i-1}\right)\\
    &=\bE\left[\1_{\{Y_j = m\}}\1_{\left\{ \left|\left< \phi(X_i),\,\phi(X_j)\right>_\cH\right|>H\norm{\phi(X_j)}_\cH\right\}} \big| \left(X_\ell\right)_{\ell=1}^{i-1}\right]\\
    &=\1_{\{Y_j = m\}}\bP\left({\left\{ \left|\left< \phi(X_i),\,\phi(X_j)\right>_\cH\right|>H\norm{\phi(X_j)}_\cH\right\}} \big| \left(X_\ell\right)_{\ell=1}^{i-1}\right)\\
    &\leq \1_{\{Y_j = m\}}\frac{\bE\left[ \left|\left< \phi(X_i),\,\phi(X_j)\right>_\cH\right|^p\big| \left(X_\ell\right)_{\ell=1}^{i-1} \right]}{H^p\norm{\phi(X_j)}_\cH^p}\leq \1_{\{Y_j = m\}}\frac{B}{H^p},
\end{align*}where we used Markov's inequality to obtain the first inequality. Hence
\begin{equation*}
    \bP\left(\left|\left<\phi(X_i),\phi(X_{j})\right>_\cH\right|>H\norm{\phi(X_j)}_\cH,\mbox{ and }Y_{j}=m \right)\leq\bP\left( Y_j = m \right)\frac{B}{H^p}.
\end{equation*}Further, by $\bE\left|\left\{j\leq N:\, Y_j=m\right\}\right|=\sum_{j\leq N}\bE\1_{\{Y_j=m\}}=\sum_{j\leq N}\bP\left(Y_j=m\right)\geq \sum_{j\leq i-1}\bP\left(Y_j=m\right)$,
\begin{equation}\label{eq:induction}
    \sum_{j=1}^{i-1} \bP\left(\left|\left<\phi(X_i),\phi(X_{j})\right>_\cH\right|>H\norm{\phi(X_j)}_\cH,\mbox{ and }Y_{j}=m \right)\leq BH^{-p}\bE\left|\left\{ j\leq N:\, Y_j = m \right\}\right|
\end{equation}

It follows from Equation~\eqref{eq:induction} that
\begin{equation*}
    \sum_{i=1}^N\bP\left( Y_i = m+1 \right) = \bE\left|\left\{ j\leq N:\, Y_j=m+1 \right\}\right| \leq BNH^{-p} \bE\left|\left\{j\leq N:\, Y_j=m\right\}\right|.
\end{equation*}We next deal with $\bE\left|\left\{j\leq N:\, Y_j=2\right\}\right|$. We simply upper bound $\bE\left|\left\{j\leq N:\, Y_j=2\right\}\right|$ by $N$. Therefore,
\begin{equation*}
    \bE\left|\left\{ j\leq N:\, Y_j=m+1 \right\}\right| \leq (BNH^{-p})^{m-1}N.
\end{equation*}Note that the set of values $\left\{Y_j:\, j\leq N\right\}$ is an interval in $\bN$, hence
\begin{equation*}
    \bP\left(\chi(\cG_H)\geq m+1\right) \leq \bP\left(\exists j\leq N:\, Y_j=m+1\right) \leq \bE\left|\left\{ j\leq N:\, Y_j=m+1 \right\}\right| \leq (BNH^{-p})^{m-1}N.
\end{equation*}
\endproof

Combining Proposition~\ref{prop:g_stochastic} and the sample coloring technique from Lemma~\ref{lemma:sample_coloring}, we obtain an upper bound for $\bE_\eta V_{I_\eta,\vu,\vv}$ uniformly over all $\vu,\vv\in S_{[N],N}^N = S_2^{N-1}$. We state this result in the following Proposition.

\begin{Proposition}\label{prop:uniform_bound_on_E_V}
There are absolute constants $\nC\label{C_DM_lambda_dominate_5}$ and $\oC{C_DM_lambda_dominate_4}$ depending only on $p$, such that the following holds. If $N\geq \oC{C_DM_lambda_dominate_5}$, then for any $\lambda> -\Tr\left(\Gamma\right)$, with probability at least $1-\bar p - N^{-2}$, where
\begin{align}\label{eq:bar_p_DM_lambda_dominate}
    \bar p := N\left(\left(\frac{ 4\kappa^2\log^2(N)\norm{\Gamma}_{op}}{\Tr(\Gamma)+\lambda}\right)^{p/2} N\right)^{\lceil(8+2p)/(p-2)\rceil-1},
\end{align}
    \begin{equation}\label{eq:upper_bound_V_lambda}
        \cM_N \leq \oC{C_DM_lambda_dominate_6} \sqrt{\Tr\left(\Gamma\right)+\lambda} + \oC{C_DM_lambda_dominate_6}\sqrt{N}B^{1/p}.
    \end{equation}
\end{Proposition}
\beginproof Let $H>0$ which will be determined later, and let $0\leq m\leq\chi$, we apply Proposition~\ref{prop:g_stochastic} to $\cC=\cC_m^H$:
\begin{align*}
    &\bE\left|\left\{I\subset[N]:\, g(N,\cC_m^H,I) > \oC{C_DM_5} \log_2(N)\max\left(\left|\left<\phi(X_i),\,\phi(X_j)\right>_\cH\right|:\, i\neq j\in \cC_m^H\right) + \oC{C_DM_5} B^{1/p} \sqrt{N}\cM_{N} \right\}\right|\\
    &\leq \frac{2^N}{N^3},
\end{align*}
then with probability at least $1-1/N^2$, there are at most $2^N/N$ subsets $I\subset[N]$, such that
\begin{equation}\label{eq:spiky_subsets_lambda}
    g(N,\cC_m^H,I) > \oC{C_DM_5} \log_2(N)\max\left(\left|\left<\phi(X_i),\,\phi(X_j)\right>_\cH\right|:\, i\neq j\in \cC_m^H\right) + \oC{C_DM_5} B^{1/p} \sqrt{N}\cM_{N}.
\end{equation}For these $2^N/N$ ``spiky'' subsets, we simply use a deterministic argument: for all $I\subset[N]$, we have
\begin{equation*}
    g(N,\cC_m^H,I)\leq N\max\left(\left|\left<\phi(X_i),\, \phi(X_j)\right>_\cH\right|:\, i\neq j\in\cC_m^H\right) \leq NH\max\left(\norm{\phi(X_i)}_\cH:\, i\leq N\right),
\end{equation*}where we use that $i,j\in\cC_m^H$ have the same color and therefore are not adjacent in $\cG_m^H$.
As a consequence, if we denote by $\cI$ the set of all subsets $I\subset[N]$ satisfying Equation~\eqref{eq:spiky_subsets_lambda}, with probability at least $1-1/N^2$, for any $(\cC_m^H)_{m\leq \chi}$, we have
\begin{eqnarray*}
    \sum_{I\subset[N]} g(N,\cC_m^H,I) &=& \sum_{I\in\cI} g(N,\cC_m^H,I) + \sum_{I\notin\cI} g(N,\cC_m^H,I) \\
    &\leq& \oC{C_DM_5} 2^N\log_2(N)H\max\left(\norm{\phi(X_i)}_\cH:\, i\leq N\right) + 2^N \oC{C_DM_5} B^{1/p} \sqrt{N}\cM_{N},
\end{eqnarray*}where we again used that $i$ and $j$ are not adjacent in $\cG_m^H$.

Let $H = \bar\delta\sqrt{\Tr(\Gamma)+\lambda}/\log{N}$ with $\bar\delta=1/2$ (unlike the case in Section~\ref{sec:proof_DM_trace_dominate}, we only need an \emph{isomorphic} upper bound, we can choose $\bar\delta$ to be an arbitrary absolute constant), then
\begin{equation*}
    BNH^{-p} \leq \left(\frac{4\kappa^2\log^2(N)\norm{\Gamma}_{op}}{(\Tr(\Gamma)+\lambda)}\right)^{p/2} N.
\end{equation*}
Let $\chi = \lceil (8+2p)/(p-2)\rceil$ and apply Lemma~\ref{lemma:sample_coloring} for $m=\chi$, with probability at least $1-\bar p$, $\chi(\cG_H)\leq \chi$.

On the other hand, by \eqref{eq:DM_lambda_dominate_bridge} and the fact that \eqref{eq:spiky_subsets_lambda} is valid uniformly over all $(\cC_m^H)_{m\leq \chi}$ (thanks to Proposition~\ref{prop:g_stochastic}), there exists an absolute constant $\nC\label{C_DM_lambda_dominate_4}$ such that with probability $1-\bar p - N^{-2}$, for all $\vlambda\in S_2^{N-1}$, we have that
\begin{align}\label{eq:DM_lambda_dominate_1}
    \begin{aligned}
        \frac{\norm{\bX_\phi^\top\vlambda}_\cH^2}{(\ell^*)^2} &\leq \left(\frac{8+2p}{p-2}\right)^2 \underset{i\in[N]}{\max}\frac{\norm{\phi(X_i)}_\cH^2}{(\ell^*)^2}\\
    &+ \frac{\left(\frac{8+2p}{p-2}\right)^2}{(\ell^*)^2}  \left(\oC{C_DM_lambda_dominate_4} \log_2(N)H\max\left(\norm{\phi(X_i)}_\cH:\, i\leq N\right) + \oC{C_DM_lambda_dominate_4} B^{1/p} \sqrt{N}\cM_{N}\right)\\
    &=\left(\frac{8+2p}{p-2}\right)^2 \underset{i\in[N]}{\max}\frac{\norm{\phi(X_i)}_\cH^2}{(\ell^*)^2}\\
    &+ \frac{\left(\frac{8+2p}{p-2}\right)^2}{(\ell^*)^2} \left(\oC{C_DM_lambda_dominate_4} \bar\delta\sqrt{\Tr(\Gamma)+\lambda} \max\left(\norm{\phi(X_i)}_\cH:\, i\leq N\right) + \oC{C_DM_lambda_dominate_4} B^{1/p} \sqrt{N}\cM_{N}\right).
    \end{aligned}
\end{align}

Solving \eqref{eq:DM_lambda_dominate_1} gives that there exists an absolute constant $\nC\label{C_DM_lambda_dominate_6}$ depending only on $p$ such that with probability at least $1-\bar p - N^{-2}$,
\begin{align*}
    \cM_N \leq \oC{C_DM_lambda_dominate_6} \sqrt{\Tr\left(\Gamma\right)+\lambda} + \oC{C_DM_lambda_dominate_6}\sqrt{N}B^{1/p}.
\end{align*}
\endproof

Recall that we have assumed that $\oC{C_comparision_trace_lambda}\Tr\left(\Gamma\right)<\lambda$ in this case, and $\Tr\left(\Gamma\right)+\lambda \geq (\kappa_{DM}/4)^{-2}N\norm{\Gamma}_{\mathrm{op}}$. Moreover, since we have $B\leq \kappa^p \norm{\Gamma}_{\mathrm{op}}^{p/2}$ for any $2< p\leq 2+\eps$, we have
\begin{align*}
    \cM_N &= \sup\left(\norm{\bX_\phi^\top\vlambda}_\cH:\, \vlambda\in S_2^{N-1}\right) \leq \oC{C_DM_lambda_dominate_6}\sqrt{1+\oC{C_comparision_trace_lambda}^{-1}\lambda} + \oC{C_DM_lambda_dominate_6}\kappa\sqrt{N\norm{\Gamma}_{\mathrm{op}}}\\
    &\leq \oC{C_DM_lambda_dominate_6}\sqrt{1+\oC{C_comparision_trace_lambda}^{-1}\lambda} + \frac{\oC{C_DM_lambda_dominate_6}}{4}\kappa\kappa_{DM} \sqrt{\Tr\left(\Gamma\right)+\lambda}\\
    &\leq \left(\oC{C_DM_lambda_dominate_6}\sqrt{1+\oC{C_comparision_trace_lambda}^{-1}} + \frac{\oC{C_DM_lambda_dominate_6}}{4}\kappa\kappa_{DM}\sqrt{1+\oC{C_comparision_trace_lambda}^{-1}}\right)\sqrt{\lambda}.
\end{align*}Letting $\oC{C_DM_lambda_dominate_2} = \oC{C_DM_lambda_dominate_6}\sqrt{1+\oC{C_comparision_trace_lambda}^{-1}} + \oC{C_DM_lambda_dominate_6}\kappa\kappa_{DM}\sqrt{1+\oC{C_comparision_trace_lambda}^{-1}}/4$, this is precisely our initial objective. As a result, we may let $\oC{C_DM_lambda_dominate_1}=\oC{C_DM_lambda_dominate_2}^2+1$, and $\oC{C_distortion_9}=\oC{C_DM_lambda_dominate_1}$.

\subsection{Proof of Proposition~\ref{prop:KRR_smooth_upper}}\label{sec:proof_smooth_kernel}



In this section, we are dealing with an inner product kernel, that is, a kernel that satisfies \eqref{eq:polynomial_inner_product_kernel}, and that possesses a smooth kernel function $h$, for which the associated RKHS may not satisfy the $L_{2+\epsilon}/L_2$ norm equivalence assumption \eqref{eq:norm_equivalence}. Consequently, we introduce a new technique, involving further truncation of $\Gamma_{k+1:\infty}$. This technique stems from the following observation that RKHS with eigenvalues exhibiting power decay, that is,  $\sigma_j\sim j^{-\alpha}$ for $\alpha>1$, possess the following property: for any $k \in \bN_+$, $\Tr\left(\Gamma_{k+1:2k}\right) \sim \Tr\left(\Gamma_{k+1:\infty}\right)$. In other words, due to the specific condition of power decay, the contribution of eigenvalues $\sigma_j$ for $j>2k$ to the trace is negligible. Readers can refer to \cite[pp. 37]{bartlett_benign_2020} for the proof of this property. In fact, for a general $\cK > k$, we have $F(\cK) - F(k+1) \leq \sum_{j>k}^\cK j^{-\alpha} \leq F(\cK) - F(k)$, where $F(x) = \frac{1}{1-\alpha}x^{1-\alpha}$. As a result, $\Tr\left(\Gamma_{k+1:2k}\right)\sim \frac{1}{\alpha-1}\left(k^{1-\alpha}-(2k)^{1-\alpha}\right)\sim \frac{1}{\alpha-1}k^{1-\alpha}$. Moreover, by the same argument, $\Tr\left(\Gamma_{k+1:\infty}\right)\sim \frac{1}{\alpha-1}k^{1-\alpha}$.


The proof strategy is almost the same as that of Theorem~\ref{theo:main_upper_k>N}.

\subsubsection{Stochastic Argument}

Like in Section~\ref{sec:proof_main_upper_k>N}, that is, like the proof of Theorem~\ref{theo:main_upper_k>N} ($k\gtrsim N$), we prove that there exist absolute constants $\nc\label{c_RIP_lower_smooth}$,  $\nC\label{C_RIP_upper_smooth}$, $\nC\label{C_DMU_smooth}$, $\nC\label{C_bX_f_star_smooth}$ and $\nC\label{C_sum_Gamma_phi_smooth}$ such that the following random event (denoted as $\Omega''$) happens with probability greater than $1/2$:

\begin{itemize}
    \item for any $\vlambda\in\bR^N$,
    \begin{align}\label{eq:DM_smooth_applied}
        \left(\frac{1}{2}\Tr(\Gamma_{k+1:2k})+\lambda\right)\norm{\blambda}_2\leq \norm{\left(\bX_{\phi,k+1:\infty}\bX_{\phi,k+1:\infty}^\top+\lambda I\right)\blambda}_2\leq \oC{C_DMU_smooth}\log^2(N)\left(\Tr(\Gamma_{k+1:\infty})+\lambda\right)\norm{\vlambda}_2,
    \end{align}
    \item for any $f_{1:k}\in\mathrm{cone}\left(\cC(R_N(\oc{c_kappa_RIP}))\right)$,
    \begin{align}\label{eq:RIP_smooth_applied}
        \oc{c_RIP_lower_smooth}\norm{\Gamma_{1:k}^{1/2}f}_\cH \leq \frac{1}{\sqrt{N}}\norm{\bX_{\phi,1:k}f}_2 \leq \oC{C_RIP_upper_smooth}\norm{\Gamma_{1:k}^{1/2}f}_\cH,
    \end{align}
      \item for all $\blambda\in\bR^N$,
  \begin{align}\label{eq:DM_upper_smooth_applied}
      \norm{\Gamma_{k+1:\infty}^{1/2}\bX_{\phi,k+1:\infty}^\top\blambda}_\cH\leq \oC{C_DMU_smooth}\left(\log(N)\sqrt{\Tr\left(\Gamma_{k+1:\infty}^2\right)} + \sqrt{N}\norm{\Gamma_{k+1:\infty}}_{\text{op}}\right)\norm{\blambda}_2,
  \end{align}
    \item \begin{align}\label{eq:bX_f_star_applied_smooth}
        \norm{\bX_{\phi,k+1:\infty}f_{k+1:\infty}^*}_2\leq \oC{C_bX_f_star_smooth}\sqrt{N}\norm{\Gamma_{k+1:\infty}^{1/2}f_{k+1:\infty}^*}_\cH,
    \end{align}
    \item
    \begin{align}\label{eq:Gamma_phi_applied_smooth}
        \sum_{i=1}^N \norm{\left(\Gamma_{k+1:\infty}^{1/2}\phi_{k+1:\infty}\right)(X_i)}_\cH^2\leq \oC{C_sum_Gamma_phi_smooth}N\Tr\left(\Gamma_{k+1:\infty}^2\right).
    \end{align}
\end{itemize}

\begin{Proposition}\label{prop:stochastic_argument_smooth_kernel}
    Suppose Assumption~\ref{assumption:smooth} holds. Then for every $k\in\bN$ such that $(2k)^{\frac{1}{\log{d}}}\lesssim N^2/\log^8(N)$, we have $\bP(\Omega'')\geq 995/1000 - N^{-2} - 2^{\frac{\log{(2k)}}{\log{d}}}\frac{\log^8(N) }{N}$.
\end{Proposition}
\beginproof
\begin{itemize}
    \item By Proposition~\ref{prop:upper_dvoretzky_infty}, the upper bound in \eqref{eq:DM_smooth_applied} and \eqref{eq:DM_upper_smooth_applied} hold with probability at least $1 -3/1000$, by setting $\odelta{delta_P_upper_dvoretzky} = \odelta{delta_P_upper_dvoretzky_2} = \odelta{delta_P_upper_dvoretzky_3} = 1000^{-1}$, and taking $\oC{C_DMU_smooth} = 10^6 \oC{C_DMU_pre}$. Notice that due to the uniform distribution condition in Assumption~\ref{assumption:smooth}, Assumption~\ref{assumption:upper_dvoretzky_infty} is verified with $\odelta{delta_DMU_infty_1} = \odelta{delta_DMU_infty_2} = \odelta{delta_DMU_infty_3}=1$, and $\ogamma{gamma_DMU_infty_1} = \ogamma{gamma_DMU_infty_2} = \ogamma{gamma_DMU_infty_3} = 0$. See also Section~\ref{sec:diagonal_terms_multiple_descent}.
    \item By Proposition~\ref{prop:RIP}, \eqref{eq:RIP_smooth_applied} holds with probability at least $1-\odelta{delta_P_RIP}^{-1}$. One may take $\odelta{delta_P_RIP} = 1000$, $\oc{c_RIP_lower_smooth}=1/2$, $\oC{C_RIP_upper_smooth}=3/2$ provided that $\oc{c_kappa_RIP}<10^6/(8\oC{C_Rudelson}^2\oC{C_estimate_gamma_infty}^2)$.
    \item By Markov's inequality, \eqref{eq:bX_f_star_applied_smooth} holds with probability at least $1-\oC{C_bX_f_star_smooth}^{-1}$. One may take $\oC{C_bX_f_star_smooth}=1000$ for simplicity. For \eqref{eq:Gamma_phi_applied_smooth}, we observe that they hold with probability $1$ and $\oC{C_sum_Gamma_phi_smooth}=1$. This is, again, because $X$ is uniform over $\sqrt{d}S_2^{d-1}$, thus $K$ is a translation-invariant kernel.
    \item 
    For any $k\in\bN_+$, we let $\cK=2k$. Notice that: the trace of $\Gamma_{k+1:\cK}$ is of the same order as that of $\Gamma_{k+1:\infty}$. We then separate the spectrum of $\Gamma_{k+1:\infty}$ into two pieces: $(\sigma_{k+1},\sigma_{k+2},\cdots,\sigma_\cK, \sigma_{\cK+1},\sigma_{\cK+2},\cdots)$. From Section~\ref{sec:proof_multiple_descent}, we know that $L_2(\sqrt{d}S_2^{d-1},\mu)$ is spanned by spherical harmonics, which are polynomials defined on a sphere, that is, $L_2(\sqrt{d}S_2^{d-1},\mu)=\oplus_{l\in\bN}  V_{d,l}$ where $V_{d,l}$ is the space of homogeneous harmonic polynomials of degree $l$ restricted to $\sqrt{d}S_2^{d-1}$. By \eqref{eq:multicity_spherical_harmonics}, $\dim(V_{d,l})\sim d^l$. As $\dim\left(\oplus_{l=1}^\iota V_{d,l}\right)\sim d^\iota$ for any $\iota>1$, the eigenfunction of $\Gamma$ associated with eigenvalue $\sigma_\cK$ for $\cK=2k$ belongs to $\oplus_{l=1}^\iota V_{d,l}$ with $\iota\sim \lfloor \frac{\log{\cK}}{\log{d}}\rfloor$, and thus $\cH_{k+1:\cK}\subset\oplus_{\lfloor\frac{\log{k}}{\log{d}}\rfloor\lesssim l\lesssim \lfloor \frac{\log{\cK}}{\log{d}}\rfloor }^\perp V_{d,l}$ contains polynomials of degree less than $\lfloor \frac{\log{\cK}}{\log{d}} \rfloor$ (up to universal constants).
Notice that $\bX_{\phi,k+1:\infty}\bX_{\phi,k+1:\infty}^\top\succeq \bX_{\phi,k+1:\cK}\bX_{\phi,k+1:\cK}^\top$. Due to the uniform distribution condition in Assumption~\ref{assumption:smooth}, \eqref{eq:diagonal_term_assumption} is verified with $\delta = \gamma = 0$, thus Assumption~\ref{assumption:DM_L4_L2} is verified with $\delta=\gamma=0$, $\eps=6$ and $\kappa\lesssim 2^{\lfloor \frac{\log{(2k)}}{\log{d}} \rfloor}$. 
Applying Theorem~\ref{theo:DM_RKHS} to $\bX_{\phi,k+1:\cK}^\top$ implies the lower bound in \eqref{eq:DM_smooth_applied} with probability at least $1-N^{-2} - 2^{\frac{\log{(2k)}}{\log{d}}}\frac{\log^8(N) }{N}$.
\end{itemize}

\endproof

Up to the end of the proof, we place ourselves on the event $\Omega''$.

\subsubsection{Deterministic Argument}\label{sec:deterministic_argument_smooth_kernel}

The deterministic argument is almost the same as in Section~\ref{sec:proof_main_upper_k>N}, except that the extra logarithmic factors appeared in the equations from $\Omega''$. Another distinction between the proof in this section and that in Section~\ref{sec:proof_main_upper_k>N} is that, in the specific case of power decay, we only need to consider the situation where $\sigma(\square,\triangle) = \square$. This fact can be derived through the following calculation: As $\sigma_j\sim j^{-\alpha}$ for every $j\in\bN_+$, $N\sigma_1\sim N\gtrsim \kappa_{DM}\left(\lambda + k^{1-\alpha}\right)\sim\kappa_{DM}(4\lambda + \Tr\left(\Gamma_{k+1:\infty}\right))$ is always satisfied. Therefore, we only consider the case $\sigma(\square,\triangle)=\square$.

\paragraph{Estimation property of the ridge estimator $\hat f_{1:k}$}
Since many of the concepts and techniques employed in this section have already been utilized in Section~\ref{sec:proof_main_upper_k>N}, we will omit excessively detailed procedures in this section.

As in Section~\ref{sec:proof_main_upper_k<N} and Section~\ref{sec:proof_main_upper_k>N}, by homogeneous argument, we only need to consider two cases:
\begin{enumerate}
    \item $\norm{\Gamma_{1:k}^{1/2}( f_{1:k}-f_{1:k}^*)}_\cH=\square$, and $\norm{ f_{1:k}-f_{1:k}^*}_\cH\leq\triangle$, or
    \item $\norm{\Gamma_{1:k}^{1/2}( f_{1:k}-f_{1:k}^*)}_\cH\leq\square$, and $\norm{ f_{1:k}-f_{1:k}^*}_\cH=\triangle$.
\end{enumerate}
\paragraph{In case [1]} We prove that $\cQ_{f_{1:k}}>\cM_{f_{1:k}}$. We first prove the analog of \eqref{eq:multiplier_bX_f_star_k<N}, \eqref{eq:multiplier_bxi_k<N} and \eqref{eq:quadratic_k<N}. This time, we have to include extra logarithmic factors provided by equations in $\Omega''$.

\begin{itemize}
    \item For $f_{1:k}-f_{1:k}^*\in\cC$, on $\Omega''$ we have
    \begin{align}
        \cQ_{f_{1:k}} &= \norm{ \left(\bX_{\phi,k+1:\infty}\bX_{\phi,k+1:\infty}^\top + \lambda I_N\right)^{-1/2}\bX_{\phi,1:k}(f_{1:k}-f_{1:k}^*) }_2^2\notag\\
        &\geq \frac{1}{\oC{C_DMU_smooth}^2\log^2(N)\left(\lambda + \Tr\left(\Gamma_{k+1:2k}\right)\right)}\norm{\bX_{\phi,1:k}(f_{1:k}-f_{1:k}^*)}_2^2\geq \frac{\oc{c_RIP_lower_smooth}^2 N\square^2}{\oC{C_DMU_smooth}^2\log^2(N)\left(\lambda + \Tr\left(\Gamma_{k+1:2k}\right)\right)}.\label{eq:lower_quadratic_smooth}
    \end{align}
    \item The analog of \eqref{eq:multiplier_bX_f_star_k<N} is:
    On $\Omega''$, we have
\begin{align}
    &\norm{\tilde\Gamma_{1:k}^{-1/2} \bX_{\phi,1:k}^\top\left( \bX_{\phi,k+1:\infty}\bX_{\phi,k+1:\infty}^\top + \lambda I_N\right)^{-1}\bX_{\phi,k+1:\infty}f_{k+1:\infty}^* }_\cH \notag\\
    &\leq \norm{\tilde\Gamma_{1:k}^{-1/2} \bX_{\phi,1:k}^\top}_{\text{op}} \norm{\left( \bX_{\phi,k+1:\infty}\bX_{\phi,k+1:\infty}^\top + \lambda I_N\right)^{-1}}_{\text{op}}\norm{\bX_{\phi,k+1:\infty}f_{k+1:\infty}^*}_\cH \notag\\
    &\leq 2 \oC{C_DMU_smooth}\log(N)\sqrt{N}\sigma(\square,\triangle)\frac{4}{4\lambda + \Tr\left(\Gamma_{k+1:2k}\right)}\oC{C_bX_f_star_smooth}\norm{\Gamma_{k+1:\infty}^{1/2}f_{k+1:\infty}^*}_\cH.\label{eq:upper_bX_f_star_smooth}
\end{align}
\item Analog of \eqref{eq:multiplier_bxi_k<N}. Let $$D = \tilde\Gamma_{1:k}^{-1/2}\bX_{\phi,1:k}^\top\left( \bX_{\phi,k+1:\infty}\bX_{\phi,k+1:\infty}^\top + \lambda I_N \right)^{-1},$$ thus we have
\begin{align*}
    \sqrt{\Tr\left(DD^\top\right)}&= \sqrt{\Tr\left(D^\top D\right)}\leq  \norm{\left(\bX_{\phi,k+1:\infty}\bX_{\phi,k+1:\infty}^\top + \lambda I_N \right)^{-1}}_{\text{op}}\sqrt{\sum_{i=1}^N \norm{\tilde\Gamma_{1:k}\phi_{1:k}(X_i)}_\cH^2} \\
    &\leq\frac{4\oC{C_sum_Gamma_phi_smooth}\sqrt{N}}{4\lambda + \Tr\left(\Gamma_{k+1:\infty}\right)}\sqrt{\left|J_1\right|\square^2 + \triangle^2\sum_{j\in J_2}\sigma_j}.
\end{align*}Moreover,
\begin{align*}
    \norm{D}_{\text{op}} &= \norm{D^\top}_{\text{op}} \leq \norm{\left(\bX_{\phi,k+1:\infty}\bX_{\phi,k+1:\infty}^\top + \lambda I_N \right)^{-1}}_{\text{op}} \norm{\bX_{\phi,1:k}\tilde\Gamma_{1:k}^{-1/2}}_{\text{op}}\\
    &\leq 2\oC{C_DMU_smooth}\log(N)\sqrt{N}\sigma(\square,\triangle)\cdot\frac{4}{4\lambda+\Tr\left(\Gamma_{k+1:\infty}\right)}\leq \frac{8\oC{C_DMU_smooth}\log(N)\sqrt{N}\sigma\left(\square,\triangle\right)}{4\lambda + \Tr\left(\Gamma_{k+1:\infty}\right)}.
\end{align*}By Borel-TIS's inequality with
\begin{align*}
    t = \frac{\oC{C_sum_Gamma_phi_smooth}^2\left(\left|J_1\right| \square^2 + \triangle^2 \sum_{j\in J_2}\sigma_j\right)}{\oC{C_DMU_smooth}^2\sigma^2(\square, \triangle)},
\end{align*}
with probability at least $1-\bP\left((\Omega')^c\right)-\exp\left(-t(\square,\triangle)/2\right)$,
\begin{align}
    \norm{\tilde\Gamma_{1:k}^{-1/2} \bX_{\phi,1:k}^\top\left( \bX_{\phi,k+1:\infty}\bX_{\phi,k+1:\infty}^\top + \lambda I_N\right)^{-1} \bxi }_\cH \leq \frac{12\oC{C_sum_Gamma_phi_smooth}\log(N)\sqrt{N}}{4\lambda + \Tr\left(\Gamma_{k+1:\infty}\right)}\sqrt{\left|J_1\right|\square^2 + \triangle^2\sum_{j\in J_2}\sigma_j}.\label{eq:upper_bX_bxi_smooth}
\end{align}
\end{itemize}

Notice that by Assumption~\ref{assumption:smooth}, $\sigma_j\sim j^{-\alpha}$. As a result, $\Tr\left(\Gamma_{k+1:2k}\right)\sim \Tr\left(\Gamma_{k+1:\infty}\right)\sim k^{1-\alpha}$. Therefore, up to universal constants, we can bound $\Tr\left(\Gamma_{k+1:2k}\right)$ from above and from below by $\Tr\left(\Gamma_{k+1:\infty}\right)$, that is, to prove that
\begin{itemize}
    \item $$\frac{\oc{c_RIP_lower_smooth}^2 N\square^2}{\oC{C_DMU_smooth}^2\log^2(N)\left(\lambda + \Tr\left(\Gamma_{k+1:\infty}\right)\right)} >  \log(N)\sqrt{N}\sigma(\square,\triangle)\frac{8\oC{C_DMU_smooth}\oC{C_bX_f_star_smooth}}{4\lambda + \Tr\left(\Gamma_{k+1:\infty}\right)}\norm{\Gamma_{k+1:\infty}^{1/2}f_{k+1:\infty}^*}_\cH.$$
    \item $$\frac{\oc{c_RIP_lower_smooth}^2 N\square^2}{\oC{C_DMU_smooth}^2\log^2(N)\left(\lambda + \Tr\left(\Gamma_{k+1:\infty}\right)\right)} > \frac{12\oC{C_sum_Gamma_phi_smooth}\log(N)\sqrt{N}}{4\lambda + \Tr\left(\Gamma_{k+1:\infty}\right)}\sqrt{\left|J_1\right|\square^2}.$$
    \item $$\frac{\oc{c_RIP_lower_smooth}^2 N\square^2}{\oC{C_DMU_smooth}^2\log^2(N)\left(\lambda + \Tr\left(\Gamma_{k+1:\infty}\right)\right)} > \frac{12\oC{C_sum_Gamma_phi_smooth}\log(N)\sqrt{N}}{4\lambda + \Tr\left(\Gamma_{k+1:\infty}\right)}\sqrt{\triangle^2\sum_{j\in J_2}\sigma_j}.$$
    \item $$\frac{\oc{c_RIP_lower_smooth}^2 N\square^2}{\oC{C_DMU_smooth}^2\log^2(N)\left(\lambda + \Tr\left(\Gamma_{k+1:\infty}\right)\right)} > \norm{\tilde\Gamma_{1:k}^{-1/2}f_{1:k}^*}_\cH.$$
\end{itemize}
\paragraph{In case [2]} we prove that $\cQ_{f_{1:k}}>\cR_{f_{1:k}}$ by showing that
\begin{itemize}
    \item $$\triangle^2 >  \log(N)\sqrt{N}\sigma(\square,\triangle)\frac{8\oC{C_DMU_smooth}\oC{C_bX_f_star_smooth}}{4\lambda + \Tr\left(\Gamma_{k+1:\infty}\right)}\norm{\Gamma_{k+1:\infty}^{1/2}f_{k+1:\infty}^*}_\cH.$$
    \item $$\triangle^2 > \frac{12\oC{C_sum_Gamma_phi_smooth}\log(N)\sqrt{N}}{4\lambda + \Tr\left(\Gamma_{k+1:\infty}\right)}\sqrt{\left|J_1\right|\square^2}.$$
    \item $$\triangle^2 > \frac{12\oC{C_sum_Gamma_phi_smooth}\log(N)\sqrt{N}}{4\lambda + \Tr\left(\Gamma_{k+1:\infty}\right)}\sqrt{\triangle^2\sum_{j\in J_2}\sigma_j}.$$
    \item $$\triangle^2 > \norm{\tilde\Gamma_{1:k}^{-1/2}f_{1:k}^*}_\cH.$$
\end{itemize}
As in Section~\ref{sec:proof_main_upper_k<N} and Section~\ref{sec:proof_main_upper_k>N}, we choose $\triangle$ such that $\square/\triangle = \sqrt{\kappa_{DM}\left(4\lambda + \Tr\left(\Gamma_{k+1:\infty}\right)\right)/N}$.
We check that there exists an absolute constant $\nC\label{C_square_smooth}$ such that we can take
\begin{align}\label{eq:square_smooth}
    \begin{aligned}
        \square &= \oC{C_square_smooth}\log^3(N)\max\bigg\{ \sigma_\xi\sqrt{\frac{\left|J_1\right|}{N}}, \sigma_\xi\sqrt{\frac{\sum_{j\in J_2}\sigma_j}{4\lambda + \Tr\left(\Gamma_{k+1:\infty}\right)}}, \norm{\Gamma_{k+1:\infty}^{1/2}f_{k+1:\infty}^*}_\cH,\\
        &{\norm{\tilde\Gamma_{1:\mathrm{thre}}^{-1/2}f_{1:k}^*}_\cH\frac{2\lambda + 3\Tr\left(\Gamma_{k+1:\infty}\right)}{N}} \bigg\}.
    \end{aligned}
\end{align}

\paragraph{Noise absorption} The price for overfitting is almost the same as that in Section~\ref{sec:proof_main_upper_k>N}, except for two differences:
\begin{itemize}
    \item When dealing with $\norm{\Gamma_{k+1:\infty}^{1/2}A\bX_{\phi,1:k}(f_{1:k}^* - \hat f_{1:k})}_\cH$.
    We no longer have \eqref{eq:lower_bound_of_quadratic_k>N}, but
    \begin{align*}
        \cQ_{\hat f_{1:k}}\geq \frac{4}{4\lambda + \Tr\left(\Gamma_{k+1:2k}\right)}\norm{\bX_{\phi,1:k}(\hat f_{1:k} - f_{1:k}^*)}_2^2.
    \end{align*}If we define $\vertiiii{\cdot}$ and $\theta$ as in \eqref{eq:def_vertiiii} and \eqref{eq:def_theta_k>N}, by \eqref{eq:square_smooth} and $\cQ_{\hat f_{1:k}}+\cR_{\hat f_{1:k}}\leq\left|\cM_{\hat f_{1:k}}\right|$ we obtain $\norm{\bX_{\phi,1:k}(\hat f_{1:k}-f_{1:k}^*)}_2\leq \frac{1}{2}\sqrt{\left(4\lambda + \Tr\left(\Gamma_{k+1:2k}\right)\right)\theta\square}$. However, since $\Tr\left(\Gamma_{k+1:2k}\right)\sim \Tr\left(\Gamma_{k+1:\infty}\right)$, we can replace the right hand side by $\nc\label{c_useless_1}\sqrt{\left(4\lambda + \Tr\left(\Gamma_{k+1:\infty}\right)\right)\theta\square}$ up to some small absolute constant $\oc{c_useless_1}$. The reminder part of the proof remains the same as in Section~\ref{sec:proof_main_upper_k>N}.
    \item When dealing with $\norm{\Gamma_{k+1:\infty}^{1/2}A\bX_{\phi,k+1:\infty}f_{k+1:\infty}^*}_\cH$ and $\norm{\Gamma_{k+1:\infty}^{1/2}A\bxi}_\cH$. We once again use $\Tr\left(\Gamma_{k+1:2k}\right)\sim \Tr\left(\Gamma_{k+1:\infty}\right)$. Moreover, we apply \eqref{eq:DM_upper_smooth_applied}, \eqref{eq:bX_f_star_applied_smooth} and \eqref{eq:Gamma_phi_applied_smooth} by including extra logarithmic factors. The proof readily follows.
\end{itemize}

Before concluding this section, we additionally state a property that allows us to obtain an upper bound without Assumption~\ref{assumption:smooth} and when $\lambda\gtrsim \Tr(\Gamma_{k+1:\infty})$.

\begin{Proposition}\label{prop:smooth_large_regularization}
    Grant Assumption~\ref{assumption:upper_dvoretzky_infty}. For any $k\in\bN$ and $\lambda\geq 0$ such that $\lambda\gtrsim\Tr(\Gamma_{k+1:\infty})$, $N\leq\oc{c_kappa_DM}\kappa_{DM}d_\lambda^*(\Gamma_{k+1:\infty}^{-1/2}B_\cH)$, $\sum_{j\in J_2}\sigma_j \leq \kappa_{DM}(4\lambda + \Tr\left(\Gamma_{k+1:\infty}\right))\left(1-\frac{\left|J_1\right|}{N}\right)$, and $\kappa_{DM}\left(4\lambda + \Tr\left(\Gamma_{k+1:\infty}\right)\right)\geq N\left(R_N^*(\oc{c_kappa_RIP})\right)^2$. Then with constant probability, $\norm{\hat f_\lambda - f^*}_{L_2}\lesssim \log^3(N)r_{\lambda,k}^*$.
\end{Proposition}
\beginproof 
We only need to prove that $\bP(\Omega'')>0$. This is because, as seen in Section~\ref{sec:deterministic_argument_smooth_kernel}, we have seen that condition on $\Omega''$, we obtain the conclusion of Proposition~\ref{prop:smooth_large_regularization}. To prove $\bP(\Omega'')>0$, we only need to establish the left-hand side of \eqref{eq:DM_smooth_applied}, as the other inequalities have already been proven in Proposition~\ref{prop:stochastic_argument_smooth_kernel}. The left-hand side of \eqref{eq:DM_smooth_applied} is trivial, as it always holds when $\lambda\gtrsim \Tr(\Gamma_{k+1:\infty})$ (since $\bX_{\phi,k+1:\infty}\bX_{\phi,k+1:\infty}^\top\succeq 0$). We have provided the proof of this in Section~\ref{sec:proof_DM_lambda_dominating}.

\endproof

\subsection{Auxiliary proofs}\label{sec:aux_proofs}

\subsubsection{Proof of Proposition~\ref{prop:upper_dvoretzky} and Proposition~\ref{prop:upper_dvoretzky_infty}}\label{sec:proof_upper_dvoretzky}

\paragraph{Proof of Proposition~\ref{prop:upper_dvoretzky}} 
It follows from \cite[Theorem 4]{tsigler_benign_2023} that there exists an absolute constant $\nC\label{C_DMU_2}$, such that with probability at least $1-\frac{\oc{c_P_DMU}}{N^\eps}-\gamma'$,
\begin{align}
    &\norm{\bX_{\phi,k+1:\infty}\Gamma_{k+1:\infty}\bX_{\phi,k+1:\infty}^\top}_{\text{op}} \label{eq:proof_upper_DM} \\
    &\leq \Tr\left(\Gamma_{k+1:\infty}^2\right)\left( 1+\delta' + \oC{C_DMU_2}(1+\delta')\left(\frac{1}{N} + \sqrt{\frac{N\norm{\Gamma_{k+1:\infty}^2}_{\text{op}}}{(1+\delta')^2 {\Tr\left(\Gamma_{k+1:\infty}^2\right)}}} + \frac{N\norm{\Gamma_{k+1:\infty}^2}_{\text{op}}}{(1+\delta')^2 \Tr\left(\Gamma_{k+1:\infty}^2\right)} \right) \right)\notag\\
    &\leq (\oC{C_DMU_2}+1)(1+\delta')\Tr\left(\Gamma_{k+1:\infty}^2\right) + \oC{C_DMU_2}\sqrt{N\norm{\Gamma_{k+1:\infty}^2}_{\text{op}}\Tr\left(\Gamma_{k+1:\infty}^2\right)} + \frac{\oC{C_DMU_2}}{1+\delta'}N\norm{\Gamma_{k+1:\infty}^2}_{\text{op}}\notag\\
    &\leq \left((\oC{C_DMU_2}+1)(1+\delta')+\frac{\oC{C_DMU_2}}{2}\right)\Tr\left(\Gamma_{k+1:\infty}^2\right) + \oC{C_DMU_2}\left(\frac{1}{2}+\frac{1}{\delta'}\right)N\norm{\Gamma_{k+1:\infty}^2}_{\text{op}}.\notag
\end{align}This indicates that there exists an absolute constant $\oC{C_DMU}>0$ such that for any $\vlambda\in\bR^N$,
\begin{align*}
    \norm{\Gamma_{k+1:\infty}^{1/2}\bX_{\phi,k+1:\infty}^\top\vlambda}_{\cH} \leq \oC{C_DMU}\left( \sqrt{\Tr\left(\Gamma_{k+1:\infty}^2\right)} + \sqrt{N}\norm{\Gamma_{k+1:\infty}}_{\text{op}}\right)\norm{\vlambda}_2.
\end{align*}

\paragraph{Proof of Proposition~\ref{prop:upper_dvoretzky_infty}}
Apply Lemma~\ref{lemma:Rudelson_method} below to $\cF = B_{\cH_{k+1:\infty}}$. Then $\sup_{f\in\cF}\norm{f}_{L_2}=\norm{\Gamma_{k+1:\infty}}_{\text{op}}^{1/2}$, and by Lemma~\ref{lemma:estimate_gamma_infty} below applied to $\bar\Gamma_{1:k}$ is the identity operator (recall the definition of $\bar\Gamma_{1:k}$ in Proposition~\ref{prop:RIP}), and $R=1$, $\norm{\gamma_2\left(P_\sigma B_{\cH_{k+1:\infty}},\ell_\infty\right)}_{L_2}\leq\oC{C_estimate_gamma_infty}Q\left(B_{\cH_{k+1:\infty}}\right)\log{N}$. By Assumption~\ref{assumption:upper_dvoretzky_infty}, $Q(B_\cH)\leq \odelta{delta_DMU_infty_1}\sqrt{\Tr\left(\Gamma_{k+1:\infty}\right)}$. Therefore, there exist absolute constants $\oC{C_Rudelson},$ and $\oC{C_DMU_pre}$, such that for any $\odelta{delta_P_upper_dvoretzky}$ with probability at least $1-\odelta{delta_P_upper_dvoretzky}$,
\begin{align}\label{eq:upper_dvoretzky_infty}
    \begin{aligned}
        \norm{\bX_{\phi,k+1:\infty}}_{\text{op}}^2 &= {\underset{\norm{f}_\cH\leq 1}{\sup}\sum_{i=1}^N f^2(X_i)}\\ 
    &\leq \odelta{delta_P_upper_dvoretzky}^{-1} N\norm{\Gamma_{k+1:\infty}}_{\text{op}}\\
    &+ \odelta{delta_P_upper_dvoretzky}^{-1} \oC{C_Rudelson}\norm{\Gamma_{k+1:\infty}}_{\text{op}}^{1/2}\sqrt{N\odelta{delta_DMU_infty_1}\Tr\left(\Gamma_{k+1:\infty}\right)}\log{N} + \odelta{delta_P_upper_dvoretzky}^{-1} \odelta{delta_DMU_infty_1}\Tr\left(\Gamma_{k+1:\infty}\right)\log^2{N}\\
    &\leq \odelta{delta_P_upper_dvoretzky}^{-1}\oC{C_DMU_pre}^2\left({N\norm{\Gamma_{k+1:\infty}}_{\text{op}}} + {\odelta{delta_DMU_infty_1}\Tr\left(\Gamma_{k+1:\infty}\right)\log^2{N}}\right).
    \end{aligned}
\end{align}
When $N\norm{\Gamma_{k+1:\infty}}_{\text{op}}\leq \Tr\left(\Gamma_{k+1:\infty}\right)\log^2(N)$, with probability at least $1-\odelta{delta_P_upper_dvoretzky}-\ogamma{gamma_DMU_infty_1}$, for any $\vlambda\in\bR^N$,
    \begin{align*}
        \norm{\bX_{\phi,k+1:\infty}^\top\vlambda}_\cH\leq \odelta{delta_P_upper_dvoretzky}^{-1/2}\oC{C_DMU_pre}\log{N}\sqrt{\odelta{delta_DMU_infty_1}\Tr\left(\Gamma_{k+1:\infty}\right)}\norm{\vlambda}_2.
    \end{align*}
Similarly, we have: with probability at least $1-\odelta{delta_P_upper_dvoretzky_2} - \odelta{delta_P_upper_dvoretzky_3}- \ogamma{gamma_DMU_infty_2} - \ogamma{gamma_DMU_infty_3}$, for any $\vlambda\in\bR^N$,
    \begin{align*}
    \norm{\Gamma_{k+1:\infty}^{1/2}\bX_{\phi,k+1:\infty}^\top}_{\text{op}}&\leq \odelta{delta_P_upper_dvoretzky_2}^{-1/2}\oC{C_DMU_pre}\left(\sqrt{N}\norm{\Gamma_{k+1:\infty}}_{\text{op}} + \log{N}\sqrt{\odelta{gamma_DMU_infty_2}\Tr\left(\Gamma_{k+1:\infty}^2\right)}\right)\norm{\vlambda}_2,\\
    \norm{\tilde\Gamma_{1:k}^{-1/2} \bX_{\phi,1:k}^\top \vlambda }_\cH &\leq \odelta{delta_P_upper_dvoretzky_3}^{-1/2}\oC{C_DMU_pre}\left(\log(N)\sqrt{\odelta{gamma_DMU_infty_3}}\sqrt{ \square^2 \left|J_1\right| + \triangle^2\sum_{j\in J_2}\sigma_j } + \sqrt{N}\sigma(\square,\triangle)\right)\norm{\vlambda}_2.
\end{align*}


\subsubsection{Proof of Proposition~\ref{prop:IP} and Proposition~\ref{prop:RIP}}\label{sec:proof_RIP}

By homogeneity, it suffices to prove the same result for all $f\in \cH_{1:k}$ such that $\norm{\Gamma_{1:k}^{1/2}f}_\cH=1$, that is, to prove that $\oc{c_RIP_lower}^2\leq (1/N)\sum_{i=1}^N  \left(\Gamma_{1:k}^{-1/2}g\right)^2(X_i)\leq \oC{C_RIP_upper}^2$ where $g=\Gamma_{1:k}^{1/2}f$ for any $f\in \cH_{1:k}$ such that $\norm{\Gamma_{1:k}^{1/2}f}_\cH=1$.

\beginproof We separate the proof into two parts: the upper bound and the lower bound.
\paragraph{Upper bound of Proposition~\ref{prop:IP}} For the upper bound of RIP, we employ Proposition~\ref{prop:upper_dvoretzky} applied to $$\left[\Gamma_{1:k}^{-1/2}\phi_{1:k}(X_1) | \cdots | \Gamma_{1:k}^{-1/2}\phi_{1:k}(X_N)\right].$$ Notice that $\bX_{\phi,1:k}^\top$ and $\bX_{\phi,1:k}$ have the same operator norm. By Proposition~\ref{prop:upper_dvoretzky},
\begin{align*}
    \bP\left(\sup\left(\norm{\Gamma_{1:k}^{-1/2}\bX_{\phi,1:k}^\top\vlambda}_\cH:\, \norm{\vlambda}_2 = 1\right)\leq \oC{C_DMU}\left( \sqrt{k} +\sqrt{N} \right)\right) \geq 1-\ogamma{gamma_RIP}-\frac{\oc{c_P_DMU} }{N^\eps}.
\end{align*}
As a result, with the same probability, $\norm{\bX_{\phi,1:k}\Gamma_{1:k}^{-1/2}}_{\text{op}}^2 = \sup\left(\sum_{i=1}^N\left(\Gamma_{1:k}^{-1/2}g\right)^2(X_i):\, \norm{g}_\cH=1 \right)\leq 2\oC{C_DMU}^2\left( k+N\right)$, which is precisely the upper bound of Proposition~\ref{prop:IP} with $\oC{C_RIP_upper}=\sqrt{2\oC{C_DMU}^2(1+\oc{c_RIP})}$.
\paragraph{Lower bound of Proposition~\ref{prop:IP}} For the lower bound, we use the following lemma taken from \cite[Proposition 4]{zhivotovskiy_dimension-free_2024}:
\begin{Lemma}\label{lemma:nikita}
    Assume that $M_1,\cdots,M_N$ are independent copies of a positive semi-definite symmetric random matrix $M$ with mean $\bE M = \Sigma$. Let $M$ satisfy that for some $\kappa\geq 1$, $\sqrt{\bE (\vv^\top M\vv)^2}\leq\kappa^2 \vv^\top\Sigma\vv$ for all $\vv\in\bR^d$. Then for any $t>\log{2}$, with probability at least $1-2\exp(-t)$,
    \begin{align*}
        \inf\left(\frac{1}{N}\sum_{i=1}^N\vv^\top M_i\vv - \norm{\Sigma^{1/2}\vv}_2^2\right) \geq -7\kappa^2\sqrt{\frac{ \frac{\Tr\left(\Sigma\right)}{\norm{\Sigma}_{\text{op}}} +t}{N}}.
    \end{align*}
\end{Lemma}As we can embed $\Gamma_{1:k}^{-1/2}\cH_{1:k}$ into $\ell_2^k$, we apply Lemma~\ref{lemma:nikita} to $\Sigma=I_k$, $M_i = \Gamma_{1:k}^{-1/2}\phi_{1:k}(X_i)\otimes\Gamma_{1:k}^{-1/2}\phi_{1:k}(X_i)$, $\kappa = \kappa''$. By Assumption~\ref{assumption:RIP}, with probability at least $1-2\exp(-t)$, for any $g\in\cH_{1:k}$ such that $\norm{g}_\cH = 1$, we have $\sum_{i=1}^N \left(\Gamma_{1:k}^{-1/2}g\right)^2(X_i)\geq N - 7(\kappa'')^2\sqrt{N}\sqrt{ k + t }$. Set $t=k$, if $\oc{c_RIP}$ is taken such that $7\sqrt{2}(\kappa'')^2\sqrt{\oc{c_RIP}}\leq 3/4$, that is $\oc{c_RIP}\leq \frac{9}{1568(\kappa'')^4}$.

In summary, Proposition~\ref{prop:IP} is verified with $\oc{c_RIP_lower}=1/2$, $\oC{C_RIP_upper}=\sqrt{2\oC{C_DMU}^2(1+\oc{c_RIP})}$ and $\bar p_{RIP} = \ogamma{gamma_RIP} + \frac{\oc{c_P_DMU} }{N^\eps} + 2\exp(-k)$.

\endproof

\paragraph{Proof of Proposition~\ref{prop:RIP}} We first introduce a notation that will be used only in this paragraph. Let $(T,d)$ be a metric space and $B\subset T$. For every $\alpha>0$, we define the Talagrand's $\gamma_\alpha$-functional as
\begin{align*}
    \gamma_\alpha\left( B,d \right) := \underset{(\cA_n)_{n\geq 0}}{\inf}\underset{\vv\in B}{\sup}\sum_{n\geq 0}2^{n/\alpha}d\left(\vv,\cA_n\right),
\end{align*}where the infimum is taken over all admissible sequences, that is, $(\cA_n)_{n\geq 0}$ is an increasing partition of $B$ such that $\left|\cA_0\right|=1$ and for any $n\geq 1$, $\left|\cA_n\right|\leq 2^{2^n}$, see \cite[Chapter 2]{talagrand_upper_2021}. Here $d(\vv,\cA_n)=\inf\left(d(\vv,\vu): \vu\in \cA_n\right)$ is the usual distance between a point and a set.

The following lemma is taken from \cite[Theorem 1.2]{guedon_subspaces_2007}
\begin{Lemma}\label{lemma:Rudelson_method}
    There exists an absolute constant $\nC\label{C_Rudelson}$ such that the following holds. Let $(\Omega,\mu)$ be a measured space, $N\in\bN_+$ and $(X_i)_{i\in[N]}$ be $N$ i.i.d. random variables with values in $\Omega$. Let $\cF$ be a class of real-valued functions defined on $\Omega$. We have
    \begin{align*}
        \bE\underset{f\in\cF}{\sup}\left|\frac{1}{N}\sum_{i=1}^N f^2(X_i) - \bE f^2(X) \right| \leq \oC{C_Rudelson}\max\left( \underset{f\in\cF}{\sup}\norm{f}_{L_2}\sqrt{\frac{\bE \gamma_2^2\left(P_\sigma\cF,\norm{\cdot}_{\ell_\infty}\right)}{N}}, \frac{\bE \gamma_2^2\left(P_\sigma\cF,\norm{\cdot}_{\ell_\infty}\right)}{N}\right),
    \end{align*}and
    \begin{align*}
        \bE\underset{f\in\cF}{\sup}\sqrt{\frac{1}{N}\sum_{i=1}^N f^2(X_i)}\leq\oC{C_Rudelson}\max\left(\underset{f\in\cF}{\sup}\norm{f}_{L_2}, \sqrt{\frac{\bE \gamma_2^2\left(P_\sigma\cF,\norm{\cdot}_{\ell_\infty}\right)}{N}}\right),
    \end{align*}where $P_\sigma\cF = \left\{\left(f(X_1),\cdots,f(X_N)\right)^\top:\, f\in\cF \right\}$.

\end{Lemma}

We apply Lemma~\ref{lemma:Rudelson_method} to $\cF = \cC(R)$. By the definition of $\cC(R)$, see \eqref{eq:def_cone_RIP}, we have $\bE f^2(X)=1$ for any $f\in\cC(R)$ and thus $\sup_{f\in\cF}\norm{f}_{L_2}=1$. In order to apply Lemma~\ref{lemma:Rudelson_method}, we need to estimate $\bE\gamma_2^2\left(P_\sigma\cC(R),\norm{\cdot}_{\ell_\infty}\right)$.

If we let $\norm{f} = \max\left(\norm{f}_{L_2}, R\norm{f}_\cH\right)$ for any $f\in\cH_{1:k}$, one can check that there exists an ellipsoid norm $\norm{\bar\Gamma_{1:k}^{1/2}\cdot}_\cH$, as defined in Proposition~\ref{prop:RIP}, such that for any $f\in\cH_{1:k}$, $\norm{f}\leq \norm{\bar\Gamma_{1:k}^{1/2}f}_\cH\leq \sqrt{2}\norm{f}$. The unit ball $\bar\Gamma_{1:k}^{-1/2}B_\cH$ has principle lengths $\left(\min\left(\frac{1}{R}, \frac{1}{\sqrt{\sigma_j}}\right)\right)_{j\in[k]}$.
The relationship between $\norm{\cdot}$ and $\cC(R)$ is the unit ball of $\norm{\cdot}$ is $B_{\norm{\cdot}} = R^{-1}B_{\cH_{1:k}}\cap \Gamma_{1:k}^{-1/2}B_{\cH_{1:k}}$ which contains $\cC(R)$. Furthermore, $B_{\norm{\cdot}}\subset \sqrt{2}\bar\Gamma_{1:k}^{-1/2}B_\cH$. Therefore, $\mu$-almost surely, $\gamma_2\left(P_\sigma\cC(R),\norm{\cdot}_{\ell_\infty}\right)\leq \gamma_2\left(P_\sigma B_{\norm{\cdot}},\norm{\cdot}_{\ell_\infty}\right)\leq \sqrt{2}\gamma_2\left(P_\sigma\bar\Gamma_{1:k}^{-1/2}B_\cH,\norm{\cdot}_{\ell_\infty}\right)$.

The following lemma is a rewrite of \cite[Theorem 4.7]{mendelson_regularization_2010}:
\begin{Lemma}\label{lemma:estimate_gamma_infty}
    There exists an absolute constant $\nC\label{C_estimate_gamma_infty}$ such that
    \begin{align*}
        \norm{\gamma_2\left(P_\sigma\bar\Gamma_{1:k}^{-1/2}B_\cH,\ell_\infty\right)}_{L_2}\leq\oC{C_estimate_gamma_infty}\norm{\max_{i\in[N]} \norm{ \bar\Gamma_{1:k}^{-1/2}\phi(X_i) }_\cH }_{L_\infty}\log{N}.
    \end{align*}
\end{Lemma}

The proof of this Lemma is almost identical to that found in \cite{mendelson_regularization_2010}. For the sake of completeness, we include it here. Given two convex bodies $K,L\subset\bR^N$, we denote the covering number of $K$ by $L$ as $N(K,L)$, see \cite[section 4.2]{vershynin_high-dimensional_2018}.

\paragraph{Proof of Lemma~\ref{lemma:estimate_gamma_infty}}
Notice that with probability $1$, we have
\begin{align}\label{eq:condition_Q}
    \max_{i\in[N]}\norm{\bar\Gamma_{1:k}^{-1/2}\phi(X_i)}_\cH\leq \norm{\max_{i\in[N]} \norm{ \bar\Gamma_{1:k}^{-1/2}\phi(X_i) }_\cH }_{L_\infty}.
\end{align}Denote $Q(\bar\Gamma_{1:k}^{-1/2}B_\cH\cap RB_\cH)$ as $\norm{\max_{i\in[N]} \norm{ \bar\Gamma_{1:k}^{-1/2}\phi(X_i) }_\cH }_{L_\infty}$.
Set $\norm{f}_E = \max_{i\in[N]}\left|f(X_i)\right|$, and let $B_E$ be the unit ball of $\norm{\cdot}_E$. Recall that $\bar \Gamma_{1:k}^{-1/2}=\sum_{j=1}^k \min\{\frac{1}{R},\frac{1}{\sqrt{\sigma_j}}\}\varphi_j\otimes\varphi_j$, then for any $i\in[k]$, $\bar\Gamma_{1:k}^{-1/2}\varphi_i = \min\{\frac{1}{R},\frac{1}{\sqrt{\sigma_i}}\}\varphi_i$. For every $\eps>0$, the covering number of $\bar \Gamma_{1:k}^{-1/2}B_\cH$ by $\eps B_E$ satisfies
\begin{align*}
    N(\bar \Gamma_{1:k}^{-1/2}B_\cH,\eps B_E) = N(B_\cH,\eps \bar \Gamma_{1:k}^{1/2}B_E),
\end{align*}and $f\in \eps \bar\Gamma_{1:k}^{1/2}B_E$ if and only if $\max_{j\in[N]}\left|\left\langle f, \bar\Gamma_{1:k}^{-1/2}\phi(X_i)\right\rangle_\cH\right|\leq \eps$. Set $\psi(X_i) = \bar\Gamma_{1:k}^{-1/2}\phi(X_i)$, and $\norm{f}_{\bar E}=\max_{j\in[N]}\left|\left\langle f,\psi(X_i)\right\rangle_\cH\right|$ with corresponding unit ball $B_{\bar E}$. Then $N(\bar \Gamma_{1:k}^{-1/2}B_\cH,\eps B_E) = N(B_\cH, \eps B_{\bar E})=N(B_\cH^N, \eps B_{\bar E})$, where $B_\cH^N$ is the unit ball in the subspace of $\cH$ spanned by $\psi(X_i)$'s.  Let $G\in\cH$ be a standard Gaussian random vector, as $B_{\bar E}$ is a convex body in $\bR^N$, $\bE\norm{G}_{\bar E}=\bE\max_{j\in[N]}\left|\left\langle G,\psi(X_i)\right\rangle_\cH\right|\lesssim \sqrt{\log{N}}\max_{i\in[N]}\norm{\psi(X_i)}_\cH$, see, for example, \cite[Exercise 7.5.10]{vershynin_high-dimensional_2018}. Recall that $\psi(X_i)=\bar\Gamma_{1:k}^{-1/2}\phi(X_i)$, thus from \eqref{eq:condition_Q}, we know that $\bE\norm{G}_{\bar E}\lesssim\sqrt{\log{N}}Q(\bar\Gamma_{1:k}^{-1/2}B_\cH\cap RB_\cH)$. By Sudakov's inequality, see, for example, \cite[Theorem 7.4.1]{vershynin_high-dimensional_2018}, $\log N(B_\cH^N, \eps B_{\bar E})\lesssim \log(N)Q^2(\bar\Gamma_{1:k}^{-1/2}B_\cH\cap RB_\cH)/\eps^2$. In particular, the diameter of $B_\cH^N$ with respect to the norm $\norm{\cdot}_{\bar E}$ is at most $Q(\bar\Gamma_{1:k}^{-1/2}B_\cH\cap RB_\cH)\sqrt{\log{N}}$ (up to universal constant), and we denote this diameter by $D_2$. For small $\eps$, we use volumetric estimate \cite[Corollary 4.1.15]{artstein-avidan_asymptotic_2015}. For any norm $\norm{\cdot}_X$ on $\bR^N$ with unit ball $B_X$ and every $\eps>0$, $N(B_X,\eps B_X)\leq (3/\eps)^N$. We see that $B_\cH^N$ and $B_{\bar E}$ are convex bodies in $\bR^N$ (up to isometric), thus there exists an absolute constant $\nc\label{c_volume_estimate}$ such that for any $0<\eps<\delta$,
\begin{align*}
    \log N(B_\cH^N, \eps B_{\bar E})&\leq \log N(B_\cH^N, \delta B_{\bar E}) + \log N(\delta B_{\bar E}, \eps B_{\bar E})\\
    &\leq \oc{c_volume_estimate} \frac{Q^2(\bar\Gamma_{1:k}^{-1/2}B_\cH\cap RB_\cH)\log(N)}{\delta^2} + N\log\left(\frac{3\delta}{\eps}\right).
\end{align*}Take $\delta^2 = \oc{c_volume_estimate}Q^2(\bar\Gamma_{1:k}^{-1/2}B_\cH\cap RB_\cH)\frac{\log{N}}{N}$. There then exists an absolute constant $\nc\label{c_upper_eps}$ such that $\eps\leq \oc{c_upper_eps}Q(\bar\Gamma_{1:k}^{-1/2}B_\cH\cap RB_\cH)\sqrt{\log(N)/N}=:\eps_0$, and $\log N(B_\cH^N, \eps B_{\bar E})\lesssim N\log(\eps_0/\eps)$. By Dudley's integral, see, for example, \cite[Exercise 8.5.7]{vershynin_high-dimensional_2018}, condition on the event such that \eqref{eq:condition_Q} holds,
\begin{align*}
    \gamma_2^2\left(P_\sigma\bar\Gamma_{1:k}^{-1/2}B_\cH,\ell_\infty\right)&\lesssim \int_0^\infty \eps\log N(\bar\Gamma_{1:k}^{-1/2}B_\cH,\eps B_E)d\eps = \int_0^\infty \eps\log N(B_\cH^N, \eps B_{\bar E}) d\eps \\
    &\lesssim \int_0^{\eps_0} N\eps \log\left(\frac{\eps_0}{\eps}\right) d\eps + \int_{\eps_0}^{D_2}\frac{Q^2(\bar\Gamma_{1:k}^{-1/2}B_\cH\cap RB_\cH)\log(N)}{\eps}d\eps.
\end{align*}Using the change of variables for $\eta=\eps_0/\eps$, there exist absolute constants $\nC\label{C_upper_first_integral_1}$ and $\nC\label{C_upper_first_integral_2}$ such that the first integral is bounded by $\oC{C_upper_first_integral_1}N\eps_0^2\int_0^1\eta\log(\eta^{-1})d\eta \leq \oC{C_upper_first_integral_2}Q^2(\bar\Gamma_{1:k}^{-1/2}B_\cH\cap RB_\cH)\log(N)$. For the second integral, notice that there exists an absolute constant $\nc\label{c_eps_0}$ such that $\eps_0 = \oc{c_eps_0}D_2/\sqrt{N}$, so the second integral is bounded from above by
\begin{align*}
    &\oC{C_upper_first_integral_2}Q^2(\bar\Gamma_{1:k}^{-1/2}B_\cH\cap RB_\cH)\log(N)\left(\log(D_2) - \log(\eps_0)\right)\\
    &= \oC{C_upper_first_integral_2}Q^2(\bar\Gamma_{1:k}^{-1/2}B_\cH\cap RB_\cH)\log(N)\left(\frac{1}{2}\log(N) - \log(\oc{c_eps_0})\right)\\
    &\lesssim Q^2(\bar\Gamma_{1:k}^{-1/2}B_\cH\cap RB_\cH)\log^2(N).
\end{align*}By Fubini's theorem, there exists an absolute constant $\oC{C_estimate_gamma_infty}$ such that $\sqrt{\bE\left(\gamma_2^2\left(P_\sigma\bar\Gamma_{1:k}^{-1/2}B_\cH,\ell_\infty\right)\right)}\leq \oC{C_estimate_gamma_infty}Q(\bar\Gamma_{1:k}^{-1/2}B_\cH\cap RB_\cH)\log(N).$
\endproof

By Lemma~\ref{lemma:estimate_gamma_infty}, we have $\bE\gamma_2^2\left(P_\sigma\cC(R),\norm{\cdot}_{\ell_\infty}\right)\leq 2\oC{C_estimate_gamma_infty}^2Q^2\left(\Gamma_{1:k}^{-1/2}B_\cH\cap RB_\cH\right)\log^2{N}$. Recall the definition of $R_N(\oc{c_kappa_RIP})$ from \eqref{eq:def_fixed_point}, for $R>R_N(\oc{c_kappa_RIP})$, by Lemma~\ref{lemma:Rudelson_method} and Lemma~\ref{lemma:estimate_gamma_infty}, for any $0<\odelta{delta_P_RIP}<1$, by Markov's inequality, with probability at least $1-\odelta{delta_P_RIP}$, for any $f\in\cC(R_N(\oc{c_kappa_RIP}))$,
\begin{align*}
    \oc{c_RIP_lower}^2:= 1-\odelta{delta_P_RIP}^{-1}\sqrt{2\oc{c_kappa_RIP}}\oC{C_Rudelson}\oC{C_estimate_gamma_infty}\leq \frac{1}{N}\sum_{i=1}^N f^2(X_i) \leq 1+ \odelta{delta_P_RIP}^{-1}\sqrt{2\oc{c_kappa_RIP}}\oC{C_Rudelson}\oC{C_estimate_gamma_infty} =: \oC{C_RIP_upper}^2.
\end{align*}Therefore, Proposition~\ref{prop:RIP} is valid with probability at least $1-\odelta{delta_P_RIP}$.



\subsubsection{Proof of Proposition~\ref{prop:linearization}}\label{sec:proof_linearization}

The main approach for this proof is to make use of Theorem~\ref{theo:DM_RKHS}. To achieve this, we need to verify Assumption~\ref{assumption:DM_L4_L2} and check the Dvoretzky-Milman condition. As in Section~\ref{sec:multiple_descent}, the $L_8-L_2$ equivalence is verified with $\kappa\sim 2^L$. We will now proceed to verify the Dvoretzky-Milman condition. As in Section~\ref{sec:multiple_descent}, since $L<\infty$ but $d\to\infty$, we let $k=\sum_{0\leq i\leq\iota}d^i$, thus $k\sim d^\iota$, and $\Tr\left(\Gamma_{k+1:\infty}\right)\sim 1$, $\norm{\Gamma_{k+1:\infty}}_{\text{op}}\sim d^{-(\iota+1)}$. As a result, $N \lesssim \kappa_{DM}d_0^*\left(\Gamma_{k+1:\infty}^{-1/2}B_\cH\right)$.

For \eqref{eq:diagonal_term_assumption}, since $h\in C^\infty$, for every $X_i$, there exists $\xi_{i,i}\in(1,\norm{X_i}_2^2/d)$ or $(\norm{X_i}_2^2/d,1)$ such that $h\left(\norm{X_i}_2^2/d\right) = \norm{\phi(X_i)}_\cH^2 = h(1)+h'(\xi_{i,i})\left(\norm{X_i}_2^2/d-1\right)$.
By \cite[Theorem 3.1.1]{vershynin_high-dimensional_2018} together with a union bound over $i\in[N]$, there exists an absolute constant $c_0>0$ such that for any $t\geq 0$,
\begin{align*}
    \bP\left(\underset{i\in[N]}{\max}\left|\frac{\norm{X_i}_2^2}{d} - 1\right|\geq t\right)\leq N\exp\left(-c_0\frac{t^2 d}{\norm{x_1}_{\psi_2}^2}\right).
\end{align*}Let $t=1/\sqrt{\log\log{N}}$. We have $\lim_{N\to\infty}t = 0$, and $\lim_{N,d\to\infty}N\exp\left(-c_0(t^2 d)/\norm{x_1}_{\psi_2}^2\right)\to 0$. This indicates that $\max_{i\in[N]}\left|\norm{X_i}_2^2/d-1\right|\to 0$ in probability as $N,d\to\infty$. Moreover, as $h\in C^\infty$, $h'(\xi_{i,i})$ is bounded in a neighborhood of $1$ uniformly over $i\in[N]$. These observations imply that $h(\norm{X_i}_2^2/d)\to h(1)$ in probability uniformly over $i\in[N]$ as $N,d\to\infty$. Consequently, we may take $\delta=\gamma=0$ in \eqref{eq:diagonal_term_assumption}. 

Now that we have checked all the conditions of Theorem~\ref{theo:DM_RKHS}, we can apply this theorem. The conclusion of Proposition~\ref{prop:linearization} is a straightforward consequence of Theorem~\ref{theo:DM_RKHS}.
\endproof

\subsubsection{Estimating the fixed point of RIP}\label{sec:estimate_RIP}

In the following, we estimate $R_N(\oc{c_kappa_RIP})$ in some cases. 
We now separate into three regimes depending the following condition: For some absolute constant $\ntheta\label{theta_embedding}$ such that $0\leq \otheta{theta_embedding}\leq 1$ and an absolute constant $\nC\label{C_embedding}$ such that for any $j\in[k]$, $\norm{f_j}_{L_\infty}\leq\oC{C_embedding}\norm{f_j}_\cH^\otheta{theta_embedding}\norm{f_j}_{L_2}^{1-\otheta{theta_embedding}}$. This condition is widely used in RKHS literature, for instance, \cite{mendelson_regularization_2010,steinwart_optimal_2009,fischer_sobolev_2020,li_asymptotic_2023}.

\paragraph{$\otheta{theta_embedding}=0$} In this case,
    \begin{align*}
        R_N(\oc{c_kappa_RIP}) = \inf\left\{R>0:\, \sqrt{\oC{C_embedding} \sum_{j=1}^k \sigma_j \wedge R^2 } \leq \oc{c_kappa_RIP}R\frac{\sqrt{N}}{\log{N}} \right\}.
    \end{align*}Therefore, up to a logarithmic factor, we recover the fixed point in the sub-Gaussian case. More precisely,
    \begin{enumerate}
        \item When $k\leq \frac{\oc{c_kappa_RIP}^2}{\oC{C_embedding}}\frac{N}{\log^2(N)}$. Then $R_N(\oc{c_kappa_RIP})=0$.
        \item When there exists $k_0\in\left[\lfloor\frac{\oc{c_kappa_RIP}^2}{\oC{C_embedding}} \rfloor N\right]$ such that $\sum_{j\geq k_0}\sigma_j \leq \left(\frac{\oc{c_kappa_RIP}^2}{\oC{C_embedding}}\frac{N}{\log{N}}-k_0+1\right)\sigma_{k_0}$, let
        \begin{align*}
            k^{**} = \max\left(k_0\in\left[\lfloor\frac{\oc{c_kappa_RIP}^2}{\oC{C_embedding}} \rfloor N\right]:\,\sum_{j\geq k_0}\sigma_j \leq \left(\frac{\oc{c_kappa_RIP}^2}{\oC{C_embedding}}\frac{N}{\log{N}}-k_0+1\right)\sigma_{k_0} \right),
        \end{align*}then $R_N(\oc{c_kappa_RIP})\leq\sigma_{k^{**}}$.
        \item When for all $k_0\in\left[\lfloor\frac{\oc{c_kappa_RIP}^2}{\oC{C_embedding}} \rfloor N\right]$, we have $\sum_{j\geq k_0}\sigma_j > \left(\frac{\oc{c_kappa_RIP}^2}{\oC{C_embedding}}\frac{N}{\log{N}}-k_0+1\right)\sigma_{k_0}$. Then $R_N(\oc{c_kappa_RIP})\leq \frac{\sqrt{\oC{C_embedding}}}{\oc{c_kappa_RIP}}{\frac{\sqrt{\Tr\left(\Gamma_{1:k}\right)}\log{N}}{\sqrt{N}}}$.
    \end{enumerate}
\paragraph{$0<\otheta{theta_embedding}<1$} We have 
\begin{align}\label{eq:estimate_fixed_point_interpolation_space}
    Q\left(\Gamma_{1:k}^{-1/2}B_\cH\cap RB_\cH\right)\leq\sqrt{\oC{C_embedding}\sum_{j=1}^k \frac{\sigma_j^{1-\otheta{theta_embedding}}}{R^2}\wedge \frac{1}{\sigma_j^{\otheta{theta_embedding}}}}.
\end{align}
The estimation of $R_N(\oc{c_kappa_RIP})$ in this case is somewhat complicated, we thus only provide two trivial estimates:
    \begin{enumerate}
        \item If $\sum_{j=1}^k \sigma_j^{-\otheta{theta_embedding}}\leq \frac{\oc{c_kappa_RIP}^2}{\oC{C_embedding}}\frac{N}{\log^2(N)}$, then $R_N(\oc{c_kappa_RIP})=0$.
        \item We always have $R_N(\oc{c_kappa_RIP})\leq \sqrt{\frac{\oC{C_embedding}}{\oc{c_kappa_RIP}^2}\frac{\sum_{j=1}^k \sigma_j^{1-\otheta{theta_embedding}}}{N}}\log{N}$.
    \end{enumerate}
As an example, when $\sigma_j\sim j^{-\alpha}$, one may take $\otheta{theta_embedding}=\alpha^{-1}$, see \cite[Corollary 3]{steinwart_optimal_2009}, see also \cite[Lemma 5.1]{mendelson_regularization_2010}.
\begin{enumerate}
    \item When $1<\alpha<2$, we have $\sum_{j=1}^k \sigma_j^{1-\otheta{theta_embedding}}\sim k^{-\alpha}$. Hence $R_N(\oc{c_kappa_RIP})\lesssim\sqrt{\frac{k^{-\alpha}}{N}}\log{N}$.
    \item When $\alpha=2$, we have $\sum_{j=1}^k \sigma_j^{1-\otheta{theta_embedding}}\sim \log{k}$. Hence $R_N(\oc{c_kappa_RIP})\lesssim\frac{\log(N)\sqrt{\log(k)}}{\sqrt{N}}$.
    \item When $\alpha>2$, we have $\sum_{j=1}^k \sigma_j^{1-\otheta{theta_embedding}}\sim 1$. Hence $R_N(\oc{c_kappa_RIP})\lesssim \frac{\log(N)}{\sqrt{N}}$.
\end{enumerate}
Apart from the above estimate, when $\alpha>2$, if there exists $k_0< k$ such that $\sum_{j=k_0+1}^k \frac{j^{1-\alpha}}{R^2}\lesssim \frac{N}{\log^2(N)}-\sum_{j=1}^{k_0}j$, we can set $R_N^2(\oc{c_kappa_RIP})\lesssim \frac{k_0^{2-\alpha}}{\frac{N}{\log^2{N}}-k_0^2}$.

\paragraph{$\otheta{theta_embedding}=1$} We have $Q\left(\Gamma_{1:k}^{-1/2}B_\cH\cap RB_\cH\right)\leq\sqrt{\oC{C_embedding}\sum_{j=1}^k \frac{1}{R^2}\wedge \frac{1}{\sigma_j}}.$ The estimation of $R_N(\oc{c_kappa_RIP})$ is also difficult in this case. We provide two trivial estimates.
    \begin{enumerate}
        \item When $\sum_{j=1}^k\sigma_j^{-1}\leq\frac{\oc{c_kappa_RIP}^2}{\oC{C_embedding}}\frac{N}{\log^2{N}}$, $R_N(\oc{c_kappa_RIP})=0$.
        \item We always have $R_N(\oc{c_kappa_RIP})\leq \frac{\sqrt{\oC{C_embedding}}}{\oc{c_kappa_RIP}}\frac{\log{N}\sqrt{k}}{\sqrt{N}}$.
    \end{enumerate}

\subsubsection{Proof of Proposition~\ref{prop:KRR_smooth_applied}}\label{sec:proof_KRR_smooth_applied}

By \eqref{eq:def_DM_dimension}, we have:
\begin{align*}
    &d_\lambda^*\left(\Gamma_{k+1:2k}^{-1/2}B_\cH\right)\gtrsim \frac{\Tr\left(\Gamma_{k+1:\infty}\right)+\lambda}{\norm{\Gamma_{k+1:\infty}}_{\mathrm{op}}} \sim \frac{k^{1-\alpha}+N^{\frac{-\alpha}{1+(s\wedge 2)\alpha} +1}}{(k+1)^{-\alpha}}\\
    &\gtrsim_\alpha N^{\frac{\alpha}{1+(s\wedge 2)\alpha}}N^{\frac{-\alpha}{1+(s\wedge 2)\alpha} +1}\geq \frac{N}{\oc{c_kappa_DM}\kappa_{DM}}.
\end{align*}By \eqref{eq:def_J1}, we know that $J_1=[k]$ and $J_2=\emptyset$, hence $\tilde\Gamma_{1:\mathrm{thre}}=\Gamma_{1:k}$. Moreover, $k^{\frac{1}{\log{d}}} = N^{\frac{1}{(1+(s\wedge 2)\alpha)\log{d}}} \lesssim \frac{N}{\log^8(N)}$ and by \eqref{eq:estimate_fixed_point_interpolation_space} and $\otheta{theta_embedding}=\frac{1}{\alpha}$, we check that we can take $R\to 0$, that is, there exists an absolute constant $\nC\label{C_RIP_equation_1}$ such that
\begin{align*}
    \sqrt{\sum_{j=1}^k \frac{\sigma_j^{1-\otheta{theta_embedding}}}{R^2}\wedge \frac{1}{\sigma_j^{\otheta{theta_embedding}}}} \leq \sqrt{\sum_{j=1}^k\frac{1}{\sigma_j^{\otheta{theta_embedding}}}}\leq \oC{C_RIP_equation_1}k \sim N^{\frac{1}{1+(s\wedge 2)\alpha}}<\oc{c_kappa_RIP}\frac{\sqrt{N}}{\log^8(N)},
\end{align*}where we have used the fact  that for any $j\in\bN_+$, $\sigma_j\sim j^{-\alpha}$ and $s\alpha>1$. Hence we can take $R_N(\oc{c_kappa_RIP})=0$, and thus by \eqref{eq:def_cone_RIP}, $\mathrm{cone}\left(\cC(R_N(\oc{c_kappa_RIP}))\right)=\cH_{1:k}$. Up to now, we have checked all the assumptions in Proposition~\ref{prop:KRR_smooth_upper}. We now compute the terms appeared in the upper bound in Proposition~\ref{prop:KRR_smooth_upper}.
\begin{enumerate}
    \item $\sigma_\xi\frac{\sqrt{Nk^{1-2\alpha}}}{\lambda + k^{1-\alpha}}\lesssim \sigma_\xi N^{-\frac{(s\wedge 2)\alpha}{2(1+(s\wedge 2)\alpha)}}$,
    \item $\sigma_\xi \sqrt{\frac{\left|J\right|}{N}}\lesssim \sigma_\xi N^{-\frac{(s\wedge 2)\alpha}{2(1+(s\wedge 2)\alpha)}}$.
    \item We compute the bias term (referred to the terms without the $\sigma_\xi$ factor) at once. We first consider the case when $0<s<2$.
    \begin{align}
        \norm{\Gamma_{k+1:\infty}^{1/2}f_{k+1:\infty}^*}_\cH^2 + \norm{\Gamma_{1,k}^{-1/2} f_{1:k}^*}_\cH^2 \frac{\lambda^2 + k^{2(1-\alpha)}}{N^2} \lesssim \sum_{j=k+1}^\infty \sigma_j a_j^2 + N^{-\frac{2\alpha}{1+s\alpha}}\sum_{j=1}^k \sigma_j^{-1}a_j^2.\label{eq:bias_smooth_compute_1}
    \end{align}
Recall that $\sigma_j\sim j^{-\alpha}$ is a decreasing function of $j$. When $j>k$, by $\sigma_j^2\vee\sigma_k^2 = \sigma_k^2$, we have: $\sigma_j\leq \frac{2\sigma_k^2}{\sigma_j^2+\sigma_k^2}\sigma_j$; When $j\leq k$, we know that $\sigma_j^2\vee \sigma_k^2 = \sigma_j^2$, hence $\sigma_j^{-1}\leq \frac{2\sigma_j}{\sigma_j^2+\sigma_k^2}$. Plugging these facts into \eqref{eq:bias_smooth_compute_1}, we obtain that
\begin{align*}
    \sum_{j=k+1}^\infty \sigma_j a_j^2 + N^{-\frac{2\alpha}{1+s\alpha}}\sum_{j=1}^k \sigma_j^{-1}a_j^2\leq 2\sigma_k^2 \sum_{j=k+1}^\infty  \frac{a_j^2\sigma_j}{\sigma_j^2+\sigma_k^2} + 2N^{-\frac{2\alpha}{1+s\alpha}}\sum_{j=1}^k \frac{a_j^2\sigma_j}{\sigma_j^2+\sigma_k^2}.
\end{align*}Recall that by our choice of $k$, $\sigma_k^2\sim N^{-\frac{2\alpha}{1+s\alpha}}$, we further derive that
\begin{align}
    2\sigma_k^2 \sum_{j=k+1}^\infty  \frac{a_j^2\sigma_j}{\sigma_j^2+\sigma_k^2} + 2N^{-\frac{2\alpha}{1+s\alpha}}\sum_{j=1}^k \frac{a_j^2\sigma_j}{\sigma_j^2+\sigma_k^2}\lesssim \sigma_k^2\sum_{j=1}^\infty a_j^2 \frac{\sigma_j}{\sigma_j^2+\sigma_k^2} \lesssim \sum_{j=1}^\infty a_j^2\sigma_j^{1-s} \left(\frac{\sigma_j^{\frac{s}{2}}\sigma_k}{\sigma_j+\sigma_k}\right)^2.\label{eq:bias_smooth_compute_2}
\end{align}
For $0<s<2$, let $g(x)=\frac{\lambda x^{s/2}}{\lambda+x}$. Then $\sup_{x>0}g(x)=\frac{2-s}{2}\left(\frac{s}{s-2}\right)^{s/2}\lambda^{s/2}$. Let $x=\sigma_j$ and $\lambda = \sigma_k$. Then \eqref{eq:bias_smooth_compute_2} together with \eqref{eq:equivalent_source_condition} indicate that there exists an absolute constant $\nC\label{C_bias_smooth_compute_1}$ depending only on $s$ such that
\begin{align*}
    \eqref{eq:bias_smooth_compute_2} < \oC{C_bias_smooth_compute_1} \sigma_k^{s}\sum_{j=1}^\infty a_j^2 \sigma_j^{1-s}\leq \oC{C_bias_smooth_compute_1} \sigma_k^{s}\sim N^{-\frac{\alpha s}{1+\alpha s}}.
\end{align*}
\item When $s\geq 2$, recall that the choice of $k = N^{\frac{1}{1+2\alpha}}$. Then
\begin{align*}
    &\norm{\Gamma_{k+1:\infty}^{1/2}f_{k+1:\infty}^*}_\cH^2 = \sum_{j=k+1}^\infty \sigma_j a_j^2 = \sum_{j=k+1}^\infty \sigma_j^s \sigma_j^{1-s}a_j^2<\sigma_{k+1}^s\sum_{j=k+1}^\infty \sigma_j^{1-s}a_j^2\\
    &\lesssim k^{-\alpha s} = N^{-\frac{\alpha s}{1+2\alpha}} \leq N^{-\frac{2\alpha}{1+2\alpha}}.
\end{align*}
Moreover,
\begin{align*}
    &\frac{\lambda^2 + \Tr^2(\Gamma_{k+1:\infty})}{N^2}\norm{\Gamma_{1,k}^{-1/2}f_{1:k}^*}_\cH^2 \lesssim \sigma_k^2\sum_{j=1}^k \sigma_j^{-1}a_j^2 = \sigma_k^2\sum_{j=1}^k \sigma_j^{1-s}a_j^2 \sigma_j^{s-2}\\
    &\leq \sigma_k^2 \sigma_1^{s-2}\sum_{j=1}^k \sigma_j^{1-s}a_j^2\lesssim N^{-\frac{2\alpha}{1+2\alpha}}.
\end{align*}
\end{enumerate}
Combining the above three terms together with Proposition~\ref{prop:KRR_smooth_upper}, we know that with constant probability,
\begin{align*}
    \norm{\hat f_\lambda - f^*}_{L_2} \lesssim \log^3(N)N^{-\frac{\alpha (s\wedge 2)}{2(1+(s\wedge 2)\alpha)}}.
\end{align*}
\endproof

\subsubsection{Proof of Proposition~\ref{prop:jaouad}}\label{sec:proof_prop_jaouad}

By the choice of $k$ we obtain $\sigma_{k+1}\lesssim \frac{\lambda}{N}\lesssim \sigma_k$, and $\Tr(\Gamma_{k+1:\infty})+\lambda\gtrsim N\norm{\Gamma_{k+1:\infty}}_{\text{op}}$ and we have Dvoretzky-Milman theorem.
\begin{enumerate}
    \item The bias term. We first consider the case when $f^*\in\cH$.
    By the choice of $k$, $\sigma_j\vee \frac{\lambda+\Tr(\Gamma_{k+1:\infty})}{N}\gtrsim\sigma_j$ for any $j\leq k$.
We have
\begin{align*}
    &\norm{\tilde\Gamma_{1,\mathrm{thre}}^{-1/2}f_{1:k}^*}_\cH\frac{\lambda+\Tr(\Gamma_{k+1:\infty})}{N} \lesssim\norm{\Gamma_{1:k}^{-1/2}f_{1:k}^*}_\cH\frac{\lambda+\Tr(\Gamma_{k+1:\infty})}{N}\\
    &\leq\norm{f_{1:k}^*}_\cH\left(\frac{\lambda}{N\sqrt{\sigma_k}} + \frac{\Tr(\Gamma_{k+1:\infty})}{N\sqrt{\sigma_k}}\right)
\end{align*}
Recall that we have chosen $k$, such that $\sigma_k\gtrsim \frac{\lambda}{N}$. Therefore
\begin{align*}
    &\norm{\tilde\Gamma_{1,\mathrm{thre}}^{-1/2}f_{1:k}^*}_\cH\frac{\lambda+\Tr(\Gamma_{k+1:\infty})}{N}\lesssim \sqrt{\frac{\lambda}{N}}\left(1+\frac{\Tr(\Gamma)}{\lambda}\right)\norm{f_{1:k}^*}_\cH.
\end{align*}
On the other hand,
\begin{align*}
    &\left(\frac{\lambda}{N}+\frac{\Tr^2(\Gamma)}{N\lambda}\right)\norm{\Gamma^{1/2}\left(\Gamma+\frac{\lambda}{N}I\right)^{-1/2}f^*}_\cH^2\\
    &\sim \left(\frac{\lambda}{N}+\frac{\Tr^2(\Gamma)}{N\lambda}\right)\left(\norm{f_{1:k}^*}_\cH^2 + \frac{N}{\lambda}\norm{\Gamma_{k+1:\infty}^{1/2}f_{k+1:\infty}^*}_\cH^2\right)\\
    &\gtrsim \frac{\lambda}{N}\left(1+\frac{\Tr^2(\Gamma)}{\lambda^2}\right)\norm{f_{1:k}^*}_\cH^2 + \norm{\Gamma_{k+1:\infty}^{1/2}f_{k+1:\infty}^*}_\cH^2.
\end{align*}
The case when $f^*$ may not belong to $\cH$ follows directly from \cite[Lemma 7.2, Equation 7.20 and Equation 7.21]{bach_learning_2024} (notice that their $\lambda$ corresponds to our $\lambda/N$).

\item Variance term. 
     By the definition of $\Gamma(\Gamma+(\lambda/N)I)^{-1}$, we have 
    \begin{align*}
        \sigma_j\left(\Gamma(\Gamma+(\lambda/N)I)^{-1}\right)\sim &\begin{cases}
            1, & j\leq k\\
            \frac{N\sigma_j}{\lambda}, &j\geq k+1.
        \end{cases}
    \end{align*}
    Then
\begin{align*}
    \frac{\Tr\left(\Gamma\left(\Gamma+\frac{\lambda}{N}I\right)^{-1}\right)}{N} \sim \frac{k}{N}+\frac{\Tr(\Gamma_{k+1:\infty})}{\lambda}.
\end{align*}Therefore,
\begin{align*}
    \sigma_\xi\sqrt{\frac{k}{N}}\lesssim \sigma_\xi \sqrt{ \frac{\Tr\left(\Gamma\left(\Gamma+\frac{\lambda}{N}I\right)^{-1}\right)}{N}}.
\end{align*}Moreover, recall that $\sigma_{k+1}\lesssim \frac{\lambda}{N}$, hence
\begin{align*}
    &\sigma_\xi^2\frac{N\Tr(\Gamma_{k+1:\infty}^2)}{(\lambda+\Tr(\Gamma_{k+1:\infty}))^2}\leq \sigma_\xi^2\frac{N\Tr(\Gamma_{k+1:\infty}^2)}{\lambda^2} \leq \sigma_\xi^2\frac{N\sigma_{k+1}}{\lambda}\frac{\Tr(\Gamma_{k+1:\infty})}{\lambda}\lesssim \sigma_\xi^2\frac{\Tr(\Gamma_{k+1:\infty})}{\lambda}.
\end{align*}
\end{enumerate}
\endproof

\subsubsection{Proof of Proposition~\ref{prop:linear}}\label{sec:proof_linear_case}

The proof strategy for this Proposition involves verifying the stochastic argument of Theorem~\ref{theo:upper_KRR} in the Gaussian case.

We prove \eqref{eq:DM_applied} by the Gaussian Dvoretzky-Milman theorem, see \cite[Theorem 5]{lecue_geometrical_2022}, \eqref{eq:RIP_applied}, \eqref{eq:RIP_k>M_applied} by the Gaussian isomorphy property and restricted isomorphy property, see \cite[Theorem 6]{lecue_geometrical_2022}, \eqref{eq:upper_dvoretzky} by \cite[Proposition 5]{lecue_geometrical_2022}. \eqref{eq:bX_f_star_applied} and \eqref{eq:Gamma_phi_applied} follow from Bernstein's inequality. The deterministic arguments from Section~\ref{sec:proof_main_upper_k<N} and Section~\ref{sec:proof_main_upper_k>N} then still hold. Since the proof is exactly the same, we have omitted it here.

The variance term in the lower bound comes from \cite{tsigler_benign_2023}. Given that the lower bound on the bias term proved in \cite{tsigler_benign_2023} takes the form of a Bayesian lower bound (where they assume $f^*$ to be a random vector), we present here a general form of the lower bound. The proof methodology is inspired by the work in \cite{lecue_geometrical_2022}. The proof presented here only requires that the marginal distribution of $X$ satisfies the equivalence of $L_{2+\epsilon}$ and $L_2$, and that $X$ is a centered symmetric random vector, where $\epsilon > 2$. This encompasses the case of a Gaussian random vector as a special instance.

\paragraph{Stochastic Argument} The underlying principle of the stochastic arguments employed  is based on the lemma presented in \cite[Lemma 9]{bartlett_benign_2020}:
\begin{Lemma}\label{lemma:BLTT_lemma_9}
    Suppose $p\in\bN\cup\{\infty\}$, $\{\eta_i\}_{i=1}^p$ is a sequence of non-negative random variables, and $\{t_i\}_{i=1}^p$ is a sequence of non-negative real numbers(at least one of which is strictly positive) such that, for some $0<\delta<1$ and any $i\leq p$, $\bP\left(\eta_i>t_i\right)\geq 1-\delta$, then
    \begin{align*}
        \bP\left(\sum_{i=1}^p \eta_i \geq \frac{1}{2}\sum_{i=1}^p t_i\right)\geq 1-2\delta.
    \end{align*}
\end{Lemma}
The purpose of Lemma~\ref{lemma:BLTT_lemma_9} is to present an alternate approach to the union bound.
Consider the case where $\phi: \mathbf{x} \in \mathbb{R}^d \mapsto \mathbf{x} \in \mathbb{R}^d$, specifically for the scenario of ridge regression. In this context, the bias term of the lower bound involves $\bX_{\phi, k+1:\infty}$, which is equivalent to $\bX P_{k+1:\infty}$. Here, we employ $P_{k+1:\infty}$ to accommodate the possibility that $d$ may be infinite. That is to say, when $j > d$, we straightforwardly set $\varphi_j = \mathbf{0}$ (zero vector) and $\sigma_j = 0$. In this case, we view $\cH$ as $\bR^d$ by viewing $f\in\cH$ as some $\vv\in\bR^d$, that is, $f(\cdot)=\left<\vv,\cdot\right>$, hence $\left<\phi(\vx),f\right>_\cH=f(\vx)=\left<\vv,\vx\right>$.

Let $\nc\label{c_lower_1}$, $\nc\label{c_lower_2}$, $\nC\label{C_lower_2}$, $\nC\label{C_lower_3}$ be absolute constants, where $\oc{c_lower_2}$ depends on $\kappa$.
For each $j\in\bN$, define $\Omega_j$ as the random event on which:
\begin{itemize}
    \item $\bX_{\phi,k+1:\infty}\bX_{\phi,k+1:\infty}^\top + \lambda I_N$ is well-conditioned, that is,
    \begin{align}\label{eq:DM}
        \oc{c_lower_1}\left(\lambda + \Tr\left(\Gamma_{k+1:\infty}\right)\right)\leq \sigma_N\left(\bX_{\phi,k+1:\infty}\bX_{\phi,k+1:\infty}^\top + \lambda I_N\right).
    \end{align}
    \item  Let $z_j = (1/\sqrt{\sigma_j})\bX_{\phi,k+1:\infty}\varphi_j\in\bR^N$, and $Q$ be a (random) projection matrix that is independent with $z_j$. Condition on $Q$,
    \begin{align}\label{eq:quadratic_form_markov}
        z_j^\top Q z_j\leq \oC{C_lower_2}\Tr(Q) + \oC{C_lower_2}N\norm{Q}_{op}.
    \end{align}
    \item For $z_j$ defined as above,
    \begin{align}\label{eq:z_j_markov}
        \oc{c_lower_2} N \leq \norm{z_j}_2^2 \leq \oC{C_lower_3}N.
    \end{align}
\end{itemize}We emphasize that all the absolute constants $\oc{c_lower_1}$, $\oC{C_lower_2}$, $\oc{c_lower_2}$ and $\oC{C_lower_3}$ are independent with $j$.

In the subsequent portion of this paragraph, we establish the inequality $2\bP\left(\Omega_j\right)\geq 99/100$ for every $j\in\bN$, assuming the conditions of the Dvoretzky-Milman theorem are satisfied. Because $b>\oC{C_lower_kappa_DM}/\kappa_{DM}$, \eqref{eq:DM} is satisfied with a probability of at least $1-\bar p_{DM}$, hence our focus now shifts to establishing the validity of \eqref{eq:quadratic_form_markov} and \eqref{eq:z_j_markov} with a high level of probability.

\begin{Lemma}\label{lemma:quadratic_form_markov}
    For each $j\in\bN$, define $z_j = (1/\sqrt{\sigma_j})\bX_{\phi,k+1:\infty}\varphi_j\in\bR^N$. Let $Q$ be a (random) projection matrix that is independent with $z_j$. Then there exists an absolute constant $\nc\label{c_p_lower_1}<10^{-3}$ such that condition on $Q$, for each $j\in\bN$, with probability at least $1-\oc{c_p_lower_1}$,
    \begin{equation*}
        z_j^\top Q z_j\leq \oC{C_lower_2}\Tr(Q) + \oC{C_lower_2}N\norm{Q}_{op}.
    \end{equation*}
\end{Lemma}

\beginproof Recall that for each $j\in\bN$, $\bX_{\phi,k+1:\infty}\varphi_j = \left(\varphi_j(X_i)\right)_{i\in[N]}$, thus $z_j$ has i.i.d. coordinates $(z_{ji})_{i\in[N]}$. Furthermore, for each $i\in[N]$, $\norm{\varphi_j(X_i)}_{L_2}^2 = \sigma_j$, thus $\norm{z_{ji}}_{L_2}=1$. We write $\bE z_j = (\bE z_{ji})_{i\in[N]}$. By the linearity,
\begin{eqnarray*}
    \bE z_j^\top Q z_j &=& \bE \left[(z_j-\bE z_j + \bE z_j)^\top Q (z_j - \bE z_j + \bE z_j)\right] \\
    &=& \bE\left[ (z_j - \bE z_j)^\top Q (z_j - \bE z_j + \bE z_j) \right] + \bE\left[(\bE z_j)^\top Q (z_j - \bE z_j + \bE z_j) \right]\\
    &=&\bE\left[ (z_j - \bE z_j)^\top Q (z_j - \bE z_j )\right] + \bE\left[ (z_j - \bE z_j)^\top Q (\bE z_j) \right]\\
    &+& \bE\left[(\bE z_j)^\top Q (z_j - \bE z_j ) \right] + \bE\left[(\bE z_j)^\top Q (\bE z_j) \right]\\
    &=&2\bE[z_{j1}^2]\Tr(Q) + 2\bE\left[ (z_j - \bE z_j)^\top Q (\bE z_j) \right] + (\bE z_j)^\top Q (\bE z_j)\\
    &=& 2\bE[z_{j1}^2]\Tr(Q) + (\bE z_j)^\top Q (\bE z_j)\\
    &\leq&2\Tr(Q) + \norm{Q}_{op}\norm{\bE z_j}_2^2\leq 2\Tr(Q) + N\norm{Q}_{op}.
\end{eqnarray*}
By Markov's inequality, there exists an absolute constant $\oc{c_p_lower_1}<10^{-3}$ such that with probability at least $1-\oc{c_p_lower_1}$,
\begin{equation*}
    z_j^\top Q z_j \leq \oC{C_lower_2}\Tr(Q) + \oC{C_lower_2}N\norm{Q}_{op}.
\end{equation*}
\endproof

\begin{Lemma}\label{lemma:z_j_markov}
    Grant the following assumption: there exists an constant $\kappa\geq 1$ such that for any $f\in\cH$, we have $\norm{f}_{L_4}\leq\kappa\norm{f}_{L_2}$. Then there exist absolute constants $\nc\label{c_p_lower_2}$, $\nc\label{c_p_lower_3}$, $\oc{c_lower_2}$ and $\oC{C_lower_3}$ where $\oc{c_lower_2}$ and $\oc{c_p_lower_2}$ depend on $\kappa$, such that $\oc{c_p_lower_3}<10^{-3}$ and with probability at least $1-\oc{c_p_lower_3}-\exp(-\oc{c_p_lower_2}N)$,
    \begin{equation*}
        \oc{c_lower_2}N \leq \norm{z_j}_2^2 \leq \oC{C_lower_3}N.
    \end{equation*}
\end{Lemma}
\beginproof $\bE\norm{z_j}_2^2= N$. Therefore there exists an absolute constant $\oc{c_p_lower_3}<10^{-3}$ such that by Markov's inequality, with probability at least $1-\oc{c_p_lower_3}$,
\begin{equation*}
    \norm{z_j}_2^2 \leq \oc{c_p_lower_3}^{-1}N.
\end{equation*}For the lower side, we use Paley-Zygmund inequality(see, for example, \cite[section 3.3]{de_la_pena_decoupling_1999}) together with a Bernstein's inequality(see, for example, \cite[section 2.8]{vershynin_high-dimensional_2018}) for selectors. Hence there exists an absolute constant $\nc\label{c_lower_3}$ depending on $\kappa$ such that for one $i\in[N]$,
\begin{equation*}
    \bP\left( z_{ji}^2 \geq \frac{1}{2}\bE z_{ji}^2 \right) \geq \frac{1}{4}\frac{\left(\bE z_{ji}^2\right)^2}{\bE z_{ji}^4}\geq \oc{c_lower_3},
\end{equation*}where we used $\bE z_{ji}^4 = (1/\sigma_j^2)\bE\varphi_j^4(X_i)\leq \kappa^4(1/\sigma_j^2) \left(\bE\varphi_j^2(X_i)\right)^2 = \kappa^4 \left(\bE z_{ji}^2\right)^2$ since $\varphi_j\in\cH$.

By Bernstein's inequality for $(\delta_i)_{i\in[N]}$, where $\delta_i = \1_{\{z_{ji}^2 \geq \frac{1}{2}\bE z_{ji}^2\}}$, there exists an absolute constant $\oc{c_p_lower_2}$ such that with probability at least $1-\exp(-\oc{c_p_lower_2}N)$, we have $\left|\left\{i\in[N]:\, z_{ji}^2 \geq \frac{1}{2}\bE z_{ji}^2 \right\}\right|\geq \oc{c_lower_3}N$. Therefore, there exist absolute constants $\oc{c_lower_2}$ depending on $\kappa$ and $\oC{C_lower_3}$ such that with probability at least $1-\oc{c_p_lower_3}-\exp(-\oc{c_p_lower_2}N)$,
\begin{equation*}
     \oc{c_lower_2}N \leq \norm{z_j}_2^2 \leq \oC{C_lower_3}N.
\end{equation*}
\endproof

\paragraph{Deterministic Argument}
\begin{Proposition}\label{prop:lower_bias_term}
Assume that $\phi:\vx\in\bR^d\mapsto\vx\in\bR^d$ and $X$ is centered and symmetric with covariance matrix $\Gamma$, and satisfies Assumption~\ref{assumption:DM_L4_L2} with $\eps>2$. Moreover, assume that $Z=\Gamma^{-1/2}X$ is a random vector with independent coordinates. There exist absolute constants $b$, $\nc\label{c_lower_4}$, $\nc\label{c_lower_5}$, $\nC\label{C_lower_4}$ , where $\oc{c_lower_5}$ depends on $b$ and $b$ depends on $\kappa_{DM}$, such that if $N\geq \oC{C_lower_4}$, then
    \begin{align*}
        \bE\norm{\Gamma^{1/2}\left(\bX_{\phi}^\top\left(\bX_{\phi}\bX_{\phi}^\top + \lambda I_N \right)^{-1} - I \right)f^*}_{\cH}^2 &\geq \oc{c_lower_4} \norm{\Gamma_{k_{b,\lambda}^*+1:\infty}^{1/2}f_{k_{b,\lambda}^*+1:\infty}^*}_{\cH}^2\\
        &+ \oc{c_lower_5}\norm{\Gamma_{1:k_{b,\lambda}^*}^{-1/2}f_{1:k_{b,\lambda}^*}^*}_{\cH}^2\left( \frac{\lambda + \Tr(\Gamma_{k_{b,\lambda}^* +1:\infty})}{N}\right)^2.
    \end{align*}
\end{Proposition}
\beginproof Let $U = \sum_{j\in\bN}\eps_j\varphi_j\otimes\varphi_j:\cH\to\cH$ where $(\eps_j)_{j\in\bN}$ are i.i.d. Rademacher random variables. Then
\begin{align*}
    UU^\top = \left(\sum_{j\in\bN}\eps_j \varphi_j\otimes \varphi_j\right)\left(\sum_{i\in\bN}\eps_i \varphi_i\otimes \varphi_i\right)=\sum_{j\in\bN}\varphi_j\otimes \varphi_j=I = U^\top U
\end{align*}implies that $U$ is a unitary operator, and the absolute convergence of $\sum_{j\in\bN}\eps_j \varphi_j\otimes \varphi_j$ implies that
\begin{align*}
    U\Gamma &= \left(\sum_{j\in\bN}\eps_j \varphi_j\otimes \varphi_j\right) \left(\sum_{j\in\bN}\sigma_j \varphi_j \otimes \varphi_j\right) = \sum_{j\in\bN}\eps_j\sigma_j \varphi_j\otimes \varphi_j\\
    &= \left(\sum_{j\in\bN}\eps_j \varphi_j\otimes \varphi_j\right) \left(\sum_{j\in\bN}\sigma_j \varphi_j \otimes \varphi_j\right)=\Gamma U,
\end{align*}
thus $U$ commutes with $\Gamma$. As $\phi:\vx\in\bR^d\mapsto\vx\in\bR^d$ and $X$ is symmetric, we have
\begin{align*}
    &\bE\norm{\Gamma^{1/2}(\bX_{\phi}^\top(\bX_{\phi}\bX_{\phi}^\top + \lambda I_N)^{-1}\bX_{\phi}-I)Uf^*}_\cH \\
    &=\bE\norm{\Gamma^{1/2}\left((\bX_{\phi} U)^\top\left((\bX_{\phi} U)(\bX_{\phi} U)^\top + \lambda I_N\right)^{-1}\bX_{\phi} U-I\right)f^*}_\cH\\ 
    &= \bE\norm{\Gamma^{1/2}\left(\bX_{\phi}^\top\left(\bX_{\phi}\bX_{\phi}^\top + \lambda I_N \right)^{-1}\bX_\phi - I \right)f^*}_{\cH}.
\end{align*}
In the following, we obtain lower bound for $\bE\norm{\Gamma^{1/2}(\bX_{\phi}^\top(\bX_{\phi}\bX_{\phi}^\top + \lambda I_N)^{-1}\bX_{\phi}-I)Uf^*}_\cH$.
For all $j\in\bN$, denote $f_j^* := \left< f^*, \varphi_j\right>_{\cH}$, then
\begin{eqnarray*}
    &&\bE_U\norm{\Gamma^{1/2}(\bX_{\phi}^\top(\bX_{\phi}\bX_{\phi}^\top + \lambda I_N)^{-1}\bX_{\phi}-I)Uf^*}_\cH^2\\
    &=&\sum_{j\in\bN}(f_j^*)^2\norm{\Gamma^{1/2}(\bX_{\phi}^\top(\bX_{\phi}\bX_{\phi}^\top + \lambda I_N)^{-1}\bX_{\phi}-I)\varphi_j}_\cH^2 \\
    &=& \sum_{j\in\bN}(f_j^*)^2\left<\varphi_j, \Gamma^{1/2}(\bX_{\phi}^\top(\bX_{\phi}\bX_{\phi}^\top+ \lambda I_N)^{-1}\bX_{\phi}-I)\varphi_j\right>_\cH^2\\
    &=&\sum_{j\in\bN}(f_j^*)^2 \sigma_j\left( 1- \sigma_j\norm{ 
\left(\bX_{\phi}\bX_{\phi}^\top+ \lambda I_N\right)^{-1/2} z_j }_2^2 \right)^2.
\end{eqnarray*}Let $j\geq k_{b,\lambda}^*+1$. On $\Omega_j$, there exist absolute constants $\nC\label{C_lower_6}$ and $\nC\label{C_lower_7}$ that do not depend on $j$, such that
\begin{align*}
    \norm{ 
\left(\bX_{\phi}\bX_{\phi}^\top+ \lambda I_N\right)^{-1/2} z_j }_2 &\leq \norm{\left(\bX_{\phi}\bX_{\phi}^\top+ \lambda I_N\right)^{-1/2}}_{op}\norm{z_j}_2\\
&\leq \norm{\left(\bX_{\phi,k_{b,\lambda}^*+1:\infty}\bX_{\phi,k_{b,\lambda}^*+1:\infty}^\top+ \lambda I_N\right)^{-1/2}}_{op}\norm{z_j}_2\\
&\leq \frac{\oC{C_lower_6}\norm{z_j}_2}{\sqrt{\lambda + \Tr(\Gamma_{k_{b,\lambda}^*+1:\infty})}}\leq \frac{\oC{C_lower_7}\sqrt{N}}{\sqrt{\lambda + \Tr(\Gamma_{k_{b,\lambda}^*+1:\infty})}},
\end{align*}where we used that (by \eqref{eq:DM})
\begin{align*}
    &\bX_\phi\bX_\phi^\top +\lambda I_N = \bX_{\phi,1:k_{b,\lambda}^*}\bX_{\phi,1:k_{b,\lambda}^*}^\top + \bX_{\phi,k_{b,\lambda}^*+1:\infty}\bX_{\phi,k_{b,\lambda}^*+1:\infty}^\top+\lambda I_N\succeq\bX_{\phi,k_{b,\lambda}^*+1:\infty}\bX_{\phi,k_{b,\lambda}^*+1:\infty}^\top+\lambda I_N \\
    &\succeq \oc{c_lower_1}\left(\lambda + \Tr\left(\Gamma_{k_{b,\lambda}^*+1:\infty}\right)\right)I_N.
\end{align*}
Therefore, there exist absolute constants $\nc\label{c_lower_6}$ and $\nc\label{c_lower_7}$ such that if $b>2\oc{c_lower_7}^{-1}$ such that
\begin{align*}
    &\left( 1- \sigma_j\norm{ 
\left(\bX_{\phi}\bX_{\phi}^\top + \lambda I_N\right)^{-1/2} z_j }_2^2 \right)^2 \geq 1-2\sigma_j\norm{\left(\bX_{\phi}\bX_{\phi}^\top+ \lambda I_N\right)^{-1/2} z_j}_2^2\\
&\geq 1-\frac{\oc{c_lower_6}\sigma_{k_{b,\lambda}^*+1}N}{\lambda + \Tr(\Gamma_{k_{b,\lambda}^*+1:\infty})} > 1-\frac{\oc{c_lower_7}}{b}\geq \frac{1}{2}.
\end{align*}

Therefore, by Lemma~\ref{lemma:BLTT_lemma_9}, there exists an absolute constant $\oc{c_lower_4}$ such that
\begin{equation*}
    \bE_{\bX_{\phi}}\left[\sum_{j>k_{b,\lambda}^*}^\infty (f_j^*)^2 \sigma_j\left( 1- \sigma_j\norm{ 
\left(\bX_{\phi}\bX_{\phi}^\top + \lambda I_N\right)^{-1/2} z_j }_2^2 \right)^2 \right]\geq \oc{c_lower_4}\norm{\Gamma_{k_{b,\lambda}^*+1:\infty}^{1/2}f^*}_\cH^2.
\end{equation*}

Now, let $1\leq j\leq k_{b,\lambda}^*$, define $A = \sum_{j\in\bN}\sigma_j z_j\otimes z_j + \lambda I_N$, and $A_{-i} = A- \sigma_i z_i\otimes z_i$. By Sherman-Morrison-Woodbury formula
\begin{equation*}
    \sigma_j\left(1-\sigma_j\norm{\left(\bX_{\phi}\bX_{\phi}^\top + \lambda I_N\right)^{-1/2}z_j}_2^2 \right)^2 = \frac{\sigma_j}{\left(1 + \sigma_j z_j^\top A_{-j}^{-1}z_j \right)^2}.
\end{equation*}By Cauchy-Schwartz inequality, $z_j^\top A_{-j}^{-1}z_j\leq \norm{z_j}_2\norm{A_{-j}^{-1}z_j}\leq \norm{z_j}_2^2\norm{A_{-j}}_{op}=\norm{z_j}_2^2/s_N(A_{-j})$, where $s_N(A_{-j})$ is the $N$-th largest singular value of $A_{-j}$. Since $j\leq k_{b,\lambda}^*$, $A_{-j}\succeq \bX_{\phi,k_{b,\lambda}^*+1:\infty}\bX_{\phi,k_{b,\lambda}^*+1:\infty}^\top + \lambda I_N$. Condition on $\Omega_j$, there exists an absolute constant $\nC\label{C_lower_8}$ that is independent with $j$ such that
\begin{align*}
    z_j^\top A_{-j}^{-1}z_j\leq \frac{\norm{z_j}_2^2}{s_N(A_{-j})}\leq \frac{\oC{C_lower_8}N}{\lambda + \Tr\left(\Gamma_{k_{b,\lambda}^*+1:\infty}\right)}
\end{align*}
Further, since by definition of $k_{b,\lambda}^*$ and since $j\leq k_{b,\lambda}^*$, $bN\sigma_j > \Tr(\Gamma_{j:\infty}) + \lambda \geq \Tr(\Gamma_{k_{b,\lambda}^*+1:\infty}) + \lambda $ for $j\leq k_{b,\lambda}^*$. Therefore, there exist absolute constants $\nC\label{C_lower_9}$ and $\nC\label{C_lower_10}$ depend on $b$ such that
\begin{equation*}
    \frac{\sigma_j}{\left(1 + \sigma_j z_j^\top A_{-j}^{-1}z_j \right)^2}\geq \frac{\sigma_j}{ \left( \frac{\oC{C_lower_9}\sigma_j N}{\Tr\left(\Gamma_{k_{b,\lambda}^*+1:\infty}\right)} \right)^2} \geq \frac{1}{\oC{C_lower_10}}\frac{\left(\Tr\left(\Gamma_{k_{b,\lambda}^*+1:\infty}\right) + \lambda \right)^2}{\sigma_j N^2}.
\end{equation*} As we have $2\bP\left(\Omega_j\right) >99/100$ and by Lemma~\ref{lemma:BLTT_lemma_9}, there exists an absolute constant $\nc\label{c_lower_8}$ depending on $b$ such that with probability at least $9/10$,
\begin{align*}
    &\sum_{j>k_{b,\lambda}^*}^\infty (f_j^*)^2 \sigma_j\left( 1- \sigma_j\norm{ 
\left(\bX_{\phi}\bX_{\phi}^\top + \lambda I_N\right)^{-1/2} z_j }_\cH^2 \right)^2\\
&\geq \sum_{j>k_{b,\lambda}^*}^\infty \frac{1}{\oC{C_lower_10}}\frac{\left(\Tr\left(\Gamma_{k_{b,\lambda}^*+1:\infty}\right) + \lambda \right)^2}{N^2}\left(\frac{f_j^*}{\sigma_j}\right)^2\\
&\geq \oc{c_lower_8}\frac{\left(\Tr\left(\Gamma_{k_{b,\lambda}^*+1:\infty}\right)+\lambda\right)^2}{N^2}\sum_{j=1}^{k_{b,\lambda}^*}\left(\frac{f_j^*}{\sigma_j}\right)^2.
\end{align*}As a result,
\begin{align*}
    &\bE_{\bX_{\phi}}\left[\sum_{j>k_{b,\lambda}^*}^\infty (f_j^*)^2 \sigma_j\left( 1- \sigma_j\norm{ 
\left(\bX_{\phi}\bX_{\phi}^\top + \lambda I_N\right)^{-1/2} z_j }_\cH^2 \right)^2 \right]\\
&\geq \oc{c_lower_5}\norm{\Gamma_{1:k_{b,\lambda}^*}^{-1}f_{1:k_{b,\lambda}^*}}_\cH^2\frac{\Tr\left(\Gamma_{k_{b,\lambda}^*+1:\infty}\right)+\lambda}{N}.
\end{align*}

\endproof

Proposition~\ref{prop:gaussian_equivalent_conjecture} follows from Proposition~\ref{prop:lower_bias_term} and the lower bound of the variance term from \cite{tsigler_benign_2023}.

\endproof

\subsubsection{Proof of Proposition~\ref{prop:conjugate_kernel}}\label{sec:proof_conjugate_kernel}

In order to achieve a more accurate probability deviation, we will proceed to provide a new stochastic argument for the proof of Theorem~\ref{theo:upper_KRR} and Theorem~\ref{theo:main_upper_k>N}. Since this section focuses only on the estimation error of KRR defined by the data-dependent conjugate kernel, we replace $N_2$ with $N$.

\paragraph{Proof of diagonal concentration assumption}

The following lemma is taken from \cite[Corollary 3.5]{wang_deformed_2023}, with their $\vw$ replaced by our design vector $X$, their matrix $X$ replaced by $\bW^\top$ defined in the following lemma.
\begin{Lemma}[\cite{wang_deformed_2023}]\label{lemma:wang_deformed}
    Let $X\in\bR^d$ be a Gaussian random vector, $\bW\in\bR^{m\times d}$ and $A\in\bR^{m\times m}$ be deterministic matrices. Define $\vphi = \sigma(\bW X)$ a random vector in $\bR^m$. Suppose $\sigma$ is Lipschitz with Lipschitz constant $\lambda_\sigma>0$. There exists an absolute constant $\nC\label{C_wang_deformed}$ such that for any $t>0$,
    \begin{align*}
        \bP\left(\left| \vphi^\top A\vphi - \Tr\left(A\bE\left[\vphi\otimes\vphi\right]\right) \right|>t \right)&\leq 2\exp\left( -\frac{1}{\oC{C_wang_deformed}}\min\left\{ \frac{t^2  }{8\lambda_\sigma^4 \norm{\bW}_{\text{op}}^4 \norm{A}_{HS}^2 }, \frac{t}{\lambda_\sigma^2\norm{\bW}_{\text{op}}^2 \norm{A}_{\text{op}} } \right\} \right)\\
        &+ 2\exp\left( 
-\frac{t^2 }{32\lambda_\sigma^2\norm{\bW}_{\text{op}}^2 \norm{A}_{\text{op}}^2 t_0 } \right),
    \end{align*}where $t_0 = 2\lambda_\sigma^2\sum_{i=1}^m\left(\norm{W_i}_2 
 - 1\right)^2 + 2m\left|\bE\sigma(g)\right|^2$, $W_i$ is the $i$-th row of $\bW$ and $g$ is a standard Gaussian random variable.
\end{Lemma}

Let $\bW=W_{J^c}$, $A=\frac{1}{m}I_m$, and $t=\frac{1}{2}\bE\norm{\phi_{k+1:m}(X)}_2^2$ in Lemma~\ref{lemma:wang_deformed}. Take independent copies of $\vphi$ as $\vphi_1,\cdots,\vphi_N$. A union bound together with Lemma~\ref{lemma:wang_deformed} leads to \eqref{eq:diagonal_term_assumption} with $\delta=1/2$ and
\begin{align}\label{eq:def_gamma_conjugate}
\begin{split}
    \gamma &= 2N\exp\left(-\frac{1}{\oC{C_wang_deformed}}\min\left\{\frac{ m (\bE\norm{\phi_{k+1:\infty}(X)}_2^2)^2 }{32\lambda_\sigma^4 \norm{W_{J^c}}_{\text{op}}^4}, \frac{ m\bE\norm{\phi_{k+1:\infty}(X)}_2^2 }{ 2\lambda_\sigma^2 \norm{W_{J^c}}_{\text{op}}^2 } \right\}\right) \\
    &+2N\exp\left(-\frac{ m^2 (\bE\norm{\phi_{k+1:\infty}(X)}_2^2)^2 }{256\lambda_\sigma^2 \norm{W_{J^c}}_{\text{op}}^2 (\sum_{i\in J^c}\left(\norm{W_i}_2 
 - 1\right)^2 + 2(m-k)\left|\bE\sigma(g)\right|^2) }\right).
\end{split}
\end{align}
Let $\bW = W_{J^c}$, $A=\frac{1}{m}\Gamma_{k+1:\infty}$ where $\Gamma_{k+1:\infty}=\bE\left[\phi_{k+1:\infty}(X)\otimes \phi_{k+1:\infty}(X)\right]$ and $t=\frac{1}{2}\bE\norm{\Gamma_{k+1:m}^{1/2}\phi_{k+1:m}(X)}_2^2$ in Lemma~\ref{lemma:wang_deformed}. Take independent copies of $\vphi$ as $\vphi_1,\cdots,\vphi_N$. A union bound  together with Lemma~\ref{lemma:wang_deformed} implies \eqref{eq:diagonal_upper_dvoretzky} with $\odelta{delta_DMU_L2}=1/2$ and
\begin{align}\label{eq:def_gamma_DMU_conjugate}
    \begin{split}
        \ogamma{gamma_DMU_L2} &=  2N\exp\left(-\frac{1}{\oC{C_wang_deformed}}\min\left\{ \frac{m^2\left(\bE\norm{\Gamma_{k+1:m}^{1/2}\phi_{k+1:m}(X)}_2^2\right)^2}{32\lambda_\sigma^4\norm{W_{J^c}}_{\mathrm{op}}^4 \norm{\Gamma_{k+1:m}}_{HS}^2}, \frac{m \bE\norm{\Gamma_{k+1:m}^{1/2}\phi_{k+1:m}(X)}_2^2}{2\lambda_\sigma^2 \norm{W_{J^c}}_{\text{op}}^2 \norm{\Gamma_{k+1:m}}_{\text{op}} } \right\}\right)\\
    &+ 2N\exp\left(-\frac{m^2 \left(\bE\norm{\Gamma_{k+1:m}^{1/2}\phi_{k+1:m}(X)}_2^2\right)^2 }{ 256\lambda_\sigma^2 \norm{W_{J^c}}_{\text{op}}^2 (\sum_{i\in J^c}(\norm{W_j}_2-1)^2 + 2(m-k)\left|\bE\sigma(g)\right|^2) \norm{\Gamma_{k+1:m}}_{\text{op}}^2 }\right).
    \end{split}
\end{align}

\paragraph{Proof of \eqref{eq:DM_applied}} By \cite[Equation 5.2]{wang_deformed_2023}, $(\Gamma_{ij}) = \frac{1}{m}\bE[\sigma(\left<W_i,G\right>)\sigma\left<W_j,G\right>] =\frac{1}{m}\|\sigma\|_{L_2(\gamma_d)}^2\1_{\{i=j\}}$ for any $i,j\in J^c$, since $\left<W_i,W_j\right>=\delta_{ij}$ by Assumption~\ref{assumption:conjugate_kernel}, item \emph{3}. Then for any $\va\in\bR^m\cap\ell_2^{J^c}$ where $\ell_2^{J^c}=\Span(\ve_j:j\in J^c)$,
        \begin{align*}
            &\norm{\left<\va,\phi_{k+1:m}(G)\right>}_{L_2(\gamma_d)}^2 = \frac{1}{m}\sum_{i\in J^c} a_i^2\bE[\sigma^2(\left<W_i,G\right>)]\\
            &+ \frac{1}{m}\sum_{\substack{j\neq l\\ j,l\in J^c}}a_ja_l\bE[\sigma(\left<W_j,G\right>)\sigma(\left<W_l,G\right>)]= \frac{\|\va\|_2^2}{m}\|\sigma\|_{L_2(\gamma_d)}^2.
        \end{align*}On the other hand, by Borel-TIS inequality, $\|\left<\va,\phi_{k+1:m}(G)\right>\|_{\psi_2}^2\lesssim \lambda_\sigma^2\|\va\|_2^2\|W_{J^c}\|_{\text{op}}^2/m$. Hence
        \begin{align*}
            &\|\left<\va,\phi_{k+1:m}(G)\right>\|_{\psi_2}^2\lesssim \frac{\|\sigma\|_{L_2(\gamma_d)}^2\|\va\|_2^2}{m}\frac{\norm{W_{J^c}}_{\text{op}}^2\lambda_\sigma^2}{\|\sigma\|_{L_2(\gamma_d)}^2}\\
            &= \frac{\lambda_\sigma^2}{\|\sigma\|_{L_2(\gamma_d)}^2}\|W_{J^c}\|_{\text{op}}^2\norm{\left<\va,\phi_{k+1:m}(G)\right>}_{L_2(\gamma_d)}^2.
        \end{align*}Since $\|W_{J^c}\|_{\text{op}}=1$, we omit it in the following.
        \begin{enumerate}
            \item We first deal with the case when $\Tr(\Gamma_{k+1:\infty})$ dominates. By \cite[Equation 2.14]{vershynin_high-dimensional_2018}, there exists an absolute constant $\nc\label{c_P_sub_gaussian}$ such that for any $t>0$, $\bP\left(\left|f_{\vw}(G)\right|\geq t\right)\leq 2\exp\left(-\oc{c_P_sub_gaussian}t^2/\norm{f_{\vw}}_{\psi_2}^2\right)$. As a result, for any $\vw\in S_2^{m-1}\cap \ell_2^{J^c}$ and any $t>0$, we have $\bP\left(\sqrt{\oc{c_P_sub_gaussian}}\|\sigma\|_{L_2(\gamma_d)}\left|f_\vw(G)\right|/(\lambda_\sigma \sqrt{\|\Gamma_{k+1:m}\|_{\text{op}}})\geq t\right)\leq 2\exp\left(- t^2\right)$. We have so far verified all the conditions of \cite[Theorem 4]{tsigler_benign_2023}. Therefore we use \cite[Theorem 4]{tsigler_benign_2023} but with the use of \cite[Theorem 2.1, case 1]{guedon_interval_2017} replaced by \cite[Theorem 2.1, case 2]{guedon_interval_2017}. Moreover, because we have assumed that
            \begin{align*}
                \frac{\Tr(\Gamma_{k+1:m})+\lambda}{\norm{\Gamma_{k+1:m}}_{\text{op}}}\gtrsim \frac{\lambda_\sigma^2}{\norm{\sigma}_{L_2(\gamma_d)}^2}N,
            \end{align*}which is Assumption~\ref{assumption:conjugate_kernel}, \emph{4.(b)}, we thus have that \eqref{eq:DM_applied} holds with probability at least $1-\gamma-\sqrt{\beta}=:1-\bar p_{DM}$, where $\gamma$ is defined in \eqref{eq:def_gamma_conjugate} and $\beta$ is defined as
            \begin{align}\label{eq:def_beta_guedon_et_al}
                \beta = \frac{1}{\left(10N\right)^2}\exp\left(-\frac{4N}{\left(3.5\log(2N)\right)^4}\right) + \frac{N^2}{2\exp\left(4N\right)}.
            \end{align}
            \item When $\lambda$ dominates, the proof idea is the same as in Section~\ref{sec:proof_DM_lambda_dominating}. The lower side is trivial, so we omit it. The upper side follows by applying the upper bound of $A_k$ from \cite[Theorem 2.1, case 2]{guedon_interval_2017}, with their $k=N$, $t=\sqrt{N}$, $\tau=1$, $\alpha=\lambda=2$, and their $X_i$ replaced by our $\sqrt{\oc{c_P_sub_gaussian}}\|\sigma\|_{L_2(\gamma_d)}/(\lambda_\sigma \sqrt{\|\Gamma_{k+1:m}\|_{\text{op}}})\phi_{k+1:\infty}(X_i)$. We obtain that with probability at least $1-\bar p_{DM}$ (with $\bar p_{DM}$ defined above), $\norm{\bX_{\phi,k+1:\infty}^\top}_{\text{op}}\lesssim \sqrt{\Tr(\Gamma_{k+1:m})}+\lambda_\sigma \sqrt{N\|\Gamma_{k+1:m}\|_{\text{op}}}/(\sqrt{\oc{c_P_sub_gaussian}}\|\sigma\|_{L_2(\gamma_d)})$. Due to Assumption~\ref{assumption:conjugate_kernel}, \emph{4.(b)}, we have $\norm{\bX_{\phi,k+1:\infty}^\top}_{\text{op}}\lesssim \sqrt{\Tr(\Gamma_{k+1:m})+\lambda}$. \eqref{eq:DM_applied} then follows.
        \end{enumerate}



\paragraph{Proof of \eqref{eq:RIP_applied}}
The following Lemma is taken from \cite[Theorem 5.39]{vershynin_introduction_2011}:
\begin{Lemma}[\cite{vershynin_introduction_2011}]\label{lemma:covariance_vershynin}
    Assume $A\in\bR^{N\times n}$ is a matrix whose rows $A_i$ are i.i.d. sub-Gaussian isotropic random vectors in $\bR^n$. There exist absolute constants $\nc\label{c_P_vershynin}$ and $\nC\label{C_vershynin}\sim \norm{A_i}_{\psi_2}^2$ such that for every $t\geq 0$, the following inequality holds with probability at least $1-2\exp(-\oc{c_P_vershynin}t^2\norm{A_i}_{\psi_2}^{-4})$:
    \begin{align*}
        \sqrt{N}-\oC{C_vershynin}\sqrt{n}-t\leq\sigma_{\text{min}}(A)\leq\sigma_1(A)\leq \sqrt{N}+\oC{C_vershynin}\sqrt{n}+t.
    \end{align*}
\end{Lemma}
Recall that we have assumed $\norm{W_{J}}_{\text{op}}\leq B$, $\sum_{j\in J}\left(\norm{W_j}_2 ^2 - 1\right)^2\leq B^2$ and $\norm{\sigma'}_{L_\infty}\leq\lambda_\sigma$. As in \cite[Proof of Lemma D.4]{fan_spectra_2020}, we know that $\norm{\sigma(W_J G)}_{\psi_2}\lesssim \lambda_\sigma B$. Hence $\norm{\Gamma_{1:k}^{-1/2}\phi_{1:k}(X)}_{\psi_2}\lesssim \lambda_\sigma B\norm{\Gamma_{1:k}^{-1/2}}_{\text{op}}/\sqrt{m}$ thus $\oC{C_vershynin}\sim \frac{\lambda_\sigma^2 B^2}{m\sigma_k }$.
Let $A$ in Lemma~\ref{lemma:covariance_vershynin} be $\bX_{\phi,1:k}\Gamma_{1:k}^{-1/2}$, $n=k$ with $\oc{c_kappa_RIP}\leq \frac{1}{4}\min\{1,\oC{C_vershynin}^{-2}\}$, and $t=\frac{\oC{C_vershynin}}{2}\sqrt{k}$. Then $t^2\norm{A_i}_{\psi_2}^{-4} \sim k$, and $\sqrt{N}-\oC{C_vershynin}\sqrt{n}-t = \sqrt{N} - \frac{1}{2}\oC{C_vershynin}\sqrt{k}\geq \sqrt{N}-\frac{1}{2}\oC{C_vershynin}\sqrt{\frac{1}{4}\oC{C_vershynin}^{-2}N}=\frac{1}{4}\sqrt{N}$.
Then by Lemma~\ref{lemma:covariance_vershynin}, with probability at least $1-2\exp(-\oc{c_P_vershynin} k)=:1-\bar p_{RIP}$. By the assumption that $2\lambda_\sigma^2 B^2 \leq m\sigma_k\sqrt{N}$ (this is Assumption~\ref{assumption:conjugate_kernel}, 4.(c)), we have $k\geq 1$, and thus
\begin{align*}
    \frac{1}{4}\sqrt{N}\leq \sigma_{\text{min}}(\bX_{\phi,1:k}\Gamma_{1:k}^{-1/2})\leq \sigma_1(\bX_{\phi,1:k}\Gamma_{1:k}^{-1/2})\leq \frac{7}{4}\sqrt{N}.
\end{align*}By the homogeneity argument as in Section~\ref{sec:proof_RIP} together with the isometry of $\cH_{1:k}\hookrightarrow\ell_2^k$, there exists some $\ndelta\label{delta_approx_isometry}=\frac{3}{4}$ such that
\begin{align*}
    (1-\odelta{delta_approx_isometry})\norm{\Gamma_{1:k}^{1/2}f}_\cH\leq \frac{1}{\sqrt{N}}\norm{\bX_{\phi,1:k}f}_2\leq (1+\odelta{delta_approx_isometry})\norm{\Gamma_{1:k}^{1/2}f}_\cH.
\end{align*}Therefore, we can take $\bar p_{RIP}=2\exp(-\oc{c_P_vershynin}k)$, $\oc{c_RIP_lower}=1-\odelta{delta_approx_isometry}$ and $\oC{C_RIP_upper}=1+\odelta{delta_approx_isometry}$.

\paragraph{Proof of \eqref{eq:DM_upper_applied}} A routine modification of \cite[Lemma D.4]{fan_spectra_2020} shows that for any $\vw\in\Gamma_{k+1:m}^{1/2}\vu$ with $\vu\in S_2^{m-1}\cap \ell_2^{J^c}$, we have $\norm{f_\vw}_{\psi_2}\lesssim\frac{1}{\sqrt{m}}\lambda_\sigma  \sqrt{\sigma_{k+1}(\Gamma)}$. Therefore, for any $\vw\in S_2^{m-1}\cap \ell_2^{J^c}$ and any $t>0$, we have
\begin{align*}
    \bP\left( 
\frac{\sqrt{\oc{c_P_sub_gaussian}m\Tr\left(\Gamma_{k+1:m}^2\right)}}{\lambda_\sigma  \sqrt{\sigma_{k+1}(\Gamma)}}\frac{\left|\left\langle \vw,\Gamma_{k+1:m}^{1/2}\phi_{k+1:m}(X)\right\rangle\right|}{\sqrt{\Tr\left(\Gamma_{k+1:m}^2\right)}} \geq t\right)\leq 2\exp(-t^2).
\end{align*}
Therefore, by \cite[Theorem 4]{tsigler_benign_2023} but with \cite[Theorem 2.1, case 1]{guedon_interval_2017} replaced by \cite[Theorem 2.1, case 2]{guedon_interval_2017} again, because  $m\Tr\left(\Gamma_{k+1:m}^2\right)\gtrsim \lambda_\sigma^2  \norm{\Gamma_{k+1:m}}_{\text{op}}$ (this is Assumption~\ref{assumption:conjugate_kernel}, 4.(a)), \eqref{eq:DM_upper_applied} holds with probability at least $1-\ogamma{gamma_DMU_L2}-\sqrt{\beta}=:1-\bar p_{DMU}$, where $\ogamma{gamma_DMU_L2}$ is defined in \eqref{eq:def_gamma_DMU_conjugate} and $\beta$ in \eqref{eq:def_beta_guedon_et_al}.

\paragraph{Proof of \eqref{eq:bX_f_star_applied}} As in the proof of \eqref{eq:DM_applied}, $f_{k+1:m}^*=f_{\vw^*}$ for some $\vw^*\in \ell_2^{J^c}$, and $\norm{f_{\vw^*}}_{\psi_2}\lesssim \frac{1}{\sqrt{m}}\lambda_\sigma \norm{f_{k+1:m}^*}_\cH$ due to the isometry between $\cH_{k+1:m}$ and $\ell_2^{J^c}$. By Bernstein's inequality for $\psi_1$ variables, there exists an absolute constant $\nc\label{c_P_Bernstein}$ such that with probability at least
\begin{align*}
    1-\exp\left(-\oc{c_P_Bernstein}N\min\left\{\left(\frac{m\norm{f_{k+1:m}^*}_{L_2}^2}{\lambda_\sigma^2 \norm{f_{k+1:m}^*}_\cH^2}\right)^2, \frac{m\norm{f_{k+1:m}^*}_{L_2}^2}{\lambda_\sigma^2 \norm{f_{k+1:m}^*}_\cH^2} \right\}\right),
\end{align*}\eqref{eq:bX_f_star_applied} holds with $\oC{C_bX_f_star}\kappa=\sqrt{2}$.



\paragraph{Proof of \eqref{eq:Gamma_phi_applied}}

\eqref{eq:Gamma_phi_applied} holds with probability $1-\ogamma{gamma_DMU_L2}$.

So far, we have successfully reestablished the proofs for Proposition~\ref{prop:stochastic_argument_upper_bound}. Since the deterministic argument remains unchanged, we can now conclude our proof.

\endproof

\begin{Remark}\label{remark:general_distribution}
    Our Assumption~\ref{assumption:conjugate_kernel} [1] requires the design vector $X$ to be a Gaussian random vector (when $X$ is anisotropic, replace the hidden layer matrix $W$ with $W\Sigma^{1/2}$, where $\Sigma=\bE[XX^\top]$). When $X=(x_1,\cdots,x_d)$ has independent (but not necessarily identically distributed) coordinates, we need to generalize Lemma~\ref{lemma:wang_deformed} and \cite[Lemma D.4]{fan_spectra_2020}. The proof of Lemma~\ref{lemma:wang_deformed} only utilizes the convex concentration property of $X$, as defined in \cite[Definition 2.2]{adamczak_note_2015}. This property has been proved in \cite{talagrand_new_1996,klochkov_uniform_2020,huang_dimension-dependent_2023,sambale_notes_2023,adamczak_orlicz_2023} to the case where $x_1,\cdots,x_d$ are Orlicz random variables satisfying the Hoffmann-Jorgensen condition, as stated in \cite[Proposition 4.7]{adamczak_orlicz_2023}. Therefore, Lemma~\ref{lemma:wang_deformed} still holds for a broad class of probability measures, but with a logarithmic factor in $m$ as a cost, which is unavoidable (see \cite{huang_dimension-dependent_2023}). 
    The proof of Proposition~\ref{prop:conjugate_kernel} also relies on showing that $\sigma(W_J X)$ and $\sigma(W_{J^c}X)$ are sub-Gaussian random vectors (we conjecture that sub-Weibull random vectors would be sufficient, see \cite[Section 4.2]{kuchibhotla_moving_2022}). We utilize \cite[Lemma D.4]{fan_spectra_2020} to establish this result, which in turn relies on Gaussian Poincaré inequality and Gaussian Lipschitz concentration inequality.
    When $X=(x_1,\cdots,x_d)$ satisfies that $x_1,\cdots,x_d$ are independent, and for any $i\in[d]$, $\|x_i\|_{L_\infty}<\infty$, then by the bounded difference inequality, see for example \cite[Theorem 6.2]{boucheron_concentration_2013}, for any $\va\in S_2^{m-1}$, $\|\left<\va,\sigma(W X)\right> - \bE \left<\va,\sigma(W X)\right>\|_{\psi_2}^2\leq \|\sigma\|_{Lip}^2\|W\|_{\text{op}}^2\sum_{i=1}^d \|x_i\|_{L_\infty}^2$. Furthermore, by the Efron-Stein inequality, see \cite[Corollary 3.2]{boucheron_concentration_2013}, combined with \cite[Equation 34]{fan_spectra_2020}, we obtain $\norm{\left<\va,\sigma(W X)\right>}_{\psi_2}^2\lesssim \|\sigma\|_{Lip}^2\|W\|_{\text{op}}^2\sum_{i=1}^d \|x_i\|_{L_\infty}^2$. Therefore, when $X=(x_1,\cdots,x_d)$ satisfies that all coordinates have $\|x_i\|_{L_\infty}<\infty$, our results still hold, but with some additional multiplicative factor $\sum_{i=1}^d \|x_i\|_{L_\infty}^2$, respectively, and in $\bar p_{DM}$, $\bar p_{DMU}$, $\bar p_{RIP}$ incurring some logarithmic factors. However, for a general design vector, we may need to examine more carefully whether the $L_2$ norm is equivalent to the $\psi_2$ norm. If this equivalence does not hold, then we need to rely on the boundedness of the $\|\cdot\|_{\psi_2}$ norm rather than the $L_p$-$L_2$ norm equivalence to establish \eqref{eq:DM_applied}.
\end{Remark}

\subsubsection{Examples when the spectrum and eigenvectors of the integral operator of the data-dependent conjugate kernel}\label{sec:example_Gamma}

We provide cases where the spectrum of $\Gamma$ can be conveniently estimated. Note that when the spectrum of $\Gamma$ is difficult to estimate, we can still use $\Tr(\Gamma) = \sum_{j=1}^m (\Gamma)_{j,j}$ as an upper bound, and by choosing an appropriate $\lambda$, we can compute $r_{\lambda,k}^*$, following the approach of \cite{bietti_learning_2022}.
\begin{enumerate}
    \item When $W = [W(T)|\cdots|W(T)]^\top\in\bR^{m\times d}$ for some $W(T)\in\bR^d$, and $\sigma$ satisfies item \emph{2.} of Assumption~\ref{assumption:conjugate_kernel}, by \cite[Equation 5.2]{wang_deformed_2023}, $\Gamma = (\Gamma_{ij})_{i,j\in[m]} = \frac{1}{m}\|\sigma(\|W(T)\|_2\cdot)\|_{L_2(\gamma_d)}^2$, where $\|\sigma(\|W(T)\|_2\cdot)\|_{L_2(\gamma_d)}$ is the $L_2$ norm (under standard Gaussian probability measure $\gamma_d$) of function $x\in\bR\mapsto \sigma(\|W(T)\|_2 x)$. In particular, when $\|W(T)\|_2=1$, $\Gamma_{ij}=\frac{1}{m}\|\sigma\|_{L_2(\gamma_d)}^2$ for all $i,j\in[m]$. In this case, $\Gamma$ is of rank $1$, and the eigenvector associated to the only non-zero eigenvalue $\sigma_1(\Gamma) = \|\sigma\|_{L_2(\gamma_d)}^2$, is $\frac{1}{\sqrt{m}}(1,1,\cdots,1)\in\bR^m$.

    \item When $W = [W_1|\cdots|W_m]^\top\in\bR^{m\times d}$ with $\left<W_i,W_j\right>=\delta_{i,j}$ for any $i,j\in[m]$. When $\sigma$ satisfies item \emph{2.} of Assumption~\ref{assumption:conjugate_kernel}, by \cite[Equation 5.2]{wang_deformed_2023}, $\Gamma = (\Gamma_{i,j})_{i,j\in[m]} = \frac{1}{m}\|\sigma\|_{L_2(\gamma_d)}^2\1_{\{i=j\} } = \frac{\|\sigma\|_{L_2(\gamma_d)}^2}{m}I_m$. The spectrum of $\Gamma = \{\frac{\|\sigma\|_{L_2(\gamma_d)}^2}{m},\cdots,\frac{\|\sigma\|_{L_2(\gamma_d)}^2}{m}\}$ and the eigenvectors alone are $\ve_1,\cdots,\ve_m$.
    
    \item When there exists $J\subset[m]$ such that for any $i,j\in J$, $W_i=W_j=\vw$ for some $\vw\in\bR^d$; for any $i,j\in J^c$, $\left<W_i,W_j\right>=\eps_i\delta_{ij}$ for some $\eps_i>0$ and for any $i\in J^c, j\in J$, $\left<W_i,W_j\right>=0$. This is a simplified case of the training dynamics studied in \cite{boursier_early_2024}. This situation implies that there exists a subset of neurons, that is, neurons $W_j$ corresponding to $j \in J$, that are well-aligned; the remaining neurons are pairwise orthogonal and orthogonal to the direction of the neurons in $J$. For the case of almost orthogonal neurons, we conjecture that the method of \cite[Lemma 5.2]{wang_deformed_2023} together with certain matrix perturbation technique could be applied, but this is beyond the scope of this paper. For the sake of simplicity, we assume that $J=[k]$ for some $k\in[m]$. Then by \cite[Equation 5.2]{wang_deformed_2023}, $\Gamma$ is a block diagonal matrix with diagonal matrices $\frac{1}{m}\|\sigma(\|W(T)\|_2\cdot)\|_{L_2(\gamma_d)}^2\1_{k\times k}$, and $\diag(\frac{1}{m}\|\sigma(\eps_i\cdot)\|_{L_2(\gamma_d)}^2: i\in J^c)$, where $\1_{k\times k}$ is the $k\times k$ matrix with all $1$'s. In this case, the spectrum of $\Gamma$ is given by
    \begin{align*}
        \bigg\{\frac{k}{m}\|\sigma(\|W(T)\|_2\cdot)\|_{L_2(\gamma_d)}^2,\frac{1}{m}\|\sigma(\eps_i\cdot)\|_{L_2(\gamma_d)}^2,\cdots,\frac{1}{m}\|\sigma(\eps_i\cdot)\|_{L_2(\gamma_d)}^2,0,0,\cdots,0\bigg\}.
    \end{align*}Moreover, the eigenvector associated with eigenvalue $\frac{k}{m}\|\sigma(\|W(T)\|_2\cdot)\|_{L_2(\gamma_d)}^2$ is $(1/\sqrt{k})(\1_k|\vzero)^\top$ where $\1_k$ is the vector on $1$'s, and other eigenvectors are $(\vzero|\ve_i)^\top$.

\end{enumerate}

\subsubsection{Proof of Proposition~\ref{prop:feature_learning}}\label{sec:proof_feature_learning}
\paragraph{Spectrum and eigenfunctions of $\Gamma(T)$} In this paragraph, we investigate the spectrum and eigenfunctions of $\Gamma(T)$. We omit $T$ from $\Gamma(T)$ in the following. Since $X\sim\cN(\vzero,I_d)$, by \cite[Proof of Lemma 5.2]{wang_deformed_2023}, $\Gamma$ is a diagonal matrix under the basis $\{\ve_1,\cdots,\ve_m\}$, that is, $\Gamma = (\Gamma_{j,j})_{j\in[m]} = \frac{1}{m}\bE[\sigma^2(\|W_j\|_2g)]$, where $g\sim\cN(0,1)$. Since we have that $\delta\in(0,1)\mapsto \bE[\sigma^2(\delta g)]$ is increasing and since we have $\|W_1\|_2=\|\vw_2\|_2=\cdots=\|\vw_m\|_2=1$ by construction, we have: $\bE[\sigma^2(g)]=\bE[\sigma^2(\|W_1\|_2g)]=\sigma_1(\Gamma)> \sigma_2(\Gamma)=\cdots=\sigma_m(\Gamma) = \frac{1}{m}\bE[\sigma^2(\|W_2\|_2g)] = \frac{1}{m}\bE[\sigma^2(\odelta{delta_feature_learning}\|\vw_2\|_2g)]$. Therefore, if we take $k=1$, we have:
\begin{align*}
    &\Tr(\Gamma_{k+1:m}) = \frac{m-1}{m}\bE[\sigma^2(\odelta{delta_feature_learning}g)],\, \Tr(\Gamma_{k+1:m}^2) = \frac{m-1}{m^2}(\bE[\sigma^2(\odelta{delta_feature_learning}g)])^2,\\
    &\sigma_k = \frac{1}{m}\|\sigma\|_{L_2(\gamma_d)}^2,\mbox{ and } \norm{\Gamma_{k+1:m}}_{\text{op}} = \frac{1}{m}\bE[\sigma^2(\odelta{delta_feature_learning}g)].
\end{align*}
Since $\Gamma$ is diagonal along the basis $\{\ve_1,\cdots,\ve_m\}$, its eigenvectors are given by $\varphi_1 = \ve_1, \varphi_2 = \ve_2,\cdots, \varphi_m = \ve_m$.

\begin{Remark}
    Here we consider a simplified model of the training dynamics of gradient descent. The activation functions we consider are bounded, such as $\tanh$ and sigmoid. In the case of an unbounded activation function like ReLU, we can consider $\|W_1\|_2$ to be very large, while $\|W_2\|_2, \cdots, \|W_m\|_2$ are of constant magnitude.
\end{Remark}

\paragraph{Structure of $W$} As we choose $k=1$ and $J=\{1\}$, we have $\|W_J\|_{\text{op}}\leq1$, $\sum_{j\in J}(\|W_j\|_2^2-1)^2\leq (1-\odelta{delta_feature_learning})^2$ and hence we can take $B^2=1$.

\paragraph{Checking the conditions}In the following we check Assumption~\ref{assumption:conjugate_kernel}. Item \emph{1.} in Assumption~\ref{assumption:conjugate_kernel} is valid; Item \emph{2.} in Assumption~\ref{assumption:conjugate_kernel} is also valid, since in Assumption~\ref{assumption:Bietti}, we already assume that $\sigma$ is Lipschitz, and from previous analysis, we can take $D=1$. Now we check Item \emph{3.}. We have
\begin{enumerate}
    \item $\frac{m\Tr(\Gamma_{k+1:m}^2)}{\lambda_\sigma^2  \|\Gamma_{k+1:m}\|_{\text{op}}} = (m-1)\frac{\|\sigma(\odelta{delta_feature_learning}\cdot)\|_{L_2(\gamma_d)}^4}{\lambda_\sigma^2\|\sigma(\odelta{delta_feature_learning}\cdot)\|_{L_2(\gamma_d)}^2}\gtrsim 1$ if $m\gtrsim \|\sigma(\odelta{delta_feature_learning})\|_{L_2(\gamma_d)}^{-2}$;
    \item $\frac{(\Tr\left(\Gamma_{k+1:m}\right)+\lambda)\norm{\sigma}_{L_2(\gamma_d)}^2}{\lambda_\sigma^2N_2\|\Gamma_{k+1:m}\|_{\text{op}}} \gtrsim \frac{\|\sigma\|_{L_2(\gamma_d)}^2}{\lambda_\sigma^2}(\frac{m-1}{N_2}+\frac{m\lambda}{N_2\|\sigma(\odelta{delta_feature_learning}\cdot)\|_{L_2(\gamma_d)}^2})\gtrsim 1$ if either $m\gtrsim N_2$, or $\frac{\lambda}{N_2}\gtrsim \frac{\|\sigma(\odelta{delta_feature_learning}\cdot)\|_{L_2(\gamma_d)}^2}{m}$;
    \item $\frac{m\sigma_k(\Gamma)\sqrt{N_2}}{2\lambda_\sigma^2 B^2} \gtrsim \frac{m\|\sigma\|_{L_2(\gamma)}^2\sqrt{N_2} }{\lambda_\sigma^2m}\gtrsim 1$.
\end{enumerate}

\paragraph{Approximation error} 
\begin{Lemma}[Theorem 6.1 of \cite{bietti_learning_2022}]\label{lemma:bietti}
Grant Assumption~\ref{assumption:Bietti}, we then have: for any $0<\delta<1/4$, with probability at least $\frac{1}{2}-\delta$,
\begin{align*}
    1-\left|\left<\vw(T),\vw^*\right>\right| = \tilde O\left(\eta^{-4}\max\left\{\frac{d+m}{N_1}, \frac{d^4}{N_1^2}\right\}\right).
\end{align*}
\end{Lemma}
Let $f^{**}:\vx\in\bR^d\mapsto a^*\sigma(\left<W(T),\vx\right>) = \left<\phi(\vx),\va^*\right>$ where $\va^*= \sqrt{m}(a^*,0,0,\cdots,0)\in\bR^m$. By the definition of $f^{**}$ and $f^*$ in Assumption~\ref{assumption:Bietti}, $(f^*-f^{**})(\vx)=a^*(\sigma(\left<W(T),\vx\right>) - \sigma( \left<\vw^*,\vx\right> ) )$. Because $\sigma$ is Lipschitz, we have
\begin{align*}
    \norm{f^*-f^{**}}_{L_4}^4 &= (a^*)^4\bE\left[\left|\sigma\left(\left<W(T),G\right>\right) - \sigma\left(\left<\vw^*,G\right>\right)  \right|^4\right]\\
    &\leq (a^*)^4\norm{\sigma}_{Lip}^4\bE\left[\left|\left(\left<W(T),G\right>\right) - \left(\left<\vw^*,G\right>\right)  \right|^4\right]\\
    &=(a^*)^4\norm{\sigma}_{Lip}^4\bE\left[\left|\left<W(T)-\vw^*,G\right>\right|^4\right] = 3(a^*)^4\norm{\sigma}_{Lip}^4\norm{W(T)-\vw^*}_2^4.
\end{align*}Because $\Pi_{\vw^\perp}$ is the retraction map at $\vw^\perp$ of the Euclidean sphere (see, for example, \cite[Equation 3.40]{boumal_introduction_2023}), we have $W(T)\in S_2^{d-1}$. Hence by Lemma~\ref{lemma:bietti}, $\norm{W(T)-\vw^*}_2^2 = 2( 1 - \left|\left<\vw(T),\vw^*\right>\right|) = \tilde O\left(\eta^{-4}\max\left\{\frac{d+m}{N_1}, \frac{d^4}{N_1^2}\right\}\right)$ with probability at least $1/2-\delta$. As a result, with probability at least $1/2-\delta$,
\begin{align}\label{eq:model_noise_feature_learning}
    \norm{f^*-f^{**}}_{L_4}\lesssim \left|a^*\right|\norm{\sigma}_{Lip}\tilde O\left(\eta^{-2}\max\left\{\sqrt{\frac{d+m}{N_1}},\frac{d^2}{N_1}\right\}\right).
\end{align}

\paragraph{Estimation error}In the following, we compute $r_{\lambda,k}^*$. By our construction of $\va^*$ and $\Gamma$, we know that $\va^*\in\Span(\varphi_1)$ and hence $\norm{\Gamma_{1,\mathrm{thre}}^{-1/2}\va_{1:k}^* }_\cH^2=\frac{1}{\sigma_1(\Gamma)}\left<\va^*,\ve_1\right>^2 = \frac{m^2(a^*)^2}{\|\sigma\|_{L_2(\gamma_d)}^2}$, and $\norm{\Gamma_{k+1:m}^{1/2}\va_{k+1:m}^*}_\cH^2=0$. Moreover,
\begin{align*}
    \frac{\lambda+\Tr(\Gamma_{k+1:m})}{N_2} = \frac{\lambda+\frac{m-1}{m}\|\sigma(\odelta{delta_feature_learning}\cdot)\|_{L_2(\gamma_d)}^2}{N_2}\lesssim \sigma_1(\Gamma) = \frac{\|\sigma\|_{L_2(\gamma_d)}^2}{m},
\end{align*}if $\lambda\vee \|\sigma(\odelta{delta_feature_learning}\cdot)\|_{L_2(\gamma_d)}^2\lesssim N_2/m$. At this time, $\left|J_1\right|=1$ and
\begin{align*}
    \sigma_\xi^2\frac{N_2\Tr(\Gamma_{k+1:m}^2)}{\left(\lambda+\Tr(\Gamma_{k+1:m})\right)^2}= \sigma_\xi^2\frac{N_2(m-1)\|\sigma(\odelta{delta_feature_learning}\cdot)\|_{L_2(\gamma_d)}^4}{m^2\left(\lambda + \frac{m-1}{m}\|\sigma(\odelta{delta_feature_learning}\cdot)\|_{L_2(\gamma_d)}^2\right)^2}\lesssim \sigma_\xi^2\frac{N_2\|\sigma(\odelta{delta_feature_learning}\cdot)\|_{L_2(\gamma_d)}^4}{m\lambda} \wedge \sigma_\xi^2\frac{N_2}{m-1}.
\end{align*}
Hence if
\begin{enumerate}
    \item $\lambda\vee \|\sigma(\odelta{delta_feature_learning}\cdot)\|_{L_2(\gamma_d)}^2\lesssim N_2/m$,
    \item $m\gtrsim \|\sigma(\odelta{delta_feature_learning}\cdot)\|_{L_2(\gamma_d)}^{-2}$ and 
    \item $\frac{\lambda}{N_2}\gtrsim\frac{\|\sigma(\odelta{delta_feature_learning}\cdot)\|_{L_2(\gamma_d)}^2}{m}$ or $m\gtrsim N_2$,
\end{enumerate}we have
\begin{align*}
    (r_{\lambda,k}^*)^2\lesssim&\frac{\sigma_\xi^2}{N_2} + \sigma_\xi^2\frac{N_2\|\sigma(\odelta{delta_feature_learning}\cdot)\|_{L_2(\gamma_d)}^4}{m\lambda} \wedge \sigma_\xi^2\frac{N_2}{m-1}\\
    & + \frac{m^2(a^*)^2}{\|\sigma\|_{L_2(\gamma_d)}^2}\frac{\lambda^2}{N_2^2} +\frac{m^2(a^*)^2}{\|\sigma\|_{L_2(\gamma_d)}^2}\frac{ \|\sigma(\odelta{delta_feature_learning}\cdot)\|_{L_2(\gamma_d)}^4}{N_2^2}.
\end{align*}
We check that the choice for $\lambda = 0$, $\odelta{delta_feature_learning}$ such that $\|\sigma(\odelta{delta_feature_learning}\cdot)\|_{L_2(\gamma_d)}^{-2}\sim m$, and $m\gtrsim N_2$ satisfy the above constraints and in this case, $(r_{0,1}^*)^2\lesssim \frac{\sigma_\xi^2}{N_2} + \sigma_\xi^2\frac{N_2}{m} + \frac{(a^*)^2}{\|\sigma\|_{L_2(\gamma_d)}^2}\frac{1}{N_2^2}$.

\paragraph{Probability deviation}
\begin{enumerate}
    \item By \eqref{eq:def_gamma_conjugate}, the fact that $\|W_{J^c}\|_{\text{op}}^2 = \|\sigma(\odelta{delta_feature_learning}\cdot)\|_{L_2(\gamma_d)}^2$, and $\bE\|\phi_{k+1:m}(X)\|_\cH^2\sim \|\sigma(\odelta{delta_feature_learning}\cdot)\|_{L_2(\gamma_d)}^2$, we know that $\gamma<1/100$ provided that $m\gtrsim N_2$.
    
    \item We know that $m^2\Tr(\Gamma_{k+1:m}^2) \sim m\|\sigma(\odelta{delta_feature_learning}\cdot)\|_{L_2(\gamma_d)}^4$ and $\|W_{J^c}\|_{\text{op}}^2 = \|\sigma(\odelta{delta_feature_learning}\cdot)\|_{L_2(\gamma_d)}^2$, hence $\ogamma{gamma_DMU_L2}<1/100$ provided that that $m\gtrsim N_2$ by \eqref{eq:def_gamma_DMU_conjugate}.
    
    \item By \eqref{eq:def_beta_guedon_et_al}, we know that $\bar p_{DM}$, $\bar p_{DMU}<1/90$ provided that $m\gtrsim N_2$.
    
    \item For $\bar p_{RIP}$, since $k=1$ is too small in this case, we can set a relatively large $t$ by adjusting constant $\oC{C_vershynin}$ such that $\bar p_{RIP}<1/90$, see the discussion after Lemma~\ref{lemma:covariance_vershynin}.

    \item Similarly, for $\bar p_\xi$, one can modify the constants in Proposition~\ref{prop:noise_concentration_dependent} to satisfy $\bar p_\xi<1/90$. Here, one can use Borel-TIS in \eqref{eq:noise_concentration_dependent_2} to get a smaller probability deviation since the noise is assumed to be Gaussian.
\end{enumerate}

\subsubsection{Concentration of noise}

The following lemma is taken from \cite[Lemma 10]{lecue_geometrical_2022}.

\begin{Lemma}\label{lemma:latala_diagonal_noise}
    Let $\bxi = (\xi_i)_{i=1}^N$ be a random vector with independent mean zero and variance $\sigma_\xi$ real-valued coordinates. We assume that for all $i$'s,  $\norm{\xi_i}_{L_r}\leq \okappa{kappa_noise}\sigma_\xi$ for some $\nkappa\label{kappa_noise}>0$ and $r>4$. There then exists some absolute constant $C_\okappa{kappa_noise}$ (depending only on $\okappa{kappa_noise}$) such that for any matrix $D \in\bR^{p\times N}$ the following holds: if for some integer $k$ for which  $\sqrt{k} \norm{D}_{\text{op}}\leq \sqrt{\Tr(DD^\top)}$ then with probability at least $1-(C_\okappa{kappa_noise}/k)^{r/4}$,
    \begin{equation*}
        \norm{D\bxi}_2 \leq  (3/2) \sigma_\xi \sqrt{\Tr(DD^\top)}.
    \end{equation*}
\end{Lemma}

We emphasize that Lemma~\ref{lemma:latala_diagonal_noise} does not depend on $p$, so we can set $p=\infty.$ Note that there exists an isometric embedding from $\cH$ to $\ell_2$, given by $f\in\cH\mapsto \sum_{j=1}^\infty\left\langle f,\varphi_j\right\rangle_\cH \ve_j$, where we recall that $(\ve_j)_{j=1}^\infty$ is ONB of $\ell_2$. Therefore, we extend $D:\bR^N\to\ell_2$ to $D:\bR^N\to\cH$. As a result, we have the following proposition:
\begin{Proposition}\label{prop:noise_concentration}
    Let $\bxi = (\xi_i)_{i=1}^N$ be a random vector with independent mean zero and variance $\sigma_\xi$ real-valued coordinates. We assume that for all $i$'s,  $\norm{\xi_i}_{L_r}\leq \okappa{kappa_noise}\sigma_\xi$ for some $\okappa{kappa_noise}>0$ and $r>4$. There then exists some absolute constant $\nC\label{C_noise}$ (depending only on $\okappa{kappa_noise}$) such that for any Hilbert-Schmidt operator $D:\bR^N\to\cH$ the following holds: if for some integer $k$ for which  $\sqrt{k} \norm{D}_{\text{op}}\leq \sqrt{\Tr(DD^\top)}$ then with probability at least $1-(\oC{C_noise}/k)^{r/4}$,
    \begin{equation*}
        \norm{D\bxi}_\cH \leq  (3/2) \sigma_\xi \sqrt{\Tr(DD^\top)}.
    \end{equation*}
\end{Proposition}

\begin{Proposition}\label{prop:noise_concentration_dependent}
    Suppose $f^*-f^{**}\in L_{2+\varepsilon}(\mu)$ for some $\varepsilon\geq 0$, where we recall that $\mu$ is the probability distribution of design vector $X$. Suppose $\xi,\xi_1,\cdots,\xi_N$ are i.i.d. mean zero sub-Gaussian random variables with variance $\sigma_\xi^2$ and suppose $\xi$ is independent with $X$.
    Let $\veps = (f^*(X_i)-f^{**}(X_i)+\xi_i)_{i=1}^N$.
    
    \begin{enumerate}
        \item When $A = (\bX_{\phi,k+1:\infty}\bX_{\phi,k+1:\infty}^\top + \lambda I_N)^{-1}\bX_{\phi,1:k}\tilde\Gamma_{1:k}^{-1}\bX_{\phi,1:k}^\top(\bX_{\phi,k+1:\infty}\bX_{\phi,k+1:\infty}^\top +\lambda I_N)^{-1}$.
        \begin{enumerate}
            \item When $\varepsilon>0$. Suppose $N\geq e^{-1}\left(\exp(2\oc{c_Men_1}/\varepsilon)-4\right)$. Recall the definition of $\sigma(\square,\triangle)$ from \eqref{eq:def_sigma_square_triangle} and the definition of $t(\square,\triangle)$ from \eqref{eq:def_t_square_triangle}. There exist absolute constants $\oC{C_noise_dependent_1}$ depending on $\varepsilon$, $\oC{C_RIP_upper}$ and $\oc{c_Men_3}$ such that for any $t_1\geq 0$ and $t_3>2$,
            \begin{align}\label{eq:noise_concentration_dependent_1}
            \begin{aligned}
                    &\bP\bigg(\veps^\top A\veps \leq 2(1+t_1)\sigma_\xi^2 \frac{16\oC{C_RIP_upper}^2 N\left(\left|J_1\right|\square^2 + \triangle^2\sum_{j\in J_2}\sigma_j\right)}{\left(4\lambda + \Tr\left(\Gamma_{k+1:\infty}\right)\right)^2} \\
                    &+t_3^2\oC{C_noise_dependent_1}^2\frac{\left(\sqrt{N}\sigma\left(\square,\triangle\right)\right)^2}{\left(4\lambda + \Tr\left(\Gamma_{k+1:\infty}\right)\right)^2}  \norm{f^*-f^{**}}_{L_{2+\varepsilon}}^2 N\bigg)\\
                    &\geq 1-\bP(\Omega_0^c) - \exp\left(-c(t_1^2\wedge t_1)\frac{\left|J_1\right|\square^2 + \triangle^2\sum_{j\in J_2}\sigma_j}{\sigma^2\left(\square,\triangle\right)}\right)  - \oc{c_Men_3}t_3^{-(2+\varepsilon)}N^{-\frac{\varepsilon}{4}}.
            \end{aligned}
            \end{align}
            \item When $\varepsilon=0$, \eqref{eq:noise_concentration_dependent_1} is still valid with not necessarily $N\geq e^{-1}\left(\exp(2\oc{c_Men_1}/\varepsilon)-4\right)$, but with $\oC{C_noise_dependent_1}$ replaced by $\nC\label{C_noise_dependent_2}$, where $\oC{C_noise_dependent_2} = 16\oC{C_RIP_upper}^2$; and $\oc{c_Men_3}$ replaced by $1$.
        \end{enumerate}

        \item When $A = (\bX_{\phi,k+1:\infty}\bX_{\phi,k+1:\infty}^\top + \lambda I_N)^{-1}\bX_{\phi,k+1:\infty}\Gamma_{k+1:\infty}\bX_{\phi,k+1:\infty}^\top(\bX_{\phi,k+1:\infty}\bX_{\phi,k+1:\infty}^\top +\lambda I_N)^{-1}$.  
        \begin{enumerate}
            \item When $\varepsilon>0$. Suppose $N\geq e^{-1}\left(\exp(2\oc{c_Men_1}/\varepsilon)-4\right)$. There exist absolute constants $\oC{C_noise_dependent_3} $ depending on $\varepsilon$, $\oC{C_sum_Gamma_phi}$, $\oC{C_DMU}$ and $\oc{c_Men_3}$, such that for any $t_1,t_2\geq 0$ and $t_3>2$,
            \begin{align}\label{eq:noise_concentration_dependent_2}
            \begin{aligned}
                &\bP\bigg( \veps^\top A\veps \leq (1+t_1)\sigma_\xi^2 \frac{16\oC{C_sum_Gamma_phi}^2 N\Tr(\Gamma_{k+1:\infty}^2)}{(4\lambda + \Tr(\Gamma_{k+1:\infty}))^2}\\
                &+ t_2t_3\sigma_\xi\frac{\oC{C_noise_dependent_3}\left(\Tr(\Gamma_{k+1:\infty}^2)+N\norm{\Gamma_{k+1:\infty}}_{\text{op}}^2 \right)}{(4\lambda+\Tr(\Gamma_{k+1:\infty}))^2}\sqrt{N}\norm{f^*-f^{**}}_{L_2}\\
                &+ t_3^2\frac{\oC{C_noise_dependent_3}\left(\Tr(\Gamma_{k+1:\infty}^2)+N\norm{\Gamma_{k+1:\infty}}_{\text{op}}^2 \right)}{(4\lambda+\Tr(\Gamma_{k+1:\infty}))^2} N \norm{f^*-f^{**}}_{L_{2+\varepsilon}}^2  \bigg)\\
                &\geq 1- \bP(\Omega_0^c) - \exp(-ct_2^2) - \oc{c_Men_3}t_3^{-(2+\varepsilon)}N^{-\frac{\varepsilon}{4}}\\
                &- \exp\left( -c(t_1^2\wedge t_1)\frac{8\oC{C_sum_Gamma_phi}N\Tr(\Gamma_{k+1:\infty}^2)}{\oC{C_DMU}^2\left(\Tr(\Gamma_{k+1:\infty}^2)+N\norm{\Gamma_{k+1:\infty}}_{\text{op}}\right)} \right)
            \end{aligned}
        \end{align}
        \item When $\varepsilon=0$, \eqref{eq:noise_concentration_dependent_2} is still valid with not necessarily $N\geq e^{-1}\left(\exp(2\oc{c_Men_1}/\varepsilon)-4\right)$, but with $\oC{C_noise_dependent_3}$ replaced by $32\oC{C_DMU}\oc{c_Men_2}$ and with $\oc{c_Men_3}$ replaced by $1$.
        \end{enumerate}
    \end{enumerate}

\end{Proposition}

\beginproof
Recall that $\veps = \vr + \bxi$. We have:
\begin{align*}
    \veps^\top A\veps = \vr^\top A\vr + \bxi^\top A\bxi + \vr^\top A\bxi + \bxi^\top A\vr \leq 2\norm{A}_{\text{op}} \norm{\vr}_2^2+ 2\bxi^\top A\bxi.
\end{align*}
By Hanson-Wright inequality, see, for example \cite[Theorem 6.2.1]{vershynin_high-dimensional_2018}, there exists some absolute constant $c>0$ such that for any $t_1\geq 0$,
\begin{align*}
    \bP\left( \bxi^\top A\bxi - \sigma_\xi^2\Tr(A) \leq t_1 \sigma_\xi^2\Tr(A) \right) \geq 1- \exp\left(-c(t_1^2\wedge t_1)\frac{\Tr(A)}{\norm{A}_{\text{op}}}\right).
\end{align*}

Let $\Omega_{\text{noise},1}$ as the random event on which
\begin{align*}
    \sqrt{\Tr(A)}\leq \frac{4\oC{C_RIP_upper}\sqrt{N}}{4\lambda + \Tr\left(\Gamma_{k+1:\infty}\right)}\sqrt{\left|J_1\right|\square^2 + \triangle^2\sum_{j\in J_2}\sigma_j},\mbox{ and } \norm{A}_{\text{op}}^{1/2}\leq \frac{4\oC{C_RIP_upper}\sqrt{N}\sigma\left(\square,\triangle\right)}{4\lambda + \Tr\left(\Gamma_{k+1:\infty}\right)}.
\end{align*}By \eqref{eq:trace_DD_top_1} and \eqref{eq:op_D_1}, $\bP(\Omega_{\text{noise},1})\geq 1-\bP(\Omega_0^c)$.
In Lemma~\ref{lem:moment-sum-variid-positive}, let $Z = f^*(X)-f^{**}(X)$, $q=2+\varepsilon$ for some $\varepsilon>0$, $r=2$, $p=1$. When $N\geq e^{-1}\left(\exp(2\oc{c_Men_1}/\varepsilon)-4\right)$, by Lemma~\ref{lem:moment-sum-variid-positive}, we have $j_0=1$ and thus for $\beta = \varepsilon/4$, for any $t_3>2$, with probability at least $1-\oc{c_Men_3}t_3^{-(2+\varepsilon)}N^{-\frac{\varepsilon}{4}}$,
\begin{align*}
    \norm{\vr}_2 \leq \oc{c_Men_2}\left(\frac{2+\varepsilon}{2+\varepsilon-2(\varepsilon/4+1)}\right)^{1/2}t_3 \norm{f^*-f^{**}}_{L_{2+\varepsilon}}\sqrt{N}.
\end{align*}
Combining the above together, we obtain that for any $t_1\geq 0$ and $t_3>2$,
\begin{eqnarray*}
    &\bP\bigg(\veps^\top A\veps \leq 2(1+t_1)\sigma_\xi^2 \frac{16\oC{C_RIP_upper}^2 N\left(\left|J_1\right|\square^2 + \triangle^2\sum_{j\in J_2}\sigma_j\right)}{\left(4\lambda + \Tr\left(\Gamma_{k+1:\infty}\right)\right)^2} + t_3^2\oC{C_noise_dependent_1}^2\frac{\left(\sqrt{N}\sigma\left(\square,\triangle\right)\right)^2}{\left(4\lambda + \Tr\left(\Gamma_{k+1:\infty}\right)\right)^2}  \norm{f^*-f^{**}}_{L_{2+\varepsilon}}^2 N\bigg)\\
    &\geq 1-\bP(\Omega_0^c) - \exp\left(-c(t_1^2\wedge t_1)\frac{\left|J_1\right|\square^2 + \triangle^2\sum_{j\in J_2}\sigma_j}{\sigma^2\left(\square,\triangle\right)}\right) - \oc{c_Men_3}t_3^{-(2+\varepsilon)}N^{-\frac{\varepsilon}{4}},
\end{eqnarray*}where $\nC\label{C_noise_dependent_1} = 16\oC{C_RIP_upper}^2\oc{c_Men_2}\left(\frac{2+\varepsilon}{2+\varepsilon-2(\varepsilon/4+1)}\right)^{1/2}$. In particular, when we use Markov's inequality to replace Lemma~\ref{lem:moment-sum-variid-positive}, with the probability deviation $t_3^{-2}$ instead of $\oc{c_Men_3}t_3^{-(2+\varepsilon)}N^{-\frac{\varepsilon}{4}}$, the term $\sqrt{N}\norm{f^*-f^{**}}_{L_{2+\varepsilon}}$ can be improved to $\sqrt{N}\norm{f^*-f^{**}}_{L_2}$.

Similarly, let $\Omega_{\text{noise},2}$ as the random event on which
\begin{align*}
\Tr(A)\leq \frac{16\oC{C_sum_Gamma_phi}N \Tr(\Gamma_{k+1:\infty}^2)}{\left(4\lambda + \Tr(\Gamma_{k+1:\infty})\right)^2}\mbox{ and }    \norm{A}_{\text{op}}^{1/2}\leq \frac{4\oC{C_DMU}}{4\lambda + \Tr(\Gamma_{k+1:\infty})}\left(\sqrt{\Tr(\Gamma_{k+1:\infty}^2)}+\sqrt{N}\norm{\Gamma_{k+1:\infty}}_{\text{op}}\right).
\end{align*}By \eqref{eq:trace_DD_top_2} and \eqref{eq:op_D_2}, $\bP(\Omega_{\text{noise},2})\geq 1-\bP(\Omega_0^c)$. Repeat the above arguments, we obtain that for any $t_1,t_2\geq 0$ and $t_3>2$,
\begin{align*}
    &\bP\bigg( \veps^\top A\veps \leq (1+t_1)\sigma_\xi^2 \frac{16\oC{C_sum_Gamma_phi}^2 N\Tr(\Gamma_{k+1:\infty}^2)}{(4\lambda + \Tr(\Gamma_{k+1:\infty}))^2} + t_2t_3\sigma_\xi\frac{\oC{C_noise_dependent_3}\left(\Tr(\Gamma_{k+1:\infty}^2)+N\norm{\Gamma_{k+1:\infty}}_{\text{op}}^2 \right)}{(4\lambda+\Tr(\Gamma_{k+1:\infty}))^2}\sqrt{N}\norm{f^*-f^{**}}_{L_2}\\
    &+ t_3^2\frac{\oC{C_noise_dependent_3}\left(\Tr(\Gamma_{k+1:\infty}^2)+N\norm{\Gamma_{k+1:\infty}}_{\text{op}}^2 \right)}{(4\lambda+\Tr(\Gamma_{k+1:\infty}))^2} N \norm{f^*-f^{**}}_{L_{2+\varepsilon}}^2  \bigg)\\
    &\geq 1- \bP(\Omega_0^c) - \exp(-ct_2^2) - \oc{c_Men_3}t_3^{-(2+\varepsilon)}N^{-\frac{\varepsilon}{4}} - \exp\left( -c(t_1^2\wedge t_1)\frac{8\oC{C_sum_Gamma_phi}N\Tr(\Gamma_{k+1:\infty}^2)}{\oC{C_DMU}^2\left(\Tr(\Gamma_{k+1:\infty}^2)+N\norm{\Gamma_{k+1:\infty}}_{\text{op}}\right)} \right)
\end{align*}where $\nC\label{C_noise_dependent_3} = 32\oC{C_DMU}\oc{c_Men_2}\left(\frac{2+\varepsilon}{2+\varepsilon-2(\varepsilon/4+1)}\right)^{1/2}$.

\endproof

For simplicity, we only consider the case where $\sigma_1 N > \kappa_{DM}(4\lambda+\Tr(\Gamma_{k+1:\infty}))$ and $k\lesssim N$. 
Replace the usage of Proposition~\ref{prop:noise_concentration} in Section~\ref{sec:proof_main_upper_k<N} with Proposition~\ref{prop:noise_concentration_dependent} and repeat the proof. We can conclude as follows:
\begin{Proposition}\label{prop:dependent_noise_result}
    Grant the assumptions of Theorem~\ref{theo:upper_KRR} and Proposition~\ref{prop:noise_concentration_dependent}. Assume that $\sigma_1 N > \kappa_{DM}(4\lambda+\Tr(\Gamma_{k+1:\infty}))$. Then with the same probability deviation as in Theorem~\ref{theo:upper_KRR} but with $\bar p_{\xi}$ replaced by \eqref{eq:noise_concentration_dependent_1} and
    \begin{align*}
        \left(\frac{\oC{C_noise}\Tr(\Gamma_{k+1:\infty})}{|J_1|\Tr(\Gamma_{k+1:\infty}) + N \left(\sum_{j\in J_2}\sigma_j\right)}\right)^{\frac{r}{4}}
    \end{align*}
    replaced by \eqref{eq:noise_concentration_dependent_2}, we have:
    \begin{align}
    \begin{aligned}\label{eq:new_term_data_dependent_noise}
        \norm{\hat f_\lambda - f^*}_{L_2}&\lesssim \norm{f^*-f^{**}}_{L_2} + \sigma_\xi\sqrt{\frac{\left|J_1\right|}{N}}+ \sigma_\xi\sqrt{\frac{\sum_{j\in J_2}\sigma_j}{  4\lambda + \Tr\left(\Gamma_{k+1:\infty}\right)}} +  \norm{\Gamma_{k+1:\infty}^{1/2}f_{k+1:\infty}^{**}}_\cH\\
        &+ {\norm{\tilde\Gamma_{1,\mathrm{thre}}^{-1/2}f_{1:k}^{**}}_\cH \frac{2\lambda + 3\Tr\left(\Gamma_{k+1:\infty}\right)}{N} }
        +4\oC{C_sum_Gamma_phi}\sigma_\xi \frac{\sqrt{ N\Tr(\Gamma_{k+1:\infty}^2)}}{4\lambda + \Tr(\Gamma_{k+1:\infty})} \\
    &+ \frac{\oC{C_noise_dependent_3}^{1/2}\sqrt{\Tr(\Gamma_{k+1:\infty}^2)+N\norm{\Gamma_{k+1:\infty}}_{\text{op}}^2 }}{4\lambda+\Tr(\Gamma_{k+1:\infty})} \sqrt{N} \norm{f^*-f^{**}}_{L_{2+\varepsilon}}.
    \end{aligned}
    \end{align}
\end{Proposition}
\cite[Section 7.5.2]{bach_learning_2024} uses $\inf\{\norm{f^*-f}_{L_2}^2+\lambda\norm{f}_\cH^2 : f\in\cH\}$ to characterize the approximation error of $\cH$, that is, the trade-off between the approximation error and $\norm{f}_\cH$. Our Proposition~\ref{prop:dependent_noise_result} shows that the approximation error is actually traded off against ${\norm{\tilde\Gamma_{1,\mathrm{thre}}^{-1/2}f_{1:k}^{**}}_\cH \frac{2\lambda + 3\Tr\left(\Gamma_{k+1:\infty}\right)}{N} }$ and $\norm{\Gamma_{k+1:\infty}^{1/2}f_{k+1:\infty}^{**}}_\cH$ instead of $\norm{f^{**}}_\cH$.

We observe that when the following holds, we still have benign overfitting, even though there is model-misspecification, when $\sigma_\xi\sim 1$: for $k\lesssim N$ such that $N\norm{\Gamma_{k+1:\infty}}_{\text{op}}\lesssim\Tr(\Gamma_{k+1:\infty})$,
\begin{align}\label{eq:conditions_BO}
\begin{aligned}
    &\left|J_1\right|=o(N),\, \sum_{j\in J_2}\sigma_j = o(\Tr(\Gamma_{k+1:\infty})),\, \norm{\Gamma_{k+1:\infty}^{1/2}f_{k+1:\infty}^{**}}_\cH = o(1),\, {\norm{\tilde\Gamma_{1,\mathrm{thre}}^{-1/2}f_{1:k}^{**}}_\cH \frac{\Tr\left(\Gamma_{k+1:\infty}\right)}{N} }=o(1),\\
    &\frac{\sqrt{ N\Tr(\Gamma_{k+1:\infty}^2)}}{\Tr(\Gamma_{k+1:\infty})} = o(1),\, \sqrt{N}\norm{f^*-f^{**}}_{L_{2+\varepsilon}} = o\left(\frac{\Tr(\Gamma_{k+1:\infty})}{\sqrt{ \Tr(\Gamma_{k+1:\infty}^2)}}\wedge \frac{\Tr(\Gamma_{k+1:\infty})}{\sqrt{N}\norm{\Gamma_{k+1:\infty}}_{\text{op}}}\wedge \sqrt{N}\right).
\end{aligned}
\end{align}

\begin{Remark}\label{remark:no_contradiction}
    We emphasize that this conclusion does not contradict the counterexamples presented in \cite{chinot_robustness_2022} and \cite{shamir_implicit_2022} in the setting of adversarial noise and model-misspecification. This is because in \cite{chinot_robustness_2022}, $\|f^* - f^{**}\|_{L_2} = \Theta(\sigma_\xi)$ (see \cite[Appendix D]{chinot_robustness_2022}); in \cite{shamir_implicit_2022}, $f^*: \vx \in \mathbb{R}^d \mapsto \exp(x_1)$, where $x_1$ is the first coordinate of $\vx$, and $f^{**} = \left<\vx,\bbeta^*\right>$ for some $\bbeta^* \in \bR^d$. Under the probability measure $\mu$ assumed in \cite[Example 1]{shamir_implicit_2022}, $\norm{f^*-f^{**}}_{L_2}$ is $\inf_{\mathbf{w}\in\bR^d}R_d(\mathbf{w})$ in the language of \cite[Example 1]{shamir_implicit_2022}, and it is grater than a constant depending only on $\mu$. Therefore, when $\|f^*-f^{**}\|_{L_2}=\Theta(1)$, our \eqref{eq:conditions_BO} does not necessarily yield benign overfitting, thus not conflicting with \cite{chinot_robustness_2022,shamir_implicit_2022}.
\end{Remark}

\subsubsection{Proof of Proposition~\ref{prop:multiple_descent}}\label{sec:proof_multiple_descent}

The proof of the results in the asymptotic regime is presented in Section~\ref{sec:multiple_descent_asymptotic}. In this section, we focus on the proof of the non-asymptotic conclusions.

\begin{flushleft}
The proof of Proposition~\ref{prop:multiple_descent} is rather lengthy, and we complete it in Section~\ref{sec:spectrum_Gamma_multiple_descent}, Section~\ref{sec:norm_equivalence_multiple_descent}, Section~\ref{sec:diagonal_terms_multiple_descent}, Section~\ref{sec:proof_spectrum_Gamma_sub_gaussian}, and Section~\ref{sec:proof_diagonal_concentration}. Section~\ref{sec:proof_diagonal_concentration} relies on the proof of Theorem~\ref{theo:polynomial_HS} in Section~\ref{sec:proof_polynomial_HS}.
\end{flushleft}

Before we begin, we first indicate how to obtain an upper bound for
\begin{align*}
    \left|\|\hat{f}_0 - f^*\|_{L_2(\mu)}^2 - \|f_{>\iota}^*\|_{L_2(\mu)}^2 - \| f^* - f^{**} \|_{L_2(\mu)}^2\right|.
\end{align*}For convenience, we denote $\hat f_0$ by $\hat f$. Recall that $\cH_{1:k}$ is orthogonal to $\cH_{k+1:\infty}$, we have $\|\hat f_0 - f^*\|_{L_2(\mu)}^2 = \|\Gamma^{1/2}(\hat f - f^*)\|_\cH^2 = \|\hat f_{1:k}- f_{1:k}^*\|_{L_2(\mu)}^2 + \|\hat f_{k+1:\infty} - f_{k+1:\infty}^*\|_{L_2(\mu)}^2$. By subtracting \( (\|f_{k+1:\infty}^*\|_{L_2(\mu)}^2 + \|f^* - f^{**}\|_{L_2(\mu)}^2) \) from both sides of the equation and taking the absolute value, we have by the triangle inequality:
\begin{align}\label{eq:pre_multiple_descent}
    \begin{aligned}
        &\left| \norm{\hat f_0 - f^*}_{L_2(\mu)}^2 - \norm{f_{k+1:\infty}^*}_{L_2(\mu)}^2 - \norm{f^* - f^{**}}_{L_2(\mu)}^2 \right|\\
    &\leq \|\hat f_{1:k}- f_{1:k}^*\|_{L_2(\mu)}^2 + \|\hat f_{k+1:\infty} - f_{k+1:\infty}^*\|_{L_2(\mu)}^2 + \norm{f_{k+1:\infty}^*}_{L_2(\mu)}^2 + \norm{f^* - f^{**}}_{L_2(\mu)}^2.
    \end{aligned}
\end{align}The right side can be further bounded from above by \( (r_{0,k}^*)^2 \) (see Proposition~\ref{prop:KRR_estimate} and Proposition~\ref{prop:KRR_noise_absorption}), and in the following, we will see that \( f_{k+1:\infty}^* \) is identified as \( f_{>\iota}^* \).

\subsubsection{Spectrum of $\Gamma$}\label{sec:spectrum_Gamma_multiple_descent}
The goal of this section is to use Theorem~\ref{theo:upper_KRR} and Theorem~\ref{theo:main_upper_k>N} to prove the multiple descent phenomenon. For applying Theorem~\ref{theo:upper_KRR}, we need to proceed with the following steps:
\begin{enumerate}
    \item Verify Assumption~\ref{assumption:DM_L4_L2}, Assumption~\ref{assumption:upper_dvoretzky} and Assumption~\ref{assumption:RIP}.

     To achieve diagonal concentration, it suffices to verify \eqref{eq:diagonal_term_assumption}, \eqref{eq:diagonal_upper_dvoretzky}, and \eqref{eq:diagonal_RIP}. This is done in Section~\ref{sec:diagonal_terms_multiple_descent}. To validate the norm-equivalence condition stated in Assumption~\ref{assumption:DM_L4_L2}, Assumption~\ref{assumption:upper_dvoretzky}, and Assumption~\ref{assumption:RIP}, we will examine the scenario where $\eps=6$ over the entire space $\cH$. This shows the norm-equivalence on $\cH_{1:k}$ and $\cH_{k+1:\infty}$, as these subspaces are contained within $\cH$. This is done in Section~\ref{sec:norm_equivalence_multiple_descent}.
    
    \item Compute the appropriate $k$ such that $N\leq\oc{c_kappa_DM}\kappa_{DM}d_0^*\left(\Gamma_{k+1:\infty}^{1/2}B_\cH\right)$, and $k\leq\oc{c_RIP}N$. In the context of multiple descent, we will make use of Theorem~\ref{theo:upper_KRR}, meaning that we will observe multiple descents occurring when $k\lesssim N$. Given that $\Tr(\Gamma_{k+1:\infty})\sim 1$, our findings remain valid even when $\lambda\lesssim 1$.
    
    \item Compute $\Tr\left(\Gamma_{k+1:\infty}\right)$, $\norm{\Gamma_{k+1:\infty}}_{\text{op}}$ and $\Tr\left(\Gamma_{k+1:\infty}^2\right)$.
    
    The last two steps are done in the following two paragraphs.
\end{enumerate}

In this subsection, the determination of the appropriate value for $k$ is discussed. The selection of the value of $k$ can be divided into two distinct scenarios: under a sub-Gaussian (the setup in \cite{liang_multiple_2020}) or a uniform distribution assumption (the setup in \cite{ghorbani_linearized_2021,mei_generalization_2022,mei_generalization_2020,misiakiewicz_spectrum_2022}) and for a kernel $K$ such that \eqref{eq:polynomial_inner_product_kernel} holds.

\paragraph{Sub-Gaussian design} 


We have the following proposition, whose proof can be found in Section~\ref{sec:proof_spectrum_Gamma_sub_gaussian}

\begin{Proposition}\label{prop:spectrum_Gamma_sub_gaussian}
    Suppose that Assumption~\ref{assumption:polynomial_feature_sub_gaussian} holds. For a given $0\leq \iota\leq L-1$, there exist absolute constants $\nC\label{C_multiple_lower}$ and $\nc\label{c_multiple_upper}$ (where $\oc{c_multiple_upper}$ depending on $L$) such that $\oC{C_multiple_lower} d^{\iota} \leq N\leq \oc{c_multiple_upper} d^{\iota+1} $. We consider $k  =\sum_{0\leq l\leq\iota} d^{l}$ and decompose the kernel function into two terms:
    $$K(X,X) = h\left(\frac{\norm{X}_2^2}{d}\right) = \sum_{0\leq i\leq \iota}\alpha_i\left(\frac{\norm{X}_2^2}{d}\right)^i + \sum_{i>\iota}\alpha_i \left(\frac{\norm{X}_2^2}{d}\right)^i.$$
    Then,
\begin{align}\label{eq:quantities_sub_gaussian}
    \begin{aligned}
        &\frac{k}{N}\sim\frac{d^\iota}{N},\quad \frac{\Tr\left(\Gamma_{k+1:\infty}\right)}{N}\sim \frac{1}{N},\quad \frac{\Tr\left(\Gamma_{k+1:\infty}\right)}{\norm{\Gamma_{k+1:\infty}}_{\text{op}}} \gtrsim d^{\iota+1}>\oc{c_kappa_DM}\kappa_{DM} N,\mbox{ and }\\
        & \frac{\sqrt{N\Tr\left(\Gamma_{k+1:\infty}^2\right)}}{\Tr\left(\Gamma_{k+1:\infty}\right)} \sim \sqrt{\frac{N}{d^{\iota+1}}}.
    \end{aligned}
\end{align}
\end{Proposition}

It is worth noting that the selection of $k=\sum_{l\leq\iota}d^l$ means that $\hat f_0$ learns $f_{\leq\iota}^*$ while considering $f_{>\iota}^*$ as noise. By observing that $f_{\leq\iota}^*$ represents the projection onto the subspace formed by the orthogonal polynomials of degree $\iota$ in the Hilbert space $L_2(\mu)$, it becomes evident that $\hat f_0$ acquires knowledge of a polynomial approximation of degree $\iota$ for $f^*$. By increasing the number of samples inside the range of $\oC{C_multiple_lower} d^{\iota}\leq N\leq \oc{c_multiple_upper} d^{\iota+1}$ to $\oC{C_multiple_lower} d^{\iota+1}\leq N \leq \oc{c_multiple_upper} d^{\iota+2}$, a polynomial approximation may be achieved with an additional degree. Proposition~\ref{prop:spectrum_Gamma_sub_gaussian} is still valid when $\Tr\left(\Gamma_{k+1:\infty}\right)$ is replaced by $\lambda + \Tr\left(\Gamma_{k+1:\infty}\right)$, where $\lambda\lesssim 1$.

\paragraph{Uniform distribution design}


When $\Omega=\sqrt{d}S_2^{d-1}$ and $\mu$ (the distribution of $X$) is uniform over $\sqrt{d}S_2^{d-1}$, $L_2(\sqrt{d}S_2^{d-1}, \mu)$ admits an orthonormal decomposition. The following statements are taken from \cite{mei_generalization_2022,mei_generalization_2020,ghorbani_linearized_2021,misiakiewicz_spectrum_2022}. In fact, for $l\in\bN$, let $\tilde V_{d,l}$ be the space of homogeneous harmonic polynomials of degree $l$ on $\bR^d$, that is, the space of degree-$l$ polynomials $q(\vx)$ such that $\Delta q(\vx)=0$, where $\Delta$ is the Laplacian. Let $V_{d,l}$ be the linear space of functions obtained by restricting the polynomials in $\tilde V_{d,l}$ to $\sqrt{d}S_2^{d-1}$. Then $L_2(\sqrt{d}S_2^{d-1},\mu) = \oplus^\perp_{l\in\bN} V_{d,l}$. Furthermore, by \cite[section 2]{xiao_precise_2022}, when $l< L$ and $L!=O(d)$, we have
\begin{align}\label{eq:multicity_spherical_harmonics}
    \dim\left(V_{d,l}\right) = \frac{d^l}{l!}+O(d^{l-1})\sim d^l.
\end{align}
When Assumption~\ref{assumption:polynomial_feature_gaussian} is true and for $0\leq\iota\leq L-1$ such that $\oC{C_multiple_lower} d^{\iota}\leq N\leq \oc{c_multiple_upper} d^{\iota+1}$, we define $k=\sum_{l=0}^\iota\dim(V_{d,l})$. In the work of \cite[section 2]{xiao_precise_2022}, it is demonstrated that the eigenvalues with multiplicity (denoted by $\lambda_l$) of the matrix $\Gamma$ exhibit a decay rate on the order of $d^{-l}$. Under Assumption~\ref{assumption:polynomial_feature_gaussian}, $\sigma_l$ exhibits a multi-plateau trend as in the sub-Gaussian case. Consequently, the spectrum of the variable $\Gamma$ exhibits similarity to that of the sub-Gaussian design scenario, but with a constant factor, and we decompose
\begin{align}\label{eq:decomposition_inner_product_kernel}
    K(X,X)=\sum_{i\leq\iota}\alpha_i \left(\frac{\norm{X}_2^2}{d}\right)^i + \sum_{i>\iota}\alpha_i \left(\frac{\norm{X}_2^2}{d}\right)^i =: K_{1:k}(X,X)+K_{k+1:\infty}(X,X)
\end{align}
as well. Therefore,
\begin{align}\label{eq:quantities_uniform}
    \begin{aligned}
        &\frac{k}{N}\sim \frac{d^\iota}{N},\quad \frac{\Tr\left(\Gamma_{k+1:\infty}\right)}{N}\sim\frac{1}{N},\quad \frac{\Tr\left(\Gamma_{k+1:\infty}\right)}{\norm{\Gamma_{k+1:\infty}}_{\text{op}}}\gtrsim d^{\iota+1}>\oc{c_kappa_DM}\kappa_{DM} N,\mbox{ and }\\
        &\frac{\sqrt{N\Tr\left(\Gamma_{k+1:\infty}^2\right)}}{\Tr\left(\Gamma_{k+1:\infty}\right)}\sim\sqrt{\frac{N}{d^{\iota+1}}}.
    \end{aligned}
\end{align}

In this way, we have completed the steps 2) and 3) outlined at the beginning of this subsection. In the following two subsections, we will complete step 1).

\subsubsection{Norm equivalence.}\label{sec:norm_equivalence_multiple_descent}
\paragraph{Under Assumption~\ref{assumption:polynomial_feature_sub_gaussian}} We strengthen \cite[Lemma 10]{liang_multiple_2020} from $L_4-L_2$ norm equivalence to $L_8-L_2$ norm equivalence. The proof of the following Proposition can be found in Section~\ref{sec:proof_norm_equivalence_sub_gaussian}.
\begin{Proposition}\label{prop:norm_equivalence_sub_gaussian}
    Under Assumption~\ref{assumption:polynomial_feature_sub_gaussian}, there exists $\kappa\geq 1$ depending only on $L$ (with $L^L$ dependence), such that for all $f\in\cH$,
    \begin{align*}
        \norm{f}_{L_8}\leq \kappa \norm{f}_{L_2}.
    \end{align*}
\end{Proposition}

The proof of Proposition~\ref{prop:norm_equivalence_sub_gaussian} may be found in Section~\ref{sec:proof_norm_equivalence_sub_gaussian}.

\paragraph{Under Assumption~\ref{assumption:polynomial_feature_gaussian}}  For every $l\in\bN$, $(V_{d,l},L_2(\mu))$ admits an ONB $(Y_{l,m})_{m\leq \dim\left(V_{d,l}\right)}$ where $(Y_{l,m})_{m\leq \dim\left(V_{d,l}\right)}$ are homogeneous polynomials of degree $l$. For further details, please refer to section 4.2 of \cite{frye_spherical_2012}. Given that $\mu$ follows a uniform distribution over $\sqrt{d}S_2^{d-1}$ and $\cH$ is a subset of $\oplus_{l>\iota}V_{d,l}$ that consists of polynomials with degrees lower than $L$, \eqref{eq:norm_equivalence} can be derived directly from the work of \cite{beckner_sobolev_1992} by setting the values of the parameters $\eps=6$ and $\kappa<7^L$, respectively.

\subsubsection{Diagonal terms.}\label{sec:diagonal_terms_multiple_descent}

In this subsection, we check \eqref{eq:diagonal_term_assumption} under Assumption~\ref{assumption:polynomial_feature_gaussian} and Assumption~\ref{assumption:polynomial_feature_sub_gaussian}. The concentration property of the diagonal terms appeared in Assumption~\ref{assumption:upper_dvoretzky} and Assumption~\ref{assumption:RIP} followed by the same idea. The scenario described in Assumption~\ref{assumption:polynomial_feature_gaussian} is comparatively simpler than the case outlined in Assumption~\ref{assumption:polynomial_feature_sub_gaussian}. We therefore prioritize addressing the former case initially.

\paragraph{Under Assumption~\ref{assumption:polynomial_feature_gaussian}} Because $\mu$ is a uniform distribution over $\sqrt{d}S_2^{d-1}$, $K$ is a translation-invariant kernel. Recall that we have chosen $k=\sum_{l=0}^\iota\dim(V_{d,l})$. In fact, $\norm{\phi_{k+1:\infty}(X)}_{\cH} = \sqrt{K_{k+1:\infty}(X,X)} = \sqrt{\sum_{i>\iota}\alpha_i\left(\norm{X}_2^2/d\right)^i} = \sqrt{\sum_{i=\iota+1}^L \alpha_i} = \sqrt{\Tr\left(\Gamma_{k+1:\infty}\right)}$ which is a constant. As a result, \eqref{eq:diagonal_term_assumption} follows by taking $\delta = \gamma = 0$. The validity of the diagonal terms appeared in Assumption~\ref{assumption:RIP} and Assumption~\ref{assumption:upper_dvoretzky} follows by noticing that $\sigma_k^{-1}\sim k$ and $\sigma_{k+1}\sim \Tr(\Gamma_{k+1:\infty}^2)$ for our choice of $k$. Hence $\norm{\Gamma_{1:k}^{-1/2}\phi_{1:k}(X)}_\cH^2\leq\sigma_k^{-1}\norm{\phi_{1:k}(X)}_\cH^2\sim k$, and $\norm{\Gamma_{k+1:\infty}^{1/2}\phi_{k+1:\infty}(X)}_\cH^2\leq \sigma_{k+1}\norm{\phi_{k+1:\infty}(X)}_\cH^2\sim\Tr(\Gamma_{k+1:\infty}^2)$. As a result, we can take $\ogamma{gamma_RIP_k>N}=\ogamma{gamma_DMU_L2}=0$ and $\odelta{delta_DMU_L2}, \odelta{delta_RIP_k>N}\sim 1$.

\paragraph{Under Assumption~\ref{assumption:polynomial_feature_sub_gaussian}} 

Under Assumption~\ref{assumption:polynomial_feature_sub_gaussian}, we need to develop a new concentration inequality similar to the Hanson-Wright inequality for non-asymptotic verification of \eqref{eq:diagonal_term_assumption}. For the sake of generality, we will establish a concentration inequality for the general kernel function $K$ in the following theorem, although we will only utilize the special case of $K = K_{k+1:\infty}$.

\begin{Theorem}\label{theo:polynomial_HS}
    Let $X\in\bR^d$ be a random vector with i.i.d. mean zero coordinates $(x_i)_{i=1}^d$. Assume that $x_1$ is sub-Gaussian with the sub-Gaussian norm denoted by $\norm{x_1}_{\psi_2}$, that is to say, for any $q\geq 2$, we have $\bE x_1^q\leq\oC{C_subgaussian}q^{q/2}\norm{x_1}_{\psi_2}^q$. Given $L\in\bN_+$, satisfying $L\geq 2$, and $2^{2L}L^{3L}<\sqrt{d}$. Given $\gamma_1,\cdots,\gamma_L\in\bR$, we let $h:t\in\bR\mapsto \sum_{\iota=1}^L \gamma_\iota t^\iota$. We suppose $h(1)\neq 0$ for the sake of simplicity. Define $K:(\vx,\vy)\in\bR^d\times\bR^d\mapsto h(\left<\vx,\vy\right>/d)\in\bR$. There exists an absolute constant $\nC\label{C_non_linear_HS}$ depending on $L$ and $\norm{x_1}_{\psi_2}$, such that for any $t>0$, we have
    \begin{align*}
        \bP\left( \left| K(X,X)- \bE K(X,X) \right|\geq t \right)\leq \exp\left(-\oC{C_non_linear_HS}\frac{t^{\frac{1}{L}} d^{\frac{1}{2L}} }{ 4L^3 \left(\sum_{\iota=1}^L \gamma_\iota\right)^{\frac{1}{L}} }\right).
    \end{align*}
\end{Theorem} The proof of Theorem~\ref{theo:polynomial_HS} is postponed to Section~\ref{sec:proof_polynomial_HS}.

Applying Theorem~\ref{theo:polynomial_HS} to $K=K_{k+1:\infty}$, $h:u\mapsto \sum_{i>\iota}^L \alpha_i u^i$ and $t=\delta(\sum_{i>\iota}^L\alpha_i)$ together with a union bound, we can take 
\begin{align}\label{eq:def_gamma_multiple_descent}
    \delta = 1/(200\sqrt{\oC{C_distortion_1}}), \mbox{ and } \gamma = N \exp\left(-\oC{C_non_linear_HS_1}\frac{ d^{\frac{1}{2L}} }{ 4L^3 }\right)
\end{align}
for an absolute constant $\nC\label{C_non_linear_HS_1}$. Here, the item [4] from Assumption~\ref{assumption:polynomial_feature_sub_gaussian} ensures that when $d\gtrsim_{\mu_1}(\log(N))^{2L}L^{6L}$, we have $\gamma<1$.

We postpone the verification under Assumption~\ref{assumption:polynomial_feature_sub_gaussian} of the diagonal concentration condition appeared in Assumption~\ref{assumption:RIP} and Assumption~\ref{assumption:upper_dvoretzky} to Section~\ref{sec:proof_diagonal_concentration}.


\subsubsection{Proof of Proposition~\ref{prop:spectrum_Gamma_sub_gaussian}}\label{sec:proof_spectrum_Gamma_sub_gaussian}

We write the feature map explicitly
\begin{equation}\label{eq:feature_map_polynomial}
    \phi:\vx\in\bR^d\mapsto \sum_{\left|\cI\right|=0}^L{\sqrt{\frac{\alpha_{\left|\cI\right|}}{d^{\left|\cI\right|}}}\sqrt{C_{\cI}}}{P_{\cI}(\vx)}e_\cI ,\mbox{ s.t. }K(\vx,\vy)=\left\langle \phi(\vx),\phi(\vy)\right\rangle_{\ell_2}
\end{equation}where $C_{\cI}=\left|\cI\right|!/\left|\cI_1\right|!\cdots\left|\cI_d\right|!$, $P_{\cI}(\vx)=x_1^{\left|\cI_1\right|}\cdot\cdots\cdot x_d^{\left|\cI_d\right|}$ and $e_\cI$ are ONB of $\ell_2^{\binom{d+L}{L}}$. Notice that here we embed $\cH$ into $\ell_2$. The notation $\cI$ is defined in Section~\ref{sec:notations}.

Notice that $(P_\cI(X))_{\left|\cI\right|\leq L}$ is not an ONB of $L_2(\mu)$, because they are not orthogonal in $L_2(\mu)$. To obtain an ONB, we utilize the Gram-Schmidt procedure, employing the identical rationale as in \cite[Appendix A.2]{liang_multiple_2020}:

\begin{Lemma}\label{lemma:gram_schmidt_LRZ}
    Define $q_0, q_1, \cdots$ as the ONB of $L_2(\bR, \mu_1)$ produced by applying the Gram-Schmidt process on the basis $\{1, t, \cdots, t^k, \cdots\}$. 
    Then we have the following properties:
    \begin{itemize}
        \item $(Q_\cI)_{\left|\cI\right|=1}^\infty$ is an ONB of $L_2(\bR^d,\mu)$ when $Q_\cI:\vx\in\bR^d\mapsto\prod_{j=1}^d q_{\left|\cI_j\right|}(x_j)$ .
        \item Let
        \begin{equation}\label{eq:feature_map_polynomial_ONB}
            \psi:\vx\in\bR^d\mapsto \sum_{\left|\cI\right|=0}^L \sqrt{\frac{\alpha_{\left|\cI\right|}}{d^{\left|\cI\right|}}}\sqrt{C_\cI}Q_\cI(\vx)e_\cI.
        \end{equation}
        There exists an upper triangular non-singular matrix $\Lambda\in\bR^{\binom{d+L}{L}\times\binom{d+L}{L}}$ with $\max\{\norm{\Lambda}_{\text{op}},\, \norm{\Lambda^{-1}}_{\text{op}}\}\leq C_L$, such that $\phi=\Lambda^\top\psi$.
        
    \end{itemize}
\end{Lemma}The matrix $\Lambda^\top$ functions as the transformation matrix that converts the orthonormal basis $(Q_\cI)_{\left|\cI\right|\leq L}$ into a linearly independent basis $(P_\cI)_{\left|\cI\right|\leq L}$, which may be found in \cite[section 4.8, section 5.5]{meyer_matrix_2023}. As a result, $\Gamma = \bE\left[ \phi(X)\otimes\phi(X) \right]=\bE\left[\left(\Lambda^\top\psi(X)\right)\otimes\left(\Lambda^\top\psi(X)\right)\right] = \Lambda^\top\bE\left[\psi(X)\otimes\psi(X)\right]\Lambda$. By Courant-Fisher's max-min theorem, see for example, \cite[section 4.1]{vershynin_high-dimensional_2018}, together with Lemma~\ref{lemma:gram_schmidt_LRZ}, the eigenvalue of $\Gamma$ is equivalent to that of $\bE\left[\psi(X)\otimes\psi(X)\right]$ up to absolute constants depending on $L$, because $\max\{\norm{\Lambda}_{\text{op}},\norm{\Lambda^{-1}}_{\text{op}}\}\leq C_L$.

We now study the spectrum of $\bE[\psi(X)\otimes\psi(X)]$, which is up to constant, the spectrum of $\Gamma$. The spectrum of $\Gamma$ performs a ``multi-stage'' pattern up to $\alpha$. More precisely, by \eqref{eq:feature_map_polynomial_ONB}, we have the following.
\begin{itemize}
    \item When $\left|\cI\right|=0$. By \eqref{eq:feature_map_polynomial_ONB}, there is a unique eigenvalue $\alpha_0$, which is of constant order.
    \item When $\left|\cI\right|=1$. By \eqref{eq:feature_map_polynomial_ONB}, for every such $\cI$, $C_\cI=1$, thus the first stage of $\mathrm{spec}(\Gamma)$ is a plateau of height $1/d$ (up to $\alpha$), and of length $d$.
    \item When $\left|\cI\right|=2$. By \eqref{eq:feature_map_polynomial_ONB}, we know that $C_\cI=1$ for diagonal $\cI$ and $C_\cI=2$ for off-diagonal. Therefore the second stage is of height $1/d^2$ (up to constants depending on $\alpha$) and of length $d^2$.
    \item Generally, if $\left|\cI\right|=qd+r$ for $q\in\bN$ and $0\leq r<d$, the multi-nomial coefficient $C_\cI$ has maximum $\left|\cI\right|!/\left((q!)^{d-r}\left((q+1)!\right)^r\right)$. When $L<d$, we have $1\leq C_\cI\leq \left|\cI\right|!/2^{\left|\cI\right|}<\exp(\left|\cI\right|)$. As in \cite{liang_multiple_2020}, we ignore $C_\cI$ (this happens when $e^L<d$, which is guaranteed by Assumption~\ref{assumption:polynomial_feature_sub_gaussian}). The constraint mentioned above may be removed through careful computation; however, we have chosen to avoid doing so in order to give priority to simplicity. In fact, as $\left|\cI\right|$ increases, each stage's slope rises rapidly (with an exponential rate), because of $C_\cI$. For instance, denoting $\left|\cI\right|=l$, then the $l$-th stage of $\mathrm{spec}(\Gamma)$ has length $d^l$ and drops from $l!/(2^l\cdot d^l)$ to $1/d^l$. 
\end{itemize}

\subsubsection{Proof of Theorem~\ref{theo:polynomial_HS}}\label{sec:proof_polynomial_HS}

In this section, we prove Theorem~\ref{theo:polynomial_HS}. We will introduce the notation for tensor products and differential calculus used in the proof of Theorem~\ref{theo:polynomial_HS} in Section~\ref{sec:notation_tensor} and Section~\ref{sec:differential_calculus}, respectively. Subsequently, in Section~\ref{sec:proof_polynomial_HS_2}, we will employ these notations and lemmas to prove Theorem~\ref{theo:polynomial_HS}.

\subsubsection{Tensors}\label{sec:notation_tensor}

In algebra, the tensor product of multiple vector spaces is used to deal with multilinear mappings on these vector spaces. Given vectors\footnote{In this paper, we only consider real vector spaces.} spaces $U,V,W$, we write $U\otimes V$ as the tensor product between $U$ and $V$, with the correspondence given by $(\vu,\vv)\mapsto\vu\otimes\vv$. We also write $L(U;V)$ as the space of linear map from $U$ to $V$, which is itself a vector space. The tensor product provides the following isomorphism: $L(U;L(V;W))\cong L(U\otimes V; W)$, where ``$\cong$'' is algebraic isomorphism, see \cite[Chapter 16, Section 2, pp.607]{lang_algebra_2002}. Tensor product has associativity, that is, $U\otimes V\otimes W \cong (U\otimes V)\otimes W \cong U\otimes (V\otimes W)$, see \cite[Proposition 1.7, Proposition 1.8]{yokonuma_tensor_1992}. The following property of tensor product is frequently utilized in the remaining part of this section. Let $V_1,V_2$, $W_1,W_2$ and $U_1, U_2$ be real vector spaces, let $F_1, G_1\in L(V_1;W_1)$, $F_2, G_2 \in L(V_2;W_2)$, $H_1\in L(W_1;U_1)$ and $H_2\in L(W_2;U_2)$. Then $(H_1\circ F_1) \otimes (H_2\circ F_2) = (H_1\otimes H_2)\circ (F_1\otimes F_2)$, see \cite[Proposition 1.9]{yokonuma_tensor_1992}. Given a tensor $T$, and some $k\geq 1$, we write $T^{\otimes k}$ as $T\otimes\cdots\otimes T$ ($k$ times).

In statistics and computer science, it is common to select a basis and define tensors by enumerating their elements. Given $k\in\bN_+$ and $d_1,d_2,\cdots,d_k \in\bN_+$, a tensor of degree $k$, dimension $d_1, d_2,\cdots,d_k=d$ is an array $T = (t_{i_1,\cdots,i_k})$, where $\{t_{i_1,\cdots,i_k}\}'s$ are in $\bR$ (we only consider real field in this paper). If $d_1=d_2=\cdots=d_k$, we call $T$ as a tensor of degree $k$, dimension $d$. A degree $k$ dimension $d$ tensor $T$ is called symmetric, if it is unchanged upon permuting the indices. The space of symmetric degree $k$, dimension $d$ tensors, denoted as $\mathrm{Sym}_k(\bR^d)$. Given a tensor $T\in\mathrm{Sym}_k(\bR^d)$, we can always unfold $T$ into $\bR^{d^k}$. Given two matrices $A\in\bR^{m\times n}$, $B\in\bR^{p\times q}$ for some $m,n,p,q\in\bN_+$, $(A\otimes B)^\top = (A^\top)\otimes (B^\top)$, \cite[Proposition 1.9']{yokonuma_tensor_1992}.

Given two tensors $A=(a_{i_1,\cdots,i_k})$ and $B=(b_{i_1,\cdots,i_k})$ of degree $k$ and dimension $d$, the inner product between $A$ and $B$ is defined as:
\begin{align*}
    \left<A,B\right> = \sum_{i_1,\cdots,i_k=1}^d a_{i_1,\cdots,i_k} b_{i_1,\cdots,i_k},
\end{align*}which is equivalent to the inner product of the two vectors obtained by unfolding $A$ and $B$ in $\bR^{d^k}$. We let $\norm{A}_{HS}=\sqrt{\left<A,A\right>}$ as the Hilbert-Schmidt norm of $A$.

Given a tensor $T\in\mathrm{Sym}_k(\bR^d)$, whose entries are random variables on $\bR$, the expectation of $T$ is defined as a tensor $\bE T=(\bE t_{i_1,\cdots,i_k})\in\mathrm{Sym}_k(\bR^d)$.

Given $r,d\in\bN_+$ and $A$ be a degree $r$, dimension $d$ tensor. Let $P_r$ denote the set of partitions of $[r]$ into non-empty, pairwise disjoint sets, and for each partition $\cJ\in P_r$, $\cJ=\left\{J_1,\cdots,J_k\right\}$, we define
\begin{equation*}
    \norm{A}_\cJ = \sup\left( \sum_{\cI\in[d]^r}a_\cI \prod_{l=1}^k x_{\cI_{J_l}}^{(l)}:\, \sqrt{\sum_{\left|\cI_{J_l}\right|\leq d}\left(x_{\cI_{J_l}}^{(l)}\right)^2}\leq 1,\, 1\leq l\leq k \right),
\end{equation*}see also \cite{lehec_moments_2011} for an alternative notation by using injective tensor products.
For example, when $r=2$, $P_2 = \{\{1,2\},\{\{1\},\{2\}\}\}$, $a_\cI = a_{ij}$. For $J_1=\{1,2\}$, we have $k=l=1$, $\cI\in [d]^2$ and $\sqrt{\sum_{i,j\leq d}x_{ij}^2}=\norm{x}_{HS}$ for $x\in\bR^{d^2}$, thus $\norm{A}_{\left\{1,2\right\}}=\norm{A}_{HS}$. For $J_1 = \{1\}, J_2=\{2\}$, we have $k=2$ and for the sake of simplicity, we write $x_{\cI_1}^{(1)}$ and $x_{\cI_2}^{(2)}$ as $x,y$ , thus $\sqrt{\sum_{i\leq d}x_{i}^2}, \sqrt{\sum_{j\leq d}y_{j}^2}\leq 1$ , and $\norm{A}_{\left\{1\right\},\left\{2\right\}}=\norm{A}_{\text{op}}$. For $r=3$, we can write $[3]=\{1,2,3\}$, corresponding to view $A\in\bR^{d^3}$ as a linear functional with norm $\norm{A}_{\{1,2,3\}}=\sqrt{\sum_{i,j,k}a_{i,j,k}^2}$; or $\{1,2\},\{3\}$, corresponding to view $A\in\bR^{d^3}$ as a bi-linear functional on $\bR^{d^2}\times \bR$, with the norm $\norm{A}_{\{1,2\},\{3\}}=\sup\left(\left|\sum_{i,j,k}a_{i,j,k}x_{i,j}y_k\right|:\, \sum x_{i,j}^2\leq 1,\, \sum y_k^2\leq 1\right)$ ($\{1,3\},\{2\}$ or $\{2,3\},\{1\}$ are similar); or $\{1\},\{2\},\{3\}$, corresponding to view $A$ as a tri-linear functional on $\left(\bR^{d}\right)^3$. 
We have
\begin{equation*}
    \norm{A}_{\{1\},\{2\},\{3\}} \leq\norm{A}_{\{1,2\},\{3\}},\, \norm{A}_{\{1,3\},\{2\}} ,\, \norm{A}_{\{2,3\},\{1\}} \leq \norm{A}_{\{1,2,3\}}.
\end{equation*}If $\cK$ is a finer partition than $\cJ$ (any element in $\cK$ is contained in an element in $\cJ$), then $\norm{A}_\cK\leq\norm{A}_\cJ$, see \cite[Section 2.1]{lehec_moments_2011}. 

\subsubsection{Differential Calculus}\label{sec:differential_calculus}

Let $\cX,\cY$ be two normed vector spaces with norm $\norm{\cdot}_\cX$ and $\norm{\cdot}_\cY$, $\Omega\subseteq\cX$ be an open set, $f:\Omega\to\cY$ is Fréchet differentiable at $\vx\in\Omega$, if there exists a bounded linear operator $A_\vx\in L(\cX,\cY)$, such that
\begin{align}\label{eq:def_Frechet_derivative}
    \lim_{\vh\to \vzero} \frac{\norm{ f(\vx+\vh) - f(\vx) - A_\vx\vh }_\cY}{\norm{\vh}_\cX} = 0,
\end{align}where the limit is in the sense of the topology generated by $\norm{\cdot}_\cX$, see \cite[Definition 14.2.1]{albiac_topics_2016}. Denote $\D  f(\vx) = A_\vx$. In the case of $\cX=(\bR^d,\norm{\cdot}_2)$ and $\cY=(\bR,\norm{\cdot}_1)$, the Fréchet derivative at $\vx$ is simply the transpose of the gradient of $f$ at $\vx$. If we omit $\vx$, then $\D  f:\vx\in\cX\mapsto A_\vx \in L(\cX,\cY)$. Given $g:\cZ\to\Omega\subseteq\cX$, the chain rule gives that $\D  (f\circ g)(\vz) = \D   f(g(\vz))\D  g(\vz) \in L(\cZ,\cY)$. We omit the variables and write $\D  (f\circ g) = (\D   f\circ g)(\D  g)$ and $\D  (fg) = (\D  f)g + f(\D  g)$. These notations make sense because $\D  g:\cZ\to L(\cZ,\cX)$ and $\D  f\circ g: \cZ\to L(\cX,\cY)$, hence $(\D  f\circ g)(\D  g):\cZ\to L(\cZ,\cY)$.

The second order Fréchet derivative is defined as the Fréchet derivative of $\D  f$, denoted as $\D  ^2 f:\cX\to L(\cX,L(\cX,\cY))\cong L(\cX\otimes\cX,\cY)$. That is to say, the second-order Fréchet derivative is a degree 2 tensor. If $\cX$ and $\cY$ are both (finite-dimensional) Euclidean spaces, and we represent $\D  ^2 f$ in the standard orthogonal basis of Euclidean space, we will obtain that $\D  ^2 f$ is actually the Hessian, that is, $(\D  ^2 f)_{ij}=\frac{\partial^2}{\partial x_i\partial x_j}f(\vx)$. Generally, for $k\geq 1$, $\D ^k f$ for $f$ defined on $\bR^d$ can be defined iteratively as $\D ^k f = \D (\D ^{k-1}f)$, resulting in a tensor in $\mathrm{Sym}_k(\bR^d)$, and if we represent $\D ^k f$ in the standard ONB of $\bR^d$, we have $(\D ^k f)_{i_1,\cdots,i_k}=\frac{\partial^k}{\partial x_{i_1}\cdots\partial x_{i_k}}f(\vx)$.

Providing a high-order chain rule is not easy. To see this, recall that we define $\D ^2 (f\circ g) = \D  (\D  (f\circ g))$. Applying the product rule gives that $\D ^2 (f\circ g)=\D (\D  f\circ g)\D  g + (\D  f\circ g)\D ^2 g$. However, if we use the chain rule for computing $\D  (\D  f\circ g) = (\D ^2 f\circ g)\D  g$, we face the issue that $\D  g:\cZ\to L(\cZ;\cX)$, however, $\D ^2 f\circ g:\cZ\to L(\cX\otimes\cX;\cY)$. This leads to its evaluation at $\vz$, that is, $((\D ^2 f\circ g)(\vz))(\D  g(\vz))$ not being well-defined.

To address this issue, we adopt the tensor notations developed by \cite{manton_differential_2013}, that is, 
\begin{align}\label{eq:chain_rule_1}
\begin{aligned}
    \D ( \D  f\circ g ) &= (\D ^2 f\circ g)(\D  g \otimes I_{\cZ\to L(\cX;\cX)}),\\
    \D ( (\D  f\circ g)(\D  g)) &= \D  (\D  f\circ g)(I_{\cZ\to L(\cZ;\cZ)}\otimes \D  g) + (\D  f\circ g)(\D ^2 g),\\
    \mbox{ and } \D  (f\otimes g) &= (\D  f \otimes g) + (f\otimes \D  g).
\end{aligned}
\end{align}
Let's verify that this notation makes sense. Recall that $\D ^2 f\circ g: \cZ\to L(\cX;L(\cX;\cY))\cong L(\cX\otimes\cX;\cY)$, and $\D  g:\cZ\to L(\cZ;\cX)$. Hence, $\D  g \otimes I_{\cZ\to L(\cX;\cX)}:\cZ\otimes\cZ\to L(\cZ;\cX)\otimes L(\cX;\cX)\cong L(\cZ\otimes\cX;\cX\otimes\cX)$. Therefore, the right-hand-side of the modified chain rule satisfies, $(\D ^2 f\circ g)(\D  g \otimes I_{\cZ\to L(\cX;\cX)}):\cZ\otimes\cZ\to L(\cZ\otimes\cX;\cY)\cong L(\cZ;L(\cX;\cY))$. The left-hand-side of the modified chain rule $\D ( \D  f\circ g):\cZ\to L(\cZ;L(\cX;\cY))$, coincides with the right-hand-side formula. As a result, the modified chain rule notation makes sense. The modified product rule follows from the same idea.

The chain rule of higher-order Fréchet derivatives can be conveniently represented using the tensor product language introduced in \eqref{eq:chain_rule_1} as follows:
\begin{align}
    \D  ^2 (f\circ g) &= \D  (\D   (f\circ g)) = \D  \left( (\D  f\circ g) (\D  g) \right) = \D  (\D  f\circ g)(I_{\cZ\to L(\cZ;\cZ)}\otimes \D  g) + (\D  f\circ g)(\D  ^2 g)\notag\\
    &=(\D  ^2 f\circ g)(\D  g\otimes I_{\cZ\to L(\cX;\cX)})(I_{\cZ\to L(\cZ;\cZ)}\otimes \D  g) + (\D  f\circ g)(\D  ^2 g)\notag\\
    &=(\D  ^2 f\circ g)(\D  g\otimes \D  g) + (\D   f\circ g)(\D  ^2 g)\label{eq:chain_rule_second_order}.
\end{align}Since the domain and codomain of the identity map are unambiguously understood, we henceforth omit the subscripts for the identity map.

Let $(\vv_i)_{i=1}^\infty$ are an ONB of $\ell_2$, define $\leq$ be a partial order on $\ell_2$ defined by $\vx\leq\vy\Leftrightarrow \left<\vx,\vv_i\right>\leq \left<\vy,\vv_i\right>$ for any $i\in\bN_+$. We say $\vx\lesssim \vy$ if there exists an absolute constant $C$ such that $\vx\leq C\vy$. We say $\vx\lesssim_K \vy$ if the constant $C$ depends on $K$. Let us prove the lemma below.
\begin{Lemma}\label{lemma:HS_norm_higher_order_chain_rule}
    Let $d,p,q\geq 1$, $K\in\bN_+$, $K\geq 3$. Suppose $g:\bR^d\to\bR^p$ and $f:\bR^p\to\bR^q$ be $C^K$ functions. Given any $n, k\in[K]$ and $i_1,\cdots,i_k\geq 0$, we view every tensor of the form $(\D ^n f\circ g)(\D ^{i_1}g\otimes\cdots\otimes\D ^{i_k}g)$ as orthonormal vectors $\vv_{(i_1,\cdots,i_k;n)}$ in $\ell_2$, that is, given $k,l,m,n\in\bN_+$, $i_1,\cdots,i_k, i_1',\cdots,i_l'\geq 0$, then $\left<\vv_{(i_1,\cdots,i_k;n)},\vv_{(i_1',\cdots,i_l';m)}\right>=\1_{\{(i_1,\cdots,i_k;n)=(i_1',\cdots,i_l';m)\}}$. In particular, if $l\neq k$ or $n\neq m$, then $\left<\vv_{(i_1,\cdots,i_k;n)},\vv_{(i_1',\cdots,i_l';m)}\right>=0$.
    There then exists an absolute constant $\nC\label{C_concentration}$ depending only on $K$, such that
\begin{align}\label{eq:chain_rule_higher_order}
\begin{aligned}
    \oC{C_concentration}^{-1}\D ^K(f\circ g) &\leq (\D ^K f\circ g)(\D  g)^{\otimes K} + \sum_{k=2}^{K-1} (\D  ^k f\circ g)\left(\sum_{\substack{i_1,\cdots,i_k\geq 1\\ i_1+\cdots+i_k=K}} \D ^{i_1} g \otimes \cdots\otimes \D ^{i_k} g \right)\\
    &  + (\D  f\circ g)\D ^K g.
\end{aligned}
\end{align}
\end{Lemma}

The overall strategy for proving Lemma~\ref{lemma:HS_norm_higher_order_chain_rule} is based on the induction method. We first prove the following lemma. To avoid certain pathological cases, in the following text, we always assume that the functions appearing are sufficiently many times Fréchet differentiable (for our purposes, this assumption always holds).

\begin{Lemma}\label{lemma:diff_tensor_power}
     Given any $k\geq 1$, for any normed vector spaces $\cZ,\cX$ and $g:\cZ\to\cX$ be $(k+1)$-times Fréchet differentiable on some open set $\Omega\subseteq\cZ$. We have
    \begin{align}\label{eq:diff_tensor_product}
        \D  \left((\D  g)^{\otimes k}\right) &= \sum_{\substack{i_1,\cdots,i_k\geq 1\\ i_1+\cdots+i_k=k+1}} \D ^{\otimes i_1}g \otimes \cdots\otimes \D ^{\otimes i_k} g.
    \end{align}
    Moreover, given any $(i_1,\cdots,i_k)$ satisfying $i_1,\cdots,i_k\geq 1$ and $i_1+\cdots+i_k=k+1$, we denote $(i_1,\cdots,i_k)+1$ as
    \begin{align}\label{eq:def_i_1_i_k+1}
    \begin{aligned}
        (i_1,\cdots,i_k)+1 := \bigcup_{j\in[k]} (i_1,\cdots,i_{j-1}, i_j+1, i_{j+1},\cdots,i_k),
    \end{aligned}
    \end{align}then
    \begin{align}\label{eq:diff_tensor_product_2}
        \D \left(\D ^{\otimes i_1}g \otimes \cdots\otimes \D ^{\otimes i_k} g\right) = \sum_{(l_1,\cdots,l_k)\in (i_1,\cdots,i_k)+1} \D ^{\otimes l_1} g \otimes \cdots \otimes \D ^{\otimes l_k} g.
    \end{align}
\end{Lemma}
\beginproof When $k=1$, $\D (\D  g)=\D ^2 g$. Suppose \eqref{eq:diff_tensor_product} holds for some $K\geq 1$, by \eqref{eq:chain_rule_1} we have
\begin{align}
    \D ((\D  g)^{\otimes (K+1)}) &= \D \left((\D  g)^{\otimes K}\otimes \D  g\right) = \D (\D  g)^{\otimes K}\otimes \D  g + (\D  g)^{\otimes K}\otimes \D ^2 g\notag\\
    &=\left( \sum_{\substack{i_1,\cdots,i_K\geq 1\\ i_1+\cdots+i_K=K+1}}  \D ^{\otimes i_1}g \otimes \cdots\otimes \D ^{\otimes i_K} g \otimes \D  g\right) + (\D  g)^{\otimes K}\otimes \D ^2 g\label{eq:diff_tensor_1}\\
    &=\sum_{\substack{i_1,\cdots,i_{K+1}\geq 1\\ i_1+\cdots+i_{K+1}=K+2}} \D ^{\otimes i_1}g \otimes \cdots\otimes \D ^{\otimes i_{K+1}} g, \notag
\end{align}where the last equality is via the following observation: given $l\in\bN_+$ and $i_1+\cdots+i_l=l+1$ and $i_1,\cdots,i_l\geq 1$, by the Pigeonhole principle, there must be a unique $n(l)\in[l]$ such that $n(l)=2$. Now let $l=K+1$, we have:
\begin{enumerate}
    \item the first term in \eqref{eq:diff_tensor_1} corresponds to the case where $n(K+1)$ belongs to $[K]$, that is, we have: $i_1+\cdots+i_K + i_{K+1} = K+2$, $i_1,\cdots,i_K,i_{K+1}\geq 1$, and $i_{K+1}=1$. Therefore, $n(K+1)$ must necessarily appear among $[K]$, and
    \item the second term in \eqref{eq:diff_tensor_1} corresponds to the case where $n(K+1)=K+1$. That is, $i_1,\cdots,i_K,i_{K+1}\geq 1$, $i_1+\cdots+i_K+i_{K+1}=K+2$, and $i_1=\cdots=i_K=1$.
\end{enumerate}
\eqref{eq:diff_tensor_product} is therefore derived. 

Now we prove \eqref{eq:diff_tensor_product_2}. When $k=1$, $i_1=2$ and $(i_1)+1=(i_1+1)=(3)$. $\D (\D ^2 g)=\D ^3 g$. Hence \eqref{eq:diff_tensor_product_2} is valid when $k=1$. Suppose \eqref{eq:diff_tensor_product_2} is valid for some $K\geq 1$, then for any $(i_1,\cdots,i_K,i_{K+1})$ such that $i_1,\cdots,i_{K+1}\geq 1$ and $i_1+\cdots+i_{K+1}=K+1$, by \eqref{eq:chain_rule_1}, we have
\begin{align*}
    &\D \left(\D ^{\otimes i_1}g \otimes \cdots\otimes \D ^{\otimes i_K} g\otimes \D ^{i_{K+1}} g\right)\\
    &= \D \left(\D ^{\otimes i_1}g \otimes \cdots\otimes \D ^{\otimes i_K} g\right)\otimes\D ^{i_{K+1}}g + \D ^{\otimes i_1}g \otimes \cdots\otimes \D ^{\otimes i_K} g\otimes \D ^{i_{K+1}+1}g\\
    &=\left(\sum_{\substack{(l_1,\cdots,l_k)\in (i_1,\cdots,i_K)+1}}\D ^{\otimes l_1} g \otimes \cdots \otimes \D ^{\otimes l_K} g\otimes \D ^{i_{K+1}}g \right) + \D ^{\otimes i_1}g \otimes \cdots\otimes \D ^{\otimes i_K} g\otimes \D ^{i_{K+1}+1}g\\
    &=\sum_{(l_1,\cdots,l_{K+1})\in (i_1,\cdots,i_{K+1})+1} \D ^{\otimes l_1}g\otimes \cdots\otimes \D ^{\otimes l_{K+1}} g.
\end{align*}
\endproof

Having proved Lemma~\ref{lemma:diff_tensor_power}, we proceed to prove Lemma~\ref{lemma:HS_norm_higher_order_chain_rule}.

\paragraph{Proof of Lemma~\ref{lemma:HS_norm_higher_order_chain_rule}} We begin with $\D ^3(f\circ g)$.
\begin{align}\label{eq:chain_rule_third_order}
\begin{aligned}
    \D ^3 (f\circ g) &= \D \left((\D  ^2 f\circ g)(\D  g\otimes \D  g)\right) + \D  \left( (\D   f\circ g)(\D  ^2 g) \right)\\
    &= (\D (\D ^2 f\circ g))(I\otimes (\D  g)^{\otimes 2}) + (\D ^2 f\circ g)\D ((\D  g)^{\otimes 2}) + (\D (\D  f\circ g))(I\otimes \D ^2 g) + (\D  f\circ g)(\D ^3 g)\\
    &= (\D ^3 f\circ g)(\D  g\otimes I)(I\otimes (\D  g)^{\otimes 2}) + (\D ^2 f\circ g)(\D ^2 g\otimes \D  g + \D  g \otimes \D ^2 g)\\
    &+ (\D ^2 f\circ g)(\D  g \otimes I)(I\otimes \D ^2 g) + (\D  f\circ g)(\D ^3 g)\\
    &= (\D ^3 f\circ g)(\D  g)^{\otimes 3} + (\D ^2 f\circ g)( \D ^2 g\otimes \D  g + \D  g\otimes \D ^2 g + \D  g\otimes \D ^2 g ) + (\D  f \circ g)\D ^3 g.
\end{aligned}
\end{align}
All terms in \eqref{eq:chain_rule_third_order} of the form $(\D ^n f \circ g)(\D ^{i_1}g \otimes \cdots \otimes \D ^{i_k}g)$ can be regarded as orthonormal vectors in $\ell_2$: $(\D ^3 f\circ g)(\D  g)^{\otimes 3}$ is viewed as $\vv_{(1,1,1;3)}$, $(\D ^2 f\circ g)(\D ^2 g\otimes \D  g)$ is viewed as $\vv_{(2,1;2)}$, $(\D ^2 f\circ g)(\D  g\otimes \D ^2 g)$ is viewed as $\vv_{(1,2;2)}$, $(\D  f \circ g)\D ^3 g$ is viewed as $\vv_{(3;1)}$. If we consider $\D ^3 (f\circ g)$ as a vector in $\ell_2$ and expand it in this ONB, the coordinates of $\D ^3 (f\circ g)$ are (we ignore terms with coordinates equal to $0$): $(1,1,2,1)$. The coordinates of the right-hand side of \eqref{eq:chain_rule_higher_order} in this ONB when $K=3$ are $(1,1,1,1)$. Therefore, when $K=3$ and $\oC{C_concentration}=2$, \eqref{eq:chain_rule_higher_order} holds.

Suppose \eqref{eq:chain_rule_higher_order} is satisfied with some $K\geq 3$. Then
\begin{align*}
&\oC{C_concentration}^{-1}\D^{K+1}(f\circ g) \leq \D\left((\D ^K f\circ g)(\D  g)^{\otimes K}\right)\\
&+\D\left(\sum_{k=2}^{K-1} (\D  ^k f\circ g)\left(\sum_{\substack{i_1,\cdots,i_k\geq 1\\ i_1+\cdots+i_k=K}} \D ^{i_1} g \otimes \cdots\otimes \D ^{i_k} g \right)\right)\\
&+\D\left((\D  f\circ g)\D ^K g\right).
\end{align*}
By \eqref{eq:chain_rule_1} and \eqref{eq:diff_tensor_product} we have:
\begin{align}\label{eq:proof_chain_rule_higer_1}
    \begin{aligned}
        \D \left( (\D ^K f\circ g)(\D  g)^{\otimes K} \right) &= \D (\D ^K f\circ g)(I\otimes (\D  g)^{\otimes K}) + (\D ^K f\circ g)\D (\D  g)^{\otimes K}\\
    &= (\D ^{K+1} f\circ g)(\D  g\otimes I)(I\otimes (\D  g)^{\otimes K})\\
    &+ (\D ^K f\circ g)\sum_{\substack{i_1,\cdots,i_K\geq 1\\ i_1+\cdots+i_K=K+1}} \D ^{\otimes i_1}g \otimes \cdots\otimes \D ^{\otimes i_K} g\\
    &=(\D ^{K+1} f\circ g)(\D  g)^{\otimes (K+1)}\\
    &+  (\D ^K f\circ g)\sum_{\substack{i_1,\cdots,i_K\geq 1\\ i_1+\cdots+i_K=K+1}} \D ^{\otimes i_1}g \otimes \cdots\otimes \D ^{\otimes i_K} g, 
    \end{aligned}
\end{align}and
\begin{align}\label{eq:proof_chain_rule_higer_2}
    \D \left((\D  f\circ g)(\D ^K g)\right) &= \D (\D  f\circ g)(I\otimes \D ^K g) + (\D  f\circ g)\D ^{K+1}g\notag\\
    &= (\D ^2 f\circ g)(\D  g\otimes I)(I\otimes \D ^K g) + (\D  f\circ g)\D ^{K+1}g \notag\\
    &= (\D ^2 f\circ g)(\D  g\otimes \D ^K g) + (\D  f\circ g)\D ^{K+1}g.
\end{align}
We are left with the middle term in \eqref{eq:chain_rule_higher_order}. For any $2\leq k\leq K-1$, by \eqref{eq:chain_rule_1}, we have:
\begin{align}
    &\D \left( (\D  ^k f\circ g)\left(\sum_{\substack{i_1,\cdots,i_k\geq 1\\ i_1+\cdots+i_k=K}} \D ^{i_1} g \otimes \cdots\otimes \D ^{i_k} g \right)\right)\notag \\
    &= \D (\D ^k f\circ g)\left(I\otimes \left(\sum_{\substack{i_1,\cdots,i_k\geq 1\\ i_1+\cdots+i_k=K}} \D ^{i_1} g \otimes \cdots\otimes \D ^{i_k} g \right) \right)\label{eq:proof_chain_rule_higer_3}\\
    &+ (\D ^k f\circ g)\D  \left(\sum_{\substack{i_1,\cdots,i_k\geq 1\\ i_1+\cdots+i_k=K}} \D ^{i_1} g \otimes \cdots\otimes \D ^{i_k} g \right)\label{eq:proof_chain_rule_higer_4}.
\end{align}
We now separately compute \eqref{eq:proof_chain_rule_higer_3} and \eqref{eq:proof_chain_rule_higer_4}. By \eqref{eq:chain_rule_1} we have:
\begin{align}
    \eqref{eq:proof_chain_rule_higer_3} &= (\D ^{k+1} f\circ g)(\D  g\otimes I)\left(I\otimes \left(\sum_{\substack{i_1,\cdots,i_k\geq 1\\ i_1+\cdots+i_k=K}} \D ^{i_1} g \otimes \cdots\otimes \D ^{i_k} g \right) \right)\notag\\
    &=(\D ^{k+1} f\circ g) \left(\sum_{\substack{i_1,\cdots,i_k\geq 1\\ i_1+\cdots+i_k=K}} \D  g\otimes \D ^{i_1} g \otimes \cdots\otimes \D ^{i_k} g \right)\notag\\
    &=(\D ^{k+1} f\circ g) \left(\sum_{\substack{i_1,\cdots,i_{k+1}\geq 1\\ i_1+\cdots+i_{k+1}=K+1 \\ i_1=1}} \D ^{i_1} g \otimes \cdots\otimes \D ^{i_{k+1}} g \right). \label{eq:proof_chain_rule_high_6}
\end{align}By \eqref{eq:diff_tensor_product_2} we have:
\begin{align}
    \eqref{eq:proof_chain_rule_higer_4} &= (\D ^k f\circ g)\sum_{\substack{i_1,\cdots,i_k\geq 1\\ i_1+\cdots+i_k=K}}\D \left( \D ^{i_1} g \otimes \cdots\otimes \D ^{i_k} g
 \right)\notag\\
 &=(\D ^k f\circ g)\sum_{\substack{i_1,\cdots,i_k\geq 1\\ i_1+\cdots+i_k=K}}\sum_{(l_1,\cdots,l_k)\in(i_1,\cdots,i_k)+1}\D ^{\otimes l_1}g\otimes \cdots\otimes \D ^{\otimes l_k}g. \label{eq:proof_chain_rule_higer_5}
\end{align}By the definition of $(i_1,\cdots,i_k)+1$, see \eqref{eq:def_i_1_i_k+1}, $l_1+\cdots+l_k=i_1+\cdots+i_k+1=K+1$. This implies that every term $\D ^{\otimes l_1}g\otimes \cdots\otimes \D ^{\otimes l_k}g$ appeared in \eqref{eq:proof_chain_rule_higer_5} is contained in one term of the following sum
\begin{align}\label{eq:middle_term_K+1_2}
    \sum_{\substack{i_1'\cdots,i_k'\geq 1\\ i_1'+\cdots+i_k'=K+1}} \D ^{i_1'} g \otimes \cdots\otimes \D ^{i_k'} g.
\end{align}
We therefore only need to compute the repetition counts. Given any $i_1',\cdots,i_k'\geq 1$ such that $i_1'+\cdots+i_k'=K+1$, then there are $\left|\left\{ j\in[k]: i_j'\geq 2 \right\}\right|$ choices of $(i_1,\cdots,i_k)$ in \eqref{eq:proof_chain_rule_higer_5} that satisfy $l_1=i_1',\cdots,l_k=i_k'$. For example, with $K=7$ and $k=4$, consider a term $(i_1',i_2',i_3',i_4')=(3,1,3,1)$ from \eqref{eq:middle_term_K+1_2} (this term is in \eqref{eq:middle_term_K+1_2} because $3+1+3+1=K+1$). Then, only $(i_1,i_2,i_3,i_4) = (2,1,3,1)$ and $(i_1,i_2,i_3,i_4) = (3,1,2,1)$ from $\eqref{eq:proof_chain_rule_higer_5}$ may be transformed into $(3,1,3,1)$. In other words, given any $i_1',\cdots,i_k'\geq 1$ such that $i_1'+\cdots+i_k'=K+1$, if we compute the multiplicity, then $(\D ^{i_1'}g\otimes\cdots\otimes\D ^{i_k'}g)$ appears $\left|\{j\in[k]: i_j'\geq 2\}\right|$ times in the summation of \eqref{eq:proof_chain_rule_higer_5}. We have $\left|\{j\in[k]: i_j'\geq 2\}\right|\leq k< K$ uniformly over all $i_1',\cdots,i_k'\geq 1$, $i_1'+\cdots+i_k'=K+1$ and $k<K$.

Therefore, if we view each $(\D ^k f\circ g) (\D ^{\otimes l_1}g\otimes \cdots\otimes \D ^{\otimes l_k}g)$ in \eqref{eq:proof_chain_rule_higer_5} as a vector $\vv_{(l_1,\cdots,l_k;k)}$ in $\ell^2$, and define a partial order on $\ell^2$ as $\vx \leq \vy$ if $x_i \leq y_i$ for any $i \in \mathbb{N}_+$, then
\begin{align}\label{eq:proof_chain_rule_high_7}
    \eqref{eq:proof_chain_rule_higer_5} \leq K(\D ^k f\circ g)\sum_{\substack{i_1,\cdots,i_k\geq 1\\ i_1+\cdots+i_k=K+1}} \D ^{i_1} g \otimes \cdots\otimes \D ^{i_k} g.
\end{align}
Combining \eqref{eq:proof_chain_rule_higer_1}, \eqref{eq:proof_chain_rule_higer_2}, \eqref{eq:proof_chain_rule_higer_3}, \eqref{eq:proof_chain_rule_higer_4}, \eqref{eq:proof_chain_rule_high_6} and \eqref{eq:proof_chain_rule_high_7} indicates
\begin{align*}
    \D ^{K+1}(f\circ g) &\lesssim_K (\D ^{K+1} f\circ g)(\D  g)^{\otimes {K+1}}\\
    &+ \sum_{k=2}^{K} (\D  ^k f\circ g)\left(\sum_{\substack{i_1,\cdots,i_k\geq 1\\ i_1+\cdots+i_k=K+1}} \D ^{i_1} g \otimes \cdots\otimes \D ^{i_k} g \right) \\
    &+ (\D  f\circ g)\D ^{K+1} g.
\end{align*}
Lemma~\ref{lemma:HS_norm_higher_order_chain_rule} follows by induction.

\subsubsection{Continue: Proof of Theorem~\ref{theo:polynomial_HS}}\label{sec:proof_polynomial_HS_2}

The proof of Theorem~\ref{theo:polynomial_HS} relies on the following theorem, taken from \cite{adamczak_concentration_2015}.
\begin{Theorem}\label{theo:adamczak_wolff_polynomial}
    Let $X=(x_j)_{j=1}^d$ be a random vector with independent coordinates, such that $\max(\|x_j\|_{\psi_2}:\, j\in[d])<\infty$. Then for every polynomial $g:\bR^d\to\bR$ of degree $D$ and every $q\geq 2$,
    \begin{equation*}
        \norm{g(X)-\bE g(X)}_{L_q} \leq C_D\sum_{r=1}^D \norm{x_1}_{\psi_2}^r \sum_{\cJ\in P_r}q^{\left|J\right|/2}\norm{\bE\D ^r g(X)}_\cJ.
    \end{equation*}As a consequence, for any $t>0$,
    \begin{equation*}
        \bP\left(\left|g(X)-\bE g(X)\right|\geq t\right)\leq 2\exp\left(-\frac{1}{C_D}\eta_g(t)\right),
    \end{equation*}where $C_D$ is an absolute constant depending on $D$ and
    \begin{equation*}
        \eta_g(t) = \min\left(\min\left(\left(\frac{t}{\norm{x_1}_{\psi_2}^r\norm{\bE \D ^r g(X)}_\cJ}\right)^{2/\left|\cJ\right|}:\, \cJ\in P_r\right):\, 1\leq r\leq D\right).
    \end{equation*}
\end{Theorem}

Due to Theorem~\ref{theo:adamczak_wolff_polynomial}, we only need to apply Theorem~\ref{theo:adamczak_wolff_polynomial} to $g(\cdot)=K(\cdot,\cdot)$, which is a $2L$-degree polynomial. Consequently, we only need to estimate $\eta_g(t)$ from below. Given that $\cJ$ is a partition of the set $[r]$, we will employ the straightforward upper bound $\norm{\bE\D ^r g(X)}_\cJ\leq\norm{\bE\D ^rg(X)}_{HS}$ uniformly for all $\cJ$ belonging to the set $P_r$ and for all values of $r\in[2L]$. For any $r\in[2L]$, we have:
\begin{align}\label{eq:proof_polynomial_HS_1}
    \max\left( \norm{\bE\D ^r K(X,X)}_{\cJ}:\, \cJ\in P_r\right)\leq \norm{\bE\D ^r K(X,X)}_{HS}\leq \sum_{\iota=1}^L \norm{\bE\D ^r \frac{\gamma_\iota}{d^\iota}(\norm{X}_2^2)^\iota}_{HS}.
\end{align}
Given any $r\in[2L]$, it suffices to bound $\norm{\bE\D ^r (\norm{X}_2^2)^\iota}_{HS}$ from above.

\paragraph{Warm up: the cases for $r=1$ and $r=2$} We first deal with the cases where $r=1$ and $r=2$, that is, we bound $\norm{\bE\D (\norm{X}_2^2)^\iota}_{HS}$ and $\norm{\bE\D ^2(\norm{X}_2^2)^\iota}_{HS}$ from above.
\begin{enumerate}
    \item $\norm{\bE\D (\norm{X}_2^2)^\iota}_{HS}$. For any $j\in[d]$, we have
    \begin{align*}
        &(\bE\D (\norm{X}_2^2)^\iota)_j = 2\iota\bE (\norm{X}_2^2)^{\iota-1}x_j\\
        &= 2\iota\bE\left[ \sum_{\substack{n_1+\cdots+n_d=\iota-1\\ \forall i\in[d], n_i\geq 0}} \frac{(\iota-1)!}{n_1!\cdot\cdots\cdot n_d!} x_1^{2n_1}\cdot\cdots\cdot x_{j-1}^{2n_{j-1}} x_j^{2n_j +1} x_{j+1}^{2n_{j+1}}\cdot\cdots\cdot x_d^{2n_d} \right].
    \end{align*}As we have assumed that $d>\iota$, given any $n_1+\cdots+n_d=\iota-1$ and $n_1,\cdots,n_d\geq 0$, there exists $i^*\in[d]$ such that $n_{i^*}=0$. If $j=i^*$, by the independence of the coordinates and the mean zero assumption, $(\bE\D (\norm{X}_2^2)^\iota)_j=0$. Given the $j$-th component of this first-order Fréchet derivative, in the summation above, only the terms corresponding to $n_j \geq 1$ contributes.
    By \cite[Proposition 1.5]{jukna_extremal_2011} we have
\begin{align*}
    &\frac{1}{2\iota}\left(\bE\D  (\norm{X}_2^{2})^\iota\right)_j =  \bE\left[\sum_{\substack{n_1+\cdots+n_d=\iota-1\\0\leq n_1,\cdots,n_{j-1}, n_{j+1},\cdots,n_d\leq \iota-1 \\ n_j\geq 1}}^d \frac{(\iota-1)!}{n_1!\cdot\cdots\cdot n_d!} x_{1}^{2n_1}\cdot\cdots\cdot x_{d}^{2n_d} x_j\right]\\
    &=\sum_{b=1}^{\iota-1}\sum_{\substack{n_1+\cdots+n_{j-1}+n_{j+1}+\cdots+n_d=\iota-b-1\\ 0\leq n_1,\cdots,n_{j-1},n_{j+1},\cdots,n_d\leq \iota-1 }} \frac{(\iota-1)!}{n_1!\cdot\cdots\cdot n_d!} \bE [x_{1}^{2n_1}]\cdot\cdots\cdot \bE[x_{{j-1}}^{2n_{j-1}}]\cdot \bE[x_{j}^{2b+1}] \cdot\bE[x_{{j+1}}^{2n_{j+1}}] \cdot\cdots\cdot \bE[x_{d}^{2n_d}]\\
    &\leq (\iota-1)!(2(\iota-2))^{\iota-\frac{3}{2}}\norm{x_1}_{\psi_2}^{2\iota-1}\sum_{b=1}^{\iota-1}\binom{d-1+\iota-b-2}{\iota-b-1}\\
    &<(\iota-1)!\norm{x_1}_{\psi_2}^{2\iota-1}(2(\iota-2))^{\iota-\frac{3}{2}}\sum_{b=1}^{\iota-1}\binom{2d}{\iota-b-1}\\
    &<(\iota-1)!(2(\iota-2))^{\iota-\frac{3}{2}}\norm{x_1}_{\psi_2}^{2\iota-1}\left(\frac{2ed}{\iota-1}\right)^{\iota-1}.
\end{align*}We conclude that
\begin{align}\label{eq:K=1}
   \frac{ \norm{\gamma_\iota\bE\D (\norm{X}_2^2)^\iota}_{HS}^2}{d^{2\iota}} &< 4\iota^2 \gamma_\iota^2 ((\iota-1)!)^2 (2(\iota-2))^{2\iota-3} \norm{x_1}_{\psi_2}^{4\iota-2}d^{-2\iota}\left(\frac{2ed}{\iota-1}\right)^{2(\iota-1)}\notag\\
   &=\frac{4\iota^2 \gamma_\iota^2 ((\iota-1)!)^2 (2(\iota-2))^{2\iota-3} \norm{x_1}_{\psi_2}^{4\iota-2}(2e)^{2(\iota-1)}}{(\iota-1)^{2(\iota-1)}d^2}\notag\\
   &\sim \frac{4\iota^2 \gamma_\iota^2  (2(\iota-2))^{2\iota-3} \norm{x_1}_{\psi_2}^{4\iota-2}}{d^2}.
\end{align}
\item $\norm{\bE\D ^2(\norm{X}_2^2)^\iota}_{HS}$. By \eqref{eq:chain_rule_second_order}, we know that $\bE\D ^2(\norm{X}_2^2)^\iota = \iota(\iota-1)(\norm{X}_2^2)^{\iota-2}((X^\top)\otimes(X^\top)) + 2\iota(\norm{X}_2^2)^{\iota-1}I_d$. Take the $\norm{\cdot}_{HS}$ and by triangular inequality we obtain that
\begin{align*}
    \norm{\bE\D ^2(\norm{X}_2^2)^\iota}_{HS} &\leq \iota(\iota-1)\norm{\bE\left[(\norm{X}_2^2)^{\iota-2}((X^\top)\otimes(X^\top))\right]}_{HS} + 2\iota \sqrt{d} \bE(\norm{X}_2^2)^{\iota-1}\\
    &=\iota(\iota-1)\norm{\bE\left[(\norm{X}_2^2)^{\iota-2}(X\otimes X)\right]}_{HS} + 2\iota \sqrt{d} \bE(\norm{X}_2^2)^{\iota-1}.
\end{align*}We have: $\sqrt{d}\bE(\norm{X}_2^2)^{\iota-1}\lesssim (2(\iota-1))^{\iota-1} \norm{x_1}_{\psi_2}^{2(\iota-1)} d^{\iota-1/2}$.
For any $j_1,j_2\in[d]$, we discuss the following two cases separately:
\begin{enumerate}
    \item $j_1\neq j_2$. We have:
    \begin{align}\label{eq:proof_K=2_1}
        \left(\bE\left[(\norm{X}_2^2)^{\iota-2}( X\otimes X)\right]\right)_{j_1,j_2} = \bE\left[ x_{j_1}x_{j_2} \sum_{\substack{ n_1+\cdots+n_d=\iota-2\\ \forall i\in[d], n_i\geq 0 }} \frac{(\iota-2)!}{n_1!\cdot\cdots\cdot n_d!} x_1^{2n_1}\cdot\cdots\cdot x_d^{2n_d} \right].
    \end{align}Repeat the argument in case \emph{1.}, we obtain that
    \begin{align*}
        \eqref{eq:proof_K=2_1}&= \sum_{\substack{n_1+\cdots+n_d=\iota-2\\ \forall i\in[d], n_i\geq 0\\ n_{j_1}n_{j_2}\geq 1 }} \frac{(\iota-2)!}{n_1!\cdot\cdots\cdot n_d!}\prod_{\substack{\alpha\in[d]\\ \alpha\neq j_1, j_2}} \bE[x_\alpha^{2n_\alpha}] \bE [x_{j_1}^{2n_{j_1}+1}]\bE[x_{j_2}^{2n_{j_2}+1}]\\
        &\lesssim (\iota-2)!(2\iota-1)^{\iota-1}\norm{x_1}_{\psi_2}^{2\iota-2} \left(\frac{2ed}{\iota-2}\right)^{\iota-2}.
    \end{align*}
    Therefore,
    \begin{align*}
        \sum_{\substack{j_1,j_2\in[d]\\ j_1\neq j_2}}\left(\bE\left[(\norm{X}_2^2)^{\iota-2}( X\otimes X)\right]\right)_{j_1,j_2}^2 \lesssim  (2\iota-1)^{2(\iota-1)}\norm{x_1}_{\psi_2}^{4\iota-4}  d^{2\iota-2}.
    \end{align*}
    \item $j_1=j_2$. In this scenario, we cannot employ the method of $j_1 \neq j_2$ to obtain sparsity in this second-order gradient tensor. Fortunately, in this case, the Hilbert-Schmidt norm of this second-order tensor only possesses one ``degree of freedom'':
    \begin{align*}
        \sum_{j_1=1}^d\left(\bE\left[(\norm{X}_2^2)^{\iota-2}( X\otimes X)\right]\right)_{j_1,j_1}^2 &\lesssim \sum_{j_1=1}^d\left((2\iota-1)^{\iota-1}\norm{x_1}_{\psi_2}^{2\iota-2} d^{\iota-2}\right)^2\\
        &= (2\iota-1)^{2(\iota-1)}\norm{x_1}_{\psi_2}^{4\iota-4}d^{2\iota-3}.
    \end{align*}
\end{enumerate}
In summary,
\begin{align}\label{eq:K=2}
    \frac{\gamma_\iota^2\norm{\bE\D ^2(\norm{X}_2^2)^\iota}_{HS}^2}{d^{2\iota}}&\lesssim \gamma_\iota^2   (2\iota-1)^{2(\iota-1)}\norm{x_1}_{\psi_2}^{4\iota-4}  d^{-2}+\gamma_\iota^2(2\iota-1)^{2(\iota-1)}\norm{x_1}_{\psi_2}^{4\iota-4}d^{-3}\notag\\
    & + \gamma_\iota^2 (2(\iota-1))^{2(\iota-1)} \norm{x_1}_{\psi_2}^{4(\iota-1)} d^{-1} \sim \frac{\gamma_\iota^2}{d} (2\iota-1)^{2(\iota-1)} \norm{x_1}_{\psi_2}^{4\iota-4} .
\end{align}
\end{enumerate}

Before addressing the case of $r \geq 3$, let us first review the proof for $r=2$, as the main proof strategy for the case of $r \geq 3$ will be utilized from the case of $r=2$. The computation of $\mathbb{E}\D ^2 (\|X\|_2^2)^\iota$ can be broadly divided into two parts: corresponding to the first and last terms in equation \eqref{eq:chain_rule_higher_order}. For the computation of the first term, we divide it into two cases: in the off-diagonal case, we exploit the independence of coordinates of $X$ and its mean-zero property to induce sparsity; in the diagonal case, since there are fewer entries (referred to as "degrees of freedom") in the diagonal of the gradient tensor, we compute it directly. Next, in the case of $r \geq 3$, we will continue using this approach. The difference lies in the fact that for $r \geq 3$, we also need to consider the intermediate term in equation \eqref{eq:chain_rule_higher_order}.

\paragraph{The case for $3\leq r\leq 2\iota$} We temporarily substitute $r$ with $K$ for ease of use with Lemma~\ref{lemma:HS_norm_higher_order_chain_rule}.

For any $\iota\geq 1$, let $g:\vz\in\bR^d\mapsto \norm{\vz}_2^2\in\bR$ and $f:t\in\bR\mapsto t^\iota\in\bR$. Then $f\circ g:\vz\in\bR^d\mapsto \norm{\vz}_2^{2\iota}$. Moreover, $\D  g(\vz) = 2\vz^\top$, $\D ^2 g(\vz) = 2I_d$ and $\D ^3 g(\vz) = \vzero$, the zero tensor of degree $3$, dimension $d$. We also have that for any $k\in[\iota]$, $(\D ^k f\circ g)(\vz) = ((\iota!)/(\iota-k)!) \norm{\vz}_2^{2\iota-2k}$ and for $k>\iota$, $(\D ^k f\circ g)(\vz)=\vzero$. Let $\vz=X$ introduced in Theorem~\ref{theo:polynomial_HS}.

By Lemma~\ref{lemma:HS_norm_higher_order_chain_rule}, take expectation and Hilbert-Schmidt norm, and by triangular inequality, we obtain that there exists an absolute constant $\nC\label{C_Concentration_2}$ depending only on $K$ such that
\begin{align*}
    \begin{aligned}
        \norm{\bE\D ^K (\norm{X}_2^2)^\iota}_{HS}&\leq \oC{C_Concentration_2} \1_{\{K\leq \iota\}}\norm{ \bE \left[ (\norm{X}_2^2)^{\iota-K}(X^\top)^{\otimes K} \right] }_{HS} + \oC{C_Concentration_2} \norm{ \bE \left[ (\norm{X}_2^2)^{\iota-1}\D ^K \norm{X}_2^2 \right] }_{HS}\\
        &+ \oC{C_Concentration_2}\norm{ \sum_{k=2}^{K-1\wedge\iota}\bE\left[ (\norm{X}_2^2)^{\iota-k}\left( \sum_{\substack{i_1,\cdots,i_k\geq 1\\ i_1+\cdots+i_k=K}} \D ^{i_1}g\otimes \cdots \otimes \D ^{i_k}g \right)(X)   \right] }_{HS}.
    \end{aligned}
\end{align*}By triangular inequality again, we have
\begin{align}\label{eq:equivalent_k_derivative}
    \begin{aligned}
        &\norm{\bE\D ^K \norm{X}_2^{2\iota} }_{HS}\lesssim_{K} \1_{\{K\leq \iota\}}\norm{\bE \norm{X}_2^{2\iota-2K} (X^\top)^{\otimes K} }_{HS}\\
        &+ \sum_{k=2}^{K-1\wedge\iota}\norm{\bE \norm{X}_2^{2\iota-2k} \left( \sum_{\substack{i_1,\cdots,i_k\in\{1,2\}\\ i_1+\cdots+i_k = K}} \D ^{i_1} g\otimes \cdots\otimes \D ^{i_k}g \right) }_{HS}.
    \end{aligned}
\end{align}

We first deal with the last term in \eqref{eq:equivalent_k_derivative}.

\paragraph{For any $2\leq k\leq K-1\wedge\iota$}

Given $i_1,\cdots,i_k\in\{1,2\}$ such that $i_1+\cdots+i_k=K$, we can separate $[K]$ into $k$ groups, that is, $[K]=\sqcup_{a=1}^k J_a$, where $J_a$ is of the form $J_a=\{j,j+1\}$ or $J_a=\{j\}$ for $j\in[K]$ according to $i_a$ for $a\in[k]$, that is, we have $i_a = \left|J_a\right|$ equals either $1$ or $2$. Denote $\cJ = \cup_{a=1}^k\{J_a:\, \left|J_a\right|=2  \}$, and $\cJ^c=\cup_{a=1}^k\{J_a:\, \left|J_a\right|=1  \}$. As $\cJ\sqcup\cJ^c=[K]$, there exists a surjective $\pi:l\in[K]\mapsto \pi(l)\in[k]$ such that $l\in J_{\pi(l)}$. Given any $\alpha_1,\cdots,\alpha_K$, denote $\cI = \cup_{j=1}^{K-1} \{\alpha_j,\alpha_{j+1}:\, \pi(j)=\pi(j+1)\}$, and $\cI^* = \cup_{j=1}^K\{\alpha_j:\, \left|J_{\pi(j)}\right|=1\}$. We have $\cI\sqcup\cI^* = \{\alpha_1,\cdots,\alpha_K\}$. We have $\left|\cI^*\right|=\sum_{a=1}^k \1_{\{ \left|J_a\right|=1 \}}=\left|\cJ^c\right|$.

Let $\alpha_1,\cdots,\alpha_K\in[d]$. For $a^*\in[k]$ such that $i_{a^*}=2$, we have $\D ^{i_{a^*}}g(X)=\D ^2 g(X) = 2I_d$, and for $a^*\in[k]$ such that $i_{a^*}=1$, we have $\D ^{i_{a^*}}g(X)=2X^\top$. We therefore have: for any $j\in[K]$, if $\left|J_{\pi(j)}\right|=2$, then there exists an $j'=j\pm 1$ such that $\pi(j')=\pi(j)$ and $(\D ^{\left|J_{\pi(j)}\right|}g(X))_{\alpha_j,\alpha_{j'}}=\1_{\{\alpha_j=\alpha_{j'}\}}$; if $\left|J_{\pi(j)}\right|=1$, then $((\D ^{\left|J_{\pi(j)}\right|})g(X)^\top)_{\alpha_j}=x_{\alpha_j}$. Therefore,
\begin{align*}
    \left(\D ^{i_1}g\otimes\cdots\otimes\D ^{i_k}g\right)_{\alpha_1,\cdots,\alpha_K}^\top = \prod_{\substack{ j\in[K]\\ \left| J_{\pi(j)} \right|=1 }} x_{\alpha_j} \prod_{\substack{j\in[K]\\ \left|J_{\pi(j)}\right|=2 }} \1_{\{\alpha_j=\alpha_{j'}\}}.
\end{align*}

 Depending on the relation between $2k$ and $K$, there are three regimes in general:
\begin{enumerate}
    \item When $2k<K$, it is impossible to have $i_1+\cdots+i_k=K$ while $i_1,\cdots,i_k\in\{1,2\}$.
    \item When $2k= K$, $\left|\cJ^c\right|=0$, and $\left|\cJ\right|=K/2$.

    In this case,
    \begin{align*}
        \bE\left[\norm{X}_2^{2\iota-2k}\D ^{i_1}g\otimes\cdots\otimes\D ^{i_k}g\right] = 2^k\bE[\norm{X}_2^{2\iota-2k}](I_d)^{\otimes k}.
    \end{align*}Therefore,
    \begin{align*}
        \norm{ \bE\left[\norm{X}_2^{2\iota-2k}\D ^{i_1}g\otimes\cdots\otimes\D ^{i_k}g\right] }_{HS}^2& = 2^{2k}\left(\bE\norm{X}_2^{2\iota-2k}\right)^2d^k\\
        &\lesssim 2^{2k}\norm{x_1}_{\psi_2}^{4(\iota-k)} (2(\iota-k))^{\iota-k} d^{2\iota-k} .
    \end{align*}
    \item When $2k>K$, $\left|\cJ^c\right|=2k-K\geq 1$, and $\left|\cJ\right|=K-k$. Recall that
    \begin{align*}
        (\norm{X}_2^2)^{\iota-k} = \sum_{\substack{ n_1+\cdots+n_d=\iota-k\\ \forall i\in[d], n_i\geq 0 }} \frac{(\iota-k)!}{n_1!\cdot\cdots\cdot n_d!} x_1^{2n_1}\cdot\cdots\cdot x_d^{2n_d}.
    \end{align*}
As we have assumed that $d>\iota$, there exists $i^*\in[d]$ such that $n_{i^*}=0$. Given $(\alpha_j)_{j=1}^K\in[d]^K$, if there exists $j^*\in[K]$ such that $\alpha_{j^*}=i^*$, and $\left|J_{\pi(j^*)}\right|=1$, then $(\bE[(\norm{X}_2^2)^{\iota-k}\D ^{i_1}g\otimes\cdots\otimes\D ^{i_k}g])_{\alpha_1,\cdots,\alpha_K}=0$. Therefore,
\begin{align}
    &\left(\bE[(\norm{X}_2^2)^{\iota-k}\D ^{i_1}g\otimes\cdots\otimes\D ^{i_k}g]\right)_{\alpha_1,\cdots,\alpha_K}\notag\\
    &= \bE\left[ \prod_{\substack{ j\in[K]\\ \left| I_{\pi(j)} \right|=1 }} x_{\alpha_j} \sum_{\substack{ n_1+\cdots+n_d=\iota-k\\ \forall i\in[d], n_i\geq 0 }} \frac{(\iota-k)!}{n_1!\cdot\cdots\cdot n_d!} x_1^{2n_1}\cdot\cdots\cdot x_d^{2n_d}\prod_{\substack{j\in[K]\\ \left|I_{\pi(j)}\right|=2 }} \1_{\{\alpha_j=\alpha_{j'}\}} \right]\notag\\
    &= \sum_{\substack{ n_1+\cdots+n_d=\iota-k\\ \forall i\in[d], n_i\geq 0 }} \frac{(\iota-k)!}{n_1!\cdot\cdots\cdot n_d!} \prod_{i^*\in\cI^*}\bE\left[x_{i^*}^{2n_{i^*}+1}\right] \prod_{i\notin\cI^*}\bE\left[ x_i^{2n_i} \right] \prod_{(i_1,i_1+1)\in\cI} \1_{\{i_1=i_1+1\}}. \label{eq:proof_middle_1}
\end{align}

To have \eqref{eq:proof_middle_1} not equals to $0$, we necessarily need for any $i^*\in\cI^*$, $n_{i^*}\geq 1$. Therefore,
\begin{align*}
    \eqref{eq:proof_middle_1}&\leq \sum_{\substack{ n_1+\cdots+n_d=\iota-k\\ \forall i\in[d], n_i\geq 0\\ \forall i^*\in\cI^*, n_{i^*}\geq 1 }} \frac{(\iota-k)!}{n_1!\cdot\cdots\cdot n_d!} \prod_{i^*\in\cI^*}\bE\left[x_{i^*}^{2n_{i^*}+1}\right] \prod_{i\notin\cI^*}\bE\left[ x_i^{2n_i} \right] \\
    &\lesssim (2(\iota-k)+1)^{\iota+\frac{k}{2} } \norm{x_1}_{\psi_2}^{2(\iota-k)+\left|\cI^*\right|} \sum_{\substack{ n_1+\cdots+n_d=\iota-k\\ \forall i\in[d], n_i\geq 0\\ \forall i^*\in\cI^*, n_{i^*}\geq 1 }} \frac{(\iota-k)!}{n_1!\cdot\cdots\cdot n_d!}\\
    &\leq (2(\iota-k)+1)^{\iota+\frac{k}{2} }\norm{x_1}_{\psi_2}^{2\iota-k} (\iota-k)! \sum_{\substack{ n_1+\cdots+n_d=\iota-k\\ \forall i\in[d], n_i\geq 0\\ \forall i^*\in\cI^*, n_{i^*}\geq 1 }}1\\
    &= (2(\iota-k)+1)^{\iota+\frac{k}{2} }\norm{x_1}_{\psi_2}^{2\iota-k} (\iota-k)! \binom{d+(\iota-k)-1-\left|\cI^*\right|}{\iota-k-\left|\cI^*\right|}\\
    &\leq (2(\iota-k)+1)^{\iota+\frac{k}{2} }\norm{x_1}_{\psi_2}^{2\iota-k} (\iota-k)! \left(\frac{2ed}{\iota-k-1}\right)^{\iota-k-1}.
\end{align*}
Therefore,
\begin{align}
    &\gamma_\iota^2 d^{-2\iota}\norm{  \bE(\norm{X}_2^2)^{\iota-k}\left( \sum_{\substack{i_1,\cdots,i_k\in\{1,2\}\\ i_1+\cdots+i_k=K }} \D ^{i_1}g\otimes \cdots\otimes \D ^{i_k}g \right)}_{HS}^2\notag\\
    &<\gamma_\iota^2(2(\iota-k)+1)^{2\iota+k } d^{-2\iota}d^k \norm{x_1}_{\psi_2}^{2(2\iota-k)} ((\iota-k)!)^2 \left(\frac{2e}{\iota-k-1}\right)^{2(\iota-k-1)} d^{2(\iota-k-1)}\notag\\
    &= \gamma_\iota^2(2(\iota-k)+1)^{2\iota+k } \norm{x_1}_{\psi_2}^{2(2\iota-k)} ((\iota-k)!)^2 \left(\frac{2e}{\iota-k-1}\right)^{2(\iota-k-1)} d^{-k-2}\notag\\
    &\sim \gamma_\iota^2\left(2(\iota-k)+1\right)^{2\iota+k}\norm{x_1}_{\psi_2}^{2(2\iota-k)}d^{-k-2} . \label{eq:K>2_2}
\end{align}
\end{enumerate}

We then deal with the first term in \eqref{eq:equivalent_k_derivative}. 

\paragraph{The case for $\iota=K$: $\norm{\bE (X^\top)^{\otimes K}}_{HS}$} Given any $\valpha:=(\alpha_1,\cdots,\alpha_K)\in[d]^{ K}$, we construct a graph $G(\valpha)$ with vertex set $[K]$. Let $\{i\}$ is connected with $\{j\}$, if $\alpha_i=\alpha_j$. To have $(\bE (X^\top)^{\otimes K})_{\valpha}\neq 0$, we necessarily need $G(\valpha)$ has no isolated vertex. Let $E(\valpha)$ as the set of edges of $G(\valpha)$, then $\left|E(\valpha)\right|\leq \lceil K/2\rceil$. For each $\valpha$ such that $\left|E(\valpha)\right|\leq \lceil K/2\rceil$, we can bound $(\bE(X^\top)^{\otimes K})_{\valpha}$ from above as $\oC{C_subgaussian}K^{K/2}\norm{x_1}_{\psi_2}^K$. Therefore,
\begin{align*}
    \norm{\bE (X^\top)^{\otimes K}}_{HS}^2 &= \norm{\bE X^{\otimes K}}_{HS}^2 \lesssim K^K\norm{x_1}_{\psi_2}^{2K} \sum_{i=0}^{\lceil K/2\rceil}\sum_{\substack{\valpha\in [d]^{ K}\\ \left|E(\valpha)\right|=i}} d^i \lesssim K^K\norm{x_1}_{\psi_2}^{2K} \sum_{i=0}^{\lceil K/2\rceil} d^i \binom{K}{2}^i\\
    &\lesssim K^K\norm{x_1}_{\psi_2}^{2K} \sum_{i=0}^{\lceil K/2\rceil} d^i K^{2i}\lesssim \norm{x_1}_{\psi_2}^{2K} d^{\lceil K/2\rceil} K^{2K }.
\end{align*}Therefore,
\begin{align}\label{eq:K>2_2_1}
    \frac{\gamma_K^2\norm{\bE (X^\top)^{\otimes K}}_{HS}^2}{d^K}\lesssim \gamma_\iota^2\norm{x_1}_{\psi_2}^{2K} \iota^{2\iota} d^{-\iota/2+1/2}.
\end{align}

\paragraph{The case for $K<\iota$: $\1_{\{K<\iota\}}\norm{\bE \norm{X}_2^{2\iota-2K} (X^\top)^{\otimes K} }_{HS}$}
Given $\valpha :=(\alpha_1,\cdots,\alpha_K)\in[d]^{ K}$, we construct a undirected graph $G(\valpha)$ with vertex $[K]$, and $\{i\}$ is connected with $\{j\}$ for $i,j\in[K]$, if $\alpha_i=\alpha_j$. We define a equivalence relation on $[K]$, denoted by $i\sim j$, where $i,j\in[K]$, if there is a path connecting $\{i\}$ and $\{j\}$. We partition $[K]$ into equivalence classes $\cC(\valpha)=\cup_{l} \cC_l(\valpha)$, that is, for every $i,j\in\cC_l(\valpha)$ (possibly $i=j$), we have $i\sim j$, and $\cC_l(\valpha)$ is largest, in the sense that if $i\sim j$, $i\in\cC_l$, then $j\in\cC_l$; for every $i\in\cC_l(\valpha)$, $j\in\cC_m(\valpha)$, $l\neq m$, we have $i \not\sim j$. Define the map $q:j\in[K]\mapsto q(j)\in[\left|\cC(\valpha)\right|]$ such that $j\in \cC_{q(j)}(\valpha)$, that is, the map $q$ maps each vertex into its unique equivalent class.
Let $\cO^*(\valpha)=\{\alpha_j:\, j\in[K], \left|\cC_{q(j)}(\valpha)\right|=1\}$, and $\cO^{**}(\valpha)=\{\alpha_j:\, j\in[K], \left|\cC_{q(j)}(\valpha)\right|>1\}$. We have the following relation: $\left|\cC(\valpha)\right|\leq \lfloor \frac{K-\left|\cO^*(\valpha)\right|}{2}\rfloor$. This is because $\cO^*(\valpha)$ and $\{\cC_l(\valpha): l\in[\left|\cC(\valpha)\right|], \left|\cC_l(\valpha)\right|=1\}$ are one-to-one corresponded; and for $\cC_l(\valpha)$ with $\left|\cC_l(\valpha)\right|>1$, the smallest size of $\cC_l(\valpha)$ is $2$.

Notice that $\cO^*(\valpha)\cup\cO^{**}(\valpha)$ might not be $[d]$. This is illustrated in the following example: let $K=8$, $d=100$, $\{\alpha_1,\cdots,\alpha_8\}=\{1,1,2,2,3,3,3,4\}$, then $\{1\},\{2\}$ are connected, so as $\{3\},\{4\}$ and $\{5\},\{6\},\{7\}$. Moreover, $1\sim 2$, $3\sim 4$, $5\sim 6\sim 7$, thus $\cC(\valpha)=\cC_1(\valpha)\cup\cC_2(\valpha)\cup\cC_3(\valpha)\cup\cC_4(\valpha)$ where $\cC_1(\valpha)=\{1,2\}$, $\cC_2(\valpha)=\{3,4\}$, $\cC_3(\valpha)=\{5,6,7\}$ and $\cC_4(\valpha)=\{8\}$. For any $j\in[8]$, only for $j=8$ we have $\left|\cC_{q(8)}(\valpha)\right|=\left|\cC_4(\valpha)\right|=1$, hence $\cO^*(\valpha)=\{\alpha_8\}=\{4\}$. Similarly, $\cO^{**}(\valpha)=\{1,2,3\}$. Thus $\cO^*(\valpha)\cup\cO^{**}(\valpha)=[4]$, however $[d]=[100]$.

\begin{align}
    &\left( \bE \norm{X}_2^{2\iota-2K} X^{\otimes K} \right)_{\valpha}\notag\\
    &= \bE\left[ \left(\sum_{\substack{ n_1+\cdots+n_d = \iota-K\\ \forall i\in[d], n_i\geq 0 }} \frac{(\iota-K)!}{n_1!\cdot\cdots\cdot n_d!} x_1^{2n_1}\cdot\cdots\cdot x_d^{2n_d} \right) \left(\prod_{\substack{ \alpha_j\in\cO^*(\valpha) }} x_{\alpha_j} \right)\left( \prod_{\substack{ \alpha_j\in\cO^{**}(\valpha) }} x_{\alpha_j}\right) \right].\label{eq:proof_head_1}
\end{align}We now seek $n_1,\cdots,n_d$ such that \eqref{eq:proof_head_1} is non-zero. A necessary condition is as follows: for every $i\in \cO^*(\valpha)$, we must have $n_i\geq 1$. Therefore,
\begin{align*}
    \eqref{eq:proof_head_1}&< (2\iota-K)^\iota \sum_{\substack{ n_1+\cdots+n_d=\iota-K\\ \forall i\in[d], n_i\geq 0\\ \forall i\in \cO^*(\valpha), n_i\geq 1 }} \frac{(\iota-K)!}{n_1!\cdot\cdots\cdot n_d!} \norm{x_1}_{\psi_2}^{2\iota-K}\\
    &< (2\iota-K)^\iota\norm{x_1}_{\psi_2}^{2\iota-K}(\iota-K)!\binom{d+(\iota-K)-1 -\left|\cO^*(\valpha)\right| }{\iota-K-\left|\cO^*(\valpha)\right|},
\end{align*}where we emphasize that when $\iota-K-\left|\cO^*(\valpha)\right|=0$, the binomial coefficient defined as $1$, and when $\iota-K-\left|\cO^*(\valpha)\right|<0$, the binomial coefficient is $0$. Before we taking square and taking sum over all $\valpha\in[d]^k$, we need to count the ``degree of freedom''. Recall the example $K=8$, $d=100$ and $\{\alpha_1,\cdots,\alpha_8\}=\{1,1,2,2,3,3,3,4\}$. In this example, if we take sum over all $\valpha$ with the edge of associated $G(\valpha)$ are the same, it is equivalent to taking sum over $x_1=x_2$, $x_3=x_4$, $x_5=x_6=x_7$, and none of $x_1,x_3,x_5,x_8$ are equal. Therefore, the sum is over $x_1,x_3,x_5,x_8\in[d]$, thus have $d^4=d^{\left|\cC(\valpha)\right|}$ terms. Moreover, given $\left|\cO^*(\valpha)\right|=\chi$ where $\chi$ is fixed, there are at most $(K-\chi)^{K-\chi}$ choices of $G(\valpha)$.

Now, we take square and sum over all $\valpha=(\alpha_1,\cdots,\alpha_K)\in[d]^{ k}$,
\begin{align*}
    &\norm{\bE(\norm{X}_2^2)^{\iota-K}(X^\top)^{\otimes K}}_{HS}^2 = \norm{\bE(\norm{X}_2^2)^{\iota-K} X^{\otimes K}}_{HS}^2\\
    &=\sum_{\chi=1}^K \sum_{\substack{\valpha\in[d]^{\otimes K}\\ \left|\cO^*(\valpha)\right|=\chi }} \left( \bE \norm{X}_2^{2\iota-2K} X^{\otimes K} \right)_{\valpha}^2\\
    &\lesssim \sum_{\chi=1}^{\iota-K} (K-\chi)^{K-\chi} d^{\lfloor \frac{K-\chi}{2}\rfloor}  \norm{x_1}_{\psi_2}^{2(2\iota-K)} (2\iota-K)^{2\iota} ((\iota-K)!)^2 \left( \frac{2ed}{\iota-K-\chi} \right)^{2(\iota-K-\chi)}\\
    &\lesssim (K-1)^{K-1} \norm{x_1}_{\psi_2}^{2(2\iota-K)} (2\iota-K)^{2\iota}((\iota-K)!)^2 \left(\frac{2e}{\iota-K}\right)^{2(\iota-K)} \sum_{\chi=1}^{\iota-K} d^{2\iota - \frac{3}{2}K - \frac{5}{2}\chi}\\
    &\lesssim (K-1)^{K-1} \norm{x_1}_{\psi_2}^{2(2\iota-K)} (2\iota-K)^{2\iota}((\iota-K)!)^2 \left(\frac{2e}{\iota-K}\right)^{2(\iota-K)} d^{2\iota - \frac{3}{2}K - \frac{5}{2}}.
\end{align*}
As a result,
\begin{align}\label{eq:K>2_2_2}
    \frac{\gamma_\iota^2 \norm{\bE(\norm{X}_2^2)^{\iota-K}(X^\top)^{\otimes K}}_{HS}^2 }{d^{2\iota}}\lesssim \gamma_\iota^2(K-1)^{K-1} \norm{x_1}_{\psi_2}^{2(2\iota-K)} (2\iota-K)^{2\iota} d^{2\iota - \frac{3}{2}K - \frac{5}{2}}.
\end{align}
Combining \eqref{eq:K=1}, \eqref{eq:K=2}, \eqref{eq:K>2_2}, \eqref{eq:K>2_2_1} and \eqref{eq:K>2_2_2}, we know that: for any $r\in[2L]$, we have,
\begin{align}\label{eq:result_HS_exp_derivative}
    \norm{\bE\D ^r K(X,X)}_{HS}&\leq \sum_{\iota=1}^L \norm{\bE\D ^r \frac{\gamma_\iota}{d^\iota}(\norm{X}_2^2)^\iota }_{HS} \lesssim
    \frac{ 2^{2L} L^{3L}}{\sqrt{d}} (\norm{x_1}_{\psi_2}^{2L-1}\vee 1)\sum_{\iota=1}^L  \gamma_\iota.
\end{align}Plug it into \eqref{eq:proof_polynomial_HS_1}, we obtain that for any $r\in[2L]$,
\begin{align*}
    \left(\norm{x_1}_{\psi_2}^r\max\left(\norm{\bE\D ^r K(X,X)}_\cJ: \cJ\in P_r\right)\right)^{2/r}\lesssim \left(\frac{ 2^{2L} L^{3L}}{\sqrt{d}} (\norm{x_1}_{\psi_2}^{2L+r-1}\vee \norm{x_1}_{\psi_2}^r)\sum_{\iota=1}^L  \gamma_\iota \right)^{2/r}.
\end{align*}By Theorem~\ref{theo:adamczak_wolff_polynomial}, we have, for any $t>0$:
\begin{align*}
    \bP\left(\left|K(X,X)-\bE K(X,X)\right|\geq t\right) \leq \exp\left(-\oC{C_non_linear_HS}\frac{t^{\frac{1}{L}} d^{\frac{1}{2L}} }{ 4L^3 \left(\sum_{\iota=1}^L \gamma_\iota\right)^{\frac{1}{L}} }\right).
\end{align*}

\subsubsection{Verifying \eqref{eq:diagonal_term_assumption} from Assumption~\ref{assumption:RIP} and Assumption~\ref{assumption:upper_dvoretzky} under sub-Gaussian assumption.}\label{sec:proof_diagonal_concentration}

We apply Theorem~\ref{theo:adamczak_wolff_polynomial} to $g(X)=\norm{\Gamma_{1:k}^{-1/2}\phi_{1:k}(X)}_\cH^2$, and $g(X) = \norm{\Gamma_{k+1:\infty}^{1/2}\phi_{k+1:\infty}(X)}_\cH^2$.

\paragraph{Checking diagonal concentration in Assumption~\ref{assumption:upper_dvoretzky}, $g(X) = \norm{\Gamma_{k+1:\infty}^{1/2}\phi_{k+1:\infty}(X)}_\cH^2$} By the definition of $\Gamma_{k+1:\infty}$ and $\phi_{k+1:\infty}(X)$, we know that $\Gamma_{k+1:\infty}^{1/2}\phi_{k+1:\infty}(X)=\sum_{j>k}\sqrt{\sigma_j}\varphi_j(X)\varphi_j$. For each $j\in\bN_+$, let $f_j = \frac{1}{\sqrt{\sigma_j}}\varphi_j$ be the corresponding ONB in $L_2(\mu)$, then $\norm{\Gamma_{k+1:\infty}^{1/2}\phi_{k+1:\infty}(X)}_\cH^2=\sum_{j>k}\sigma_j^2 f_j^2(X)$. Compared to $\norm{\phi_{k+1:\infty}(X)}_\cH^2 = \sum_{j>k}\sigma_j f_j^2(X)$, we observe that the eigenvalues are replaced by their square. 

For any $r\geq 0$, $i_1,\cdots,i_r\in[d]$ and $X\in\bR^d$, we have
\begin{align*}
    \left(\D ^r \norm{\Gamma_{k+1:\infty}^{1/2}\phi_{k+1:\infty}(X)}_\cH^2\right)_{i_1,\cdots,i_r} &= \left(\frac{\partial^r}{\partial x_1\cdots\partial x_r}\right)\sum_{j>k}\sigma_j^2 f_j^2(X)\\
    &\leq \sigma_{k+1}\left(\frac{\partial^r}{\partial x_1\cdots\partial x_r}\right)\sum_{j>k}\sigma_j f_j^2(X)\\
    &= \sigma_{k+1}\left(\D ^r \norm{\phi_{k+1:\infty}(X)}_\cH^2\right)_{i_1,\cdots,i_r}.
\end{align*}
By the positive homogeneity of expectation and Hilbert-Schmidt norm, we know that for any $r$, we have
\begin{align*}
    \norm{\bE \D ^r \norm{\Gamma_{k+1:\infty}^{1/2}\phi_{k+1:\infty}(X)}_\cH^2}_{HS} \leq \sigma_{k+1} \norm{\bE \D ^r \norm{\phi_{k+1:\infty}(X)}_\cH^2}_{HS},
\end{align*}where upper bounds for the right-hand-side may be found in \eqref{eq:result_HS_exp_derivative}. We apply Theorem~\ref{theo:adamczak_wolff_polynomial} with $t=1000\Tr(\Gamma_{k+1:\infty}^2)$ for $k=\sum_{i=0}^\iota d^\iota$. Notice that at this time, $\Tr(\Gamma_{k+1:\infty}^2)\sim d^{-(\iota+1)}$ (see the discussion after Lemma~\ref{lemma:gram_schmidt_LRZ}), and $\sigma_{k+1}\sim d^{-(\iota+1)}$. By Theorem~\ref{theo:adamczak_wolff_polynomial} together with \eqref{eq:result_HS_exp_derivative}, we can take $\odelta{delta_DMU_L2}=999$ and $\ogamma{gamma_DMU_L2} \sim \gamma$.

\paragraph{Checking diagonal term concentration in Assumption~\ref{assumption:RIP}, $g(X)=\norm{\Gamma_{1:k}^{-1/2}\phi_{1:k}(X)}_\cH^2$} The idea is again to use Theorem~\ref{theo:adamczak_wolff_polynomial}, by noticing that for any $r\geq 0$, $i_1,\cdots,i_r\in[d]$ and $X\in\bR^d$, we have
\begin{align*}
    \left(\D ^r \norm{\Gamma_{1:k}^{-1/2}\phi_{1:k}(X)}_\cH^2\right)_{i_1,\cdots,i_r} &= \left(\frac{\partial^r}{\partial x_1\cdots\partial x_r}\right)\sum_{j\leq k} f_j^2(X)\\
    &\leq \sigma_{k}^{-1}\left(\frac{\partial^r}{\partial x_1\cdots\partial x_r}\right)\sum_{j\leq k}\sigma_j f_j^2(X)\\
    &= \sigma_k^{-1}\left(\D ^r \norm{\phi_{1:k}(X)}_\cH^2\right)_{i_1,\cdots,i_r}.
\end{align*}As a result, for any $r$, we have
\begin{align*}
    \norm{\bE \D ^r \norm{\Gamma_{1:k}^{-1/2}\phi_{1:k}(X)}_\cH^2}_{HS} \leq \sigma_k^{-1} \norm{\bE \D ^r \norm{\phi_{1:k}(X)}_\cH^2}_{HS},
\end{align*}where the upper bound for the right-hand-side may be found in \eqref{eq:result_HS_exp_derivative}. We apply Theorem~\ref{theo:adamczak_wolff_polynomial} with $t=1000k$, $k=\sum_{i=0}^\iota d^i$. Notice that $\sigma_k^{-1}\sim d^\iota\sim k$, we can take $\odelta{delta_RIP} = 999$ and $\ogamma{gamma_RIP}\sim\gamma$.

\subsubsection{Proof of Proposition~\ref{prop:norm_equivalence_sub_gaussian}}\label{sec:proof_norm_equivalence_sub_gaussian}

In this subsection, we verify that \eqref{eq:norm_equivalence}, \eqref{eq:diagonal_upper_dvoretzky} and \eqref{eq:diagonal_RIP} hold under Assumption~\ref{assumption:polynomial_feature_sub_gaussian} with $\eps=6$. The following Lemma is taken from  \cite[Lemma 10]{liang_multiple_2020}:
\begin{Lemma}\label{lemma:liang_et_al}
    Suppose [1] of Assumption~\ref{assumption:polynomial_feature_sub_gaussian} holds. There exists $\kappa_L\geq 1$ such that for any $\beta\in\ell_2^{\binom{d+L}{L}}$ and $f_\beta(\cdot) := \sum_{\left|\cI\right|\leq L}\beta_\cI Q_\cI(\cdot)$, we have
    \begin{align*}
        \norm{f_\beta}_{L_4} \leq \kappa_L \norm{f_\beta}_{L_2} = \kappa_L \sqrt{\sum_{\left|\cI\right|\leq L}\beta_\cI^2}.
    \end{align*}
\end{Lemma}

The original assumption of \cite{liang_multiple_2020} in Lemma~\ref{lemma:liang_et_al} is: there exist absolute constants $\nu>1, \nC\label{C_LRZ}>0$ such that $\bP\left(\left|x_1\right|\geq t\right)\leq \oC{C_LRZ}(1+t)^{-\nu}$ for all $t\geq 0$. This assumption is implied by Assumption~\ref{assumption:polynomial_feature_sub_gaussian}. By Assumption~\ref{assumption:polynomial_feature_sub_gaussian}, $x_1$ has a finite $\psi_2$ norm $\norm{x_1}_{\psi_2}$ such that $\bP\left(\left|x_1\right|\geq t\right)\leq\exp\left(-\nC\label{C_subgaussian}t^2/\norm{x_1}_{\psi_2}^2\right)$ for some absolute constant $\oC{C_subgaussian}>0$. In fact, one can choose $\nu = \sqrt{2\oC{C_subgaussian}\log{\oC{C_LRZ}}/\norm{x_1}_{\psi_2}^2}>1$ and $\oC{C_LRZ}>\exp\left(\norm{x_1}_{\psi_2}^2/(2\oC{C_subgaussian})\right)>1$. This is because $\exp\left(-\oC{C_subgaussian}t^2/\norm{x_1}_{\psi_2}^2\right)\leq \oC{C_LRZ}(1+t)^{-\nu}$ for all $t\geq 0$ is equivalent to $(-\oC{C_subgaussian}/\norm{x_1}_{\psi_2}^2)t^2+\nu\log(1+t)-\log{\oC{C_LRZ}}\leq 0$ for all $t\geq 0$, and the above choice of $\oC{C_LRZ}$ and $\nu$ is a sufficient condition for it.

Lemma~\ref{lemma:liang_et_al} should not be confused with Bourgain's \( \Lambda_p \) problem \cite{bourgain_bounded_1989}, as the ONB here is polynomial and not uniformly bounded.

\beginproof[Proposition~\ref{prop:norm_equivalence_sub_gaussian}]
To prove Proposition~\ref{prop:norm_equivalence_sub_gaussian}, we simply note that $(Q_\cI)_\cI$ is an ONB of $L_2(\mu)$ by Lemma~\ref{lemma:gram_schmidt_LRZ}. Notice that for any $\gamma\in\ell_2^{\binom{d+L}{L}}$, $(f_\gamma)^2\in L_2(\mu)$ due to $\mu$ being sub-Gaussian and $(Q_\cI)$ being finite-degree polynomial. There therefore exists $\zeta\in\ell_2^{\binom{d+2L}{2L}}$ such that
\begin{align*}
     \left(f_\gamma(\cdot)\right)^2 = f_\zeta (\cdot) = \sum_{\left|\cI\right|\leq 2L}\zeta_\cI Q_\cI(\cdot).
\end{align*}
Applying Lemma~\ref{lemma:liang_et_al} with $L$ replaced by $2L$ and $\beta$ replaced by $\zeta$, we have
\begin{align*}
    \norm{f_\gamma}_{L_8}^2 = \norm{f_\zeta}_{L_4}\leq \kappa_{2L} \sqrt{\sum_{\left|\cI\right|\leq 2L}\zeta_\cI^2}.
\end{align*}On the other hand, we know that $\sum_{\left|\cI\right|\leq 2L}\zeta_\cI^2 = \norm{f_\gamma}_{L_4}^2$. This is because
\begin{align*}
    f_\gamma^4(X) = \left(f_\gamma^2(X)\right)^2=\left(f_\zeta(X)\right)^2 = \sum_{\left|\cI\right|,\left|\cJ\right|\leq 2L}\zeta_\cI\zeta_\cJ Q_\cI(X) Q_\cJ(X),
\end{align*}which implies that $\norm{f_\gamma}_{L_4}^4 = \sum_{\left|\cI\right|\leq 2L}\zeta_\cI^2$. Applying Lemma~\ref{lemma:liang_et_al} again with $\beta = \gamma$, we have $\sum_{\left|\cI\right|\leq 2L}\zeta_\cI^2 = \norm{f_\gamma}_{L_4}^4\leq\kappa_L^4\norm{f_\gamma}_{L_2}^4$. As a result,
\begin{align*}
    \norm{f_\gamma}_{L_8}\leq\sqrt{\kappa_{2L}}\kappa_L\norm{f_\gamma}_{L_2}.
\end{align*}

\endproof

Thus, Proposition~\ref{prop:spectrum_Gamma_sub_gaussian} as a corollary of Theorem~\ref{theo:upper_KRR} has been proven.

\subsubsection{Multiple descent in asymptotic regime}\label{sec:multiple_descent_asymptotic}

By \eqref{eq:pre_multiple_descent}, and when $f^*=f^{**}$, we only need to prove that $r_{0,k}^*=o_d(1)(\sigma_\xi + \|f^*\|_\cH)$ and the probability deviation tends to $1$ when $N,d\to\infty$ with $\omega_d(d^\iota)\leq N\leq o_d(d^{\iota+1})$. For the probability deviation, by \eqref{eq:prob_devi_multiple_descent_sub_Gaussian}, \eqref{eq:def_gamma_multiple_descent} and the fact that $\ogamma{gamma_DMU_L2} \sim \gamma$ and $\ogamma{gamma_RIP}\sim\gamma$, we have: $\eqref{eq:prob_devi_multiple_descent_sub_Gaussian}\to 1$. Moreover, by the definition of $r_{0,k}^*$ from \eqref{eq:def_rate},
\begin{enumerate}
    \item $\sigma_\xi\sqrt{d^\iota/N}=o_d(1)\sigma_\xi$;
    \item $\sigma_\xi\sqrt{N/d^{\iota+1}}=o_d(1)\sigma_\xi$;
    \item $\|\Gamma_{>\iota}^{1/2}f_{>\iota}^*\|_\cH\leq \|\Gamma_{k+1:\infty}^{1/2}\|_{\text{op}}\|f^*\|_\cH=o_d(1)\|f^*\|_\cH$, see the discussion after Lemma~\ref{lemma:gram_schmidt_LRZ};
    \item $\|\Gamma_{\leq \iota}^{-1/2}f_{\leq \iota}^*\|N^{-1}\leq \|f^*\|_\cH \sigma_k(\Gamma)^{-1/2}N^{-1}\lesssim \|f^*\|_\cH d^{\iota}/N=o_d(1)\|f^*\|_\cH$.
\end{enumerate}Therefore,
\begin{align*}
    \left| \norm{\hat f_0 - f^*}_{L_2(\mu)} - \norm{f_{>\iota}^*}_{L_2(\mu)} \right|=o_{d,\bP}(1)\left(\sigma_\xi + \norm{f^*}_\cH\right).
\end{align*}

\bibliographystyle{alpha}
\bibliography{biblio}
\end{document}

%% file: packages_notes.tex
\usepackage{graphicx}
\usepackage{xcolor}
\usepackage{amsmath}
\usepackage{amssymb}
\usepackage{amscd}
\usepackage{bbm}
\usepackage{tikz}
\usetikzlibrary{patterns}
\usepackage{hyperref}
\hypersetup{
    bookmarks=true,         
    colorlinks=true,       
    linkcolor=red,          
    citecolor=green,        
    filecolor=magenta,      
    urlcolor=cyan           
}

\usepackage[utf8]{inputenc}
\usepackage{amsthm}
\usepackage{bbm} 
\usepackage{glossaries}
\usepackage[linesnumbered,lined,boxed,commentsnumbered]{algorithm2e}
\usepackage{xparse}
\usepackage[normalem]{ulem}
\usepackage{enumitem}

\usepackage{marginnote}

\newtheorem{Theorem}{Theorem}
\newtheorem{Assumption}{Assumption}
\newtheorem{Definition}{Definition}
\newtheorem{Lemma}{Lemma}
\newtheorem{Proposition}{Proposition}

\newtheorem{Remark}{Remark}

\newenvironment{Theorembis}[1]
     {\renewcommand{\theTheorem}{\ref{#1} (informal)}%
     \addtocounter{Theorem}{-1}%
     \begin{Theorem}}
     {\end{Theorem}} 
\newenvironment{Propositionbis}[1]
    {\renewcommand{\theProposition}{\ref{#1} (informal)}%
    \addtocounter{Proposition}{-1}%
    \begin{Proposition}}
    {\end{Proposition}}

\newcounter{BigConst}                     
\newcommand{\nC}{                   
    \refstepcounter{BigConst}             
    \ensuremath{{C_{\theBigConst}}}
    }
\newcommand{\oC}[1]{\ensuremath{C_{\ref*{#1}}}}  

\newcounter{SmallConst}                     
\newcommand{\nc}{                   
    \refstepcounter{SmallConst}             
    \ensuremath{{c_{\theSmallConst}}}
    }
\newcommand{\oc}[1]{\ensuremath{c_{\ref*{#1}}}}  

\newcounter{gamma}                     
\newcommand{\ngamma}{                   
    \refstepcounter{gamma}             
    \ensuremath{{\gamma_{\thegamma}}}
    }
\newcommand{\ogamma}[1]{\ensuremath{{\gamma_{\ref*{#1}}}}}  

\newcounter{delta}                     
\newcommand{\ndelta}{                   
    \refstepcounter{delta}             
    \ensuremath{{\delta_{\thedelta}}}
    }
\newcommand{\odelta}[1]{\ensuremath{\delta_{\ref*{#1}}}}  

\newcounter{theta}                     
\newcommand{\ntheta}{                   
    \refstepcounter{theta}             
    \ensuremath{\theta_{\thetheta}}
    }
\newcommand{\otheta}[1]{\ensuremath{{\theta_{\ref*{#1}}}}}  

\newcounter{kappa}                     
\newcommand{\nkappa}{                   
    \refstepcounter{kappa}             
    \ensuremath{\kappa_{\thekappa}}
    }
\newcommand{\okappa}[1]{\ensuremath{{\kappa_{\ref*{#1}}}}}  

%% file: commandes_guillaume.tex
\newcommand{\inr}[1]{\bigl< #1 \bigr>}

\newcommand{\norm}[1]{\left\|#1\right\|}%
\newcommand{\vertiii}[1]{{\left\vert\kern-0.25ex\left\vert\kern-0.25ex\left\vert #1 
    \right\vert\kern-0.25ex\right\vert\kern-0.25ex\right\vert}}

\newcommand{\vertiiii}[1]{{\left\vert\kern-0.25ex\left\vert\kern-0.25ex\left\vert\kern-0.25ex\left\vert #1 
    \right\vert\kern-0.25ex\right\vert\kern-0.25ex\right\vert\kern-0.25ex\right\vert}}

\newcommand\eps{\epsilon}

\newcommand{\beginproof}{{\bf Proof. {\hspace{0.2cm}}}}
\def \endproof
{{\mbox{}\nolinebreak\hfill\rule{2mm}{2mm}\par\medbreak}}

\DeclareMathOperator*{\argmin}{argmin}

\DeclareMathOperator{\diag}{diag}
\DeclareMathOperator{\Span}{span}

\def\ds1{\textrm{1\kern-0.25emI}} 

\newcommand{\D}{\ensuremath{\mathbf{D}}}

\newcommand{\1}{\ensuremath{\mathbbm{1}}}

\newcommand \cA{{\cal A}}
\newcommand \cB{{\cal B}}
\newcommand \cC{{\cal C}}

\newcommand \cF{{\cal F}}
\newcommand \cG{{\cal G}}
\newcommand \cH{{\cal H}}
\newcommand \cI{{\cal I}}
\newcommand \cJ{{\cal J}}
\newcommand \cK{{\cal K}}
\newcommand \cL{{\cal L}}
\newcommand \cM{{\cal M}}
\newcommand \cN{{\cal N}}
\newcommand \cO{{\cal O}}

\newcommand \cQ{{\cal Q}}
\newcommand \cR{{\cal R}}

\newcommand \cV{{\cal V}}

\newcommand \cX{{\cal X}}
\newcommand \cY{{\cal Y}}
\newcommand \cZ{{\cal Z}}

\newcommand \bE{{\mathbb E}}

\newcommand \bG{{\mathbb G}}

\newcommand \bK{{\mathbb K}}

\newcommand \bN{{\mathbb N}}

\newcommand \bP{{\mathbb P}}

\newcommand \bR{{\mathbb R}}

\newcommand \bW{{\mathbb W}}
\newcommand \bX{{\mathbb X}}

\newcommand{\bbeta}{{\boldsymbol{\beta}}}
\newcommand{\blambda}{{\boldsymbol{\lambda}}}

\newcommand{\vmu}{{\boldsymbol{\mu}}}

\newcommand{\valpha}{{\boldsymbol{\alpha}}}

\newcommand{\va}{{\boldsymbol{a}}}

\newcommand{\vlambda}{{\boldsymbol{\lambda}}}
\newcommand{\vv}{{\boldsymbol{v}}}
\newcommand{\vh}{{\boldsymbol{h}}}
\newcommand{\vr}{{\boldsymbol{r}}}
\newcommand{\ve}{{\boldsymbol{e}}}
\newcommand{\vzero}{{\boldsymbol{0}}}
\newcommand{\vx}{{\boldsymbol{x}}}
\newcommand{\vz}{{\boldsymbol{z}}}

\newcommand{\vu}{{\boldsymbol{u}}}
\newcommand{\vw}{{\boldsymbol{w}}}

\newcommand{\vy}{\mathbf{y}}
\newcommand{\veps}{{\boldsymbol{\eps}}}

\newcommand{\vphi}{\boldsymbol{\phi}}

\newcommand{\bxi}{{\boldsymbol{\xi}}}

\DeclareMathOperator*{\Tr}{Tr}

\newcommand{\mysymbol}[3]{%
\newglossaryentry{#1}{%
      name={\ensuremath{#2}},%
      text={\ensuremath{#2}},%
      description={#3},%
      sort={#1}%
    }%
\expandafter\newcommand\expandafter{\csname smb#1\endcsname}{\gls{#1}}%
\expandafter\newcommand\expandafter{\csname #1\endcsname}{\ensuremath{#2}}%
}

\newcommand\hl{\bgroup\markoverwith
  {\textcolor{yellow}{\rule[-.5ex]{2pt}{2.5ex}}}\ULon}

\newif\ifPrintRemark